\documentclass[11pt]{amsart}
\usepackage{amssymb}
\usepackage{graphicx}                     % Pictures
\usepackage{hyperref}
\usepackage{setspace}                     % Control line spacing with, say, \begin{spacing}{0.9} blah blah \end{spacing}.
\usepackage{ifthen}                       % For \ifthenelse{cond}{true}{false}.

\newcommand{\adfmod}[1]{~(\mathrm{mod}~#1)}%     This is like \pmod except that I think it works properly.
\newcommand{\adfhide}[1]{}%                       Hide it.

\begin{document}
% \title[short-title]{title}
\title{Steiner systems with large block size}
\author{A. D. Forbes}
\address{LSBU Business School,
London South Bank University,
103 Borough Road,
London SE1 0AA, UK.}
\email{anthony.d.forbes@gmail.com}
\author{C. G. Rutherford}
\address{LSBU Business School,
London South Bank University,
103 Borough Road,
London SE1 0AA, UK.}
\email{c.g.rutherford@lsbu.ac.uk}
% \date{\today\ version 3.25}

\begin{abstract}
We construct Steiner systems $S(2, k, v)$ with $k > 9$ and small $v$.
\end{abstract}

\maketitle

%%%%%%%%%%%%%%%%%%%%%%%%%%%%%%%%%%%%%%%%%%%%%%%%%%%%%%%%%%%%%%%%%%%%%%%%%%%%%%%%%%%%%%%%%%
%%%%%%%%%%%%%%%%%%%%%%%%%%%%%%%%%%%%%%%%%%%%%%%%%%%%%%%%%%%%%%%%%%%%%%%%%%%%%%%%%%%%%%%%%%
%%%%%%%%%%%%%%%%%%%%%%%%%%%%%%%%%%%%%%%%%%%%%%%%%%%%%%%%%%%%%%%%%%%%%%%%%%%%%%%%%%%%%%%%%%
%%%%%%%%%%%%%%%%%%%%%%%%%%%%%%%%%%%%%%%%%%%%%%%%%%%%%%%%%%%%%%%%%%%%%%%%%%%%%%%%%%%%%%%%%%
%%%%%%%%%%%%%%%%%%%%%%%%%%%%%%%%%%%%%%%%%%%%%%%%%%%%%%%%%%%%%%%%%%%%%%%%%%%%%%%%%%%%%%%%%%
%%%%%%%%%%%%%%%%%%%%%%%%%%%%%%%%%%%%%%%%%%%%%%%%%%%%%%%%%%%%%%%%%%%%%%%%%%%%%%%%%%%%%%%%%%
\section{Introduction}\label{sec:Introduction}

A Steiner system $S(2, k, v)$ is a pair $(V,\mathcal{B})$,
where $V$ is a set of cardinality $v$, and
$\mathcal{B}$ is a set of $k$-element subsets of $V$ such that
every pair $\{x,y\}$ with $x, y \in V$ is a subset of precisely one $B \in \mathcal{B}$.
The elements of $V$ are the {\em points} of the Steiner system $(V,\mathcal{B})$,
and the elements of $\mathcal{B}$ are its {\em blocks}.
The parameter $k$ is the {\em block size}.
A Steiner system $S(2,k,v)$ is sometimes called a {\em balanced incomplete block design} and
denoted by $(v, k, 1)$-BIBD.
To avoid trouble we require $k \ge 2$, and we consider only positive $v$.

The number of blocks in an $S(2,k,v)$ is $v(v - 1)/(k(k - 1))$,
and necessary conditions for its existence are
\begin{align}
v(v - 1)/(k(k - 1)) &\ge v \text{~or~} v \in \{1,k\}, \nonumber \\
v(v - 1) &\equiv 0 \adfmod{k(k - 1)}, \label{eqn:necessary-conditions} \\
v - 1    &\equiv 0 \adfmod{k - 1}. \nonumber
\end{align}
The first item is known as Fisher's inequality, \cite[Theorem 1.1.3]{Anderson1997}.
For further information about block designs, see \cite{Anderson1997} for example.

Given $k$, by a theorem of Wilson, \cite{Wilson1975}, the conditions (\ref{eqn:necessary-conditions})
are sufficient for all large $v$, and hence
the determination of the $v$ for which an $S(2,k,v)$ exists is actually a finite problem.
However, currently it seems to be impossible to resolve every case not covered by `all large'
whenever $k$ is not small.
Indeed, the existence spectrum has been computed only for
$k = 2$ (trivial), $k = 3$ by Kirkman, \cite{Kirkman1847}, and $k \in \{4,5\}$ by Hanani, \cite{Hanani1961},
where in each case the necessary conditions (\ref{eqn:necessary-conditions}) suffice.
For $k \in \{6, 7, \dots, 9\}$ the spectrum is known only up to a smallish number of possible exceptions;
see \cite{AdamsBryantBuchanan2008} and the references therein.

The objective of this paper is to exhibit many Steiner systems with block size greater than 9.
We hope that some of them might be of interest to someone.

See \cite[Theorem II.5.11]{ColbournMathon2007} or \cite[Theorem 2.1]{AdamsBryantBuchanan2008}
for some infinite families of $S(2,k,v)$ with $k \ge 10$.
In particular, if $k$ is a prime power, there exists an $S(2,k,k^2)$, an affine plane of order $k$,
and if $k-1$ is a prime power, there exists an $S(2,k,k(k-1)+ 1)$, a projective plane of order $k - 1$.
Further infinite families are given in \cite{Nechvatal2014}.
Table II.1.35 of \cite{MathonRosa2007} includes all admissible parameters for Steiner systems $S(2,k,v)$
with $10 \le k \le v/2$ and $(v - 1)/(k - 1) \le 41$, and we note that existence is confirmed only for
affine planes, projective planes and the design $S(2, 16, 496)$.
Some $S(2,10,v)$s appear in \cite{Abel1996}.

%%%%%%%%%%%%%%%%%%%%%%%%%%%%%%%%%%%%%%%%%%%%%%%%%%%%%%%%%%%%%%%%%%%%%%%%%%%%%%%%%%%%%%%%%%
%%%%%%%%%%%%%%%%%%%%%%%%%%%%%%%%%%%%%%%%%%%%%%%%%%%%%%%%%%%%%%%%%%%%%%%%%%%%%%%%%%%%%%%%%%
%%%%%%%%%%%%%%%%%%%%%%%%%%%%%%%%%%%%%%%%%%%%%%%%%%%%%%%%%%%%%%%%%%%%%%%%%%%%%%%%%%%%%%%%%%
%%%%%%%%%%%%%%%%%%%%%%%%%%%%%%%%%%%%%%%%%%%%%%%%%%%%%%%%%%%%%%%%%%%%%%%%%%%%%%%%%%%%%%%%%%
%%%%%%%%%%%%%%%%%%%%%%%%%%%%%%%%%%%%%%%%%%%%%%%%%%%%%%%%%%%%%%%%%%%%%%%%%%%%%%%%%%%%%%%%%%
%%%%%%%%%%%%%%%%%%%%%%%%%%%%%%%%%%%%%%%%%%%%%%%%%%%%%%%%%%%%%%%%%%%%%%%%%%%%%%%%%%%%%%%%%%
\section{Steiner systems from single base blocks}

%
% DESIGN SPECIFICATION
%
% Format:
%
%xxxxDesignB{k}{t}{m}{v}{w}{0}{0}{o}{B}
%xxxxDesignC{k}{t}{m}{v}{w}{R}{P}{o}{B}
%xxxxDesignE{k}{t}{m}{v}{w}{e}{a}{o}{B}
%xxxxDesignF{k}{t}{m}{v}{w}{e}{a}{o}{B}
%
% type = 11th character of the design specification, "B", "C", "E", "F".
%
% xxxx = \adf;
% type = "B" : the base block of the Steiner system is given;
% type = "C" : the base block of the Steiner system is given;
% type = "E" : the base block of the Steiner system is (e, 1, a, a^2, ..., a^{k-2}), e not = a^{k-1};
% type = "F" : the base block of the Steiner system is (e, a, a^2, ..., a^{k-1}).
% k : the block size of the system.
% t : the number of points in the system is v = k(k - 1)t + m.
% m : the number of points in the system is v = k(k - 1)t + m.
% v : the number of points in the system is v = k(k - 1)t + m.
% w : the multiplier for developing the base block.
% R : the ring specification if X = "C".
% P : the polynomial for multiplication in the Galois field R.
% e : parameter for generating the base block.
% a : parameter for generating the base block.
% o : originator (if available).
% B : the base block; this is needed only if type = "B" or "C" but is supplied anyway.
%

A well-known method for creating a Steiner system $S(2,k,v)$, $v = k(k - 1)t + 1$,
is to find by some method a single base block $B = (b_1, b_2, \dots, b_k)$ and develop it into $tv$ blocks by
$x \mapsto w^e x + d$, $e \in \{0, 1, \dots, t - 1\}$, $d \in \mathcal{R}$, where $\mathcal{R}$ is a suitable ring of order $v$.
Typically $\mathcal{R}$ is $\mathbb{Z}_v$, the ring of integers modulo $v$,
but if $v$ is a composite prime power, then $\mathcal{R} = \mathrm{GF}(v)$ is also available.

Here are some examples.
% The ring $\mathcal{R}$ is $\mathbb{Z}_v$ unless specified otherwise.
If $\mathcal{R}$ is $\mathrm{GF}(p^n)$, $n \ge 2$,
then $P(x)$ is the polynomial for multiplication in the Galois field,
and the element $g(x)$ of $\mathrm{GF}(p^n)$ is represented by the number $g(p)$.

\newcommand{\adfDesignB}[9]{$S(2,#1,#4)$, $w=#5$, $\mathcal{R} = \mathbb{Z}_{#4}$,%
\begin{center} $B=(#9)$ \end{center}}

\newcommand{\adfDesignC}[9]{$S(2,#1,#4)$, $w=#5$, $\mathcal{R} = \mathrm{GF}(#6)$, $P(x) = #7$,%
\begin{center} $B=(#9)$ \end{center}}

%%%%%%%%%%%%%%%%%%%%%%%%%%%%%%%%%%%%%%%%%%%%%%%%%%%%%%%%%%
%%%%%%%%%%%%%%%%%%%%%%%%%%%%%%%%%%%%%%%%%%%%%%%%%%%%%%%%%%
% Designs go after here.
\adfDesignB{10}{13}{1}{1171}{166}{0}{0}{"Audrey 26 Mar 2014"}{{0,1,608,339,291,302,316,355,1029,1090}}

\adfDesignB{10}{17}{1}{1531}{502}{0}{0}{"Audrey 26 Mar 2014"}{{0,1,1059,1125,163,315,825,893,997,1226}}

\adfDesignB{10}{18}{1}{1621}{1434}{0}{0}{"Mononoke 26 Mar 2014"}{{0,1,213,984,809,71,170,623,868,1276}}

\adfDesignB{10}{20}{1}{1801}{1242}{0}{0}{"Audrey 26 Mar 2014"}{{0,1,1548,208,608,502,594,947,1197,1570}}

\adfDesignB{10}{24}{1}{2161}{669}{0}{0}{"Audrey 26 Mar 2014"}{{0,1,1921,397,265,89,268,900,1410,1456}}

\adfDesignB{10}{25}{1}{2251}{1057}{0}{0}{"Rebecca 26 Mar 2014"}{{0,1,711,463,1115,321,537,851,922,1625}}

\adfDesignB{10}{26}{1}{2341}{2335}{0}{0}{"Audrey 25 Mar 2014"}{{0,1,174,679,688,745,840,887,1889,1926}}

\adfDesignB{10}{28}{1}{2521}{1331}{0}{0}{"Audrey 26 Mar 2014"}{{0,1,2376,2050,56,353,1590,1710,2057,2122}}

\adfDesignB{10}{31}{1}{2791}{1627}{0}{0}{"Audrey 19 Apr 2014"}{{0,1,2378,2340,1191,19,265,1283,1645,2488}}

\adfDesignB{10}{33}{1}{2971}{1699}{0}{0}{"Audrey 19 Apr 2014"}{{0,1,1291,2202,1788,584,1046,1671,2144,2513}}

\adfDesignB{10}{34}{1}{3061}{2061}{0}{0}{"Audrey 23 Apr 2014"}{{0,1,1842,2243,2046,104,563,1875,2263,2356}}

\adfDesignB{10}{37}{1}{3331}{42}{0}{0}{"Audrey 26 Apr 2014"}{{0,1,1648,2777,1316,482,819,1045,2348,2481}}

\adfDesignB{10}{39}{1}{3511}{437}{0}{0}{"Felicity 27 Apr 2014"}{{0,1,3183,765,331,642,939,1243,1325,2503}}

\adfDesignB{10}{41}{1}{3691}{1752}{0}{0}{"Audrey 26 Apr 2014"}{{0,1,1855,3306,733,711,2382,2630,3389,3520}}

\adfDesignB{10}{45}{1}{4051}{77}{0}{0}{"Charlotte 7 Aug 2025"}{{0,1,3373,1352,1554,181,75,2081,2633,2635}}

\adfDesignB{10}{47}{1}{4231}{44}{0}{0}{"Tanya 21 Aug 2025"}{{0,1,550,1730,527,1767,1302,3701,1167,3318}}

\adfDesignB{10}{51}{1}{4591}{125}{0}{0}{"Tanya 18 Aug 2025"}{{0,1,207,4565,1011,3475,1732,1863,3901,3904}}

\adfDesignB{10}{54}{1}{4861}{59}{0}{0}{"Tanya 22 Aug 2025"}{{0,1,4474,2964,821,2077,3466,1379,2371,4462}}

\adfDesignB{10}{55}{1}{4951}{11}{0}{0}{"Tanya 26 Aug 2025"}{{0,1,4113,2013,338,4934,4782,3947,2281,3694}}

\adfDesignC{10}{56}{1}{5041}{4981}{71^2}{x^2 + 1}{"Charlotte 28 Jan 2016"}{{0,1,1501,2640,3534,1591,3575,3839,4556,4769}}

\adfDesignB{10}{62}{1}{5581}{2}{0}{0}{"Tanya 31 Aug 2025"}{{0,1,1739,1999,4188,5086,4119,1810,2653,3239}}

\adfDesignB{10}{65}{1}{5851}{30}{0}{0}{"Tanya 30 Aug 2025"}{{0,1,710,2724,2581,4353,4606,4127,407,2427}}

\adfDesignB{10}{68}{1}{6121}{209}{0}{0}{"Tanya 30 Aug 2025"}{{0,1,772,6054,2370,1421,3542,2362,1941,2318}}

\adfDesignB{10}{69}{1}{6211}{5}{0}{0}{"Tanya 30 Aug 2025"}{{0,1,2212,1846,1734,5776,35,1755,1206,1276}}

\adfDesignB{10}{70}{1}{6301}{197}{0}{0}{"Tanya 03 Sep 2025"}{{0,1,5612,5674,3782,5560,380,4946,937,1062}}

\adfDesignB{10}{72}{1}{6481}{112}{0}{0}{"Tanya 05 Sep 2025"}{{0,1,408,280,1638,333,2054,5643,1897,2273}}

\adfDesignB{10}{73}{1}{6571}{57}{0}{0}{"Tanya 04 Sep 2025"}{{0,1,4814,553,4165,1952,559,4821,2081,4525}}

\adfDesignB{10}{74}{1}{6661}{166}{0}{0}{"Tanya 08 Sep 2025"}{{0,1,4795,3489,5393,4549,2833,2931,4686,5566}}

\adfDesignB{10}{76}{1}{6841}{201}{0}{0}{"Tanya 04 Sep 2025"}{{0,1,50,1933,5662,6380,2636,5478,4301,4602}}

\adfDesignB{10}{84}{1}{7561}{78}{0}{0}{"Tanya 04 Sep 2025"}{{0,1,2226,3792,5275,5093,2099,725,4798,5538}}

\adfDesignB{10}{86}{1}{7741}{132}{0}{0}{"Tanya 08 Sep 2025"}{{0,1,705,3836,3614,702,3212,7547,2168,7427}}

\adfDesignC{10}{88}{1}{7921}{6319}{89^2}{x^2 + 3}{"Charlotte 28 Jan 2016"}{{0,1,5612,513,4483,2172,208,1210,1571,5580}}

\adfDesignB{10}{89}{1}{8011}{6}{0}{0}{"Doris 11 Sep 2025"}{{0,1,4339,4082,823,999,685,6481,3990,6161}}

\adfDesignB{10}{90}{1}{8101}{59}{0}{0}{"Ethel 09 Sep 2025"}{{0,1,4806,3335,6752,5504,7121,4511,4751,7351}}

\adfDesignB{10}{91}{1}{8191}{107}{0}{0}{"Tanya 06 Sep 2025"}{{0,1,5001,1037,6134,4421,2067,5219,6263,7767}}

\adfDesignB{10}{94}{1}{8461}{134}{0}{0}{"Ethel 08 Sep 2025"}{{0,1,7087,3453,1323,971,118,2074,3536,4210}}

\adfDesignB{10}{96}{1}{8641}{962}{0}{0}{"Charlotte 10 Sep 2025"}{{0,1,713,8278,3195,4877,947,310,1177,4619}}

\adfDesignB{10}{97}{1}{8731}{125}{0}{0}{"Tanya 08 Sep 2025"}{{0,1,4785,5261,6995,188,685,7979,660,8313}}

\adfDesignB{10}{98}{1}{8821}{39}{0}{0}{"Ethel 13 Sep 2025"}{{0,1,4871,5631,4287,3903,5445,206,2589,6230}}

\adfDesignB{10}{100}{1}{9001}{302}{0}{0}{"Charlotte 09 Sep 2025"}{{0,1,4429,33,1813,4010,2108,6132,5664,8349}}

\adfDesignB{10}{101}{1}{9091}{21}{0}{0}{"Tanya 14 Sep 2025"}{{0,1,366,5145,3609,473,6831,4038,665,730}}

\adfDesignB{10}{102}{1}{9181}{367}{0}{0}{"Tanya 10 Sep 2025"}{{0,1,5942,4325,3783,3081,952,8887,1006,4247}}

\adfDesignB{10}{107}{1}{9631}{29}{0}{0}{"Charlotte 09 Sep 2025"}{{0,1,6120,1839,2271,7230,1554,6288,5897,7991}}

\adfDesignB{10}{108}{1}{9721}{22}{0}{0}{"Ethel 12 Sep 2025"}{{0,1,632,2300,3899,1643,3079,2746,8860,8899}}

\adfDesignB{10}{109}{1}{9811}{38}{0}{0}{"Tanya 10 Sep 2025"}{{0,1,9101,2137,6759,144,1856,115,6539,7996}}

\adfDesignB{10}{110}{1}{9901}{24}{0}{0}{"Charlotte 11 Sep 2025"}{{0,1,782,1493,9015,3458,1518,8530,1112,5162}}

\adfDesignB{10}{117}{1}{10531}{45}{0}{0}{"Ethel 13 Sep 2025"}{{0,1,5030,2612,9226,8902,1297,5480,3109,8157}}

\adfDesignB{10}{119}{1}{10711}{59}{0}{0}{"Ethel 15 Sep 2025"}{{0,1,187,2170,8605,9819,6289,8909,1896,3225}}

\adfDesignB{10}{121}{1}{10891}{24}{0}{0}{"Ethel 14 Sep 2025"}{{0,1,2676,6033,983,8341,3308,8487,7051,9732}}

\adfDesignB{10}{123}{1}{11071}{19}{0}{0}{"Tanya 13 Sep 2025"}{{0,1,9504,7062,10996,1000,10962,3850,5408,11014}}

\adfDesignB{10}{124}{1}{11161}{11}{0}{0}{"Ethel 13 Sep 2025"}{{0,1,363,1421,9030,8079,6326,7271,1614,7353}}

\adfDesignB{10}{125}{1}{11251}{131}{0}{0}{"Charlotte 13 Sep 2025"}{{0,1,1053,7108,9683,5305,9001,6284,4908,5502}}

\adfDesignC{10}{132}{1}{11881}{5651}{109^2}{x^2 + 2}{"Charlotte 29 Jan 2016"}{{0,1,3949,10375,7563,859,2452,8283,10518,11572}}

\adfDesignC{11}{72}{1}{7921}{7914}{89^2}{x^2 + 3}{"Charlotte 29 Jan 2016"}{{0,1,98,382,4035,5348,5760,5998,6700,6955,7168}}

\adfDesignC{11}{108}{1}{11881}{1962}{109^2}{x^2 + 2}{"Charlotte 30 Jan 2016"}{{0,1,1670,1793,3812,4073,4842,5975,6806,7365,8266}}

\adfDesignC{11}{156}{1}{17161}{12562}{131^2}{x^2 + 1}{"Charlotte 30 Jan 2016"}{{0,1,1578,2199,2378,4086,6472,6945,10339,12266,17118}}

For clarification, we explain by an example how multiplication works in
our encoding of elements in a Galois field.
To multiply 1501 by 4981 in $\mathrm{GF}(71^2)$ with $P(x) = x^2 + 1$, we proceed thus:
$1501 = 21\cdot 71 + 10 \rightarrow 21x + 10$;
$4981 = 70\cdot 71 + 11 \rightarrow 70x + 11$;
$(21x + 10)(70x + 11) = 1470 x^2 + 931 x + 110
\rightarrow 931x - 1360$ (using $x^2 + 1 = 0$)
$\rightarrow 8x + 60
\rightarrow 8 \cdot 71 + 60 = 628$.

% 56 designs
% Designs go before here.
%%%%%%%%%%%%%%%%%%%%%%%%%%%%%%%%%%%%%%%%%%%%%%%%%%%%%%%%%%
%%%%%%%%%%%%%%%%%%%%%%%%%%%%%%%%%%%%%%%%%%%%%%%%%%%%%%%%%%

%%%%%%%%%%%%%%%%%%%%%%%%%%%%%%%%%%%%%%%%%%%%%%%%%%%%%%%%%%%%%%%%%%%%%%%%%%%%%%%%%%%%%%%%%%
%%%%%%%%%%%%%%%%%%%%%%%%%%%%%%%%%%%%%%%%%%%%%%%%%%%%%%%%%%%%%%%%%%%%%%%%%%%%%%%%%%%%%%%%%%
%%%%%%%%%%%%%%%%%%%%%%%%%%%%%%%%%%%%%%%%%%%%%%%%%%%%%%%%%%%%%%%%%%%%%%%%%%%%%%%%%%%%%%%%%%
%%%%%%%%%%%%%%%%%%%%%%%%%%%%%%%%%%%%%%%%%%%%%%%%%%%%%%%%%%%%%%%%%%%%%%%%%%%%%%%%%%%%%%%%%%
%%%%%%%%%%%%%%%%%%%%%%%%%%%%%%%%%%%%%%%%%%%%%%%%%%%%%%%%%%%%%%%%%%%%%%%%%%%%%%%%%%%%%%%%%%
%%%%%%%%%%%%%%%%%%%%%%%%%%%%%%%%%%%%%%%%%%%%%%%%%%%%%%%%%%%%%%%%%%%%%%%%%%%%%%%%%%%%%%%%%%
\section{A 2-parameter construction for $S(2,k,k(k-1)t+1)$} \label{sec:simple-construction-m-1}

Suppose $k \ge 3$ and $t \ge 1$.
Let $K = k(k - 1)/2$ and $v = 2Kt + 1$.
Suppose $v$ is prime, let $\rho$ be a primitive root modulo $v$ and write $y = \log x$,
the {\em logarithm} of $x$, if $x = \rho^y$.
We tacitly perform arithmetic modulo $v$ for elements of $\mathbb{Z}_v$ and modulo $v-1$ for their logarithms.
Let $a$ and $e$ be parameters for a base block of the form $B = \{e, 1, a, \dots, a^{k-2}\}$.
The objective is to choose $a$, $e$ and a multiplier $w$ such that the difference set
\begin{eqnarray*}
D &=& \{\pm w^h (a^s - a^r):~ 0 \le h \le t - 1, ~ 0 \le r < s \le k - 2\}\\
  & & \cup~ \{\pm w^h (a^r - e):~ 0 \le h \le t - 1, ~ 0 \le r \le k - 2\}
\end{eqnarray*}
is equal to $\{1, 2, \dots, v - 1\}$.
Let $w = \rho^K$. Then on taking logs and observing that $\log(-1) = K t$, $D$ becomes
\begin{eqnarray*}
\Delta &=& \{K t\epsilon + K h + \log(a^s - a^r):~ \\
       & & \hspace{5mm} \epsilon \in \{0,1\},~ 0 \le h \le t - 1, ~ 0 \le r < s \le k - 2\}\\
       & & \cup~ \{K t\epsilon + K h + \log(a^r - e):~ \\
       & & \hspace{5mm} \epsilon \in \{0,1\},~ 0 \le h \le t - 1, ~ 0 \le r \le k - 2\}.
\end{eqnarray*}
Observe that
$$\{K t \epsilon + Kh:~ \epsilon \in \{0,1\},~ 0 \le h \le t - 1\}$$
covers all multiples of $K$.
So we look for an $a$ such that
$$\{\log(a^s - a^r):~ 0 \le r < s \le k - 2\}$$
covers $(k-1)(k-2)/2$ distinct residues modulo $K$.
If successful we then look for an $e$ such that $\{\log(a^r - e):~ 0 \le r \le k - 2\}$ covers the
remaining $(k-1)$ residues modulo $K$.

If $t$ is odd we can choose $w = \rho^{2K}$, for in this case the set
$$\{K t \epsilon + 2Kh:~ \epsilon \in \{0,1\},~ 0 \le h \le t - 1\}$$
covers the multiples of $K$.
With $e = a^{k-1}$ the base block has the geometric progression format, $B = \{1, a, \dots, a^{k-1}\}$.

Steiner $S(2,k,v)$ designs can be constructed with small $v$, say less than a few million,
for the following values of $k$.
\begin{eqnarray*}
B = \{0, 1, a, \dots, a^{k-2}\}:   && k = \adfhide{3, 4, 5, 6, 7, 8,} 10, 12, 14, 16, 18, 20, 22, 28, 30; \\
B = \{1, a, a^2, \dots, a^{k-1}\}: && k = \adfhide{3, 4, 5, 6, 7, 8, 9,} 10, 11, 13, 15, 17, 19, 21, 23, 25, 27.
\end{eqnarray*}

%%%%%%%%%%%%%%%%%%%%%%%%%%%%%%%%%%%%%%%%%%%%%%%%%%%%%%%%%%%%%%%%%%%%%%%%%%%%%%%%%%%%%%%%%%
%%%%%%%%%%%%%%%%%%%%%%%%%%%%%%%%%%%%%%%%%%%%%%%%%%%%%%%%%%%%%%%%%%%%%%%%%%%%%%%%%%%%%%%%%%
%%%%%%%%%%%%%%%%%%%%%%%%%%%%%%%%%%%%%%%%%%%%%%%%%%%%%%%%%%%%%%%%%%%%%%%%%%%%%%%%%%%%%%%%%%
%%%%%%%%%%%%%%%%%%%%%%%%%%%%%%%%%%%%%%%%%%%%%%%%%%%%%%%%%%%%%%%%%%%%%%%%%%%%%%%%%%%%%%%%%%
%%%%%%%%%%%%%%%%%%%%%%%%%%%%%%%%%%%%%%%%%%%%%%%%%%%%%%%%%%%%%%%%%%%%%%%%%%%%%%%%%%%%%%%%%%
%%%%%%%%%%%%%%%%%%%%%%%%%%%%%%%%%%%%%%%%%%%%%%%%%%%%%%%%%%%%%%%%%%%%%%%%%%%%%%%%%%%%%%%%%%
\section{A 2-parameter construction for $S(2,k,k(k-1)t+k)$}\label{sec:simple-construction-m-k}

Suppose $v = qk$, where $q = (k - 1)t + 1$ is prime and let $\rho$ be a primitive root modulo $q$
with the additional property $\gcd(\rho, k) = 1$.
We do arithmetic modulo $v$ for elements of $\mathbb{Z}_v$ and modulo $q-1$ for their logarithms.
Let $a$ be a generator for a base block of the form $B = \{e, 1, a, \dots, a^{k-2}\}$.
The objective is to choose $a$ and a multiplier $w$ such that the difference set
\begin{eqnarray*}
D &=& \{\pm w^h (a^s - a^r): ~ 0 \le h \le t - 1, ~ 0 \le r < s \le k - 2\}\\
  & & \cup~ \{\pm w^h (a^r - e): ~ 0 \le h \le t - 1,~ 0 \le r \le k - 2\}
\end{eqnarray*}
is equal to $\{1, 2, \dots, v - 1\}\setminus \{q, 2q, \dots, (k - 1)q\}$.
Let
\begin{eqnarray*}
D_0 &=& \{\pm(a^s - a^r): ~ 0 \le r < s \le k - 2\} \\
    & & \cup \{\pm(a^r - e): ~ 0 \le r \le k - 2\},
\end{eqnarray*}
let $w = \rho^{k-1}$ and choose $a$ and $e$ such that the pairs
$$E ~=~ \{(\log(d \mathrm{~mod~} q) \mathrm{~mod~} (k - 1),~ d \mathrm{~mod~} k):~ d \in D_0\}$$
cover all pairs of residue classes modulo $(k - 1,~ k)$.
Let $d$ be an element of $D_0$ and suppose $d \equiv d_q \pmod{q}$, $d \equiv d_k \pmod{k}$.
Then $d_q \neq 0$ and $d$ can be written $d = \kappa d_q + \varphi d_k$,
where $\kappa = k(1/k \mathrm{~mod~} q)$ and $\varphi = q(1/q \mathrm{~mod~} k)$.
Observe that $k|\kappa$, $\kappa \equiv 1 \pmod{q}$, $q|\varphi$, $\varphi \equiv 1 \pmod{k}$.
So
$$\log(d \mathrm{~mod~} q) = \log d_q = (k-1)\delta + \epsilon,$$
say.
Thus $d = \kappa \rho^{(k-1)\delta + \epsilon} + \varphi d_k$ and
$$w^h d ~=~ \kappa \rho^{(k-1)(\delta + h) + \epsilon} + w^h \varphi d_k.$$
Therefore, since $\gcd(w^h \varphi, k) = 1$, $w^h D_0$ satisfies the same covering property as $D_0$.
Moreover, if $w^h \neq \pm1$ for $0 < h < t$, then $\{w^h B:~ 0 \le h < t\}$ together with
$\{0, q, 2q, \dots, (k - 1)q\}$ will be a valid set of starter blocks for a Steiner system $S(2, k, qk)$.

The construction yields $S(2,k,v)$ designs with small $v$ for the following values of $k$.
\begin{eqnarray*}
B = \{0, 1, a, \dots, a^{k-2}\}: && k = \adfhide{3, 5, 7,} 11, 13, 17, 19, 23, 29, 31.
\end{eqnarray*}

%%%%%%%%%%%%%%%%%%%%%%%%%%%%%%%%%%%%%%%%%%%%%%%%%%%%%%%%%%%%%%%%%%%%%%%%%%%%%%%%%%%%%%%%%%
%%%%%%%%%%%%%%%%%%%%%%%%%%%%%%%%%%%%%%%%%%%%%%%%%%%%%%%%%%%%%%%%%%%%%%%%%%%%%%%%%%%%%%%%%%
%%%%%%%%%%%%%%%%%%%%%%%%%%%%%%%%%%%%%%%%%%%%%%%%%%%%%%%%%%%%%%%%%%%%%%%%%%%%%%%%%%%%%%%%%%
%%%%%%%%%%%%%%%%%%%%%%%%%%%%%%%%%%%%%%%%%%%%%%%%%%%%%%%%%%%%%%%%%%%%%%%%%%%%%%%%%%%%%%%%%%
%%%%%%%%%%%%%%%%%%%%%%%%%%%%%%%%%%%%%%%%%%%%%%%%%%%%%%%%%%%%%%%%%%%%%%%%%%%%%%%%%%%%%%%%%%
%%%%%%%%%%%%%%%%%%%%%%%%%%%%%%%%%%%%%%%%%%%%%%%%%%%%%%%%%%%%%%%%%%%%%%%%%%%%%%%%%%%%%%%%%%
\section{Results}\label{sec:Results}

We present a few examples of systems created by the methods of
Sections~\ref{sec:simple-construction-m-1} and \ref{sec:simple-construction-m-k}.
%
%adfhide%
{Many more designs are listed in Appendix~\ref{sec:More-results}.
}
\adfhide%
{Many more designs are listed in the appendix of this paper's preprint, which is available at
\begin{center}
\url{https://arxiv.org/abs/blah.blah}.
\end{center}
}

There are two parts to this list.
In Part I the base block of the Steiner system $S(2,k,v)$ is
$$B = \{e, 1, a, a^2, \dots, a^{k-2}\},~~~ e \neq a^{k-1}.$$
In Part II the base block is
$$B = \{e, a, a^2, \dots, a^{k-1}\}.$$

The parameter $m$ is either 1 or $k$.
Let $t = (v - m)/(k(k - 1))$ so that $v = k(k - 1)t + m$.
The blocks of the Steiner system are generated from the base block $B$ by
$$x \mapsto w^s x + d \mathrm{~mod~}v,~~ s = 0, 1, \dots, t - 1,~~ d = 0, 1, \dots, v - 1,$$
together with, only when $m = k$, a further $v/k$ blocks generated from
$$\{0,~ v/k,~ 2v/k,~ \dots,~ (k - 1)v/k\}$$
by $x \mapsto x + d \mathrm{~mod~} v$, $d = 0, 1, \dots, v/k - 1$.

~

\newcommand{\adfDesignE}[9]{$S(2,#1,#4)$, $m=#3$, $e=#6$, $a=#7$, $w=#5$}%
\newcommand{\adfDesignF}[9]{$S(2,#1,#4)$, $m=#3$, $e=#6$, $a=#7$, $w=#5$}%
%%%%%%%%%%%%%%%%%%%%%%%%%%%%%%%%%%%%%%%%%%%%%%%%%%%%%%%%%%
%%%%%%%%%%%%%%%%%%%%%%%%%%%%%%%%%%%%%%%%%%%%%%%%%%%%%%%%%%
%%%%%%%%%%%%%%%%%%%%%%%%%%%%%%%%%%%%%%%%%%%%%%%%%%%%%%%%%%
%%%%%%%%%%%%%%%%%%%%%%%%%%%%%%%%%%%%%%%%%%%%%%%%%%%%%%%%%%
% Designs go after here.

\noindent{\bf Part I}

~

\adfDesignE{10}{1572}{1}{141481}{31756}{49696}{41852}{"Rebecca 22 Apr 2014"}{{49696, 1, 41852, 55124, 60462, 67939, 39371, 67366, 109945, 31577}}

\adfDesignE{10}{1988}{1}{178921}{129189}{0}{130353}{"Rebecca 22 Apr 2014"}{{0, 1, 130353, 135081, 61220, 155139, 109121, 30213, 125058, 14243}}

\adfDesignE{10}{2378}{1}{214021}{8377}{135861}{5392}{"Rebecca 22 Apr 2014"}{{135861, 1, 5392, 180829, 164313, 142777, 20047, 12819, 205286, 199521}}

\adfDesignE{10}{2462}{1}{221581}{45875}{0}{121062}{"Rebecca 22 Apr 2014"}{{0, 1, 121062, 197342, 196946, 117890, 188551, 194447, 42017, 48618}}

\adfDesignE{10}{3990}{1}{359101}{349308}{234242}{267403}{"Audrey 11 May 2014"}{{234242, 1, 267403, 173289, 323629, 333699, 183510, 331981, 75335, 316208}}

\adfDesignE{10}{4178}{1}{376021}{359759}{310258}{112904}{"Audrey 11 May 2014"}{{310258, 1, 112904, 201316, 40277, 212455, 263709, 82135, 316159, 318227}}

\adfDesignE{10}{4639}{1}{417511}{221872}{98148}{319439}{"Audrey 11 May 2014"}{{98148, 1, 319439, 333788, 110730, 366061, 183965, 87363, 296606, 82760}}

\adfDesignE{10}{4755}{1}{427951}{247667}{62358}{53634}{"Audrey 11 May 2014"}{{62358, 1, 53634, 347285, 144366, 8601, 402807, 328256, 206115, 369629}}

~

\adfDesignE{11}{33}{11}{3641}{89}{0}{207}{"Audrey 8 Apr 2014"}{{0, 1, 207, 2798, 267, 654, 661, 2110, 3491, 1719, 2656}}

\adfDesignE{11}{54}{11}{5951}{91}{0}{3018}{"Audrey 9 Apr 2014"}{{0, 1, 3018, 3294, 3122, 1763, 540, 5097, 5362, 1747, 5811}}

\adfDesignE{11}{57}{11}{6281}{309}{0}{1897}{"Audrey 9 Apr 2014"}{{0, 1, 1897, 5877, 6175, 6191, 5138, 4955, 3259, 1819, 2374}}

\adfDesignE{11}{91}{11}{10021}{122}{0}{1461}{"Audrey 9 Apr 2014"}{{0, 1, 1461, 48, 10002, 2304, 9109, 361, 6329, 7307, 3162}}

\adfDesignE{11}{94}{11}{10351}{90}{0}{577}{"Audrey 9 Apr 2014"}{{0, 1, 577, 1697, 6175, 2231, 3763, 7892, 9595, 8881, 592}}

\adfDesignE{11}{123}{11}{13541}{155}{0}{460}{"Audrey 10 Apr 2014"}{{0, 1, 460, 8485, 3292, 11269, 11078, 4464, 8749, 2863, 3503}}

\adfDesignE{11}{148}{11}{16291}{38}{0}{763}{"Audrey 10 Apr 2014"}{{0, 1, 763, 11984, 4541, 11091, 7404, 12566, 8750, 13231, 11124}}

\adfDesignE{11}{162}{11}{17831}{310}{0}{1753}{"Audrey 10 Apr 2014"}{{0, 1, 1753, 6077, 7874, 1928, 9725, 1489, 6891, 8336, 9419}}

~

\adfDesignE{12}{345}{1}{45541}{528}{0}{2734}{"Audrey 10 Apr 2014"}{{0, 1, 2734, 6032, 5646, 43306, 37545, 44157, 41588, 31256, 18988, 41993}}

\adfDesignE{12}{388}{1}{51217}{153}{0}{18}{"Audrey 10 Apr 2014"}{{0, 1, 18, 324, 5832, 2542, 45756, 4136, 23231, 8422, 49162, 14227}}

\adfDesignE{12}{504}{1}{66529}{8}{0}{2821}{"Audrey 10 Apr 2014"}{{0, 1, 2821, 41090, 21372, 15138, 59209, 40799, 65338, 33168, 27154, 26555}}

\adfDesignE{12}{868}{1}{114577}{88}{0}{15726}{"Audrey 10 Apr 2014"}{{0, 1, 15726, 49910, 32210, 104120, 85790, 103942, 36410, 42391, 31880, 70505}}

\adfDesignE{12}{928}{1}{122497}{44}{0}{2040}{"Audrey 10 Apr 2014"}{{0, 1, 2040, 119199, 9415, 97068, 63568, 76894, 67600, 94875, 122237, 82085}}

\adfDesignE{12}{1091}{1}{144013}{31}{0}{8577}{"Audrey 10 Apr 2014"}{{0, 1, 8577, 118299, 78938, 46113, 51503, 53360, 139419, 56824, 39456, 127575}}

\adfDesignE{12}{1145}{1}{151141}{178}{0}{4530}{"Audrey 10 Apr 2014"}{{0, 1, 4530, 116865, 102668, 25183, 118676, 144884, 70298, 146994, 106715, 70032}}

\adfDesignE{12}{1396}{1}{184273}{3}{0}{9079}{"Audrey 10 Apr 2014"}{{0, 1, 9079, 58210, 177899, 176449, 95282, 87816, 116466, 36340, 82190, 81633}}

~

\adfDesignE{13}{13}{13}{2041}{14}{0}{179}{"Audrey 8 Apr 2014"}{{0, 1, 179, 1426, 129, 640, 264, 313, 920, 1400, 1598, 302, 992}}

\adfDesignE{13}{19}{13}{2977}{27}{0}{140}{"Audrey 8 Apr 2014"}{{0, 1, 140, 1738, 2183, 1966, 1356, 2289, 1921, 1010, 1481, 1927, 1850}}

\adfDesignE{13}{38}{13}{5941}{23}{0}{127}{"Audrey 9 Apr 2014"}{{0, 1, 127, 4247, 4679, 133, 5009, 456, 4443, 5807, 805, 1238, 2760}}

\adfDesignE{13}{50}{13}{7813}{29}{0}{3000}{"Audrey 9 Apr 2014"}{{0, 1, 3000, 7237, 6486, 3630, 6491, 3004, 3611, 4182, 6135, 5385, 5529}}

\adfDesignE{13}{56}{13}{8749}{17}{0}{1330}{"Audrey 9 Apr 2014"}{{0, 1, 1330, 1602, 4653, 2947, 8707, 5383, 2708, 5801, 7461, 1764, 1388}}

\adfDesignE{13}{59}{13}{9217}{27}{0}{2649}{"Audrey 9 Apr 2014"}{{0, 1, 2649, 3064, 5576, 5190, 5763, 2835, 7277, 4026, 805, 3318, 5581}}

\adfDesignE{13}{64}{13}{9997}{42}{0}{1557}{"Audrey 9 Apr 2014"}{{0, 1, 1557, 4975, 8397, 8050, 7609, 768, 6133, 1946, 831, 4254, 5464}}

\adfDesignE{13}{89}{13}{13897}{274}{0}{283}{"Audrey 9 Apr 2014"}{{0, 1, 283, 10604, 13077, 4189, 4242, 5344, 11476, 9707, 9372, 11846, 3241}}

~

\adfDesignE{14}{6}{1}{1093}{241}{0}{11}{"Audrey 8 Apr 2014"}{{0, 1, 11, 121, 238, 432, 380, 901, 74, 814, 210, 124, 271, 795}}

\adfDesignE{14}{90}{1}{16381}{880}{0}{1510}{"Audrey 10 Apr 2014"}{{0, 1, 1510, 3141, 8801, 4519, 9194, 8233, 15032, 10635, 5470, 3676, 13982, 14092}}

\adfDesignE{14}{921}{1}{167623}{189}{0}{10555}{"Audrey 11 Apr 2014"}{{0, 1, 10555, 106353, 152307, 95815, 57766, 75279, 36825, 137761, 105453, 39695, 90848, 97080}}

\adfDesignE{14}{6276}{1}{1142233}{80554}{0}{25474}{"Audrey 19 May 2014"}{{0, 1, 25474, 136332, 533048, 1141081, 352010, 573690, 450058, 184871, 1119428, 462027, 106966, 626179}}

\adfDesignE{14}{6281}{1}{1143143}{374310}{0}{559862}{"Audrey 19 May 2014"}{{0, 1, 559862, 221016, 88900, 428723, 1123659, 662441, 1087080, 877188, 850112, 102780, 227169, 629927}}

\adfDesignE{14}{6845}{1}{1245791}{1039908}{0}{307647}{"Audrey 19 May 2014"}{{0, 1, 307647, 196966, 724762, 427625, 772984, 902031, 956852, 954481, 311179, 476418, 11505, 186504}}

\adfDesignE{14}{6915}{1}{1258531}{1115660}{0}{124353}{"Audrey 19 May 2014"}{{0, 1, 124353, 98212, 172012, 215360, 380931, 64334, 902866, 545188, 1215456, 1060992, 699322, 813628}}

\adfDesignE{14}{7356}{1}{1338793}{500348}{0}{1295712}{"Audrey 19 May 2014"}{{0, 1, 1295712, 405463, 819561, 480348, 1190006, 1090656, 1066785, 1260312, 587636, 629114, 1004051, 880699}}

~

\adfDesignE{16}{1530}{1}{367201}{10}{0}{17744}{"Audrey 18 Apr 2014"}{{0, 1, 17744, 158279, 149328, 320817, 226946, 203658, 82511, 44797, 255004, 140254, 145799, 126411, 173076, 158581}}

\adfDesignE{16}{3723}{1}{893521}{813}{0}{312001}{"Audrey 18 Apr 2014"}{{0, 1, 312001, 872177, 62669, 764147, 887322, 378766, 70348, 196504, 501089, 6198, 202754, 844517, 634348, 521806}}

\adfDesignE{16}{5948}{1}{1427521}{741373}{0}{809452}{"Audrey 19 May 2014"}{{0, 1, 809452, 386598, 1363323, 847467, 61702, 130077, 1421407, 225779, 314604, 1365818, 487992, 1047037, 1265940, 265450}}

\adfDesignE{16}{6193}{1}{1486321}{583704}{0}{797239}{"Audrey 19 May 2014"}{{0, 1, 797239, 519175, 344708, 1339917, 323253, 1159240, 1399523, 1382396, 414749, 1163067, 1016163, 1053944, 1245538, 1217976}}

\adfDesignE{16}{8027}{1}{1926481}{46472}{0}{945050}{"Audrey 19 May 2014"}{{0, 1, 945050, 984419, 530316, 1103650, 1839657, 1584033, 987271, 665997, 1783821, 86823, 1323879, 1881272, 730168, 965491}}

\adfDesignE{16}{10751}{1}{2580241}{1129999}{0}{2485809}{"Audrey 19 May 2014"}{{0, 1, 2485809, 89728, 316948, 762064, 2274083, 2092092, 880903, 1537744, 1141831, 306157, 582581, 1609610, 725549, 676346}}

\adfDesignE{16}{13357}{1}{3205681}{65431}{0}{2218134}{"Audrey 19 May 2014"}{{0, 1, 2218134, 774984, 774735, 3042182, 2412026, 1843871, 881269, 2954823, 2243327, 793654, 2349357, 1549428, 2681804, 1125853}}

\adfDesignE{16}{14415}{1}{3459601}{927424}{0}{1440955}{"Audrey 19 May 2014"}{{0, 1, 1440955, 2579855, 2979793, 2225205, 3250758, 27920, 3223172, 658780, 2795313, 208842, 1990726, 1116175, 3280629, 1898684}}

~

\adfDesignE{17}{175}{17}{47617}{256}{0}{5421}{"Rebecca 11 Apr 2014"}{{0, 1, 5421, 7552, 36389, 35155, 12021, 25785, 24590, 22407, 44997, 34463, 22432, 37671, 32795, 27434, 11823}}

\adfDesignE{17}{237}{17}{64481}{256}{0}{4360}{"Rebecca 11 Apr 2014"}{{0, 1, 4360, 52186, 41992, 23561, 7727, 30638, 41529, 3792, 25984, 61604, 30075, 37127, 26410, 49015, 15366}}

\adfDesignE{17}{427}{17}{116161}{256}{0}{2253}{"Rebecca 11 Apr 2014"}{{0, 1, 2253, 81086, 81666, 110635, 95310, 67902, 115330, 102494, 107075, 89739, 61827, 19192, 27684, 109756, 89660}}

\adfDesignE{17}{435}{17}{118337}{749}{0}{1481}{"Rebecca 11 Apr 2014"}{{0, 1, 1481, 63295, 16991, 76227, 117026, 70138, 92829, 90492, 61168, 62003, 115268, 69954, 56999, 41238, 11586}}

\adfDesignE{17}{501}{17}{136289}{256}{0}{2134}{"Rebecca 11 Apr 2014"}{{0, 1, 2134, 56419, 54959, 73966, 20782, 54863, 5391, 56118, 94070, 127972, 105381, 6204, 19303, 33324, 106847}}

\adfDesignE{17}{510}{17}{138737}{14}{0}{8475}{"Rebecca 11 Apr 2014"}{{0, 1, 8475, 98596, 126886, 8363, 120555, 44357, 87042, 16321, 138423, 113590, 117944, 114052, 10021, 20931, 84339}}

\adfDesignE{17}{643}{17}{174913}{137}{0}{14686}{"Rebecca 11 Apr 2014"}{{0, 1, 14686, 10867, 72106, 25414, 140575, 161224, 113296, 92600, 149938, 9711, 61651, 56898, 44627, 168024, 102773}}

\adfDesignE{17}{853}{17}{232033}{171}{0}{9352}{"Audrey 18 Apr 2014"}{{0, 1, 9352, 215496, 111987, 137495, 158387, 168585, 174718, 218383, 195383, 193974, 10854, 108187, 100944, 118044, 166507}}

~

\adfDesignE{18}{4390}{1}{1343341}{819}{0}{37484}{"Audrey 10 Apr 2014"}{{0, 1, 37484, 1258911, 137276, 657554, 133468, 312628, 584409, 125269, 606401, 1005364, 319103, 148588, 180806, 176759, 276544, 756140}}

\adfDesignE{18}{8002}{1}{2448613}{1590147}{0}{2443765}{"Audrey 18 May 2014"}{{0, 1, 2443765, 1465587, 709150, 2342065, 2335974, 33173, 785754, 706436, 807859, 1280368, 9891, 1020692, 332057, 1375018, 1485935, 6566}}

\adfDesignE{18}{13656}{1}{4178737}{3287159}{0}{1283274}{"Audrey 18 May 2014"}{{0, 1, 1283274, 2052220, 3506244, 4161895, 3705793, 3023224, 1749836, 3677585, 650126, 771737, 1894149, 3531981, 3047585, 2456253, 3399537, 2875930}}

\adfDesignE{18}{30705}{1}{9395731}{5279918}{0}{582650}{"Audrey 18 May 2014"}{{0, 1, 582650, 3865739, 5837, 9069159, 5168412, 5883376, 1132629, 7724334, 9264638, 5956580, 6220220, 5260101, 4352760, 2319556, 7356360, 1003496}}

\adfDesignE{18}{33682}{1}{10306693}{4834062}{0}{8544986}{"Audrey 18 May 2014"}{{0, 1, 8544986, 8624224, 8334950, 9930283, 1810943, 3404905, 9126393, 8686122, 7007004, 7290193, 6721235, 7811598, 7283139, 3930655, 2637588, 1922711}}

\adfDesignE{18}{33701}{1}{10312507}{3752302}{0}{8844691}{"Audrey 18 May 2014"}{{0, 1, 8844691, 5159923, 7290349, 5730850, 6448751, 2152302, 1644483, 8804334, 6770027, 8731689, 4940967, 5129790, 6647354, 3439130, 3386448, 7691474}}

\adfDesignE{18}{41623}{1}{12736639}{12236571}{0}{8370551}{"Audrey 18 May 2014"}{{0, 1, 8370551, 618668, 1995058, 1355635, 9067171, 5789308, 8843375, 2085193, 685577, 11701809, 1736097, 3689534, 12365204, 12598966, 12228897, 8762068}}

\adfDesignE{18}{42860}{1}{13115161}{2666363}{0}{5544998}{"Audrey 18 May 2014"}{{0, 1, 5544998, 2983858, 2822090, 11022221, 8666721, 10450850, 172394, 12760566, 9486471, 4486165, 895911, 11558954, 7938008, 736410, 12830152, 2565518}}

~

\adfDesignE{19}{185}{19}{63289}{172}{0}{499}{"Rebecca 12 Apr 2014"}{{0, 1, 499, 59134, 15192, 49417, 39662, 45170, 8946, 33824, 43302, 26149, 10817, 18118, 53844, 33620, 4795, 51012, 12810}}

\adfDesignE{19}{914}{19}{312607}{280}{0}{100006}{"Rebecca 12 Apr 2014"}{{0, 1, 100006, 276892, 133292, 124665, 168223, 51026, 222095, 105220, 279700, 229054, 183792, 261580, 303113, 243102, 212222, 271495, 273199}}

\adfDesignE{19}{1289}{19}{440857}{39}{0}{29549}{"Rebecca 12 Apr 2014"}{{0, 1, 29549, 246541, 318941, 187520, 337704, 17301, 273986, 114366, 232029, 16857, 379940, 423555, 137322, 79950, 330744, 236480, 164070}}

\adfDesignE{19}{1406}{19}{480871}{287}{0}{127551}{"Rebecca 12 Apr 2014"}{{0, 1, 127551, 429929, 306781, 307448, 269798, 433325, 205206, 421976, 50617, 74921, 379959, 47545, 148114, 109837, 123473, 98502, 311985}}

\adfDesignE{19}{2342}{19}{800983}{643}{0}{17185}{"Rebecca 12 Apr 2014"}{{0, 1, 17185, 562481, 774124, 595276, 464167, 521181, 704562, 238942, 379412, 193600, 533601, 279801, 79236, 800543, 448430, 12107, 604198}}

\adfDesignE{19}{2356}{19}{805771}{313}{0}{223900}{"Rebecca 12 Apr 2014"}{{0, 1, 223900, 167235, 543901, 39586, 630171, 755945, 663866, 632572, 84817, 115372, 384082, 49825, 723776, 5964, 177053, 650813, 597289}}

\adfDesignE{19}{2867}{19}{980533}{58}{0}{486816}{"Rebecca 12 Apr 2014"}{{0, 1, 486816, 874954, 32330, 226097, 846836, 961255, 823028, 745987, 361248, 751752, 570042, 19277, 651222, 340125, 586955, 2684, 544188}}

\adfDesignE{19}{3002}{19}{1026703}{659}{0}{343297}{"Audrey 18 Apr 2014"}{{0, 1, 343297, 672948, 534120, 851464, 639302, 372008, 724315, 20391, 108073, 197073, 1002199, 654694, 985794, 333364, 384510, 1004869, 411905}}

~

\adfDesignE{20}{347}{1}{131861}{474}{0}{19405}{"Rebecca 11 Apr 2014"}{{0, 1, 19405, 90870, 87058, 89219, 91626, 120667, 87358, 108835, 57399, 129589, 85275, 37686, 127585, 96650, 34247, 115456, 105290, 98116}}

\adfDesignE{20}{41144}{1}{15634721}{8790850}{0}{11386168}{"Audrey 19 May 2014"}{{0, 1, 11386168, 11007635, 9848882, 3704253, 7346039, 9663354, 4786022, 5731500, 11761649, 1931714, 4072199, 6843086, 3551851, 14619572, 14002663, 11226342, 8671662, 557897}}

\adfDesignE{20}{46325}{1}{17603501}{3472089}{0}{1776065}{"Audrey 19 May 2014"}{{0, 1, 1776065, 333033, 10621545, 8922789, 7091041, 6599231, 5791201, 17408276, 3869572, 10965269, 13277670, 2958930, 8442916, 12555712, 15240001, 1220960, 17051715, 14305582}}

~

\adfDesignE{22}{9911}{1}{4578883}{1500}{0}{185929}{"Felicity 9 May 2014"}{{0, 1, 185929, 3605274, 3991644, 3283987, 3528639, 227742, 2911217, 748397, 1030326, 754783, 2442223, 1410823, 2439146, 1667665, 3644557, 4522166, 4411339, 3431556, 4218304, 1911995}}

\adfDesignE{22}{54179}{1}{25030699}{3493122}{0}{7826006}{"Audrey 19 May 2014"}{{0, 1, 7826006, 4063886, 15100615, 14343592, 13888366, 6838981, 14203136, 2505516, 9726961, 18502461, 4817977, 22857232, 20691949, 10827562, 14699177, 11239852, 16458322, 5869130, 13315800, 22530468}}

~

\adfDesignE{23}{358}{23}{181171}{669}{0}{23225}{"Rebecca 11 Apr 2014"}{{0, 1, 23225, 54558, 180747, 117005, 57296, 179776, 30734, 164581, 47967, 13096, 149662, 134315, 63597, 134333, 119305, 29351, 111673, 142560, 55975, 117450, 65674}}

\adfDesignE{23}{2434}{23}{1231627}{174}{0}{1484}{"Audrey 12 Apr 2014"}{{0, 1, 1484, 970629, 641473, 1129888, 509445, 1029029, 1093183, 230813, 134186, 840077, 267744, 748202, 635841, 161762, 1119170, 615084, 149049, 727483, 679520, 936794, 927040}}

\adfDesignE{23}{3628}{23}{1835791}{749}{0}{188603}{"Audrey 12 Apr 2014"}{{0, 1, 188603, 805193, 1512277, 474525, 191534, 1027495, 705734, 1358938, 1330522, 661603, 1596339, 928835, 711330, 1201501, 323645, 367185, 648662, 651155, 875938, 1702524, 1094371}}

\adfDesignE{23}{4758}{23}{2407571}{1903}{0}{12524}{"Audrey 14 Apr 2014"}{{0, 1, 12524, 358461, 1653220, 2224251, 923054, 1579925, 1562222, 1346382, 1868455, 1347871, 1256123, 615538, 2363141, 2115152, 2067506, 19039, 94907, 1682765, 1479897, 748470, 1164377}}

~

\adfDesignE{28}{156975}{1}{118673101}{70589532}{0}{86913349}{"Rebecca 19 Jun 2014"}{{0, 1, 86913349, 6744137, 12036068, 19146903, 15756912, 67701702, 114222787, 18940118, 86664892, 8942008, 48175579, 17271926, 82278028, 59545807, 115130915, 53323599, 26240919, 32008208, 22817744, 63753385, 89365107, 17248170, 43370382, 54513585, 69799816, 73492064}}

~

\adfDesignE{29}{4275}{29}{3471329}{1306}{0}{293718}{"Mononoko 16 Apr 2014"}{{0, 1, 293718, 795216, 881323, 3425384, 1653642, 3009934, 667550, 374993, 396133, 2858401, 2078294, 2623771, 3318591, 1556112, 2100302, 3154917, 1991501, 3397573, 1099481, 3094817, 1247666, 702316, 2596392, 1211433, 1512736, 1565764, 990645}}

~

\adfDesignE{30}{76249}{1}{66336631}{4061041}{0}{35121401}{"Rebecca 20 May 2014"}{{0, 1, 35121401, 4262075, 44246324, 58292371, 40121810, 52253278, 60320367, 62139720, 59295571, 56762194, 38342173, 9982875, 44342132, 43871904, 12702974, 27995646, 54837245, 43755643, 48485887, 49527641, 7013576, 63643665, 39615504, 7167801, 10445216, 28455800, 46265624, 12119678}}

~

\adfDesignE{31}{16217}{31}{15081841}{32}{0}{4421147}{"Mononoko 18 Apr 2014"}{{0, 1, 4421147, 7322538, 621490, 13446445, 14029075, 3872090, 13168632, 11966763, 6664140, 12083030, 5775586, 11684226, 7453185, 5866981, 9730219, 9201161, 8732325, 3757109, 13798171, 2998810, 4885108, 11156918, 7007166, 10850256, 6404388, 2008590, 8262725, 14439810, 9701371}}

~

\noindent{\bf Part II}

~

\adfDesignF{10}{9145}{1}{823051}{475317}{1}{259876}{"???"}{{1, 259876, 85571, 657278, 534345, 745653, 660741, 67139, 779666, 255389}}

~

\adfDesignF{11}{201}{1}{22111}{2}{1}{527}{"Audrey 11 Apr 2014"}{{1, 527, 12397, 10474, 14159, 10386, 12005, 2889, 18955, 17224, 11538}}

\adfDesignF{11}{320}{1}{35201}{56}{1}{3717}{"Audrey 11 Apr 2014"}{{1, 3717, 17297, 15923, 12910, 7507, 24327, 27291, 26566, 7017, 33449}}

\adfDesignF{11}{414}{1}{45541}{7}{1}{2734}{"Audrey 11 Apr 2014"}{{1, 2734, 6032, 5646, 43306, 37545, 44157, 41588, 31256, 18988, 41993}}

\adfDesignF{11}{453}{1}{49831}{10}{1}{12651}{"Audrey 11 Apr 2014"}{{1, 12651, 40460, 45259, 13419, 39383, 23995, 40124, 30158, 22722, 30814}}

\adfDesignF{11}{959}{1}{105491}{43}{1}{8082}{"Audrey 11 Apr 2014"}{{1, 8082, 19795, 58834, 48451, 103881, 68864, 93823, 8178, 57230, 60316}}

\adfDesignF{11}{1106}{1}{121661}{58}{1}{3427}{"Audrey 11 Apr 2014"}{{1, 3427, 64873, 45124, 8817, 43931, 56880, 26838, 119771, 92664, 24318}}

\adfDesignF{11}{1517}{1}{166871}{108}{1}{3471}{"Audrey 11 Apr 2014"}{{1, 3471, 33129, 16640, 20074, 91647, 50611, 122489, 138882, 135974, 54566}}

\adfDesignF{11}{1662}{1}{182821}{72}{1}{4313}{"Audrey 11 Apr 2014"}{{1, 4313, 137048, 27731, 38869, 177961, 63235, 146444, 149238, 133574, 35691}}

~

\adfDesignF{13}{135}{1}{21061}{1295}{1}{240}{"Audrey 12 Apr 2014"}{{1, 240, 15478, 7984, 20670, 11465, 13670, 16345, 5454, 3178, 4524, 11649, 15708}}

\adfDesignF{13}{342}{1}{53353}{162}{1}{3751}{"Audrey 12 Apr 2014"}{{1, 3751, 38162, 52916, 14756, 22695, 30910, 7341, 5943, 43992, 46516, 17206, 35929}}

\adfDesignF{13}{460}{1}{71761}{104}{1}{6645}{"Audrey 12 Apr 2014"}{{1, 6645, 23010, 50520, 7442, 8761, 18674, 13961, 55633, 40374, 42612, 59595, 31577}}

\adfDesignF{13}{627}{1}{97813}{136}{1}{6674}{"Audrey 12 Apr 2014"}{{1, 6674, 37361, 21977, 52811, 40375, 85748, 76102, 59652, 18538, 86980, 82178, 18481}}

\adfDesignF{13}{890}{1}{138841}{136}{1}{936}{"Audrey 12 Apr 2014"}{{1, 936, 43050, 30910, 52832, 23356, 63179, 128119, 99601, 64225, 135288, 6576, 46132}}

\adfDesignF{13}{1220}{1}{190321}{356}{1}{10240}{"Audrey 12 Apr 2014"}{{1, 10240, 181050, 35139, 116670, 55883, 136994, 152790, 130980, 43113, 122721, 163798, 182868}}

\adfDesignF{13}{1322}{1}{206233}{113}{1}{8850}{"Audrey 12 Apr 2014"}{{1, 8850, 160193, 62408, 18826, 180069, 48259, 189840, 109982, 127173, 67569, 116183, 148045}}

\adfDesignF{13}{1386}{1}{216217}{141}{1}{26055}{"Audrey 12 Apr 2014"}{{1, 26055, 157862, 214636, 104492, 150813, 121174, 204153, 51998, 208385, 46088, 169839, 58023}}

~

\adfDesignF{15}{363}{1}{76231}{68254}{1}{5734}{"Mononoko 12 Apr 2014"}{{1, 5734, 23195, 53266, 45858, 29053, 25167, 2295, 47798, 23287, 47177, 45330, 50741, 51398, 7086}}

\adfDesignF{15}{693}{1}{145531}{5495}{1}{128690}{"Mononoko 12 Apr 2014"}{{1, 128690, 124893, 36530, 103338, 89971, 67161, 8531, 114057, 29732, 55559, 95211, 12107, 140475, 12461}}

\adfDesignF{15}{1071}{1}{224911}{59838}{1}{189220}{"Mononoko 12 Apr 2014"}{{1, 189220, 176488, 49169, 89754, 222270, 22222, 135695, 143229, 15880, 2640, 13469, 137639, 32513, 119277}}

\adfDesignF{15}{1159}{1}{243391}{305}{1}{16662}{"Mononoko 12 Apr 2014"}{{1, 16662, 156504, 221865, 92122, 113118, 195603, 131696, 148887, 114122, 130272, 31126, 198582, 116030, 37147}}

\adfDesignF{15}{1661}{1}{348811}{447}{1}{87983}{"Mononoko 12 Apr 2014"}{{1, 87983, 194577, 173122, 262989, 183502, 339331, 277072, 271419, 308006, 165308, 270308, 225973, 253081, 126627}}

\adfDesignF{15}{2235}{1}{469351}{82776}{1}{228292}{"Mononoko 12 Apr 2014"}{{1, 228292, 32873, 189777, 188127, 395180, 130095, 55162, 355974, 237513, 74170, 110964, 381316, 392951, 43711}}

\adfDesignF{15}{2313}{1}{485731}{411390}{1}{341497}{"Mononoko 12 Apr 2014"}{{1, 341497, 73757, 213224, 393580, 249981, 52576, 471319, 257459, 279275, 235749, 91658, 426386, 15048, 298607}}

\adfDesignF{15}{2356}{1}{494761}{43}{1}{58574}{"Mononoko 12 Apr 2014"}{{1, 58574, 240702, 169492, 444943, 60846, 230121, 333531, 111948, 174619, 433914, 206066, 415289, 213321, 369960}}

~

\adfDesignF{17}{3716}{1}{1010753}{16}{1}{74004}{"Mononoko 16 Apr 2014"}{{1, 74004, 332262, 128817, 561725, 658269, 287488, 932808, 125591, 362529, 179237, 143329, 77334, 141850, 797495, 963063, 298716}}

\adfDesignF{17}{4143}{1}{1126897}{4}{1}{27314}{"Mononoko 16 Apr 2014"}{{1, 27314, 48782, 439294, 803957, 566556, 360980, 585867, 433838, 529177, 359656, 482835, 79599, 382773, 838253, 876093, 1073304}}

\adfDesignF{17}{10604}{1}{2884289}{2566626}{1}{1819206}{"Audrey 18 May 2014"}{{1, 1819206, 2280322, 661458, 107259, 1181215, 234487, 2467089, 1948549, 661071, 2725342, 1289035, 951962, 1337902, 536006, 2211850, 2661558}}

\adfDesignF{17}{10614}{1}{2887009}{2777952}{1}{2628521}{"Audrey 18 May 2014"}{{1, 2628521, 1996857, 1879885, 2092964, 2086114, 105388, 283580, 2018479, 2101773, 28414, 2759873, 306921, 2498881, 2867714, 1661417, 666299}}

\adfDesignF{17}{12164}{1}{3308609}{2496275}{1}{42435}{"Audrey 18 May 2014"}{{1, 42435, 845929, 1898074, 3301303, 979136, 124338, 2360284, 439892, 2953651, 1454047, 235204, 2116996, 2682301, 676117, 2076256, 974299}}

\adfDesignF{17}{16373}{1}{4453457}{1817841}{1}{1831015}{"Audrey 18 May 2014"}{{1, 1831015, 59141, 2251160, 1694136, 4210002, 3475047, 4307783, 4074448, 2133261, 3731269, 1205448, 2185579, 460512, 191671, 2250637, 1566546}}

\adfDesignF{17}{16418}{1}{4465697}{2326577}{1}{4316015}{"Audrey 18 May 2014"}{{1, 4316015, 299275, 3791754, 1506593, 3379377, 2056773, 3254994, 2462186, 117164, 3915968, 4068953, 596702, 2856933, 3699414, 1810258, 2058913}}

\adfDesignF{17}{16440}{1}{4471681}{687211}{1}{2776620}{"Audrey 18 May 2014"}{{1, 2776620, 2355662, 1235249, 1977932, 3454475, 4157819, 1460288, 3893258, 1492743, 1308165, 2172896, 2694295, 503520, 3694388, 601628, 1994509}}

~

\adfDesignF{19}{423}{1}{144667}{785}{1}{2816}{"Rebecca 12 Apr 2014"}{{1, 2816, 117838, 110377, 76916, 28957, 95391, 119104, 58758, 108147, 17917, 110156, 33248, 26819, 6130, 46707, 24609, 3451, 25327}}

\adfDesignF{19}{18121}{1}{6197383}{541}{1}{60662}{"Audrey 22 Apr 2014"}{{1, 60662, 4830125, 5169276, 3435678, 3305929, 3148501, 3418368, 604436, 2578804, 1066562, 5302907, 3582436, 300354, 5965711, 1977780, 952863, 5781448, 4294606}}

\adfDesignF{19}{22570}{1}{7718941}{5874365}{1}{2032211}{"Audrey 18 May 2014"}{{1, 2032211, 1107409, 626985, 222965, 2369974, 7681918, 5689015, 3411185, 5708873, 4433675, 5373886, 4449972, 2340722, 5010387, 437383, 3049781, 6041838, 1763348}}

\adfDesignF{19}{25841}{1}{8837623}{531419}{1}{656241}{"Audrey 18 May 2014"}{{1, 656241, 3718914, 3088447, 6952268, 2294199, 5342171, 7997079, 804441, 593999, 5460098, 5464875, 2971967, 4601915, 8120447, 7317449, 8053552, 4549195, 1551349}}

\adfDesignF{19}{39149}{1}{13388959}{5706975}{1}{12410065}{"Audrey 18 May 2014"}{{1, 12410065, 12445524, 6024906, 11771782, 1294873, 2727027, 6888323, 9476818, 7743038, 10139518, 3141747, 6594482, 273516, 9230778, 839347, 8504735, 7254151, 3740000}}

\adfDesignF{19}{43288}{1}{14804497}{6898585}{1}{6562832}{"Audrey 18 May 2014"}{{1, 6562832, 11129130, 2757786, 3710530, 8409091, 2618459, 8969174, 13520367, 4483575, 4991104, 13186693, 9666053, 4476976, 16461, 2362943, 5763052, 5327544, 10484205}}

\adfDesignF{19}{47649}{1}{16295959}{2054397}{1}{583157}{"Audrey 18 May 2014"}{{1, 583157, 8014237, 7732681, 14463273, 10505354, 488036, 8981676, 10457024, 8519296, 3260978, 5212041, 13696511, 12696721, 14980793, 4755314, 11303268, 11105207, 16204022}}

~

\adfDesignF{21}{4837}{1}{2031541}{1489191}{1}{1915986}{"Rebecca 22 Apr 2014"}{{1, 1915986, 1670573, 57428, 951907, 174060, 784141, 1409968, 735960, 511542, 512667, 531916, 751116, 348304, 621572, 1411136, 1880967, 1461446, 547718, 1106365, 899096}}

\adfDesignF{21}{6761}{1}{2839621}{179}{1}{34187}{"Rebecca 22 Apr 2014"}{{1, 34187, 1666738, 937020, 138239, 847349, 1346442, 556244, 2211412, 2312161, 2157951, 517257, 1145092, 245098, 2283376, 594022, 1700343, 2584271, 2184125, 847180, 1248081}}

\adfDesignF{21}{8139}{1}{3418381}{1048576}{1}{2357122}{"Rebecca 22 Apr 2014"}{{1, 2357122, 3003487, 2008460, 1037600, 733130, 3215216, 3139922, 1500812, 3364070, 715508, 2758863, 2485031, 2335566, 278458, 2978828, 469205, 694413, 3040299, 200220, 1285980}}

\adfDesignF{21}{70348}{1}{29546161}{17584789}{1}{7958054}{"Rebecca 18 May 2014"}{{1, 7958054, 22856110, 1762205, 21330513, 1924316, 5848603, 21998482, 26378810, 27384951, 1255887, 13959394, 5084850, 4707291, 7547678, 24755458, 8851261, 26528069, 23802415, 22536995, 18990267}}

\adfDesignF{21}{96923}{1}{40707661}{5077810}{1}{31600789}{"Rebecca 18 May 2014"}{{1, 31600789, 15808610, 11078307, 3382188, 13546087, 10159925, 38920198, 17265056, 39549703, 14997665, 40241388, 31701485, 11705045, 1787462, 32549477, 26057009, 14631689, 37791746, 36167750, 18842691}}

~

\adfDesignF{23}{26097}{1}{13205083}{9862926}{1}{6390971}{"Audrey 18 May 2014"}{{1, 6390971, 146371, 4736521, 5825015, 9252708, 676504, 1605105, 8854650, 9192302, 11488866, 7324089, 9100485, 5977830, 10752476, 12155350, 2847241, 3876845, 1092931, 8852819, 7027939, 7585225, 10493669}}

~

\adfDesignF{25}{13289}{1}{7973401}{221}{1}{305397}{"Rebecca 17 Apr 2014"}{{1, 305397, 2456112, 7484191, 2321568, 5085576, 4793085, 4930561, 1758867, 626631, 1630106, 2217246, 7170138, 3518156, 2556380, 3197346, 5296298, 5340248, 331114, 2550776, 6033773, 3334776, 4023144, 2854474, 7891447}}

\adfDesignF{25}{84407}{1}{50644201}{17091760}{1}{41647216}{"Rebecca 18 May 2014"}{{1, 41647216, 49103704, 29913875, 45036551, 37966447, 34704877, 17483706, 45102786, 2333938, 47951097, 9918608, 46975170, 22794529, 35840003, 12885452, 29877292, 45076507, 26972485, 15310960, 2485405, 35784610, 43937717, 27148928, 2766950}}

~

\adfDesignF{27}{138999}{1}{97577299}{84129266}{1}{65249917}{"Felicity 1 Jun 2014"}{{1, 65249917, 21850293, 79456466, 20949842, 43068068, 96606162, 89324971, 19333894, 38254368, 90440342, 51203549, 22477934, 19912767, 61633401, 80139957, 47462139, 23337084, 59322930, 39464339, 24523457, 93220198, 56714690, 90453129, 18408479, 58367488, 54906321}}

% 134 designs.
% Designs go before here.
%%%%%%%%%%%%%%%%%%%%%%%%%%%%%%%%%%%%%%%%%%%%%%%%%%%%%%%%%%
%%%%%%%%%%%%%%%%%%%%%%%%%%%%%%%%%%%%%%%%%%%%%%%%%%%%%%%%%%
%%%%%%%%%%%%%%%%%%%%%%%%%%%%%%%%%%%%%%%%%%%%%%%%%%%%%%%%%%
%%%%%%%%%%%%%%%%%%%%%%%%%%%%%%%%%%%%%%%%%%%%%%%%%%%%%%%%%%

%%%%%%%%%%%%%%%%%%%%%%%%%%%%%%%%%%%%%%%%%%%%%%%%%%%%%%%%%%%%%%%%%%%%%%%%%%%%%%%%%%%%%%%%%%
%%%%%%%%%%%%%%%%%%%%%%%%%%%%%%%%%%%%%%%%%%%%%%%%%%%%%%%%%%%%%%%%%%%%%%%%%%%%%%%%%%%%%%%%%%
%%%%%%%%%%%%%%%%%%%%%%%%%%%%%%%%%%%%%%%%%%%%%%%%%%%%%%%%%%%%%%%%%%%%%%%%%%%%%%%%%%%%%%%%%%
%%%%%%%%%%%%%%%%%%%%%%%%%%%%%%%%%%%%%%%%%%%%%%%%%%%%%%%%%%%%%%%%%%%%%%%%%%%%%%%%%%%%%%%%%%
%%%%%%%%%%%%%%%%%%%%%%%%%%%%%%%%%%%%%%%%%%%%%%%%%%%%%%%%%%%%%%%%%%%%%%%%%%%%%%%%%%%%%%%%%%
%%%%%%%%%%%%%%%%%%%%%%%%%%%%%%%%%%%%%%%%%%%%%%%%%%%%%%%%%%%%%%%%%%%%%%%%%%%%%%%%%%%%%%%%%%

\section*{ORCID}

\noindent A. D. Forbes     \url{https://orcid.org/0000-0003-3805-7056} \\
C. G. Rutherford           \url{https://orcid.org/0000-0003-1924-207X}

%%%%%%%%%%%%%%%%%%%%%%%%%%%%%%%%%%%%%%%%%%%%%%%%%%%%%%%%%%%%%%%%%%%%%%%%%%%%%%%%%%%%%%%%%%
%%%%%%%%%%%%%%%%%%%%%%%%%%%%%%%%%%%%%%%%%%%%%%%%%%%%%%%%%%%%%%%%%%%%%%%%%%%%%%%%%%%%%%%%%%
%%%%%%%%%%%%%%%%%%%%%%%%%%%%%%%%%%%%%%%%%%%%%%%%%%%%%%%%%%%%%%%%%%%%%%%%%%%%%%%%%%%%%%%%%%
%%%%%%%%%%%%%%%%%%%%%%%%%%%%%%%%%%%%%%%%%%%%%%%%%%%%%%%%%%%%%%%%%%%%%%%%%%%%%%%%%%%%%%%%%%
%%%%%%%%%%%%%%%%%%%%%%%%%%%%%%%%%%%%%%%%%%%%%%%%%%%%%%%%%%%%%%%%%%%%%%%%%%%%%%%%%%%%%%%%%%
%%%%%%%%%%%%%%%%%%%%%%%%%%%%%%%%%%%%%%%%%%%%%%%%%%%%%%%%%%%%%%%%%%%%%%%%%%%%%%%%%%%%%%%%%%

%%%%%%%%%%%%%%%%%%%%%%%%%%%%%%%%%%%%%%%%%%%%%%%%%%%%%%%%%%%%%%%%%%%%%%%%%%%%%%%%%%%%%%%%%%
%%%%%%%%%%%%%%%%%%%%%%%%%%%%%%%%%%%%%%%%%%%%%%%%%%%%%%%%%%%%%%%%%%%%%%%%%%%%%%%%%%%%%%%%%%
%%%%%%%%%%%%%%%%%%%%%%%%%%%%%%%%%%%%%%%%%%%%%%%%%%%%%%%%%%%%%%%%%%%%%%%%%%%%%%%%%%%%%%%%%%
%%%%%%%%%%%%%%%%%%%%%%%%%%%%%%%%%%%%%%%%%%%%%%%%%%%%%%%%%%%%%%%%%%%%%%%%%%%%%%%%%%%%%%%%%%
%%%%%%%%%%%%%%%%%%%%%%%%%%%%%%%%%%%%%%%%%%%%%%%%%%%%%%%%%%%%%%%%%%%%%%%%%%%%%%%%%%%%%%%%%%
%%%%%%%%%%%%%%%%%%%%%%%%%%%%%%%%%%%%%%%%%%%%%%%%%%%%%%%%%%%%%%%%%%%%%%%%%%%%%%%%%%%%%%%%%%

\appendix

%%%%%%%%%%%%%%%%%%%%%%%%%%%%%%%%%%%%%%%%%%%%%%%%%%%%%%%%%%%%%%%%%%%%%%%%%%%%%%%%%%%%%%%%%%
%%%%%%%%%%%%%%%%%%%%%%%%%%%%%%%%%%%%%%%%%%%%%%%%%%%%%%%%%%%%%%%%%%%%%%%%%%%%%%%%%%%%%%%%%%
%%%%%%%%%%%%%%%%%%%%%%%%%%%%%%%%%%%%%%%%%%%%%%%%%%%%%%%%%%%%%%%%%%%%%%%%%%%%%%%%%%%%%%%%%%
%%%%%%%%%%%%%%%%%%%%%%%%%%%%%%%%%%%%%%%%%%%%%%%%%%%%%%%%%%%%%%%%%%%%%%%%%%%%%%%%%%%%%%%%%%
%%%%%%%%%%%%%%%%%%%%%%%%%%%%%%%%%%%%%%%%%%%%%%%%%%%%%%%%%%%%%%%%%%%%%%%%%%%%%%%%%%%%%%%%%%
%%%%%%%%%%%%%%%%%%%%%%%%%%%%%%%%%%%%%%%%%%%%%%%%%%%%%%%%%%%%%%%%%%%%%%%%%%%%%%%%%%%%%%%%%%

\section{More results}\label{sec:More-results}

The lists of Steiner systems in this Appendix continue those of Section~\ref{sec:Results}.
As before, there are two parts.
In Part I %, page~\pageref{app:PART-I},
the base block of the $S(2,k,v)$ is
$$B = \{e, 1, a, a^2, \dots, a^{k-2}\},~~~ e \neq a^{k-1}.$$
In Part II, which starts on page~\pageref{app:PART-II}, the base block is
$$B = \{e, a, a^2, \dots, a^{k-1}\}.$$

With $t = (v - m)/(k(k - 1))$, $m \in \{1, k\}$,
$tv$ blocks of the Steiner system are generated from $B$ by
$$x \mapsto w^s x + d \mathrm{~mod~}v,~~ s = 0, 1, \dots, t - 1,~~ d = 0, 1, \dots, v - 1.$$
If $m = k$, we generate a further $v/k$ blocks from
$$\{0,~ v/k,~ 2v/k,~ \dots,~ (k - 1)v/k\}$$
by $x \mapsto x + d \mathrm{~mod~} v$, $d = 0, 1, \dots, v/k - 1$.

~

%%%%%%%%%%%%%%%%%%%%%%%%%%%%%%%%%%%%%%%%%%%%%%%%%%%%%%%%%%
%%%%%%%%%%%%%%%%%%%%%%%%%%%%%%%%%%%%%%%%%%%%%%%%%%%%%%%%%%
%%%%%%%%%%%%%%%%%%%%%%%%%%%%%%%%%%%%%%%%%%%%%%%%%%%%%%%%%%
%%%%%%%%%%%%%%%%%%%%%%%%%%%%%%%%%%%%%%%%%%%%%%%%%%%%%%%%%%
% Designs go after here.

\noindent{\bf Part I}\label{app:PART-I}

~

\adfDesignE{10}{5610}{1}{504901}{154346}{292428}{3862}{"???"}{{292428, 1, 3862, 272915, 269343, 106606, 218057, 466167, 364889, 22627}}

\adfDesignE{10}{5815}{1}{523351}{154763}{161939}{253062}{"???"}{{161939, 1, 253062, 7378, 298419, 6380, 521076, 493401, 485633, 406373}}

\adfDesignE{10}{7273}{1}{654571}{65789}{94187}{237847}{"???"}{{94187, 1, 237847, 551305, 613902, 250595, 651989, 521215, 222415, 275998}}

\adfDesignE{10}{7599}{1}{683911}{648983}{148830}{277603}{"???"}{{148830, 1, 277603, 334129, 467323, 557001, 394524, 222343, 116079, 44050}}

\adfDesignE{10}{7784}{1}{700561}{497742}{496852}{237363}{"???"}{{496852, 1, 237363, 677027, 172572, 405966, 543230, 247074, 162869, 16784}}

\adfDesignE{10}{7885}{1}{709651}{433355}{136568}{71211}{"???"}{{136568, 1, 71211, 550126, 158433, 140765, 196040, 659619, 328919, 619654}}

\adfDesignE{10}{7982}{1}{718381}{697820}{105537}{365907}{"???"}{{105537, 1, 365907, 392155, 683502, 286593, 718376, 325608, 194368, 174395}}

\adfDesignE{10}{8084}{1}{727561}{558229}{0}{215252}{"???"}{{0, 1, 215252, 156341, 106438, 96486, 575727, 195513, 253353, 405201}}

\adfDesignE{10}{8218}{1}{739621}{355889}{704091}{358464}{"???"}{{704091, 1, 358464, 603724, 215336, 397860, 328094, 533543, 322046, 372422}}

\adfDesignE{10}{8257}{1}{743131}{26397}{255803}{295476}{"???"}{{255803, 1, 295476, 64172, 298407, 356713, 374396, 322443, 314882, 72632}}

\adfDesignE{10}{8362}{1}{752581}{277303}{18607}{141962}{"???"}{{18607, 1, 141962, 595426, 225635, 243348, 443133, 553737, 268801, 660538}}

\adfDesignE{10}{8467}{1}{762031}{279479}{0}{281618}{"???"}{{0, 1, 281618, 321599, 682832, 21357, 566974, 204440, 255777, 426911}}

\adfDesignE{10}{9145}{1}{823051}{475317}{391633}{259876}{"???"}{{391633, 1, 259876, 85571, 657278, 534345, 745653, 660741, 67139, 779666}}

\adfDesignE{10}{9170}{1}{825301}{231462}{590110}{102553}{"???"}{{590110, 1, 102553, 307166, 706230, 65333, 301631, 57162, 21583, 769418}}

\adfDesignE{10}{9226}{1}{830341}{475086}{483829}{133549}{"???"}{{483829, 1, 133549, 441062, 659180, 77000, 330056, 827100, 605693, 365260}}

\adfDesignE{10}{9515}{1}{856351}{795284}{399388}{725870}{"???"}{{399388, 1, 725870, 177130, 813960, 48962, 624089, 372483, 246682, 350995}}

\adfDesignE{10}{9737}{1}{876331}{603520}{625845}{168418}{"???"}{{625845, 1, 168418, 417247, 675018, 513556, 833701, 120543, 527028, 63707}}

\adfDesignE{10}{9844}{1}{885961}{335370}{196722}{203890}{"???"}{{196722, 1, 203890, 70058, 662378, 785385, 12666, 780386, 507667, 515039}}

\adfDesignE{10}{9914}{1}{892261}{507982}{561837}{234031}{"???"}{{561837, 1, 234031, 851998, 378268, 763193, 690786, 136820, 443174, 35754}}

\adfDesignE{10}{9977}{1}{897931}{677932}{649375}{337048}{"???"}{{649375, 1, 337048, 511770, 305722, 18820, 258776, 303494, 644123, 407586}}

\adfDesignE{10}{9981}{1}{898291}{734767}{115125}{457844}{"???"}{{115125, 1, 457844, 432031, 21255, 287817, 488303, 832943, 188625, 26051}}

\adfDesignE{10}{10013}{1}{901171}{726694}{375158}{742712}{"???"}{{375158, 1, 742712, 828279, 84721, 840619, 241731, 721997, 340511, 581076}}

\adfDesignE{10}{10069}{1}{906211}{636017}{732045}{58378}{"???"}{{732045, 1, 58378, 637524, 196513, 310865, 801695, 83615, 424024, 519607}}

\adfDesignE{10}{10215}{1}{919351}{54974}{417588}{439833}{"???"}{{417588, 1, 439833, 472416, 707667, 744402, 415832, 28765, 607134, 119109}}

\adfDesignE{10}{10281}{1}{925291}{744376}{37092}{174775}{"???"}{{37092, 1, 174775, 594133, 663182, 131644, 719385, 121713, 874776, 367597}}

\adfDesignE{10}{10304}{1}{927361}{781839}{873507}{129451}{"???"}{{873507, 1, 129451, 148131, 662684, 504540, 99671, 137028, 777781, 17100}}

\adfDesignE{10}{10364}{1}{932761}{263470}{603529}{471588}{"???"}{{603529, 1, 471588, 767558, 177400, 377110, 338420, 336621, 762319, 411757}}

\adfDesignE{10}{10404}{1}{936361}{57174}{265591}{620365}{"???"}{{265591, 1, 620365, 871337, 754481, 450661, 325690, 572992, 10177, 508743}}

\adfDesignE{10}{10486}{1}{943741}{784577}{606778}{166461}{"???"}{{606778, 1, 166461, 85020, 174184, 288081, 883649, 680188, 391934, 910244}}

\adfDesignE{10}{10729}{1}{965611}{897950}{684166}{277298}{"???"}{{684166, 1, 277298, 645652, 210342, 649072, 339500, 426555, 328945, 313106}}

\adfDesignE{10}{10755}{1}{967951}{693019}{229194}{574675}{"???"}{{229194, 1, 574675, 25739, 300594, 417637, 156623, 462888, 771433, 733324}}

\adfDesignE{10}{10811}{1}{972991}{212441}{74011}{267288}{"???"}{{74011, 1, 267288, 37778, 878457, 772478, 544509, 727812, 458271, 502058}}

\adfDesignE{10}{11030}{1}{992701}{626820}{454867}{275837}{"???"}{{454867, 1, 275837, 482424, 805240, 122532, 368337, 11121, 137187, 481100}}

\adfDesignE{10}{11270}{1}{1014301}{211468}{475053}{464146}{"???"}{{475053, 1, 464146, 62722, 712411, 590006, 840789, 598648, 430066, 1005438}}

\adfDesignE{10}{11284}{1}{1015561}{946147}{0}{191792}{"???"}{{0, 1, 191792, 551844, 543711, 601071, 217878, 984470, 369120, 521291}}

\adfDesignE{10}{11373}{1}{1023571}{559123}{775241}{915699}{"???"}{{775241, 1, 915699, 413256, 836531, 770899, 585396, 351962, 474239, 969172}}

\adfDesignE{10}{11399}{1}{1025911}{860311}{438495}{114259}{"???"}{{438495, 1, 114259, 401606, 152746, 833193, 387742, 72554, 586606, 197502}}

~

\adfDesignE{11}{172}{11}{18931}{157}{0}{6015}{"Audrey 10 Apr 2014"}{{0, 1, 6015, 3084, 16811, 7694, 12046, 7753, 7242, 399, 14679}}

\adfDesignE{11}{174}{11}{19151}{37}{0}{1785}{"Audrey 10 Apr 2014"}{{0, 1, 1785, 7159, 5098, 3205, 13927, 1697, 3287, 7089, 14205}}

\adfDesignE{11}{216}{11}{23771}{141}{0}{2704}{"Audrey 10 Apr 2014"}{{0, 1, 2704, 13919, 7483, 4911, 15126, 14584, 22818, 14127, 23182}}

\adfDesignE{11}{228}{11}{25091}{161}{0}{1648}{"Audrey 10 Apr 2014"}{{0, 1, 1648, 6076, 1939, 8915, 13685, 21162, 23577, 14028, 9333}}

\adfDesignE{11}{237}{11}{26081}{210}{0}{1400}{"Audrey 10 Apr 2014"}{{0, 1, 1400, 3925, 17990, 17835, 9483, 971, 3188, 3349, 20101}}

\adfDesignE{11}{301}{11}{33121}{210}{0}{1763}{"Audrey 10 Apr 2014"}{{0, 1, 1763, 27916, 31223, 32168, 9032, 25336, 20260, 13942, 3964}}

\adfDesignE{11}{336}{11}{36971}{381}{0}{13776}{"Audrey 10 Apr 2014"}{{0, 1, 13776, 6033, 36771, 17625, 13443, 3029, 24216, 10283, 22707}}

\adfDesignE{11}{339}{11}{37301}{78}{0}{17175}{"Audrey 10 Apr 2014"}{{0, 1, 17175, 4317, 27388, 23290, 27127, 16735, 19420, 30259, 20793}}

\adfDesignE{11}{382}{11}{42031}{60}{0}{2896}{"Audrey 10 Apr 2014"}{{0, 1, 2896, 22647, 17352, 24347, 22925, 23851, 15563, 13216, 25326}}

\adfDesignE{11}{423}{11}{46541}{111}{0}{6890}{"Audrey 10 Apr 2014"}{{0, 1, 6890, 280, 21019, 31859, 21154, 31189, 12413, 29753, 31606}}

\adfDesignE{11}{493}{11}{54241}{177}{0}{13478}{"Audrey 10 Apr 2014"}{{0, 1, 13478, 3375, 34292, 15, 39447, 50625, 26211, 225, 49295}}

\adfDesignE{11}{558}{11}{61391}{126}{0}{8353}{"Audrey 10 Apr 2014"}{{0, 1, 8353, 32433, 55757, 26095, 33485, 2809, 12215, 53, 12972}}

\adfDesignE{11}{574}{11}{63151}{115}{0}{27355}{"Audrey 10 Apr 2014"}{{0, 1, 27355, 19826, 62593, 18452, 51668, 58760, 60548, 29263, 50440}}

\adfDesignE{11}{579}{11}{63701}{111}{0}{2242}{"Audrey 10 Apr 2014"}{{0, 1, 2242, 57886, 21475, 52695, 40536, 44086, 40561, 36435, 22588}}

\adfDesignE{11}{613}{11}{67441}{133}{0}{13282}{"Audrey 10 Apr 2014"}{{0, 1, 13282, 53309, 54520, 20623, 36785, 35766, 57049, 25183, 40687}}

\adfDesignE{11}{631}{11}{69421}{67}{0}{1654}{"Audrey 10 Apr 2014"}{{0, 1, 1654, 28297, 13484, 18395, 18932, 4657, 66368, 18071, 38404}}

\adfDesignE{11}{732}{11}{80531}{523}{0}{45851}{"Audrey 10 Apr 2014"}{{0, 1, 45851, 52446, 45886, 46611, 29283, 42001, 50048, 20003, 70525}}

\adfDesignE{11}{741}{11}{81521}{67}{0}{43107}{"Audrey 10 Apr 2014"}{{0, 1, 43107, 23775, 68434, 65532, 22232, 75469, 65157, 79786, 45633}}

\adfDesignE{11}{748}{11}{82291}{59}{0}{16146}{"Audrey 10 Apr 2014"}{{0, 1, 16146, 77719, 77806, 1270, 14961, 36221, 64420, 49371, 73540}}

\adfDesignE{11}{784}{11}{86251}{96}{0}{3172}{"Audrey 10 Apr 2014"}{{0, 1, 3172, 56468, 59420, 21805, 78409, 51715, 76829, 42513, 40923}}

\adfDesignE{11}{816}{11}{89771}{103}{0}{4426}{"Audrey 10 Apr 2014"}{{0, 1, 4426, 19398, 34472, 52143, 73448, 20057, 78534, 87943, 78433}}

\adfDesignE{11}{822}{11}{90431}{103}{0}{10686}{"Audrey 10 Apr 2014"}{{0, 1, 10686, 66674, 62946, 15178, 49325, 55082, 81304, 43927, 67032}}

\adfDesignE{11}{823}{11}{90541}{56}{0}{28550}{"Audrey 10 Apr 2014"}{{0, 1, 28550, 52418, 72252, 89538, 65847, 29067, 54585, 10058, 50389}}

\adfDesignE{11}{843}{11}{92741}{122}{0}{4211}{"Audrey 10 Apr 2014"}{{0, 1, 4211, 18990, 24148, 43092, 59016, 63237, 31596, 60162, 66511}}

\adfDesignE{11}{852}{11}{93731}{157}{0}{48013}{"Audrey 10 Apr 2014"}{{0, 1, 48013, 27955, 69226, 46678, 42604, 54239, 48734, 58589, 72616}}

\adfDesignE{11}{864}{11}{95051}{63}{0}{16988}{"Audrey 10 Apr 2014"}{{0, 1, 16988, 17308, 35561, 61163, 34563, 26217, 60461, 85413, 42529}}

\adfDesignE{11}{874}{11}{96151}{72}{0}{15574}{"Audrey 10 Apr 2014"}{{0, 1, 15574, 56654, 47820, 59185, 43704, 89318, 22015, 83295, 63189}}

\adfDesignE{11}{876}{11}{96371}{619}{0}{12488}{"Audrey 10 Apr 2014"}{{0, 1, 12488, 21866, 43565, 25425, 61326, 75122, 48222, 70328, 27141}}

\adfDesignE{11}{951}{11}{104621}{56}{0}{2407}{"Audrey 10 Apr 2014"}{{0, 1, 2407, 39494, 66190, 86168, 47554, 7104, 46105, 76475, 46986}}

\adfDesignE{11}{981}{11}{107921}{45}{0}{37273}{"Audrey 10 Apr 2014"}{{0, 1, 37273, 9496, 71449, 59981, 88298, 80459, 39559, 65905, 87184}}

\adfDesignE{11}{987}{11}{108581}{144}{0}{33356}{"Audrey 10 Apr 2014"}{{0, 1, 33356, 101810, 103585, 25259, 59225, 94967, 85739, 103706, 43838}}

\adfDesignE{11}{1131}{11}{124421}{320}{0}{10217}{"Audrey 10 Apr 2014"}{{0, 1, 10217, 122291, 11465, 57744, 90487, 57649, 115240, 11157, 21433}}

\adfDesignE{11}{1149}{11}{126401}{23}{0}{46854}{"Audrey 10 Apr 2014"}{{0, 1, 46854, 91149, 111060, 55273, 57454, 114020, 81216, 118760, 82619}}

\adfDesignE{11}{1183}{11}{130141}{78}{0}{22301}{"Audrey 10 Apr 2014"}{{0, 1, 22301, 65840, 47078, 39031, 47323, 36854, 40639, 118556, 102941}}

\adfDesignE{11}{1194}{11}{131351}{236}{0}{49872}{"Audrey 10 Apr 2014"}{{0, 1, 49872, 85199, 102380, 19288, 47763, 117302, 105857, 40912, 88181}}

\adfDesignE{11}{1197}{11}{131681}{45}{0}{7704}{"Audrey 10 Apr 2014"}{{0, 1, 7704, 95166, 90737, 75100, 95767, 112006, 120312, 112770, 80523}}

\adfDesignE{11}{1207}{11}{132781}{94150}{0}{87685}{"Rebecca 22 Apr 2014"}{{0, 1, 87685, 108201, 3892, 22850, 69741, 10630, 101711, 27608, 77069}}

\adfDesignE{11}{1221}{11}{134321}{23}{0}{3667}{"Audrey 10 Apr 2014"}{{0, 1, 3667, 14789, 99900, 39933, 24421, 94021, 107321, 119898, 33333}}

\adfDesignE{11}{1224}{11}{134651}{333}{0}{2094}{"Audrey 10 Apr 2014"}{{0, 1, 2094, 76004, 129545, 80116, 122409, 83593, 132093, 29588, 17812}}

\adfDesignE{11}{1312}{11}{144331}{249}{0}{41255}{"Audrey 10 Apr 2014"}{{0, 1, 41255, 23873, 110202, 101341, 131209, 37471, 81095, 125976, 69232}}

\adfDesignE{11}{1341}{11}{147521}{111}{0}{43822}{"Audrey 10 Apr 2014"}{{0, 1, 43822, 86827, 71162, 14745, 13410, 76877, 114338, 116592, 52310}}

\adfDesignE{11}{1344}{11}{147851}{201}{0}{34093}{"Audrey 10 Apr 2014"}{{0, 1, 34093, 75938, 83224, 95142, 120968, 6230, 85354, 118391, 119914}}

\adfDesignE{11}{1368}{11}{150491}{217}{0}{16934}{"Audrey 10 Apr 2014"}{{0, 1, 16934, 75001, 73385, 97403, 41042, 37790, 48128, 90787, 121493}}

\adfDesignE{11}{1372}{11}{150931}{138}{0}{115206}{"Audrey 10 Apr 2014"}{{0, 1, 115206, 3089, 126967, 33268, 82325, 131772, 134121, 133732, 145105}}

\adfDesignE{11}{1378}{11}{151591}{63}{0}{3591}{"Audrey 10 Apr 2014"}{{0, 1, 3591, 10046, 148119, 114101, 137809, 79095, 100202, 99939, 65052}}

\adfDesignE{11}{1401}{11}{154121}{67}{0}{7078}{"Audrey 10 Apr 2014"}{{0, 1, 7078, 8759, 39560, 121944, 42032, 48966, 117340, 128572, 102232}}

\adfDesignE{11}{1407}{11}{154781}{199}{0}{6273}{"Audrey 10 Apr 2014"}{{0, 1, 6273, 36155, 46150, 58480, 14070, 35940, 90484, 24205, 152585}}

\adfDesignE{11}{1434}{11}{157751}{401}{0}{18445}{"Audrey 10 Apr 2014"}{{0, 1, 18445, 106869, 99960, 126263, 43022, 53260, 65223, 29109, 88852}}

\adfDesignE{11}{1443}{11}{158741}{925}{0}{22729}{"Audrey 10 Apr 2014"}{{0, 1, 22729, 64227, 33247, 63903, 129878, 49426, 152238, 139925, 138131}}

\adfDesignE{11}{1446}{11}{159071}{173}{0}{64870}{"Audrey 10 Apr 2014"}{{0, 1, 64870, 52666, 75553, 145600, 72304, 152045, 120866, 126901, 143620}}

\adfDesignE{11}{1456}{11}{160171}{148}{0}{44467}{"Audrey 10 Apr 2014"}{{0, 1, 44467, 3094, 154180, 122747, 43682, 13777, 127955, 20552, 110229}}

\adfDesignE{11}{1474}{11}{162151}{250}{0}{5042}{"Audrey 10 Apr 2014"}{{0, 1, 5042, 126208, 60212, 42232, 29481, 112886, 21202, 42975, 46214}}

\adfDesignE{11}{1509}{11}{166001}{100}{0}{38932}{"Audrey 10 Apr 2014"}{{0, 1, 38932, 111494, 90260, 93152, 135818, 36523, 114871, 90832, 118122}}

\adfDesignE{11}{1527}{11}{167981}{6414}{0}{69837}{"Rebecca 22 Apr 2014"}{{0, 1, 69837, 46215, 98002, 115791, 58708, 78329, 129089, 152166, 2920}}

\adfDesignE{11}{1546}{11}{170071}{124}{0}{26601}{"Audrey 10 Apr 2014"}{{0, 1, 26601, 117841, 109840, 34060, 61843, 158931, 98613, 29309, 43245}}

\adfDesignE{11}{1554}{11}{170951}{301}{0}{12708}{"Audrey 10 Apr 2014"}{{0, 1, 12708, 115520, 71923, 93438, 155409, 111620, 86513, 21323, 15349}}

\adfDesignE{11}{1558}{11}{171391}{18}{0}{13275}{"Audrey 10 Apr 2014"}{{0, 1, 13275, 35677, 58842, 98763, 109066, 111373, 57809, 96968, 103790}}

\adfDesignE{11}{1567}{11}{172381}{100}{0}{1879}{"Audrey 10 Apr 2014"}{{0, 1, 1879, 83021, 164035, 4537, 78354, 13792, 58018, 71030, 42476}}

\adfDesignE{11}{1573}{11}{173041}{23}{0}{15849}{"Audrey 10 Apr 2014"}{{0, 1, 15849, 108310, 38470, 87587, 31461, 94268, 17538, 55916, 69723}}

\adfDesignE{11}{1599}{11}{175901}{188}{0}{11334}{"Audrey 10 Apr 2014"}{{0, 1, 11334, 51826, 62445, 101907, 47972, 4657, 12138, 17510, 42012}}

\adfDesignE{11}{1606}{11}{176671}{59}{0}{100610}{"Audrey 10 Apr 2014"}{{0, 1, 100610, 7155, 106896, 136106, 32121, 27878, 153455, 5531, 136931}}

\adfDesignE{11}{1641}{11}{180521}{133}{0}{51716}{"Audrey 10 Apr 2014"}{{0, 1, 51716, 126041, 84088, 124639, 147698, 145216, 136535, 145666, 121526}}

\adfDesignE{11}{1642}{11}{180631}{228}{0}{14991}{"Audrey 10 Apr 2014"}{{0, 1, 14991, 25117, 93943, 100237, 164209, 17851, 89830, 37425, 178920}}

\adfDesignE{11}{1648}{11}{181291}{222}{0}{128967}{"Audrey 10 Apr 2014"}{{0, 1, 128967, 125585, 145337, 389, 131847, 85286, 154592, 151321, 164421}}

\adfDesignE{11}{1687}{11}{185581}{45}{0}{92822}{"Audrey 10 Apr 2014"}{{0, 1, 92822, 140178, 147244, 184242, 50612, 109630, 112887, 122692, 153178}}

\adfDesignE{11}{1701}{11}{187121}{155}{0}{66047}{"Audrey 10 Apr 2014"}{{0, 1, 66047, 41457, 156007, 163585, 119076, 104063, 94631, 65136, 125602}}

\adfDesignE{11}{1743}{11}{191741}{474}{0}{19254}{"Audrey 10 Apr 2014"}{{0, 1, 19254, 81163, 23252, 170514, 87154, 137625, 162871, 185920, 90951}}

\adfDesignE{11}{1758}{11}{193391}{951}{0}{59569}{"Audrey 10 Apr 2014"}{{0, 1, 59569, 127693, 89505, 126866, 140647, 116241, 188765, 15981, 101687}}

\adfDesignE{11}{1812}{11}{199331}{184}{0}{2770}{"Audrey 10 Apr 2014"}{{0, 1, 2770, 98322, 65794, 60846, 108725, 178440, 137251, 61053, 84122}}

\adfDesignE{11}{1813}{11}{199441}{66793}{0}{6957}{"Rebecca 22 Apr 2014"}{{0, 1, 6957, 135127, 113106, 83697, 111750, 23732, 165817, 22125, 154614}}

\adfDesignE{11}{1831}{11}{201421}{12}{0}{71305}{"Audrey 18 Apr 2014"}{{0, 1, 71305, 134143, 187588, 197993, 91554, 1939, 85589, 68766, 168227}}

\adfDesignE{11}{1834}{11}{201751}{228}{0}{49911}{"Audrey 18 Apr 2014"}{{0, 1, 49911, 88324, 79814, 23059, 110045, 188522, 58404, 103596, 105328}}

\adfDesignE{11}{1854}{11}{203951}{124}{0}{24897}{"Audrey 18 Apr 2014"}{{0, 1, 24897, 53520, 75557, 102556, 74163, 67808, 113349, 184017, 119936}}

\adfDesignE{11}{1873}{11}{206041}{100}{0}{6450}{"Audrey 18 Apr 2014"}{{0, 1, 6450, 188259, 70937, 132630, 187309, 124667, 130168, 172566, 17218}}

\adfDesignE{11}{1903}{11}{209341}{12}{0}{5465}{"Audrey 18 Apr 2014"}{{0, 1, 5465, 139803, 138086, 175026, 38061, 127752, 12445, 185441, 15284}}

\adfDesignE{11}{1923}{11}{211541}{199}{0}{90369}{"Audrey 18 Apr 2014"}{{0, 1, 90369, 15856, 123671, 102028, 153847, 101941, 118761, 203256, 147975}}

\adfDesignE{11}{1942}{11}{213631}{227}{0}{46480}{"Audrey 18 Apr 2014"}{{0, 1, 46480, 153728, 175014, 9502, 77683, 128309, 79324, 135722, 48761}}

\adfDesignE{11}{1966}{11}{216271}{58}{0}{84197}{"Audrey 18 Apr 2014"}{{0, 1, 84197, 203971, 98719, 116571, 117965, 53430, 208910, 58569, 139022}}

\adfDesignE{11}{1975}{11}{217261}{35949}{0}{790}{"Rebecca 22 Apr 2014"}{{0, 1, 790, 189578, 73791, 68942, 148930, 116499, 132807, 197728, 211722}}

\adfDesignE{11}{1986}{11}{218471}{332}{0}{1967}{"Audrey 18 Apr 2014"}{{0, 1, 1967, 155082, 60778, 46689, 79443, 57616, 162694, 177554, 132060}}

\adfDesignE{11}{2034}{11}{223751}{104}{0}{18676}{"Audrey 18 Apr 2014"}{{0, 1, 18676, 188918, 126800, 159967, 20340, 164393, 116197, 157974, 165489}}

\adfDesignE{11}{2041}{11}{224521}{199}{0}{15282}{"Audrey 18 Apr 2014"}{{0, 1, 15282, 37684, 215044, 213052, 81643, 5129, 23549, 193176, 113524}}

\adfDesignE{11}{2055}{11}{226061}{59049}{0}{57820}{"Rebecca 22 Apr 2014"}{{0, 1, 57820, 162332, 209581, 199576, 200575, 91139, 175070, 214003, 204625}}

\adfDesignE{11}{2101}{11}{231121}{243}{0}{41116}{"Audrey 18 Apr 2014"}{{0, 1, 41116, 106462, 90973, 214725, 42021, 105961, 61626, 35093, 226506}}

\adfDesignE{11}{2122}{11}{233431}{138}{0}{59227}{"Audrey 18 Apr 2014"}{{0, 1, 59227, 69892, 61561, 114558, 21220, 4436, 121097, 44544, 203757}}

\adfDesignE{11}{2149}{11}{236401}{100}{0}{15854}{"Audrey 18 Apr 2014"}{{0, 1, 15854, 55053, 17770, 171989, 64472, 177565, 52402, 68194, 85903}}

\adfDesignE{11}{2152}{11}{236731}{212}{0}{845}{"Audrey 18 Apr 2014"}{{0, 1, 845, 3832, 160537, 6902, 150646, 171323, 125294, 54673, 36140}}

\adfDesignE{11}{2196}{11}{241571}{218}{0}{45901}{"Audrey 18 Apr 2014"}{{0, 1, 45901, 161110, 138658, 111292, 153726, 129787, 212227, 80952, 174201}}

\adfDesignE{11}{2209}{11}{243001}{1024}{0}{80826}{"Rebecca 22 Apr 2014"}{{0, 1, 80826, 3392, 56664, 84617, 233498, 36683, 84957, 12224, 217959}}

\adfDesignE{11}{2248}{11}{247291}{370}{0}{51049}{"Audrey 18 Apr 2014"}{{0, 1, 51049, 47843, 91391, 27153, 67442, 61356, 221929, 110938, 62771}}

\adfDesignE{11}{2253}{11}{247841}{45}{0}{32503}{"Audrey 18 Apr 2014"}{{0, 1, 32503, 146667, 143707, 97135, 180247, 102683, 78643, 149196, 60582}}

\adfDesignE{11}{2274}{11}{250151}{691}{0}{104461}{"Audrey 18 Apr 2014"}{{0, 1, 104461, 13599, 207761, 71212, 136445, 77467, 145588, 87872, 156198}}

\adfDesignE{11}{2301}{11}{253121}{34}{0}{15304}{"Audrey 18 Apr 2014"}{{0, 1, 15304, 75491, 70020, 124887, 207098, 99751, 16553, 206112, 197267}}

\adfDesignE{11}{2383}{11}{262141}{112795}{0}{162902}{"Rebecca 22 Apr 2014"}{{0, 1, 162902, 3892, 157646, 205627, 148292, 245952, 180123, 168393, 73682}}

\adfDesignE{11}{2439}{11}{268301}{45}{0}{29131}{"Audrey 18 Apr 2014"}{{0, 1, 29131, 247399, 147108, 99576, 146345, 141606, 264812, 48020, 217507}}

\adfDesignE{11}{2577}{11}{283481}{111}{0}{118121}{"Audrey 18 Apr 2014"}{{0, 1, 118121, 202783, 203648, 41672, 257709, 88047, 132240, 234459, 138725}}

\adfDesignE{11}{2584}{11}{284251}{107}{0}{109432}{"Audrey 18 Apr 2014"}{{0, 1, 109432, 152245, 239479, 144983, 25840, 278183, 261211, 277341, 216791}}

\adfDesignE{11}{2602}{11}{286231}{61}{0}{82635}{"Audrey 18 Apr 2014"}{{0, 1, 82635, 216489, 131015, 23181, 104083, 229617, 147805, 103174, 106924}}

\adfDesignE{11}{2617}{11}{287881}{45}{0}{41110}{"Audrey 18 Apr 2014"}{{0, 1, 41110, 170630, 90854, 39846, 26170, 37403, 64909, 40001, 64838}}

\adfDesignE{11}{2656}{11}{292171}{522}{0}{28539}{"Audrey 18 Apr 2014"}{{0, 1, 28539, 193944, 80392, 180596, 132804, 51144, 204471, 158657, 138136}}

\adfDesignE{11}{2664}{11}{293051}{212}{0}{64464}{"Audrey 18 Apr 2014"}{{0, 1, 64464, 144116, 284073, 17933, 239768, 15459, 175576, 115542, 115272}}

\adfDesignE{11}{2668}{11}{293491}{382}{0}{59702}{"Audrey 18 Apr 2014"}{{0, 1, 59702, 174100, 134435, 233484, 106723, 180427, 146072, 292461, 140050}}

\adfDesignE{11}{2695}{11}{296461}{244377}{0}{110984}{"Rebecca 22 Apr 2014"}{{0, 1, 110984, 86628, 91722, 93091, 242155, 251487, 119441, 82790, 149587}}

\adfDesignE{11}{2736}{11}{300971}{647}{0}{52996}{"Audrey 18 Apr 2014"}{{0, 1, 52996, 215615, 67554, 41739, 164165, 220614, 140078, 123973, 177149}}

\adfDesignE{11}{2755}{11}{303061}{3983}{0}{74319}{"Rebecca 22 Apr 2014"}{{0, 1, 74319, 27036, 297115, 265225, 169335, 199840, 101594, 205793, 53541}}

\adfDesignE{11}{2779}{11}{305701}{221}{0}{5069}{"Audrey 18 Apr 2014"}{{0, 1, 5069, 15877, 81150, 181505, 194536, 217259, 150869, 196760, 179778}}

\adfDesignE{11}{2815}{11}{309661}{64417}{0}{271199}{"Rebecca 22 Apr 2014"}{{0, 1, 271199, 74847, 153003, 305919, 242100, 165731, 25963, 67819, 129886}}

\adfDesignE{11}{2821}{11}{310321}{1024}{0}{256672}{"Rebecca 22 Apr 2014"}{{0, 1, 256672, 298247, 118099, 240927, 307990, 306777, 215604, 276279, 80173}}

\adfDesignE{11}{2866}{11}{315271}{79}{0}{34439}{"Audrey 18 Apr 2014"}{{0, 1, 34439, 310490, 233874, 158449, 114643, 51544, 148086, 110058, 99500}}

\adfDesignE{11}{2923}{11}{321541}{38886}{0}{44256}{"Rebecca 22 Apr 2014"}{{0, 1, 44256, 87305, 133424, 33620, 116513, 167852, 217930, 87785, 154598}}

\adfDesignE{11}{2941}{11}{323521}{1024}{0}{89182}{"Rebecca 22 Apr 2014"}{{0, 1, 89182, 312381, 45511, 191057, 288388, 69679, 244731, 226340, 8127}}

\adfDesignE{11}{2953}{11}{324841}{243}{0}{4651}{"Audrey 18 Apr 2014"}{{0, 1, 4651, 192295, 76772, 66313, 147654, 24880, 73484, 41352, 22280}}

\adfDesignE{11}{3021}{11}{332321}{78}{0}{152871}{"Audrey 18 Apr 2014"}{{0, 1, 152871, 65279, 331021, 327979, 211476, 28395, 327464, 242988, 306452}}

\adfDesignE{11}{3067}{11}{337381}{42483}{0}{84022}{"Audrey 12 May 2014"}{{0, 1, 84022, 336440, 219833, 210719, 289081, 93449, 241246, 120932, 44927}}

\adfDesignE{11}{3078}{11}{338591}{201}{0}{22654}{"Audrey 18 Apr 2014"}{{0, 1, 22654, 238351, 92877, 31084, 246247, 192813, 161802, 214933, 153602}}

\adfDesignE{11}{3112}{11}{342331}{63}{0}{64766}{"Audrey 18 Apr 2014"}{{0, 1, 64766, 53013, 202359, 182990, 31120, 215323, 71471, 233335, 314946}}

\adfDesignE{11}{3123}{11}{343541}{2960}{0}{244308}{"Audrey 12 May 2014"}{{0, 1, 244308, 272606, 281306, 272739, 137875, 113991, 105604, 316473, 235306}}

\adfDesignE{11}{3177}{11}{349481}{287}{0}{44388}{"Audrey 18 Apr 2014"}{{0, 1, 44388, 270147, 242445, 80227, 254167, 19154, 269960, 329433, 237483}}

\adfDesignE{11}{3277}{11}{360481}{1024}{0}{216324}{"Audrey 12 May 2014"}{{0, 1, 216324, 231961, 136645, 150980, 297958, 21668, 334470, 304846, 193407}}

\adfDesignE{11}{3283}{11}{361141}{137919}{0}{48734}{"Audrey 12 May 2014"}{{0, 1, 48734, 139540, 57330, 133444, 193909, 345800, 294717, 160708, 239946}}

\adfDesignE{11}{3319}{11}{365101}{252176}{0}{45170}{"Audrey 12 May 2014"}{{0, 1, 45170, 144512, 331362, 306045, 233487, 300304, 134227, 166384, 326296}}

\adfDesignE{11}{3339}{11}{367301}{287}{0}{62054}{"Audrey 18 Apr 2014"}{{0, 1, 62054, 282533, 291450, 104361, 133563, 338638, 184941, 9069, 62594}}

\adfDesignE{11}{3381}{11}{371921}{45}{0}{119139}{"Audrey 18 Apr 2014"}{{0, 1, 119139, 108277, 305539, 214967, 101432, 49916, 297455, 370681, 292198}}

\adfDesignE{11}{3391}{11}{373021}{1816}{0}{16482}{"Audrey 12 May 2014"}{{0, 1, 16482, 97036, 206325, 189214, 169588, 103063, 319753, 128258, 38349}}

\adfDesignE{11}{3426}{11}{376871}{146}{0}{8099}{"Audrey 18 Apr 2014"}{{0, 1, 8099, 18247, 49021, 175916, 171304, 128945, 16014, 53762, 132433}}

\adfDesignE{11}{3447}{11}{379181}{78}{0}{48976}{"Audrey 18 Apr 2014"}{{0, 1, 48976, 328751, 125354, 17933, 103412, 364676, 188314, 47001, 292306}}

\adfDesignE{11}{3459}{11}{380501}{59049}{0}{176949}{"Audrey 12 May 2014"}{{0, 1, 176949, 282313, 168250, 129507, 80917, 310104, 162985, 339971, 320379}}

\adfDesignE{11}{3472}{11}{381931}{367}{0}{63662}{"Audrey 18 Apr 2014"}{{0, 1, 63662, 180403, 150616, 138037, 243046, 5780, 166807, 57710, 139731}}

\adfDesignE{11}{3525}{11}{387761}{59049}{0}{11169}{"Audrey 12 May 2014"}{{0, 1, 11169, 275280, 45351, 109453, 257885, 28857, 74442, 83114, 432}}

\adfDesignE{11}{3546}{11}{390071}{248}{0}{89391}{"Audrey 18 Apr 2014"}{{0, 1, 89391, 146446, 171626, 327336, 106382, 52453, 172703, 253906, 240040}}

\adfDesignE{11}{3595}{11}{395461}{391788}{0}{40423}{"Audrey 12 May 2014"}{{0, 1, 40423, 369538, 86221, 113690, 38589, 184963, 173683, 168776, 334537}}

\adfDesignE{11}{3601}{11}{396121}{1024}{0}{43728}{"Audrey 12 May 2014"}{{0, 1, 43728, 61917, 19541, 55851, 166563, 386158, 71036, 277447, 204549}}

\adfDesignE{11}{3606}{11}{396671}{416}{0}{137091}{"Audrey 18 Apr 2014"}{{0, 1, 137091, 66972, 308157, 89787, 288487, 79275, 293638, 160636, 162640}}

\adfDesignE{11}{3619}{11}{398101}{221640}{0}{68038}{"Audrey 12 May 2014"}{{0, 1, 68038, 51016, 382090, 246019, 86076, 373178, 199186, 62826, 144951}}

\adfDesignE{11}{3624}{11}{398651}{410}{0}{103481}{"Audrey 18 Apr 2014"}{{0, 1, 103481, 152850, 193774, 180645, 181204, 222688, 354524, 241118, 362970}}

\adfDesignE{11}{3667}{11}{403381}{34}{0}{79693}{"Audrey 18 Apr 2014"}{{0, 1, 79693, 143785, 217319, 43213, 110012, 103662, 296267, 112720, 103471}}

\adfDesignE{11}{3679}{11}{404701}{243}{0}{48448}{"Audrey 18 Apr 2014"}{{0, 1, 48448, 347605, 349028, 86661, 183954, 282671, 167469, 92464, 60503}}

\adfDesignE{11}{3682}{11}{405031}{178}{0}{50362}{"Audrey 18 Apr 2014"}{{0, 1, 50362, 26922, 207007, 193625, 220925, 23280, 267646, 161203, 64122}}

\adfDesignE{11}{3693}{11}{406241}{133}{0}{12455}{"Audrey 18 Apr 2014"}{{0, 1, 12455, 349204, 119674, 41441, 221585, 246062, 20106, 175774, 32421}}

\adfDesignE{11}{3706}{11}{407671}{90}{0}{95412}{"Audrey 18 Apr 2014"}{{0, 1, 95412, 156314, 403175, 305211, 37060, 238137, 399601, 115679, 287765}}

\adfDesignE{11}{3717}{11}{408881}{45}{0}{105570}{"Audrey 18 Apr 2014"}{{0, 1, 105570, 155483, 221446, 283045, 37170, 5943, 179056, 373290, 274520}}

\adfDesignE{11}{3732}{11}{410531}{410}{0}{33163}{"Audrey 18 Apr 2014"}{{0, 1, 33163, 382551, 309851, 408314, 373209, 41479, 289227, 399248, 226143}}

\adfDesignE{11}{3759}{11}{413501}{199}{0}{8754}{"Audrey 18 Apr 2014"}{{0, 1, 8754, 134831, 178720, 240597, 225545, 367156, 353852, 84417, 60131}}

\adfDesignE{11}{3781}{11}{415921}{45}{0}{159681}{"Audrey 18 Apr 2014"}{{0, 1, 159681, 400777, 371551, 168265, 226865, 142807, 279821, 119992, 209845}}

\adfDesignE{11}{3799}{11}{417901}{392074}{0}{75662}{"Audrey 12 May 2014"}{{0, 1, 75662, 330346, 398143, 319982, 219451, 59030, 219873, 227918, 46951}}

\adfDesignE{11}{3801}{11}{418121}{342}{0}{17219}{"Audrey 18 Apr 2014"}{{0, 1, 17219, 46172, 187647, 272726, 152043, 172836, 295927, 344507, 183406}}

\adfDesignE{11}{3889}{11}{427801}{111}{0}{24919}{"Audrey 18 Apr 2014"}{{0, 1, 24919, 217310, 42832, 394914, 155563, 169536, 132709, 73841, 71778}}

\adfDesignE{11}{3916}{11}{430771}{127}{0}{117804}{"Audrey 18 Apr 2014"}{{0, 1, 117804, 63880, 180921, 391488, 78321, 273806, 171086, 132367, 313410}}

\adfDesignE{11}{3934}{11}{432751}{50}{0}{41083}{"Audrey 18 Apr 2014"}{{0, 1, 41083, 83989, 196364, 310821, 275386, 273645, 152057, 197046, 200612}}

\adfDesignE{11}{3946}{11}{434071}{90}{0}{26948}{"Audrey 18 Apr 2014"}{{0, 1, 26948, 427992, 261946, 58206, 236765, 367662, 85001, 14281, 257482}}

\adfDesignE{11}{3951}{11}{434621}{23}{0}{44443}{"Audrey 18 Apr 2014"}{{0, 1, 44443, 262425, 334361, 313933, 355598, 153112, 330240, 139771, 239221}}

\adfDesignE{11}{4003}{11}{440341}{23}{0}{120631}{"Audrey 18 Apr 2014"}{{0, 1, 120631, 329475, 160406, 31623, 40030, 79524, 230959, 440059, 328656}}

\adfDesignE{11}{4015}{11}{441661}{193953}{0}{40918}{"Audrey 12 May 2014"}{{0, 1, 40918, 387534, 161329, 194716, 266509, 405172, 198939, 373772, 165588}}

\adfDesignE{11}{4035}{11}{443861}{59049}{0}{22049}{"Audrey 12 May 2014"}{{0, 1, 22049, 130606, 405387, 349006, 15137, 416102, 26128, 408555, 70200}}

\adfDesignE{11}{4053}{11}{445841}{1024}{0}{258636}{"Audrey 12 May 2014"}{{0, 1, 258636, 380220, 322232, 183263, 160876, 213811, 225043, 124639, 44740}}

\adfDesignE{11}{4059}{11}{446501}{144}{0}{3140}{"Audrey 18 Apr 2014"}{{0, 1, 3140, 36578, 104163, 233088, 81181, 402770, 206968, 220565, 51049}}

\adfDesignE{11}{4101}{11}{451121}{1024}{0}{172990}{"Audrey 12 May 2014"}{{0, 1, 172990, 428565, 234210, 359769, 237271, 266105, 214868, 351646, 281416}}

\adfDesignE{11}{4105}{11}{451561}{1024}{0}{18329}{"Audrey 12 May 2014"}{{0, 1, 18329, 442418, 398645, 55664, 189157, 424856, 16179, 320875, 187411}}

\adfDesignE{11}{4135}{11}{454861}{6568}{0}{145304}{"Audrey 12 May 2014"}{{0, 1, 145304, 424240, 96518, 177120, 209100, 171044, 227197, 186091, 99658}}

\adfDesignE{11}{4138}{11}{455191}{570}{0}{96363}{"Audrey 18 Apr 2014"}{{0, 1, 96363, 386560, 436177, 352884, 372428, 110542, 234155, 60395, 226450}}

\adfDesignE{11}{4149}{11}{456401}{67}{0}{77746}{"Audrey 18 Apr 2014"}{{0, 1, 77746, 322073, 359395, 198049, 373418, 88218, 258801, 304461, 299843}}

\adfDesignE{11}{4152}{11}{456731}{354}{0}{7902}{"Audrey 18 Apr 2014"}{{0, 1, 7902, 326188, 204543, 384508, 207604, 365787, 255106, 293709, 238307}}

\adfDesignE{11}{4198}{11}{461791}{97}{0}{32057}{"Audrey 18 Apr 2014"}{{0, 1, 32057, 166274, 253896, 77697, 293866, 387753, 169574, 291857, 174189}}

\adfDesignE{11}{4233}{11}{465641}{59049}{0}{29434}{"Audrey 12 May 2014"}{{0, 1, 29434, 268096, 385278, 51738, 210222, 236740, 353236, 316176, 23358}}

\adfDesignE{11}{4326}{11}{475871}{168}{0}{111961}{"Audrey 18 Apr 2014"}{{0, 1, 111961, 347510, 354150, 464688, 432609, 234127, 214883, 381287, 314010}}

\adfDesignE{11}{4395}{11}{483461}{33551}{0}{172441}{"Audrey 12 May 2014"}{{0, 1, 172441, 146215, 2743, 180805, 278576, 272134, 400590, 265588, 483239}}

\adfDesignE{11}{4438}{11}{488191}{26}{0}{170608}{"Audrey 18 Apr 2014"}{{0, 1, 170608, 165862, 369163, 152003, 221904, 361964, 233567, 296552, 469331}}

\adfDesignE{11}{4462}{11}{490831}{57}{0}{81327}{"Audrey 18 Apr 2014"}{{0, 1, 81327, 133204, 441538, 255797, 312346, 186399, 446869, 406261, 190413}}

\adfDesignE{11}{4465}{11}{491161}{1024}{0}{50704}{"Audrey 12 May 2014"}{{0, 1, 50704, 158942, 25480, 184490, 219715, 406719, 394430, 85122, 194181}}

\adfDesignE{11}{4474}{11}{492151}{86}{0}{107222}{"Audrey 18 Apr 2014"}{{0, 1, 107222, 402075, 334503, 84390, 268445, 250706, 403263, 247130, 363020}}

\adfDesignE{11}{4477}{11}{492481}{383494}{0}{160680}{"Audrey 12 May 2014"}{{0, 1, 160680, 238456, 88280, 392638, 287816, 339056, 284898, 316728, 345943}}

\adfDesignE{11}{4485}{11}{493361}{65891}{0}{265598}{"Audrey 12 May 2014"}{{0, 1, 265598, 61741, 446561, 243995, 136577, 210521, 367706, 181516, 36370}}

\adfDesignE{11}{4516}{11}{496771}{262}{0}{82938}{"Audrey 18 Apr 2014"}{{0, 1, 82938, 420578, 128857, 107343, 180643, 52545, 301998, 413075, 299106}}

\adfDesignE{11}{4519}{11}{497101}{155}{0}{47480}{"Audrey 18 Apr 2014"}{{0, 1, 47480, 494466, 159652, 480912, 361527, 404430, 318972, 111494, 106571}}

\adfDesignE{11}{4554}{11}{500951}{456}{0}{95181}{"Audrey 18 Apr 2014"}{{0, 1, 95181, 224877, 385311, 164532, 91081, 223606, 139451, 388886, 290878}}

\adfDesignE{11}{4564}{11}{502051}{83}{0}{7973}{"Audrey 18 Apr 2014"}{{0, 1, 7973, 310303, 440542, 92570, 45640, 402796, 374312, 198432, 135635}}

\adfDesignE{11}{4641}{11}{510521}{59049}{0}{89868}{"Audrey 12 May 2014"}{{0, 1, 89868, 325725, 1202, 301405, 462364, 423762, 329421, 314680, 372487}}

\adfDesignE{11}{4644}{11}{510851}{37}{0}{128153}{"Audrey 18 Apr 2014"}{{0, 1, 128153, 353461, 440214, 447110, 417968, 104052, 343154, 117078, 203064}}

\adfDesignE{11}{4645}{11}{510961}{172778}{0}{5602}{"Audrey 12 May 2014"}{{0, 1, 5602, 213783, 430743, 264444, 139349, 395651, 399045, 506676, 10597}}

\adfDesignE{11}{4659}{11}{512501}{503559}{0}{120679}{"Audrey 12 May 2014"}{{0, 1, 120679, 192625, 284518, 343227, 472814, 446873, 269042, 268667, 114130}}

\adfDesignE{11}{4683}{11}{515141}{59049}{0}{92910}{"Audrey 12 May 2014"}{{0, 1, 92910, 50363, 200627, 392626, 202027, 135953, 135910, 261908, 156863}}

\adfDesignE{11}{4738}{11}{521191}{86}{0}{7362}{"Audrey 18 Apr 2014"}{{0, 1, 7362, 516371, 477339, 299996, 284285, 324305, 478630, 422900, 315957}}

\adfDesignE{11}{4743}{11}{521741}{122}{0}{33192}{"Audrey 18 Apr 2014"}{{0, 1, 33192, 313613, 188005, 239600, 426878, 14239, 445283, 476029, 471865}}

\adfDesignE{11}{4809}{11}{529001}{309}{0}{53469}{"Audrey 18 Apr 2014"}{{0, 1, 53469, 212557, 152749, 89842, 432818, 138895, 460717, 87706, 487250}}

\adfDesignE{11}{4822}{11}{530431}{40}{0}{76224}{"Audrey 18 Apr 2014"}{{0, 1, 76224, 287433, 370968, 449084, 144662, 116660, 146556, 207684, 322452}}

\adfDesignE{11}{4828}{11}{531091}{69}{0}{67105}{"Audrey 18 Apr 2014"}{{0, 1, 67105, 491527, 512780, 184919, 48280, 174300, 184407, 211435, 249610}}

\adfDesignE{11}{4834}{11}{531751}{46}{0}{268559}{"Audrey 18 Apr 2014"}{{0, 1, 268559, 421347, 447924, 246794, 241704, 408215, 303768, 518896, 333298}}

\adfDesignE{11}{4837}{11}{532081}{1024}{0}{265114}{"Audrey 12 May 2014"}{{0, 1, 265114, 193301, 483961, 420457, 195922, 449969, 521266, 176599, 528015}}

\adfDesignE{11}{4861}{11}{534721}{1024}{0}{354957}{"Audrey 12 May 2014"}{{0, 1, 354957, 301503, 469989, 419567, 424304, 158068, 137988, 432358, 363280}}

\adfDesignE{11}{4875}{11}{536261}{59049}{0}{288264}{"Audrey 12 May 2014"}{{0, 1, 288264, 346702, 351541, 446176, 176485, 263492, 322370, 69512, 414903}}

\adfDesignE{11}{4899}{11}{538901}{109264}{0}{83813}{"Audrey 12 May 2014"}{{0, 1, 83813, 44434, 340932, 385993, 465378, 149736, 435781, 97678, 241123}}

\adfDesignE{11}{4903}{11}{539341}{12}{0}{19442}{"Audrey 18 Apr 2014"}{{0, 1, 19442, 452664, 266391, 421540, 294185, 372806, 429894, 371012, 68770}}

\adfDesignE{11}{4912}{11}{540331}{348}{0}{215572}{"Audrey 18 Apr 2014"}{{0, 1, 215572, 119529, 341191, 289870, 196483, 226517, 469923, 444745, 397824}}

\adfDesignE{11}{4917}{11}{540881}{144}{0}{76431}{"Audrey 18 Apr 2014"}{{0, 1, 76431, 182961, 495698, 143312, 98341, 218595, 161236, 536893, 250056}}

\adfDesignE{11}{4933}{11}{542641}{56}{0}{201578}{"Audrey 18 Apr 2014"}{{0, 1, 201578, 189363, 418951, 85848, 246654, 538387, 402609, 271683, 358131}}

\adfDesignE{11}{4939}{11}{543301}{23}{0}{86332}{"Audrey 18 Apr 2014"}{{0, 1, 86332, 211106, 171147, 392109, 98781, 308796, 282804, 174590, 447538}}

\adfDesignE{11}{4941}{11}{543521}{1024}{0}{194830}{"Audrey 12 May 2014"}{{0, 1, 194830, 309302, 48348, 421910, 239823, 388804, 161550, 28911, 222007}}

\adfDesignE{11}{4971}{11}{546821}{111}{0}{36320}{"Audrey 18 Apr 2014"}{{0, 1, 36320, 210148, 47842, 371123, 49710, 411079, 535517, 101291, 424253}}

\adfDesignE{11}{4981}{11}{547921}{194866}{0}{125073}{"Audrey 12 May 2014"}{{0, 1, 125073, 110779, 183540, 200204, 125192, 200599, 216137, 124624, 388865}}

\adfDesignE{11}{4987}{11}{548581}{12}{0}{34389}{"Audrey 18 Apr 2014"}{{0, 1, 34389, 411266, 59713, 131674, 149612, 414450, 386670, 139771, 466778}}

\adfDesignE{11}{5026}{11}{552871}{50}{0}{51188}{"Audrey 18 Apr 2014"}{{0, 1, 51188, 155675, 162177, 158211, 50260, 200117, 547979, 38867, 294138}}

\adfDesignE{11}{5031}{11}{553421}{122}{0}{101348}{"Audrey 18 Apr 2014"}{{0, 1, 101348, 476765, 545131, 471579, 150932, 99896, 529455, 58601, 333397}}

\adfDesignE{11}{5044}{11}{554851}{82}{0}{47111}{"Audrey 18 Apr 2014"}{{0, 1, 47111, 42321, 204988, 8013, 201763, 104212, 209884, 400304, 445956}}

\adfDesignE{11}{5046}{11}{555071}{52}{0}{67235}{"Audrey 18 Apr 2014"}{{0, 1, 67235, 47001, 93032, 466492, 302765, 285992, 457609, 310656, 189501}}

\adfDesignE{11}{5058}{11}{556391}{445}{0}{321501}{"Audrey 18 Apr 2014"}{{0, 1, 321501, 467758, 523323, 124160, 404647, 184309, 442700, 336554, 533393}}

\adfDesignE{11}{5107}{11}{561781}{111}{0}{314704}{"Audrey 18 Apr 2014"}{{0, 1, 314704, 549783, 473290, 136068, 510709, 553503, 420566, 445988, 527855}}

\adfDesignE{11}{5113}{11}{562441}{78}{0}{103480}{"Audrey 18 Apr 2014"}{{0, 1, 103480, 358642, 167216, 14315, 409047, 561223, 510585, 190701, 496995}}

\adfDesignE{11}{5151}{11}{566621}{78}{0}{269239}{"Audrey 18 Apr 2014"}{{0, 1, 269239, 114728, 474798, 474775, 515109, 142849, 554915, 400889, 452423}}

\adfDesignE{11}{5172}{11}{568931}{168}{0}{32470}{"Audrey 18 Apr 2014"}{{0, 1, 32470, 71757, 177345, 241499, 465488, 174414, 83406, 81260, 379153}}

\adfDesignE{11}{5197}{11}{571681}{12}{0}{252080}{"Audrey 18 Apr 2014"}{{0, 1, 252080, 268207, 338776, 374619, 259854, 215659, 459387, 284876, 404946}}

\adfDesignE{11}{5202}{11}{572231}{136}{0}{13381}{"Audrey 18 Apr 2014"}{{0, 1, 13381, 515089, 455745, 58078, 52020, 246724, 213205, 324570, 410111}}

\adfDesignE{11}{5212}{11}{573331}{86}{0}{129859}{"Audrey 18 Apr 2014"}{{0, 1, 129859, 548509, 480115, 374190, 416967, 391351, 389669, 405942, 303383}}

\adfDesignE{11}{5239}{11}{576301}{45}{0}{138736}{"Audrey 18 Apr 2014"}{{0, 1, 138736, 376898, 378596, 245215, 523909, 228001, 513749, 302487, 173913}}

\adfDesignE{11}{5254}{11}{577951}{167}{0}{15438}{"Audrey 18 Apr 2014"}{{0, 1, 15438, 216032, 324746, 281774, 367786, 89644, 309378, 568451, 138554}}

\adfDesignE{11}{5269}{11}{579601}{12}{0}{199149}{"Audrey 18 Apr 2014"}{{0, 1, 199149, 546175, 542612, 406349, 105381, 327761, 349572, 458517, 542089}}

\adfDesignE{11}{5298}{11}{582791}{53}{0}{278}{"Audrey 18 Apr 2014"}{{0, 1, 278, 77284, 504476, 374488, 370866, 529532, 346564, 184277, 526189}}

\adfDesignE{11}{5323}{11}{585541}{12}{0}{359188}{"Audrey 18 Apr 2014"}{{0, 1, 359188, 257568, 442325, 64865, 53230, 492508, 487666, 356140, 414214}}

\adfDesignE{11}{5328}{11}{586091}{137}{0}{23995}{"Audrey 18 Apr 2014"}{{0, 1, 23995, 218663, 132053, 203789, 159842, 29286, 580552, 134352, 275740}}

\adfDesignE{11}{5338}{11}{587191}{112}{0}{40989}{"Audrey 18 Apr 2014"}{{0, 1, 40989, 144670, 423912, 160087, 533809, 386059, 549283, 483565, 213580}}

\adfDesignE{11}{5383}{11}{592141}{111}{0}{124041}{"Audrey 18 Apr 2014"}{{0, 1, 124041, 570078, 159119, 36067, 161492, 91283, 506542, 486853, 233088}}

\adfDesignE{11}{5401}{11}{594121}{111}{0}{151607}{"Audrey 18 Apr 2014"}{{0, 1, 151607, 517443, 244061, 94268, 108021, 388503, 400744, 188027, 283809}}

\adfDesignE{11}{5551}{11}{610621}{12}{0}{99938}{"Audrey 18 Apr 2014"}{{0, 1, 99938, 286768, 134370, 502649, 388576, 455172, 157320, 587273, 441038}}

\adfDesignE{11}{5554}{11}{610951}{50}{0}{21338}{"Audrey 18 Apr 2014"}{{0, 1, 21338, 151749, 590813, 404860, 55540, 478531, 70415, 186761, 483796}}

\adfDesignE{11}{5569}{11}{612601}{45}{0}{299347}{"Audrey 18 Apr 2014"}{{0, 1, 299347, 415134, 554244, 549838, 556909, 90490, 531613, 169739, 508291}}

\adfDesignE{11}{5572}{11}{612931}{469}{0}{16307}{"Audrey 18 Apr 2014"}{{0, 1, 16307, 519126, 197641, 140589, 222883, 485182, 149526, 80964, 26574}}

\adfDesignE{11}{5667}{11}{623381}{155}{0}{178649}{"Audrey 18 Apr 2014"}{{0, 1, 178649, 228144, 424295, 488141, 510038, 48035, 565250, 482441, 191911}}

\adfDesignE{11}{5668}{11}{623491}{47}{0}{176456}{"Audrey 18 Apr 2014"}{{0, 1, 176456, 202887, 398743, 258949, 566809, 163630, 250561, 621515, 477904}}

\adfDesignE{11}{5724}{11}{629651}{214}{0}{39752}{"Audrey 18 Apr 2014"}{{0, 1, 39752, 427145, 69523, 140057, 171722, 246453, 259747, 445646, 88907}}

\adfDesignE{11}{5733}{11}{630641}{111}{0}{26464}{"Audrey 18 Apr 2014"}{{0, 1, 26464, 331786, 600702, 410041, 515978, 202860, 470848, 316594, 277931}}

\adfDesignE{11}{5779}{11}{635701}{243}{0}{12972}{"Audrey 18 Apr 2014"}{{0, 1, 12972, 447720, 59504, 144874, 173372, 507147, 476936, 171660, 548618}}

\adfDesignE{11}{5823}{11}{640541}{89}{0}{65542}{"Audrey 18 Apr 2014"}{{0, 1, 65542, 285818, 461811, 532689, 174692, 633230, 587647, 469885, 631931}}

\adfDesignE{11}{5883}{11}{647141}{23}{0}{94461}{"Audrey 18 Apr 2014"}{{0, 1, 94461, 100413, 613897, 313789, 470647, 493849, 311404, 386230, 451014}}

\adfDesignE{11}{5892}{11}{648131}{71}{0}{218771}{"Audrey 18 Apr 2014"}{{0, 1, 218771, 164877, 519755, 514727, 412446, 370439, 306491, 246218, 486930}}

\adfDesignE{11}{5922}{11}{651431}{211}{0}{338144}{"Audrey 18 Apr 2014"}{{0, 1, 338144, 241323, 420297, 161791, 177662, 372508, 646992, 527239, 371198}}

\adfDesignE{11}{5928}{11}{652091}{167}{0}{243719}{"Audrey 18 Apr 2014"}{{0, 1, 243719, 633862, 594423, 382122, 59280, 586215, 551758, 354073, 507093}}

\adfDesignE{11}{5934}{11}{652751}{438}{0}{179558}{"Audrey 18 Apr 2014"}{{0, 1, 179558, 397972, 446153, 168397, 296704, 651216, 492143, 87916, 543695}}

\adfDesignE{11}{5967}{11}{656381}{100}{0}{143960}{"Audrey 18 Apr 2014"}{{0, 1, 143960, 564287, 387579, 205935, 298354, 94724, 151765, 447815, 331104}}

\adfDesignE{11}{5979}{11}{657701}{100}{0}{50099}{"Audrey 18 Apr 2014"}{{0, 1, 50099, 122785, 585963, 333903, 239163, 488020, 594707, 370693, 503171}}

\adfDesignE{11}{6109}{11}{672001}{34}{0}{227595}{"Audrey 18 Apr 2014"}{{0, 1, 227595, 302943, 337484, 628681, 183272, 16769, 246876, 395608, 348775}}

\adfDesignE{11}{6156}{11}{677171}{448}{0}{196386}{"Audrey 18 Apr 2014"}{{0, 1, 196386, 541033, 468154, 61945, 430926, 419224, 628626, 332139, 307421}}

\adfDesignE{11}{6187}{11}{680581}{34}{0}{78808}{"Audrey 18 Apr 2014"}{{0, 1, 78808, 399239, 648063, 389502, 309354, 478031, 466955, 94389, 538563}}

\adfDesignE{11}{6196}{11}{681571}{123}{0}{194951}{"Audrey 18 Apr 2014"}{{0, 1, 194951, 130299, 450750, 577362, 619609, 610542, 303428, 44938, 475975}}

\adfDesignE{11}{6201}{11}{682121}{89}{0}{103690}{"Audrey 18 Apr 2014"}{{0, 1, 103690, 24898, 527756, 544536, 372065, 20332, 471190, 92354, 571662}}

\adfDesignE{11}{6219}{11}{684101}{45}{0}{126240}{"Audrey 18 Apr 2014"}{{0, 1, 126240, 404805, 238500, 270889, 186572, 620052, 528060, 72455, 288830}}

\adfDesignE{11}{6258}{11}{688391}{86}{0}{11201}{"Audrey 18 Apr 2014"}{{0, 1, 11201, 175239, 249298, 273002, 62580, 176542, 387990, 63607, 665713}}

\adfDesignE{11}{6259}{11}{688501}{122}{0}{89588}{"Audrey 18 Apr 2014"}{{0, 1, 89588, 153587, 548172, 233808, 125181, 411140, 472323, 578466, 141738}}

\adfDesignE{11}{6276}{11}{690371}{414}{0}{173034}{"Audrey 18 Apr 2014"}{{0, 1, 173034, 65257, 662033, 267721, 251043, 140771, 499592, 216621, 485411}}

\adfDesignE{11}{6298}{11}{692791}{156}{0}{161780}{"Audrey 18 Apr 2014"}{{0, 1, 161780, 510002, 179415, 586964, 251923, 593992, 371732, 387414, 420732}}

\adfDesignE{11}{6313}{11}{694441}{100}{0}{56312}{"Audrey 18 Apr 2014"}{{0, 1, 56312, 223738, 585634, 607600, 63130, 133081, 344441, 424462, 339365}}

\adfDesignE{11}{6361}{11}{699721}{67}{0}{8979}{"Audrey 18 Apr 2014"}{{0, 1, 8979, 154526, 641932, 305551, 636109, 499909, 672417, 439455, 139726}}

\adfDesignE{11}{6417}{11}{705881}{89}{0}{217688}{"Audrey 18 Apr 2014"}{{0, 1, 217688, 156171, 617807, 486810, 192512, 103167, 613881, 665613, 475755}}

\adfDesignE{11}{6481}{11}{712921}{34}{0}{72103}{"Audrey 18 Apr 2014"}{{0, 1, 72103, 222677, 698811, 677858, 583298, 187141, 684677, 338165, 99874}}

\adfDesignE{11}{6487}{11}{713581}{243}{0}{70772}{"Audrey 18 Apr 2014"}{{0, 1, 70772, 50945, 468328, 98928, 389225, 577938, 78797, 699350, 420040}}

\adfDesignE{11}{6511}{11}{716221}{23}{0}{16162}{"Audrey 18 Apr 2014"}{{0, 1, 16162, 505800, 509327, 215021, 65110, 179171, 80199, 532449, 45423}}

\adfDesignE{11}{6538}{11}{719191}{93}{0}{67852}{"Audrey 18 Apr 2014"}{{0, 1, 67852, 352313, 671218, 713661, 196142, 716720, 628402, 374878, 593959}}

\adfDesignE{11}{6558}{11}{721391}{248}{0}{132367}{"Audrey 18 Apr 2014"}{{0, 1, 132367, 599472, 185788, 702397, 590228, 64376, 187500, 76536, 346899}}

\adfDesignE{11}{6588}{11}{724691}{372}{0}{28130}{"Audrey 18 Apr 2014"}{{0, 1, 28130, 659019, 608690, 175443, 65880, 169513, 658601, 445406, 88081}}

\adfDesignE{11}{6619}{11}{728101}{23}{0}{56963}{"Audrey 18 Apr 2014"}{{0, 1, 56963, 365313, 197839, 683780, 397145, 472565, 98024, 662644, 706231}}

\adfDesignE{11}{6649}{11}{731401}{221}{0}{97036}{"Audrey 18 Apr 2014"}{{0, 1, 97036, 660223, 522636, 624358, 332454, 102437, 337142, 75783, 173534}}

\adfDesignE{11}{6657}{11}{732281}{276}{0}{532051}{"Audrey 18 Apr 2014"}{{0, 1, 532051, 400431, 612122, 344315, 599138, 533085, 598134, 186730, 587679}}

\adfDesignE{11}{6774}{11}{745151}{118}{0}{109509}{"Audrey 18 Apr 2014"}{{0, 1, 109509, 506038, 325774, 335690, 541927, 567901, 713300, 80672, 544943}}

\adfDesignE{11}{6807}{11}{748781}{155}{0}{95232}{"Audrey 18 Apr 2014"}{{0, 1, 95232, 647133, 98432, 635666, 544567, 381265, 237790, 582278, 521541}}

\adfDesignE{11}{6814}{11}{749551}{271}{0}{72064}{"Audrey 18 Apr 2014"}{{0, 1, 72064, 330768, 743352, 7660, 340704, 200500, 486924, 210622, 605609}}

\adfDesignE{11}{6858}{11}{754391}{103}{0}{272237}{"Audrey 18 Apr 2014"}{{0, 1, 272237, 103547, 750533, 576317, 342904, 550735, 513682, 278182, 383817}}

\adfDesignE{11}{6861}{11}{754721}{45}{0}{264031}{"Audrey 18 Apr 2014"}{{0, 1, 264031, 299633, 281240, 588692, 411665, 422079, 592310, 199037, 714917}}

\adfDesignE{11}{6871}{11}{755821}{56}{0}{48392}{"Audrey 18 Apr 2014"}{{0, 1, 48392, 252206, 511065, 238539, 480976, 638718, 297482, 382178, 173727}}

\adfDesignE{11}{6879}{11}{756701}{45}{0}{101033}{"Audrey 18 Apr 2014"}{{0, 1, 101033, 527300, 680397, 47756, 206372, 242922, 298192, 695423, 227408}}

\adfDesignE{11}{6903}{11}{759341}{551}{0}{324824}{"Audrey 18 Apr 2014"}{{0, 1, 324824, 199026, 406707, 325411, 276123, 296455, 629346, 697789, 641023}}

\adfDesignE{11}{6943}{11}{763741}{177}{0}{21932}{"Audrey 18 Apr 2014"}{{0, 1, 21932, 619535, 689230, 230488, 624878, 255792, 352499, 421666, 602684}}

\adfDesignE{11}{6993}{11}{769241}{309}{0}{372267}{"Audrey 18 Apr 2014"}{{0, 1, 372267, 106934, 546869, 112891, 419585, 187181, 382583, 362234, 586419}}

\adfDesignE{11}{6999}{11}{769901}{100}{0}{6868}{"Audrey 18 Apr 2014"}{{0, 1, 6868, 205463, 661252, 602638, 699909, 483069, 214483, 248631, 727191}}

\adfDesignE{11}{7092}{11}{780131}{194}{0}{66599}{"Audrey 18 Apr 2014"}{{0, 1, 66599, 382066, 460838, 216291, 425525, 500769, 114381, 461135, 492919}}

\adfDesignE{11}{7098}{11}{780791}{234}{0}{105121}{"Audrey 18 Apr 2014"}{{0, 1, 105121, 670409, 649620, 723160, 709809, 320765, 678230, 628038, 199593}}

\adfDesignE{11}{7119}{11}{783101}{199}{0}{163432}{"Audrey 18 Apr 2014"}{{0, 1, 163432, 9716, 559585, 428536, 640718, 690860, 346239, 437089, 639329}}

\adfDesignE{11}{7266}{11}{799271}{433}{0}{204945}{"Audrey 18 Apr 2014"}{{0, 1, 204945, 761975, 599124, 260076, 290643, 158360, 691245, 418130, 611856}}

\adfDesignE{11}{7347}{11}{808181}{100}{0}{73682}{"Audrey 18 Apr 2014"}{{0, 1, 73682, 485347, 136585, 386158, 73470, 220202, 690189, 524654, 642436}}

\adfDesignE{11}{7377}{11}{811481}{144}{0}{218991}{"Audrey 18 Apr 2014"}{{0, 1, 218991, 153943, 776330, 767606, 516396, 518719, 436225, 182693, 486501}}

\adfDesignE{11}{7557}{11}{831281}{89}{0}{426006}{"Audrey 18 Apr 2014"}{{0, 1, 426006, 521, 828380, 271441, 151141, 102991, 604047, 456327, 484269}}

\adfDesignE{11}{7564}{11}{832051}{304}{0}{188334}{"Audrey 18 Apr 2014"}{{0, 1, 188334, 193477, 287875, 207090, 529486, 568076, 411651, 695458, 246756}}

\adfDesignE{11}{7603}{11}{836341}{12}{0}{275421}{"Audrey 18 Apr 2014"}{{0, 1, 275421, 598541, 422592, 479626, 684278, 104734, 541924, 489780, 584808}}

\adfDesignE{11}{7609}{11}{837001}{56}{0}{52231}{"Audrey 18 Apr 2014"}{{0, 1, 52231, 291102, 425397, 719162, 456545, 480406, 469808, 183331, 270021}}

\adfDesignE{11}{7642}{11}{840631}{327}{0}{159432}{"Audrey 18 Apr 2014"}{{0, 1, 159432, 403077, 494838, 633297, 458525, 604778, 590396, 40109, 818702}}

\adfDesignE{11}{7647}{11}{841181}{56}{0}{120047}{"Audrey 18 Apr 2014"}{{0, 1, 120047, 169317, 541396, 798009, 688238, 109366, 748335, 605669, 425527}}

\adfDesignE{11}{7677}{11}{844481}{89}{0}{184772}{"Audrey 18 Apr 2014"}{{0, 1, 184772, 14116, 484224, 808421, 76770, 199083, 216197, 667341, 726999}}

\adfDesignE{11}{7696}{11}{846571}{15}{0}{229255}{"Audrey 18 Apr 2014"}{{0, 1, 229255, 187632, 455079, 265818, 538726, 232511, 812861, 160609, 503792}}

\adfDesignE{11}{7704}{11}{847451}{268}{0}{190204}{"Audrey 18 Apr 2014"}{{0, 1, 190204, 725877, 534341, 692036, 231122, 503165, 506779, 621274, 232456}}

\adfDesignE{11}{7719}{11}{849101}{177}{0}{283650}{"Audrey 18 Apr 2014"}{{0, 1, 283650, 757245, 558886, 8099, 463145, 719833, 709384, 213024, 532238}}

\adfDesignE{11}{7749}{11}{852401}{331}{0}{151287}{"Audrey 18 Apr 2014"}{{0, 1, 151287, 789519, 418427, 710086, 387454, 546132, 295355, 511465, 452279}}

\adfDesignE{11}{7941}{11}{873521}{45}{0}{67362}{"Audrey 18 Apr 2014"}{{0, 1, 67362, 570970, 551510, 842011, 79410, 647337, 620195, 560244, 428565}}

\adfDesignE{11}{7969}{11}{876601}{34}{0}{215011}{"Audrey 18 Apr 2014"}{{0, 1, 215011, 423184, 661027, 373162, 398454, 24062, 772181, 56192, 583130}}

\adfDesignE{11}{8014}{11}{881551}{145}{0}{240960}{"Audrey 18 Apr 2014"}{{0, 1, 240960, 128087, 743010, 615459, 320563, 480309, 833605, 539246, 506515}}

\adfDesignE{11}{8083}{11}{889141}{45}{0}{37629}{"Audrey 18 Apr 2014"}{{0, 1, 37629, 429169, 621459, 472411, 646647, 447357, 379141, 429344, 93406}}

\adfDesignE{11}{8091}{11}{890021}{144}{0}{117859}{"Audrey 18 Apr 2014"}{{0, 1, 117859, 186134, 329498, 18489, 323643, 610340, 784798, 75057, 224244}}

\adfDesignE{11}{8103}{11}{891341}{34}{0}{146657}{"Audrey 18 Apr 2014"}{{0, 1, 146657, 217319, 563787, 736217, 567216, 15405, 592991, 813740, 803372}}

\adfDesignE{11}{8107}{11}{891781}{12}{0}{246508}{"Audrey 18 Apr 2014"}{{0, 1, 246508, 236724, 670057, 517698, 162141, 320989, 411844, 708150, 293012}}

\adfDesignE{11}{8142}{11}{895631}{211}{0}{20706}{"Audrey 18 Apr 2014"}{{0, 1, 20706, 626818, 304687, 24258, 732788, 223557, 350234, 20997, 382847}}

\adfDesignE{11}{8167}{11}{898381}{12}{0}{30485}{"Audrey 18 Apr 2014"}{{0, 1, 30485, 409271, 809488, 512372, 408354, 704554, 734123, 170564, 712693}}

\adfDesignE{11}{8193}{11}{901241}{45}{0}{26372}{"Audrey 18 Apr 2014"}{{0, 1, 26372, 625573, 394651, 205104, 655447, 547145, 439530, 424659, 286482}}

\adfDesignE{11}{8217}{11}{903881}{23}{0}{94165}{"Audrey 18 Apr 2014"}{{0, 1, 94165, 878496, 386720, 834953, 164341, 727545, 518411, 270648, 644125}}

\adfDesignE{11}{8281}{11}{910921}{23}{0}{158409}{"Audrey 18 Apr 2014"}{{0, 1, 158409, 270494, 782048, 7474, 662487, 338457, 557616, 294495, 572203}}

\adfDesignE{11}{8298}{11}{912791}{367}{0}{326885}{"Audrey 18 Apr 2014"}{{0, 1, 326885, 663183, 363619, 691168, 248942, 88020, 332589, 383210, 553547}}

\adfDesignE{11}{8401}{11}{924121}{56}{0}{25326}{"Audrey 18 Apr 2014"}{{0, 1, 25326, 66302, 36595, 835728, 504065, 142696, 605786, 803515, 676470}}

\adfDesignE{11}{8413}{11}{925441}{12}{0}{125196}{"Audrey 18 Apr 2014"}{{0, 1, 125196, 769640, 783402, 559612, 673047, 463721, 324063, 57908, 870615}}

\adfDesignE{11}{8422}{11}{926431}{108}{0}{131839}{"Audrey 18 Apr 2014"}{{0, 1, 131839, 749930, 378519, 434195, 589546, 373487, 344943, 295249, 408015}}

\adfDesignE{11}{8452}{11}{929731}{63}{0}{102276}{"Audrey 18 Apr 2014"}{{0, 1, 102276, 906426, 288104, 160121, 253562, 320329, 107826, 472585, 177963}}

\adfDesignE{11}{8473}{11}{932041}{23}{0}{97936}{"Audrey 18 Apr 2014"}{{0, 1, 97936, 758206, 888387, 905964, 847309, 579912, 343297, 552040, 619194}}

\adfDesignE{11}{8506}{11}{935671}{60}{0}{29078}{"Audrey 18 Apr 2014"}{{0, 1, 29078, 619171, 72956, 248411, 850609, 481288, 61317, 522471, 857382}}

\adfDesignE{11}{8509}{11}{936001}{23}{0}{467652}{"Audrey 18 Apr 2014"}{{0, 1, 467652, 823453, 846937, 150771, 340363, 723622, 538003, 174155, 615048}}

\adfDesignE{11}{8538}{11}{939191}{536}{0}{32015}{"Audrey 18 Apr 2014"}{{0, 1, 32015, 302844, 281967, 608804, 768428, 53366, 124061, 913367, 671911}}

\adfDesignE{11}{8611}{11}{947221}{12}{0}{39273}{"Audrey 18 Apr 2014"}{{0, 1, 39273, 292741, 396016, 314769, 688887, 132949, 223925, 206761, 546341}}

\adfDesignE{11}{8613}{11}{947441}{254}{0}{436885}{"Audrey 18 Apr 2014"}{{0, 1, 436885, 829129, 800517, 236010, 172261, 166032, 807360, 663710, 625300}}

\adfDesignE{11}{8749}{11}{962401}{89}{0}{151924}{"Audrey 18 Apr 2014"}{{0, 1, 151924, 600994, 504784, 843132, 262472, 635495, 798862, 807581, 206760}}

\adfDesignE{11}{8781}{11}{965921}{122}{0}{39977}{"Audrey 18 Apr 2014"}{{0, 1, 39977, 527195, 244216, 459485, 878109, 662511, 614348, 282650, 155192}}

\adfDesignE{11}{8826}{11}{970871}{59}{0}{58942}{"Audrey 18 Apr 2014"}{{0, 1, 58942, 382926, 586155, 703475, 264782, 29319, 940989, 826021, 90874}}

\adfDesignE{11}{8859}{11}{974501}{56}{0}{246585}{"Audrey 18 Apr 2014"}{{0, 1, 246585, 172330, 876945, 685426, 265772, 196370, 890762, 894875, 643439}}

\adfDesignE{11}{8872}{11}{975931}{156}{0}{229585}{"Audrey 18 Apr 2014"}{{0, 1, 229585, 214846, 890239, 195209, 354883, 214020, 583643, 351855, 869443}}

\adfDesignE{11}{8959}{11}{985501}{12}{0}{148613}{"Audrey 18 Apr 2014"}{{0, 1, 148613, 746359, 512517, 273134, 447954, 209751, 328733, 741757, 533185}}

\adfDesignE{11}{8967}{11}{986381}{23}{0}{9343}{"Audrey 18 Apr 2014"}{{0, 1, 9343, 490121, 419901, 297806, 807038, 259670, 585931, 925164, 150549}}

\adfDesignE{11}{8989}{11}{988801}{12}{0}{132366}{"Audrey 18 Apr 2014"}{{0, 1, 132366, 193037, 917702, 317684, 809018, 317089, 166527, 160990, 940790}}

\adfDesignE{11}{9007}{11}{990781}{12}{0}{234349}{"Audrey 18 Apr 2014"}{{0, 1, 234349, 462971, 326693, 548425, 720567, 396148, 707952, 574017, 982782}}

\adfDesignE{11}{9012}{11}{991331}{156}{0}{25161}{"Audrey 18 Apr 2014"}{{0, 1, 25161, 606743, 754554, 353213, 901209, 605686, 925314, 417019, 367755}}

\adfDesignE{11}{9097}{11}{1000681}{67}{0}{493623}{"Audrey 18 Apr 2014"}{{0, 1, 493623, 844672, 777791, 244799, 181941, 143174, 883777, 768716, 664911}}

\adfDesignE{11}{9258}{11}{1018391}{137}{0}{317772}{"Audrey 18 Apr 2014"}{{0, 1, 317772, 484379, 431066, 1005106, 648066, 237714, 719174, 309582, 939095}}

\adfDesignE{11}{9267}{11}{1019381}{34}{0}{20432}{"Audrey 18 Apr 2014"}{{0, 1, 20432, 539795, 408401, 815747, 463354, 257581, 850270, 425638, 296305}}

\adfDesignE{11}{9313}{11}{1024441}{232}{0}{460464}{"Audrey 18 Apr 2014"}{{0, 1, 460464, 590408, 598937, 189599, 651916, 98322, 620295, 146111, 741711}}

\adfDesignE{11}{9358}{11}{1029391}{185}{0}{282714}{"Audrey 18 Apr 2014"}{{0, 1, 282714, 141601, 598515, 365303, 561485, 372353, 794209, 150133, 851250}}

\adfDesignE{11}{9429}{11}{1037201}{56}{0}{570656}{"Audrey 18 Apr 2014"}{{0, 1, 570656, 346768, 772621, 147889, 848618, 843709, 973306, 736035, 983603}}

\adfDesignE{11}{9462}{11}{1040831}{356}{0}{80745}{"Audrey 18 Apr 2014"}{{0, 1, 80745, 1030472, 390669, 103288, 851588, 13876, 483464, 934025, 275196}}

\adfDesignE{11}{9477}{11}{1042481}{111}{0}{256402}{"Audrey 18 Apr 2014"}{{0, 1, 256402, 6301, 785933, 88323, 379083, 880850, 278012, 67006, 385532}}

\adfDesignE{11}{9496}{11}{1044571}{27}{0}{333646}{"Audrey 18 Apr 2014"}{{0, 1, 333646, 766417, 985582, 362888, 949609, 236120, 1037842, 723716, 871605}}

\adfDesignE{11}{9531}{11}{1048421}{56}{0}{276631}{"Audrey 18 Apr 2014"}{{0, 1, 276631, 461371, 1039087, 187169, 476554, 104613, 682361, 295267, 670630}}

\adfDesignE{11}{9544}{11}{1049851}{148}{0}{5289}{"Audrey 18 Apr 2014"}{{0, 1, 5289, 677395, 650543, 360200, 668086, 758239, 945102, 303867, 880533}}

\adfDesignE{11}{9546}{11}{1050071}{345}{0}{333788}{"Audrey 18 Apr 2014"}{{0, 1, 333788, 845773, 439987, 500767, 763687, 620822, 872525, 132779, 740226}}

\adfDesignE{11}{9553}{11}{1050841}{100}{0}{49146}{"Audrey 18 Apr 2014"}{{0, 1, 49146, 496698, 734319, 859952, 477654, 46385, 363081, 698646, 477482}}

\adfDesignE{11}{9562}{11}{1051831}{468}{0}{2786}{"Audrey 18 Apr 2014"}{{0, 1, 2786, 398979, 821958, 138901, 956209, 762182, 844094, 803599, 530446}}

\adfDesignE{11}{9622}{11}{1058431}{360}{0}{99889}{"Audrey 18 Apr 2014"}{{0, 1, 99889, 1041715, 459594, 1057303, 577325, 862321, 209158, 213953, 770896}}

\adfDesignE{11}{9646}{11}{1061071}{205}{0}{115207}{"Audrey 18 Apr 2014"}{{0, 1, 115207, 776781, 941598, 87101, 96460, 270637, 766595, 987622, 202282}}

\adfDesignE{11}{9691}{11}{1066021}{210}{0}{128692}{"Audrey 18 Apr 2014"}{{0, 1, 128692, 994629, 455735, 171263, 193821, 452774, 749769, 513375, 604025}}

\adfDesignE{11}{9708}{11}{1067891}{141}{0}{392824}{"Audrey 18 Apr 2014"}{{0, 1, 392824, 445476, 501836, 546264, 388323, 772148, 913658, 371893, 1007032}}

\adfDesignE{11}{9717}{11}{1068881}{177}{0}{124041}{"Audrey 18 Apr 2014"}{{0, 1, 124041, 696567, 940493, 950992, 291512, 264643, 177972, 225559, 603744}}

\adfDesignE{11}{9723}{11}{1069541}{45}{0}{276587}{"Audrey 18 Apr 2014"}{{0, 1, 276587, 379003, 519810, 709086, 97230, 15106, 496076, 1035886, 749379}}

\adfDesignE{11}{9724}{11}{1069651}{453}{0}{41024}{"Audrey 18 Apr 2014"}{{0, 1, 41024, 407553, 809142, 831576, 194481, 931386, 175893, 1038437, 918762}}

\adfDesignE{11}{9771}{11}{1074821}{188}{0}{67978}{"Audrey 18 Apr 2014"}{{0, 1, 67978, 353005, 120244, 1007748, 977109, 127444, 330972, 661444, 653339}}

\adfDesignE{11}{9777}{11}{1075481}{23}{0}{226983}{"Audrey 18 Apr 2014"}{{0, 1, 226983, 364984, 861842, 1017153, 782167, 359643, 612726, 409181, 742725}}

\adfDesignE{11}{9784}{11}{1076251}{379}{0}{291498}{"Audrey 18 Apr 2014"}{{0, 1, 291498, 1067554, 489250, 300239, 489204, 882594, 889266, 978365, 1069535}}

\adfDesignE{11}{9787}{11}{1076581}{45}{0}{164178}{"Audrey 18 Apr 2014"}{{0, 1, 164178, 57187, 1060966, 776472, 587225, 520919, 921523, 798583, 495851}}

\adfDesignE{11}{9796}{11}{1077571}{146}{0}{175783}{"Audrey 18 Apr 2014"}{{0, 1, 175783, 314664, 862482, 821561, 391843, 999749, 1056790, 19167, 745815}}

\adfDesignE{11}{9841}{11}{1082521}{23}{0}{90566}{"Audrey 18 Apr 2014"}{{0, 1, 90566, 1021260, 838920, 892335, 688876, 893544, 848549, 440423, 780652}}

\adfDesignE{11}{9862}{11}{1084831}{135}{0}{93465}{"Audrey 18 Apr 2014"}{{0, 1, 93465, 647013, 250781, 387579, 394483, 202398, 930923, 932671, 500010}}

\adfDesignE{11}{9999}{11}{1099901}{133}{0}{51346}{"Audrey 18 Apr 2014"}{{0, 1, 51346, 1048920, 93954, 1096199, 199981, 648591, 850909, 505992, 1003612}}

\adfDesignE{11}{100008}{11}{11000891}{203}{0}{379427}{"Audrey 27 Apr 2014"}{{0, 1, 379427, 7188703, 5096859, 6288230, 7000566, 5621059, 5812350, 2903789, 3712580}}

\adfDesignE{11}{100021}{11}{11002321}{166}{0}{660444}{"Audrey 27 Apr 2014"}{{0, 1, 660444, 10263412, 10934680, 7332977, 8001687, 9341666, 7740386, 6068907, 5669766}}

\adfDesignE{11}{100029}{11}{11003201}{100}{0}{1624940}{"Audrey 27 Apr 2014"}{{0, 1, 1624940, 2862831, 6289561, 1050505, 4001163, 1375933, 2138624, 5713931, 7952114}}

\adfDesignE{11}{100054}{11}{11005951}{162}{0}{199181}{"Audrey 27 Apr 2014"}{{0, 1, 199181, 7623357, 2846853, 1425922, 8004327, 9806229, 10386381, 3156393, 575160}}

~

\adfDesignE{12}{1615}{1}{213181}{5}{0}{12833}{"Audrey 18 Apr 2014"}{{0, 1, 12833, 110157, 41570, 88948, 98610, 19714, 157096, 173432, 43216, 107147}}

\adfDesignE{12}{1986}{1}{262153}{474}{0}{38128}{"Audrey 18 Apr 2014"}{{0, 1, 38128, 105999, 179224, 172574, 123325, 159392, 67330, 156064, 59398, 249330}}

\adfDesignE{12}{2083}{1}{274957}{82}{0}{20911}{"Audrey 18 Apr 2014"}{{0, 1, 20911, 88291, 191803, 269731, 152000, 244037, 130744, 90333, 273730, 188161}}

\adfDesignE{12}{2284}{1}{301489}{106}{0}{3002}{"Audrey 18 Apr 2014"}{{0, 1, 3002, 268823, 222082, 97985, 199195, 130703, 133217, 143020, 25704, 283713}}

\adfDesignE{12}{2531}{1}{334093}{17}{0}{9963}{"Audrey 18 Apr 2014"}{{0, 1, 9963, 35748, 14186, 13779, 302047, 118610, 24489, 96017, 109112, 278327}}

\adfDesignE{12}{2589}{1}{341749}{99}{0}{64358}{"Audrey 18 Apr 2014"}{{0, 1, 64358, 296033, 268562, 157521, 94182, 104892, 71339, 179296, 318732, 153829}}

\adfDesignE{12}{2694}{1}{355609}{704}{0}{79931}{"Audrey 18 Apr 2014"}{{0, 1, 79931, 93467, 276905, 189395, 256615, 322054, 273982, 186195, 160286, 294823}}

\adfDesignE{12}{2778}{1}{366697}{410}{0}{9117}{"Audrey 18 Apr 2014"}{{0, 1, 9117, 246167, 118899, 45851, 355684, 69457, 320447, 40300, 351403, 276159}}

\adfDesignE{12}{3019}{1}{398509}{86}{0}{30504}{"Audrey 18 Apr 2014"}{{0, 1, 30504, 374010, 285388, 46447, 119793, 236651, 210078, 194592, 42813, 53759}}

\adfDesignE{12}{3044}{1}{401809}{11}{0}{78055}{"Audrey 18 Apr 2014"}{{0, 1, 78055, 354967, 209590, 295824, 186326, 199175, 212606, 249630, 347622, 277058}}

\adfDesignE{12}{3145}{1}{415141}{507}{0}{101311}{"Audrey 18 Apr 2014"}{{0, 1, 101311, 387778, 138705, 234546, 249248, 197662, 178465, 246783, 380929, 375418}}

\adfDesignE{12}{3254}{1}{429529}{134}{0}{58133}{"Audrey 18 Apr 2014"}{{0, 1, 58133, 341046, 257065, 216206, 255329, 236633, 90435, 252424, 165165, 275208}}

\adfDesignE{12}{3335}{1}{440221}{75}{0}{36297}{"Audrey 18 Apr 2014"}{{0, 1, 36297, 330977, 281300, 300447, 170147, 405471, 352636, 203317, 372526, 188207}}

\adfDesignE{12}{4379}{1}{578029}{39}{0}{3773}{"Audrey 18 Apr 2014"}{{0, 1, 3773, 362833, 196237, 525081, 225230, 90160, 292628, 50054, 416288, 149831}}

\adfDesignE{12}{4611}{1}{608653}{104}{0}{24374}{"Audrey 18 Apr 2014"}{{0, 1, 24374, 46548, 31760, 520277, 554996, 159579, 285876, 82080, 584162, 144959}}

\adfDesignE{12}{4945}{1}{652741}{85}{0}{93893}{"Audrey 18 Apr 2014"}{{0, 1, 93893, 628244, 162463, 234030, 558507, 644034, 358122, 501813, 577147, 158192}}

\adfDesignE{12}{6123}{1}{808237}{627401}{0}{719274}{"Audrey 19 May 2014"}{{0, 1, 719274, 158665, 544810, 424386, 523263, 171983, 611018, 793338, 759294, 143390}}

\adfDesignE{12}{6491}{1}{856813}{671134}{0}{688466}{"Audrey 19 May 2014"}{{0, 1, 688466, 765621, 381103, 610699, 504130, 268166, 535368, 562774, 782897, 41653}}

\adfDesignE{12}{6730}{1}{888361}{98014}{0}{276685}{"Audrey 19 May 2014"}{{0, 1, 276685, 80050, 17798, 254607, 687217, 512288, 854286, 133918, 452881, 283713}}

\adfDesignE{12}{7016}{1}{926113}{449554}{0}{900241}{"Audrey 19 May 2014"}{{0, 1, 900241, 706798, 749442, 464457, 810784, 781915, 307492, 803759, 88454, 869448}}

\adfDesignE{12}{7334}{1}{968089}{420873}{0}{541847}{"Audrey 19 May 2014"}{{0, 1, 541847, 11845, 715734, 899209, 313857, 215427, 174405, 818300, 893388, 256432}}

\adfDesignE{12}{7574}{1}{999769}{817381}{0}{775318}{"Audrey 19 May 2014"}{{0, 1, 775318, 891260, 580719, 923337, 202561, 476033, 229916, 200357, 280382, 438961}}

\adfDesignE{12}{7604}{1}{1003729}{773271}{0}{176550}{"Audrey 19 May 2014"}{{0, 1, 176550, 102134, 769944, 602188, 311991, 374717, 507960, 163037, 245817, 760577}}

\adfDesignE{12}{7698}{1}{1016137}{509059}{0}{781838}{"Audrey 19 May 2014"}{{0, 1, 781838, 236113, 506904, 8401, 927607, 85889, 886474, 463348, 71754, 96219}}

\adfDesignE{12}{8018}{1}{1058377}{6573}{0}{606701}{"Audrey 19 May 2014"}{{0, 1, 606701, 575210, 775623, 959868, 1040381, 25936, 506277, 822745, 1044866, 1051031}}

\adfDesignE{12}{9243}{1}{1220077}{852559}{0}{444699}{"Audrey 19 May 2014"}{{0, 1, 444699, 1020056, 575006, 855534, 943510, 793652, 916727, 812009, 721063, 238205}}

\adfDesignE{12}{9271}{1}{1223773}{827392}{0}{782220}{"Audrey 19 May 2014"}{{0, 1, 782220, 1208768, 1219516, 1199566, 240089, 989227, 252267, 1015355, 1087327, 572075}}

\adfDesignE{12}{9455}{1}{1248061}{422656}{0}{1138477}{"Audrey 19 May 2014"}{{0, 1, 1138477, 1058175, 794432, 210706, 370257, 220022, 399611, 992544, 326393, 729687}}

\adfDesignE{12}{10038}{1}{1325017}{206460}{0}{1112514}{"Audrey 19 May 2014"}{{0, 1, 1112514, 945649, 153790, 707935, 143824, 1135667, 651811, 184179, 1087126, 502589}}

\adfDesignE{12}{10954}{1}{1445929}{46222}{0}{933079}{"Audrey 19 May 2014"}{{0, 1, 933079, 637400, 700533, 183651, 953781, 1073347, 607279, 1395876, 103513, 641185}}

\adfDesignE{12}{11889}{1}{1569349}{571412}{0}{184504}{"Audrey 19 May 2014"}{{0, 1, 184504, 976857, 568674, 661603, 1295994, 647242, 695162, 414576, 860044, 1542088}}

\adfDesignE{12}{11951}{1}{1577533}{392712}{0}{351674}{"Audrey 19 May 2014"}{{0, 1, 351674, 747675, 967642, 156679, 1435555, 525811, 371953, 318028, 1399304, 35810}}

\adfDesignE{12}{12498}{1}{1649737}{1152999}{0}{271632}{"Audrey 19 May 2014"}{{0, 1, 271632, 1105836, 1280603, 758435, 1008571, 1532178, 1172821, 1600750, 341858, 825737}}

\adfDesignE{12}{12823}{1}{1692637}{820154}{0}{1186884}{"Audrey 19 May 2014"}{{0, 1, 1186884, 1564117, 422123, 638554, 1153164, 544865, 965803, 35527, 1147561, 804586}}

\adfDesignE{12}{13425}{1}{1772101}{1374524}{0}{1302466}{"Audrey 19 May 2014"}{{0, 1, 1302466, 1342765, 1560580, 771179, 492210, 899494, 1343791, 1726542, 1525592, 1440087}}

\adfDesignE{12}{14116}{1}{1863313}{807017}{0}{443391}{"Audrey 19 May 2014"}{{0, 1, 443391, 1150877, 1605727, 456209, 1432065, 1698092, 509497, 179520, 547586, 1293600}}

\adfDesignE{12}{14985}{1}{1978021}{1048579}{0}{850294}{"Audrey 19 May 2014"}{{0, 1, 850294, 1562600, 1050364, 187075, 257272, 1561515, 239160, 1908093, 1930449, 343282}}

\adfDesignE{12}{15649}{1}{2065669}{1317372}{0}{1402832}{"Audrey 19 May 2014"}{{0, 1, 1402832, 1617621, 1560046, 1165546, 1454674, 1225682, 237866, 1997590, 740818, 991338}}

\adfDesignE{12}{15984}{1}{2109889}{1720608}{0}{227112}{"Audrey 19 May 2014"}{{0, 1, 227112, 1514050, 1873714, 1531447, 1119081, 1605021, 336489, 510188, 1042843, 789499}}

\adfDesignE{12}{16579}{1}{2188429}{905838}{0}{326712}{"Audrey 19 May 2014"}{{0, 1, 326712, 106469, 1809402, 1774170, 5097, 2045024, 2130530, 482188, 355862, 1906690}}

\adfDesignE{12}{16585}{1}{2189221}{2104105}{0}{1023225}{"Audrey 19 May 2014"}{{0, 1, 1023225, 1025038, 1861776, 1607041, 1228147, 940329, 1133083, 546401, 1336582, 244482}}

\adfDesignE{12}{16859}{1}{2225389}{815441}{0}{1851597}{"Audrey 19 May 2014"}{{0, 1, 1851597, 1636288, 1224831, 1700596, 2060873, 590435, 849166, 226376, 653544, 333238}}

\adfDesignE{12}{17241}{1}{2275813}{739175}{0}{586039}{"Audrey 19 May 2014"}{{0, 1, 586039, 1045504, 363731, 1078490, 690563, 403232, 635193, 2241369, 933994, 1125136}}

\adfDesignE{12}{17321}{1}{2286373}{2203480}{0}{2123580}{"Audrey 19 May 2014"}{{0, 1, 2123580, 211406, 1355811, 763805, 1867, 153078, 1438846, 284226, 1613556, 1199316}}

\adfDesignE{12}{17526}{1}{2313433}{1394207}{0}{378938}{"Audrey 19 May 2014"}{{0, 1, 378938, 1534967, 119588, 911940, 1978778, 1960371, 1849100, 1668760, 487227, 677495}}

\adfDesignE{12}{17936}{1}{2367553}{2031603}{0}{1044618}{"Audrey 19 May 2014"}{{0, 1, 1044618, 280247, 764643, 1912893, 1684791, 2233887, 1263693, 2295617, 491772, 2233156}}

\adfDesignE{12}{18871}{1}{2490973}{1668676}{0}{1893678}{"Audrey 19 May 2014"}{{0, 1, 1893678, 1672992, 1499121, 131750, 1202766, 459122, 2433553, 962636, 1173105, 1026141}}

\adfDesignE{12}{19455}{1}{2568061}{1441491}{0}{2182924}{"Audrey 19 May 2014"}{{0, 1, 2182924, 1873470, 347658, 173333, 2362135, 293999, 1052749, 583250, 2276542, 1694244}}

\adfDesignE{12}{19698}{1}{2600137}{1034418}{0}{1699921}{"Audrey 19 May 2014"}{{0, 1, 1699921, 1547729, 802397, 441533, 71651, 221943, 487529, 918240, 414104, 1195363}}

\adfDesignE{12}{20163}{1}{2661517}{4096}{0}{444483}{"Audrey 19 May 2014"}{{0, 1, 444483, 730379, 2512982, 308297, 1911189, 1331812, 1166607, 1606622, 1878639, 1416574}}

\adfDesignE{12}{20295}{1}{2678941}{4096}{0}{324943}{"Audrey 19 May 2014"}{{0, 1, 324943, 172675, 1792221, 42295, 496855, 495959, 1351600, 2013378, 869021, 477875}}

\adfDesignE{12}{20410}{1}{2694121}{2396457}{0}{339411}{"Audrey 19 May 2014"}{{0, 1, 339411, 1907082, 485484, 781322, 1563070, 2232692, 563653, 496173, 2458635, 2643082}}

\adfDesignE{12}{20771}{1}{2741773}{786771}{0}{757115}{"Audrey 19 May 2014"}{{0, 1, 757115, 642115, 160503, 1107712, 878548, 2255674, 828951, 204254, 2286464, 1845755}}

\adfDesignE{12}{20873}{1}{2755237}{1148131}{0}{2053103}{"Audrey 19 May 2014"}{{0, 1, 2053103, 352783, 377852, 1789799, 1396056, 1263038, 699624, 1407114, 1323132, 358446}}

\adfDesignE{12}{20984}{1}{2769889}{702799}{0}{171167}{"Audrey 19 May 2014"}{{0, 1, 171167, 1025936, 964490, 705541, 945836, 1438340, 298793, 270935, 1649507, 839121}}

\adfDesignE{12}{21226}{1}{2801833}{1779094}{0}{1794063}{"Audrey 19 May 2014"}{{0, 1, 1794063, 352559, 2059300, 131102, 2574408, 2152850, 2126219, 1290782, 1234436, 293445}}

\adfDesignE{12}{21548}{1}{2844337}{582814}{0}{2088287}{"Audrey 19 May 2014"}{{0, 1, 2088287, 2261632, 445994, 2387650, 1293583, 472952, 466355, 2005444, 566705, 2233082}}

\adfDesignE{12}{21558}{1}{2845657}{413295}{0}{1230343}{"Audrey 19 May 2014"}{{0, 1, 1230343, 2347813, 251787, 460607, 543622, 1146723, 103674, 850814, 1047810, 553777}}

\adfDesignE{12}{23706}{1}{3129193}{2930009}{0}{2604102}{"Audrey 19 May 2014"}{{0, 1, 2604102, 1104665, 2598509, 2755594, 1234146, 1267049, 472429, 1688229, 1489517, 431724}}

\adfDesignE{12}{23725}{1}{3131701}{2988393}{0}{390423}{"Audrey 19 May 2014"}{{0, 1, 390423, 836156, 2890047, 1476385, 38197, 2958870, 1563934, 1796710, 937938, 2769844}}

\adfDesignE{12}{23890}{1}{3153481}{3107527}{0}{3111985}{"Audrey 19 May 2014"}{{0, 1, 3111985, 117390, 912705, 2853611, 2922975, 569103, 921121, 506185, 684181, 14667}}

\adfDesignE{12}{23974}{1}{3164569}{76867}{0}{930810}{"Audrey 19 May 2014"}{{0, 1, 930810, 2061573, 1413910, 611380, 503668, 972006, 2627760, 2436965, 990726, 109477}}

\adfDesignE{12}{24146}{1}{3187273}{715585}{0}{2086030}{"Audrey 19 May 2014"}{{0, 1, 2086030, 1079460, 3067484, 1942803, 231670, 301475, 1871347, 68381, 1601588, 498853}}

\adfDesignE{12}{24411}{1}{3222253}{1865495}{0}{1672410}{"Audrey 19 May 2014"}{{0, 1, 1672410, 937064, 2790931, 2441825, 541941, 1890729, 3106171, 823377, 553526, 1353290}}

\adfDesignE{12}{24533}{1}{3238357}{1524451}{0}{2935969}{"Audrey 19 May 2014"}{{0, 1, 2935969, 254292, 3176026, 908688, 1482863, 2363518, 2670559, 851841, 2134543, 2651842}}

\adfDesignE{12}{24534}{1}{3238489}{2284533}{0}{1922507}{"Audrey 19 May 2014"}{{0, 1, 1922507, 723662, 101701, 349521, 2486537, 2598024, 979957, 2647383, 3036781, 1346371}}

\adfDesignE{12}{24770}{1}{3269641}{1587623}{0}{1474247}{"Audrey 19 May 2014"}{{0, 1, 1474247, 1912207, 580775, 1260960, 2297647, 165424, 3122461, 242582, 2262097, 3029804}}

\adfDesignE{12}{24805}{1}{3274261}{4096}{0}{1309081}{"Audrey 19 May 2014"}{{0, 1, 1309081, 519598, 889298, 2890849, 652840, 2302169, 918720, 265627, 740587, 2608274}}

\adfDesignE{12}{24881}{1}{3284293}{251770}{0}{2374827}{"Audrey 19 May 2014"}{{0, 1, 2374827, 2203157, 1510743, 1446140, 2960661, 530438, 1335490, 1351748, 1140606, 1567654}}

\adfDesignE{12}{25744}{1}{3398209}{1667600}{0}{986161}{"Audrey 19 May 2014"}{{0, 1, 986161, 473465, 2199474, 2851331, 3003987, 2340694, 3106513, 2571003, 2163538, 693087}}

\adfDesignE{12}{25868}{1}{3414577}{2750714}{0}{2937721}{"Audrey 19 May 2014"}{{0, 1, 2937721, 1303998, 2518028, 2656659, 1843243, 1854447, 2972828, 1991437, 2741975, 3081702}}

\adfDesignE{12}{25910}{1}{3420121}{3027327}{0}{2761610}{"Audrey 19 May 2014"}{{0, 1, 2761610, 3015652, 1942663, 524168, 1897956, 173877, 2313812, 91810, 3004128, 1374928}}

\adfDesignE{12}{27136}{1}{3581953}{1377802}{0}{950519}{"Audrey 19 May 2014"}{{0, 1, 950519, 3200265, 3577439, 536928, 20239, 2466431, 1289189, 1771932, 2822343, 1274479}}

\adfDesignE{12}{27185}{1}{3588421}{136638}{0}{710303}{"Audrey 19 May 2014"}{{0, 1, 710303, 1947630, 2956391, 192536, 585677, 2483601, 2004872, 1722366, 1365368, 384939}}

\adfDesignE{12}{27769}{1}{3665509}{419865}{0}{2171501}{"Audrey 19 May 2014"}{{0, 1, 2171501, 3181149, 259118, 1037073, 633680, 1009571, 2151315, 1418585, 1637575, 1162468}}

\adfDesignE{12}{28741}{1}{3793813}{855783}{0}{1016192}{"Audrey 19 May 2014"}{{0, 1, 1016192, 632768, 414086, 111617, 675203, 2043048, 2807096, 3266610, 38006, 376812}}

\adfDesignE{12}{28889}{1}{3813349}{452683}{0}{2073413}{"Audrey 19 May 2014"}{{0, 1, 2073413, 1459835, 2724803, 1445132, 1245370, 2832648, 1511455, 3037480, 2538243, 198365}}

\adfDesignE{12}{29215}{1}{3856381}{1838133}{0}{3226855}{"Audrey 19 May 2014"}{{0, 1, 3226855, 1991211, 1404445, 2403276, 2953982, 147764, 2342618, 2857428, 2980227, 2830479}}

\adfDesignE{12}{29733}{1}{3924757}{805691}{0}{2664462}{"Audrey 19 May 2014"}{{0, 1, 2664462, 2178639, 760639, 218259, 3717454, 3635166, 3011158, 2215372, 2777976, 2369416}}

\adfDesignE{12}{29943}{1}{3952477}{3007710}{0}{304791}{"Audrey 19 May 2014"}{{0, 1, 304791, 2486750, 172299, 2575087, 221542, 3843131, 3561855, 2394669, 1251405, 2950855}}

\adfDesignE{12}{30475}{1}{4022701}{855890}{0}{3336224}{"Audrey 19 May 2014"}{{0, 1, 3336224, 3317482, 648917, 3290630, 1913339, 954910, 3534787, 3608316, 470430, 3061170}}

\adfDesignE{12}{30704}{1}{4052929}{2195063}{0}{2200772}{"Audrey 19 May 2014"}{{0, 1, 2200772, 1335540, 2502648, 1409203, 810484, 493677, 1841614, 2996718, 827407, 1793652}}

\adfDesignE{12}{30805}{1}{4066261}{375114}{0}{3301030}{"Audrey 19 May 2014"}{{0, 1, 3301030, 303012, 291492, 98764, 2318723, 3062469, 3854269, 3433818, 2671374, 996353}}

\adfDesignE{12}{31315}{1}{4133581}{3071647}{0}{1835664}{"Audrey 19 May 2014"}{{0, 1, 1835664, 158344, 1231658, 2653571, 3302772, 2671355, 1728610, 3692390, 1688763, 1704777}}

\adfDesignE{12}{31580}{1}{4168561}{965065}{0}{2081985}{"Audrey 19 May 2014"}{{0, 1, 2081985, 58619, 918318, 1292897, 2254649, 3890263, 1043226, 2197292, 75024, 2861970}}

\adfDesignE{12}{31723}{1}{4187437}{3761973}{0}{189}{"Audrey 19 May 2014"}{{0, 1, 189, 35721, 2563832, 3008993, 3395682, 1106037, 3856580, 279582, 2591754, 4098814}}

\adfDesignE{12}{32048}{1}{4230337}{1068745}{0}{905644}{"Audrey 19 May 2014"}{{0, 1, 905644, 626165, 1670073, 2283054, 1953645, 2264826, 649787, 1978432, 64195, 295189}}

\adfDesignE{12}{32101}{1}{4237333}{3262001}{0}{4202332}{"Audrey 19 May 2014"}{{0, 1, 4202332, 480764, 3465912, 220545, 1125181, 3550054, 112838, 3988851, 2111166, 1929021}}

\adfDesignE{12}{32433}{1}{4281157}{3252793}{0}{2134177}{"Audrey 19 May 2014"}{{0, 1, 2134177, 1378500, 3277984, 2217038, 3947855, 1616567, 1243397, 1205546, 1623352, 2761368}}

\adfDesignE{12}{32784}{1}{4327489}{4255600}{0}{67617}{"Audrey 19 May 2014"}{{0, 1, 67617, 2230305, 2196513, 2197041, 3278905, 3802937, 3794749, 4065345, 4096, 4327425}}

\adfDesignE{12}{33686}{1}{4446553}{265479}{0}{1685566}{"Audrey 19 May 2014"}{{0, 1, 1685566, 3254453, 2946229, 1857412, 3673763, 2472892, 1661907, 1135763, 2354450, 1190329}}

\adfDesignE{12}{34689}{1}{4578949}{513573}{0}{3282019}{"Audrey 19 May 2014"}{{0, 1, 3282019, 878189, 4377643, 1855447, 3927107, 406986, 1414046, 764159, 2410741, 3578407}}

\adfDesignE{12}{34943}{1}{4612477}{3164696}{0}{4382263}{"Audrey 19 May 2014"}{{0, 1, 4382263, 1125066, 3089334, 3728585, 277556, 3979374, 4125796, 3713404, 3511201, 4354759}}

\adfDesignE{12}{35261}{1}{4654453}{2112412}{0}{44784}{"Audrey 19 May 2014"}{{0, 1, 44784, 4191866, 474095, 2910347, 2987142, 2533655, 950286, 1944445, 4518156, 2717488}}

\adfDesignE{12}{35329}{1}{4663429}{4468270}{0}{213010}{"Audrey 19 May 2014"}{{0, 1, 213010, 2759359, 1796288, 2284288, 3331878, 736699, 4531569, 334267, 979698, 1686659}}

\adfDesignE{12}{35434}{1}{4677289}{2164794}{0}{1411254}{"Audrey 19 May 2014"}{{0, 1, 1411254, 1423426, 4201906, 787433, 232050, 1101365, 531409, 1235915, 855576, 2251532}}

\adfDesignE{12}{35555}{1}{4693261}{432014}{0}{2069870}{"Audrey 19 May 2014"}{{0, 1, 2069870, 1181525, 1163782, 227697, 1226509, 2590883, 3460733, 4162803, 4617229, 2765273}}

\adfDesignE{12}{36455}{1}{4812061}{2340723}{0}{1504129}{"Audrey 19 May 2014"}{{0, 1, 1504129, 3945369, 4438059, 1924886, 4112424, 506856, 1984194, 196277, 1173322, 4274788}}

\adfDesignE{12}{36635}{1}{4835821}{2985448}{0}{1826890}{"Audrey 19 May 2014"}{{0, 1, 1826890, 2999993, 3660525, 1865307, 4362950, 2269713, 2579552, 1830570, 4164003, 3455422}}

\adfDesignE{12}{36743}{1}{4850077}{1463521}{0}{4080912}{"Audrey 19 May 2014"}{{0, 1, 4080912, 2404765, 3094111, 2395815, 11521, 4390791, 1657741, 4036958, 396408, 1781362}}

\adfDesignE{12}{36776}{1}{4854433}{2393741}{0}{2431642}{"Audrey 19 May 2014"}{{0, 1, 2431642, 3808143, 772687, 4432703, 2590590, 3339732, 3089914, 4331079, 4315847, 14260}}

\adfDesignE{12}{36825}{1}{4860901}{3177641}{0}{4025792}{"Audrey 19 May 2014"}{{0, 1, 4025792, 3853609, 253174, 1863930, 2165256, 1233690, 2344740, 3273170, 1038305, 4252838}}

\adfDesignE{12}{36944}{1}{4876609}{2623780}{0}{1195015}{"Audrey 19 May 2014"}{{0, 1, 1195015, 4423883, 1518179, 1831415, 4796333, 1628108, 2542108, 2874724, 1337592, 2235687}}

\adfDesignE{12}{37665}{1}{4971781}{2492901}{0}{4091049}{"Audrey 19 May 2014"}{{0, 1, 4091049, 1527766, 461566, 1727153, 2226983, 2356506, 4820715, 3800752, 1547245, 2299169}}

\adfDesignE{12}{37728}{1}{4980097}{1036442}{0}{2311640}{"Audrey 19 May 2014"}{{0, 1, 2311640, 547921, 3050333, 2234790, 1034105, 2041618, 2289627, 1047747, 3460294, 1902312}}

\adfDesignE{12}{38563}{1}{5090317}{2589624}{0}{822943}{"Audrey 19 May 2014"}{{0, 1, 822943, 4136618, 429854, 3941041, 2439966, 1044533, 4559780, 4959333, 186880, 2930636}}

\adfDesignE{12}{39186}{1}{5172553}{925937}{0}{3483925}{"Audrey 19 May 2014"}{{0, 1, 3483925, 1575180, 4393809, 1119148, 3797477, 1758710, 4009964, 5089931, 3723100, 2084073}}

\adfDesignE{12}{39465}{1}{5209381}{497505}{0}{4281983}{"Audrey 19 May 2014"}{{0, 1, 4281983, 3456685, 77245, 2549002, 4809670, 2048780, 3895214, 324992, 2225701, 2377632}}

\adfDesignE{12}{39569}{1}{5223109}{4678369}{0}{26145}{"Audrey 19 May 2014"}{{0, 1, 26145, 4556855, 5080794, 3251042, 2840333, 3565632, 1399208, 4860833, 3013706, 2744105}}

\adfDesignE{12}{39755}{1}{5247661}{4627135}{0}{65257}{"Audrey 19 May 2014"}{{0, 1, 65257, 2622978, 4716509, 4662502, 1508234, 2944083, 5155121, 1175031, 175435, 3213154}}

\adfDesignE{12}{40508}{1}{5347057}{2316224}{0}{3594280}{"Audrey 19 May 2014"}{{0, 1, 3594280, 753581, 2658045, 134876, 1880489, 3131500, 2346741, 847462, 532626, 161570}}

\adfDesignE{12}{40808}{1}{5386657}{1357821}{0}{2919281}{"Audrey 19 May 2014"}{{0, 1, 2919281, 3837203, 2766780, 1606501, 4955272, 1732531, 2493288, 1134475, 5309107, 178846}}

\adfDesignE{12}{41449}{1}{5471269}{4661259}{0}{4535006}{"Audrey 19 May 2014"}{{0, 1, 4535006, 3571065, 3628191, 322528, 4717053, 672592, 3644397, 4544556, 2312961, 4879264}}

\adfDesignE{12}{41831}{1}{5521693}{4169863}{0}{124299}{"Audrey 19 May 2014"}{{0, 1, 124299, 544387, 3933691, 2420766, 5176385, 4202290, 4851989, 1506172, 2672263, 2176222}}

\adfDesignE{12}{43094}{1}{5688409}{2164275}{0}{2485857}{"Audrey 19 May 2014"}{{0, 1, 2485857, 1363888, 2246200, 1095818, 4297742, 439724, 621619, 5326042, 1152085, 3724660}}

\adfDesignE{12}{43178}{1}{5699497}{2128864}{0}{3228948}{"Audrey 19 May 2014"}{{0, 1, 3228948, 3925610, 422232, 5593057, 1842974, 5295161, 3859262, 4573061, 4834683, 621993}}

\adfDesignE{12}{43615}{1}{5757181}{689024}{0}{2442731}{"Audrey 19 May 2014"}{{0, 1, 2442731, 2362988, 126809, 908855, 4145785, 715948, 3744437, 109869, 3662743, 4808558}}

\adfDesignE{12}{43775}{1}{5778301}{1087854}{0}{3359285}{"Audrey 19 May 2014"}{{0, 1, 3359285, 4990265, 4414375, 950525, 1020326, 4796332, 2742715, 2631265, 3774708, 211708}}

\adfDesignE{12}{44488}{1}{5872417}{1240231}{0}{1501659}{"Audrey 19 May 2014"}{{0, 1, 1501659, 986366, 3258535, 5199481, 5116536, 3975351, 4361125, 3167975, 2520910, 5146563}}

\adfDesignE{12}{44735}{1}{5905021}{873265}{0}{2767152}{"Audrey 19 May 2014"}{{0, 1, 2767152, 885089, 3381547, 4736998, 4253812, 5432486, 2154815, 5296773, 5275976, 1245477}}

\adfDesignE{12}{44763}{1}{5908717}{5357795}{0}{1801547}{"Audrey 19 May 2014"}{{0, 1, 1801547, 1975864, 5753147, 1506671, 4022011, 5828785, 87203, 5244162, 3080672, 243522}}

\adfDesignE{12}{44846}{1}{5919673}{1586339}{0}{1976101}{"Audrey 19 May 2014"}{{0, 1, 1976101, 3671021, 3553560, 5052002, 4074987, 1426730, 5839693, 791447, 1601547, 1611276}}

\adfDesignE{12}{45065}{1}{5948581}{3868244}{0}{4170805}{"Audrey 19 May 2014"}{{0, 1, 4170805, 472295, 1543649, 2676687, 5565676, 5349707, 4642587, 3882139, 5824822, 1563118}}

\adfDesignE{12}{45093}{1}{5952277}{4955202}{0}{3613417}{"Audrey 19 May 2014"}{{0, 1, 3613417, 4491060, 4773469, 1933865, 561968, 5426106, 1083033, 342294, 3511660, 4024512}}

\adfDesignE{12}{45119}{1}{5955709}{3125979}{0}{292658}{"Audrey 19 May 2014"}{{0, 1, 292658, 5609544, 4609229, 1342145, 4507151, 1235165, 5116524, 1368303, 813341, 4884484}}

\adfDesignE{12}{45919}{1}{6061309}{3387213}{0}{1869628}{"Audrey 19 May 2014"}{{0, 1, 1869628, 448556, 2266546, 3394190, 3387697, 4695020, 1649232, 1824306, 2268160, 2441900}}

\adfDesignE{12}{45959}{1}{6066589}{5858042}{0}{935901}{"Audrey 19 May 2014"}{{0, 1, 935901, 4428803, 3087910, 655446, 3353522, 2640994, 4635913, 3958881, 6054332, 581242}}

\adfDesignE{12}{46225}{1}{6101701}{5592201}{0}{1197581}{"Audrey 19 May 2014"}{{0, 1, 1197581, 1533212, 3390149, 3811385, 3809625, 1153910, 5754333, 453970, 4287470, 5116869}}

\adfDesignE{12}{46259}{1}{6106189}{2050508}{0}{4736074}{"Audrey 19 May 2014"}{{0, 1, 4736074, 1641333, 2852570, 2810546, 5267225, 5400177, 5700945, 1223479, 711329, 187266}}

\adfDesignE{12}{47151}{1}{6223933}{2751603}{0}{5336846}{"Audrey 19 May 2014"}{{0, 1, 5336846, 376714, 3542451, 1541463, 4963618, 3967615, 1976129, 3502892, 3281842, 2630094}}

\adfDesignE{12}{47203}{1}{6230797}{3439215}{0}{1462123}{"Audrey 19 May 2014"}{{0, 1, 1462123, 4754835, 4784421, 2731928, 345572, 979032, 1902156, 1257471, 6153767, 551882}}

\adfDesignE{12}{47334}{1}{6248089}{4317146}{0}{3937226}{"Audrey 19 May 2014"}{{0, 1, 3937226, 2338694, 1861319, 6236371, 5722997, 5550151, 1142657, 6101655, 4908480, 5937072}}

\adfDesignE{12}{47480}{1}{6267361}{2064834}{0}{5269509}{"Audrey 19 May 2014"}{{0, 1, 5269509, 437112, 3905371, 133098, 5808816, 5088174, 950101, 3515418, 139008, 6090197}}

\adfDesignE{12}{47616}{1}{6285313}{865017}{0}{3264415}{"Audrey 19 May 2014"}{{0, 1, 3264415, 2792940, 3610438, 438690, 1646491, 1073632, 1882098, 5202666, 3537456, 1692425}}

\adfDesignE{12}{48564}{1}{6410449}{2796454}{0}{5092272}{"Audrey 19 May 2014"}{{0, 1, 5092272, 6349634, 2269510, 6045601, 3126769, 1667131, 5102401, 1089619, 1677479, 6042277}}

\adfDesignE{12}{48629}{1}{6419029}{5174740}{0}{4999882}{"Audrey 19 May 2014"}{{0, 1, 4999882, 1439830, 562386, 1565973, 5402146, 6035137, 3350836, 3368830, 2135074, 3235050}}

\adfDesignE{12}{49110}{1}{6482521}{3603038}{0}{504483}{"Audrey 19 May 2014"}{{0, 1, 504483, 5805350, 1116586, 6475264, 1591234, 479029, 86648, 803881, 4267284, 4800324}}

\adfDesignE{12}{49298}{1}{6507337}{5696619}{0}{1022126}{"Audrey 19 May 2014"}{{0, 1, 1022126, 1619200, 2385316, 2562700, 3937590, 3269884, 96814, 5540142, 5987807, 6115105}}

\adfDesignE{12}{49326}{1}{6511033}{2997143}{0}{5999618}{"Audrey 19 May 2014"}{{0, 1, 5999618, 3617648, 3597096, 5390914, 4975045, 3726735, 1399735, 3537127, 960586, 5861693}}

\adfDesignE{12}{49594}{1}{6546409}{2770700}{0}{4562881}{"Audrey 19 May 2014"}{{0, 1, 4562881, 4610602, 3198645, 930197, 849589, 6127606, 906929, 1395643, 2644553, 958990}}

~

\adfDesignE{13}{100}{13}{15613}{17}{0}{1557}{"Audrey 9 Apr 2014"}{{0, 1, 1557, 4234, 3652, 3032, 5698, 3602, 3247, 12580, 8358, 7777, 8714}}

\adfDesignE{13}{115}{13}{17953}{27}{0}{4027}{"Audrey 10 Apr 2014"}{{0, 1, 4027, 5170, 12063, 14836, 14941, 6904, 11164, 3116, 16938, 5879, 12679}}

\adfDesignE{13}{119}{13}{18577}{66}{0}{10131}{"Audrey 10 Apr 2014"}{{0, 1, 10131, 17813, 6525, 7809, 12113, 15718, 15591, 10767, 14910, 3623, 15038}}

\adfDesignE{13}{124}{13}{19357}{29}{0}{4489}{"Audrey 10 Apr 2014"}{{0, 1, 4489, 484, 4692, 1972, 6159, 5955, 19335, 17384, 8709, 12918, 14687}}

\adfDesignE{13}{148}{13}{23101}{14}{0}{1200}{"Audrey 10 Apr 2014"}{{0, 1, 1200, 7738, 22099, 21953, 8460, 10661, 18347, 1147, 13441, 4702, 5756}}

\adfDesignE{13}{149}{13}{23257}{222}{0}{2708}{"Audrey 10 Apr 2014"}{{0, 1, 2708, 7309, 1065, 152, 16247, 17889, 22338, 23104, 4302, 21316, 23111}}

\adfDesignE{13}{166}{13}{25909}{22}{0}{875}{"Audrey 10 Apr 2014"}{{0, 1, 875, 14264, 18771, 24228, 5938, 13950, 3111, 1680, 19096, 23604, 4027}}

\adfDesignE{13}{176}{13}{27469}{55}{0}{13725}{"Audrey 10 Apr 2014"}{{0, 1, 13725, 20692, 23178, 27030, 17905, 8451, 15857, 438, 23308, 25795, 15903}}

\adfDesignE{13}{191}{13}{29809}{92}{0}{4099}{"Audrey 10 Apr 2014"}{{0, 1, 4099, 19334, 17744, 28505, 20524, 6878, 23417, 1303, 5186, 3597, 18457}}

\adfDesignE{13}{203}{13}{31681}{144}{0}{2149}{"Audrey 10 Apr 2014"}{{0, 1, 2149, 24456, 28846, 22018, 16949, 21932, 22221, 9662, 12583, 16974, 12295}}

\adfDesignE{13}{224}{13}{34957}{53}{0}{5230}{"Audrey 10 Apr 2014"}{{0, 1, 5230, 16526, 17276, 24592, 9357, 32267, 18971, 10364, 20370, 21121, 33667}}

\adfDesignE{13}{229}{13}{35737}{196}{0}{4060}{"Audrey 10 Apr 2014"}{{0, 1, 4060, 8843, 22632, 6093, 7576, 24740, 23430, 29643, 24101, 2154, 25412}}

\adfDesignE{13}{233}{13}{36361}{27}{0}{12971}{"Audrey 10 Apr 2014"}{{0, 1, 12971, 4494, 4991, 15681, 31178, 2796, 14999, 20679, 28573, 29071, 16371}}

\adfDesignE{13}{236}{13}{36829}{14}{0}{20329}{"Audrey 10 Apr 2014"}{{0, 1, 20329, 10032, 18355, 24196, 29189, 31162, 33498, 12632, 24140, 32464, 21805}}

\adfDesignE{13}{260}{13}{40573}{17}{0}{10534}{"Audrey 10 Apr 2014"}{{0, 1, 10534, 38574, 40494, 19847, 36202, 6241, 14434, 20725, 34410, 36331, 26218}}

\adfDesignE{13}{275}{13}{42913}{92}{0}{23397}{"Audrey 10 Apr 2014"}{{0, 1, 23397, 21381, 14416, 37885, 27330, 36310, 39322, 5027, 35099, 28135, 32088}}

\adfDesignE{13}{314}{13}{48997}{137}{0}{5691}{"Audrey 10 Apr 2014"}{{0, 1, 5691, 464, 43783, 19308, 30554, 41458, 16923, 29688, 12752, 7075, 37288}}

\adfDesignE{13}{316}{13}{49309}{14}{0}{2604}{"Audrey 10 Apr 2014"}{{0, 1, 2604, 25483, 37127, 33068, 15558, 30343, 20154, 16240, 31147, 42792, 41337}}

\adfDesignE{13}{344}{13}{53677}{14}{0}{8122}{"Audrey 10 Apr 2014"}{{0, 1, 8122, 51528, 44524, 1979, 24015, 41289, 29039, 51697, 21540, 14537, 33791}}

\adfDesignE{13}{346}{13}{53989}{10}{0}{11008}{"Audrey 10 Apr 2014"}{{0, 1, 11008, 24748, 51479, 12288, 23859, 37376, 38828, 41700, 19122, 45854, 17671}}

\adfDesignE{13}{356}{13}{55549}{14}{0}{11015}{"Audrey 10 Apr 2014"}{{0, 1, 11015, 11209, 37257, 45392, 51880, 25637, 35988, 10156, 48203, 18703, 37853}}

\adfDesignE{13}{385}{13}{60073}{131}{0}{18594}{"Audrey 10 Apr 2014"}{{0, 1, 18594, 16721, 32499, 12099, 55494, 41588, 27616, 47973, 46058, 1764, 60031}}

\adfDesignE{13}{388}{13}{60541}{29}{0}{75}{"Audrey 10 Apr 2014"}{{0, 1, 75, 5625, 58629, 38223, 21298, 23284, 51152, 22317, 39168, 31632, 11301}}

\adfDesignE{13}{409}{13}{63817}{27}{0}{6848}{"Audrey 10 Apr 2014"}{{0, 1, 6848, 53426, 62204, 58334, 40629, 49089, 37333, 5482, 16340, 25119, 28097}}

\adfDesignE{13}{425}{13}{66313}{27}{0}{218}{"Audrey 10 Apr 2014"}{{0, 1, 218, 47524, 15404, 42422, 30489, 15302, 20186, 23890, 35606, 3487, 30723}}

\adfDesignE{13}{440}{13}{68653}{14}{0}{9409}{"Audrey 10 Apr 2014"}{{0, 1, 9409, 35564, 6954, 3877, 23950, 26404, 48682, 64775, 35294, 6685, 13017}}

\adfDesignE{13}{454}{13}{70837}{29}{0}{33297}{"Audrey 10 Apr 2014"}{{0, 1, 33297, 20322, 26610, 3974, 69599, 5448, 59336, 66862, 38778, 45067, 55728}}

\adfDesignE{13}{470}{13}{73333}{51}{0}{6803}{"Audrey 10 Apr 2014"}{{0, 1, 6803, 7686, 1429, 41531, 56677, 62050, 21402, 31801, 9853, 3597, 50502}}

\adfDesignE{13}{475}{13}{74113}{27}{0}{29260}{"Audrey 10 Apr 2014"}{{0, 1, 29260, 68337, 45993, 11326, 39537, 22803, 50554, 62786, 5316, 57086, 51679}}

\adfDesignE{13}{485}{13}{75673}{27}{0}{1798}{"Audrey 10 Apr 2014"}{{0, 1, 1798, 54538, 62789, 66179, 31886, 46567, 33128, 9493, 41989, 50241, 55429}}

\adfDesignE{13}{496}{13}{77389}{17}{0}{40375}{"Audrey 10 Apr 2014"}{{0, 1, 40375, 18729, 15456, 48493, 40564, 65482, 72732, 28895, 73839, 70567, 66590}}

\adfDesignE{13}{503}{13}{78481}{14}{0}{20531}{"Audrey 10 Apr 2014"}{{0, 1, 20531, 510, 32837, 24657, 30417, 18110, 51913, 53823, 27533, 59861, 72212}}

\adfDesignE{13}{511}{13}{79729}{170}{0}{6205}{"Audrey 10 Apr 2014"}{{0, 1, 6205, 72647, 66598, 5183, 29728, 49063, 30593, 74545, 43796, 37748, 62267}}

\adfDesignE{13}{518}{13}{80821}{29}{0}{15077}{"Audrey 10 Apr 2014"}{{0, 1, 15077, 47277, 34930, 9974, 50938, 31084, 53310, 70846, 14806, 2460, 73402}}

\adfDesignE{13}{523}{13}{81601}{144}{0}{1096}{"Audrey 10 Apr 2014"}{{0, 1, 1096, 58802, 63803, 77632, 56430, 75323, 55397, 3968, 24075, 29077, 44002}}

\adfDesignE{13}{530}{13}{82693}{14}{0}{35279}{"Audrey 10 Apr 2014"}{{0, 1, 35279, 78191, 27195, 8219, 36443, 44526, 79219, 74473, 10971, 42669, 58972}}

\adfDesignE{13}{535}{13}{83473}{209}{0}{38172}{"Audrey 10 Apr 2014"}{{0, 1, 38172, 80369, 45772, 35421, 78231, 70630, 77406, 48051, 50543, 15947, 43768}}

\adfDesignE{13}{544}{13}{84877}{14}{0}{14024}{"Audrey 10 Apr 2014"}{{0, 1, 14024, 12567, 34956, 58269, 53577, 32644, 57795, 26607, 17276, 39666, 77003}}

\adfDesignE{13}{566}{13}{88309}{40}{0}{634}{"Audrey 10 Apr 2014"}{{0, 1, 634, 48720, 68639, 69098, 6868, 27171, 6159, 19210, 80807, 12418, 13511}}

\adfDesignE{13}{569}{13}{88777}{105}{0}{39446}{"Audrey 10 Apr 2014"}{{0, 1, 39446, 81214, 49399, 26581, 57756, 47802, 62989, 62195, 80352, 48538, 65166}}

\adfDesignE{13}{583}{13}{90961}{27}{0}{30183}{"Audrey 10 Apr 2014"}{{0, 1, 30183, 39074, 61177, 88052, 65979, 34984, 46784, 2908, 85760, 16903, 73961}}

\adfDesignE{13}{608}{13}{94861}{29}{0}{12568}{"Audrey 10 Apr 2014"}{{0, 1, 12568, 11059, 18147, 25652, 56658, 51078, 23917, 69208, 25635, 32724, 52797}}

\adfDesignE{13}{610}{13}{95173}{110}{0}{13979}{"Audrey 10 Apr 2014"}{{0, 1, 13979, 22272, 29405, 308, 22747, 7320, 15305, 94864, 58447, 65581, 50463}}

\adfDesignE{13}{635}{13}{99073}{27}{0}{10027}{"Audrey 10 Apr 2014"}{{0, 1, 10027, 80707, 20825, 65464, 48903, 38104, 43320, 33608, 40143, 79335, 34928}}

\adfDesignE{13}{646}{13}{100789}{95}{0}{37281}{"Audrey 10 Apr 2014"}{{0, 1, 37281, 93440, 67222, 85686, 53200, 23258, 94520, 15102, 10308, 84880, 39836}}

\adfDesignE{13}{656}{13}{102349}{92}{0}{12978}{"Audrey 10 Apr 2014"}{{0, 1, 12978, 64379, 35775, 32886, 101527, 78729, 97244, 69462, 90193, 61590, 71679}}

\adfDesignE{13}{661}{13}{103129}{40}{0}{36397}{"Audrey 10 Apr 2014"}{{0, 1, 36397, 49604, 60514, 2005, 63782, 39664, 50866, 101123, 2950, 13861, 94878}}

\adfDesignE{13}{691}{13}{107809}{144}{0}{13888}{"Audrey 10 Apr 2014"}{{0, 1, 13888, 6243, 24348, 56000, 101683, 91222, 27577, 51808, 100047, 10344, 55884}}

\adfDesignE{13}{724}{13}{112957}{22}{0}{16312}{"Audrey 10 Apr 2014"}{{0, 1, 16312, 67609, 38817, 58919, 48572, 26066, 18444, 54037, 48073, 19282, 55696}}

\adfDesignE{13}{730}{13}{113893}{29}{0}{12022}{"Audrey 10 Apr 2014"}{{0, 1, 12022, 112160, 8293, 42071, 92642, 96370, 40544, 71821, 9229, 19256, 65056}}

\adfDesignE{13}{751}{13}{117169}{27}{0}{28487}{"Audrey 10 Apr 2014"}{{0, 1, 28487, 113844, 70446, 41739, 105050, 63090, 106708, 75429, 100801, 57404, 57184}}

\adfDesignE{13}{754}{13}{117637}{17}{0}{10716}{"Audrey 10 Apr 2014"}{{0, 1, 10716, 18944, 80079, 82286, 87461, 18097, 61676, 35350, 19460, 80596, 93519}}

\adfDesignE{13}{761}{13}{118729}{27}{0}{40544}{"Audrey 10 Apr 2014"}{{0, 1, 40544, 12931, 85929, 40329, 81917, 36531, 87318, 78399, 114997, 69267, 64211}}

\adfDesignE{13}{778}{13}{121381}{30}{0}{6985}{"Audrey 10 Apr 2014"}{{0, 1, 6985, 116444, 108640, 97769, 26959, 46684, 58374, 23611, 87437, 79634, 75748}}

\adfDesignE{13}{805}{13}{125593}{144}{0}{2097}{"Audrey 10 Apr 2014"}{{0, 1, 2097, 1654, 77427, 98263, 84991, 9660, 36547, 27329, 38505, 114279, 11619}}

\adfDesignE{13}{896}{13}{139789}{22}{0}{14642}{"Audrey 10 Apr 2014"}{{0, 1, 14642, 91627, 47501, 59367, 43612, 10752, 28370, 80421, 81535, 37410, 63918}}

\adfDesignE{13}{899}{13}{140257}{144}{0}{74897}{"Audrey 10 Apr 2014"}{{0, 1, 74897, 122151, 59851, 46627, 103633, 118678, 119305, 93629, 101984, 39685, 101358}}

\adfDesignE{13}{944}{13}{147277}{17}{0}{20368}{"Audrey 10 Apr 2014"}{{0, 1, 20368, 123392, 113528, 89404, 47844, 101960, 115580, 57872, 79065, 69202, 65446}}

\adfDesignE{13}{1009}{13}{157417}{27}{0}{41077}{"Audrey 10 Apr 2014"}{{0, 1, 41077, 124523, 80690, 88195, 148594, 108980, 104231, 69221, 125163, 81331, 129913}}

\adfDesignE{13}{1021}{13}{159289}{144}{0}{48448}{"Audrey 10 Apr 2014"}{{0, 1, 48448, 85289, 124812, 122047, 125376, 49011, 123094, 37241, 144754, 24989, 70672}}

\adfDesignE{13}{1043}{13}{162721}{27}{0}{45426}{"Audrey 10 Apr 2014"}{{0, 1, 45426, 56475, 136785, 94025, 78842, 150203, 67227, 68695, 38453, 118764, 121430}}

\adfDesignE{13}{1051}{13}{163969}{53}{0}{25152}{"Audrey 10 Apr 2014"}{{0, 1, 25152, 30702, 86683, 118992, 124596, 63064, 113591, 44976, 14221, 70203, 127664}}

\adfDesignE{13}{1084}{13}{169117}{35}{0}{5834}{"Audrey 10 Apr 2014"}{{0, 1, 5834, 43039, 119898, 17020, 23001, 78053, 98238, 152096, 140282, 48025, 120098}}

\adfDesignE{13}{1109}{13}{173017}{105}{0}{77022}{"Audrey 10 Apr 2014"}{{0, 1, 77022, 154605, 91285, 61441, 120735, 106471, 122613, 111575, 148277, 84958, 132136}}

\adfDesignE{13}{1120}{13}{174733}{29}{0}{17658}{"Audrey 21 Apr 2014"}{{0, 1, 17658, 81292, 22541, 161937, 152734, 147850, 49547, 12795, 4341, 120324, 102645}}

\adfDesignE{13}{1133}{13}{176761}{27}{0}{34213}{"Audrey 21 Apr 2014"}{{0, 1, 34213, 18027, 38622, 86011, 153976, 149566, 47369, 90749, 165333, 9168, 90770}}

\adfDesignE{13}{1160}{13}{180973}{22}{0}{5334}{"Audrey 21 Apr 2014"}{{0, 1, 5334, 38795, 80391, 80557, 61136, 167051, 119955, 100415, 114503, 156100, 161600}}

\adfDesignE{13}{1190}{13}{185653}{202}{0}{79687}{"Audrey 21 Apr 2014"}{{0, 1, 79687, 128410, 156922, 171252, 134359, 57123, 120247, 14400, 157260, 120, 94137}}

\adfDesignE{13}{1199}{13}{187057}{27}{0}{61565}{"Audrey 21 Apr 2014"}{{0, 1, 61565, 100291, 37959, 42734, 149062, 172667, 168659, 144322, 163487, 101156, 167496}}

\adfDesignE{13}{1204}{13}{187837}{44}{0}{14577}{"Audrey 21 Apr 2014"}{{0, 1, 14577, 45282, 16496, 30832, 131960, 130040, 129913, 157004, 41300, 12515, 41428}}

\adfDesignE{13}{1228}{13}{191581}{5}{0}{97595}{"Audrey 21 Apr 2014"}{{0, 1, 97595, 143029, 132014, 84080, 181789, 147369, 108723, 107500, 103778, 92764, 142425}}

\adfDesignE{13}{1255}{13}{195793}{144}{0}{78641}{"Audrey 21 Apr 2014"}{{0, 1, 78641, 89183, 135043, 104243, 116646, 60243, 162335, 91549, 506, 46367, 94208}}

\adfDesignE{13}{1256}{13}{195949}{10}{0}{5600}{"Audrey 21 Apr 2014"}{{0, 1, 5600, 8160, 39883, 158889, 169940, 135656, 175276, 37059, 20409, 52133, 176739}}

\adfDesignE{13}{1279}{13}{199537}{27}{0}{23677}{"Audrey 21 Apr 2014"}{{0, 1, 23677, 100896, 57628, 24150, 126045, 92093, 145162, 175386, 49815, 6548, 196284}}

\adfDesignE{13}{1291}{13}{201409}{456}{0}{14752}{"Audrey 21 Apr 2014"}{{0, 1, 14752, 99784, 116596, 192741, 24379, 123943, 16234, 8667, 162278, 179091, 68579}}

\adfDesignE{13}{1336}{13}{208429}{10}{0}{98004}{"Audrey 21 Apr 2014"}{{0, 1, 98004, 167267, 102647, 199332, 116874, 112230, 190590, 9096, 201980, 137361, 123621}}

\adfDesignE{13}{1361}{13}{212329}{144}{0}{13082}{"Audrey 21 Apr 2014"}{{0, 1, 13082, 1550, 105845, 66881, 141762, 48998, 182914, 145447, 57485, 161781, 135899}}

\adfDesignE{13}{1364}{13}{212797}{146}{0}{65296}{"Audrey 21 Apr 2014"}{{0, 1, 65296, 179721, 159054, 32399, 110127, 16368, 98394, 180397, 37374, 16708, 168146}}

\adfDesignE{13}{1394}{13}{217477}{14}{0}{43840}{"Audrey 21 Apr 2014"}{{0, 1, 43840, 101351, 172730, 151537, 112161, 200747, 106721, 65939, 61476, 132856, 155503}}

\adfDesignE{13}{1420}{13}{221533}{17}{0}{1343}{"Audrey 21 Apr 2014"}{{0, 1, 1343, 31385, 58785, 82507, 40401, 204491, 152026, 139025, 179789, 207190, 10722}}

\adfDesignE{13}{1441}{13}{224809}{508}{0}{30580}{"Audrey 21 Apr 2014"}{{0, 1, 30580, 155769, 162928, 121182, 218813, 86464, 90471, 103626, 200225, 207385, 196219}}

\adfDesignE{13}{1445}{13}{225433}{27}{0}{18295}{"Audrey 21 Apr 2014"}{{0, 1, 18295, 164453, 38817, 43065, 211273, 190750, 68410, 182367, 221298, 95663, 118206}}

\adfDesignE{13}{1454}{13}{226837}{63}{0}{10209}{"Audrey 21 Apr 2014"}{{0, 1, 10209, 105498, 7006, 70599, 84042, 87244, 111934, 156237, 132586, 34095, 107897}}

\adfDesignE{13}{1478}{13}{230581}{17}{0}{104394}{"Audrey 21 Apr 2014"}{{0, 1, 104394, 157433, 169246, 228380, 118063, 53210, 108450, 2200, 8124, 19938, 183466}}

\adfDesignE{13}{1505}{13}{234793}{170}{0}{8597}{"Audrey 21 Apr 2014"}{{0, 1, 8597, 183407, 114984, 38918, 232814, 126426, 27525, 195874, 228175, 159753, 92284}}

\adfDesignE{13}{1508}{13}{235261}{23}{0}{67799}{"Audrey 21 Apr 2014"}{{0, 1, 67799, 174983, 165970, 66400, 134365, 36193, 76977, 168860, 33097, 24085, 227575}}

\adfDesignE{13}{1510}{13}{235573}{14}{0}{13972}{"Audrey 21 Apr 2014"}{{0, 1, 13972, 162340, 117636, 17371, 67422, 199330, 94754, 218201, 154179, 109476, 23183}}

\adfDesignE{13}{1514}{13}{236197}{55}{0}{15389}{"Audrey 21 Apr 2014"}{{0, 1, 15389, 151927, 126697, 170095, 56801, 181689, 148132, 66101, 164007, 138778, 197565}}

\adfDesignE{13}{1543}{13}{240721}{40}{0}{24574}{"Audrey 21 Apr 2014"}{{0, 1, 24574, 153208, 56952, 227275, 87929, 55550, 197630, 13445, 128218, 31963, 226860}}

\adfDesignE{13}{1595}{13}{248833}{27}{0}{27024}{"Audrey 21 Apr 2014"}{{0, 1, 27024, 220554, 203280, 201412, 233679, 57422, 49540, 47420, 236963, 219690, 244846}}

\adfDesignE{13}{1601}{13}{249769}{105}{0}{148425}{"Audrey 21 Apr 2014"}{{0, 1, 148425, 105056, 107899, 220333, 170817, 211342, 197409, 29435, 180296, 183140, 194230}}

\adfDesignE{13}{1603}{13}{250081}{92}{0}{14219}{"Audrey 21 Apr 2014"}{{0, 1, 14219, 114513, 233037, 229934, 122633, 153895, 24255, 20146, 113229, 231754, 242870}}

\adfDesignE{13}{1609}{13}{251017}{131}{0}{41129}{"Audrey 21 Apr 2014"}{{0, 1, 41129, 242095, 33916, 29695, 127950, 135162, 55416, 221321, 81938, 124777, 161685}}

\adfDesignE{13}{1618}{13}{252421}{10}{0}{66713}{"Audrey 21 Apr 2014"}{{0, 1, 66713, 189718, 15573, 209134, 143030, 194169, 108040, 43286, 42678, 120955, 128808}}

\adfDesignE{13}{1619}{13}{252577}{261}{0}{3592}{"Audrey 21 Apr 2014"}{{0, 1, 3592, 21037, 44381, 40465, 118505, 77715, 54695, 212111, 130480, 153825, 153501}}

\adfDesignE{13}{1625}{13}{253513}{14}{0}{181107}{"Audrey 21 Apr 2014"}{{0, 1, 181107, 233509, 89855, 116502, 201263, 39001, 228414, 137010, 124656, 234516, 188757}}

\adfDesignE{13}{1646}{13}{256789}{22}{0}{35682}{"Audrey 21 Apr 2014"}{{0, 1, 35682, 45262, 92663, 242791, 234758, 177776, 201354, 13997, 243138, 33751, 219561}}

\adfDesignE{13}{1655}{13}{258193}{27}{0}{41467}{"Audrey 21 Apr 2014"}{{0, 1, 41467, 204902, 55990, 65874, 173411, 158887, 256448, 192318, 43315, 152597, 203948}}

\adfDesignE{13}{1666}{13}{259909}{22}{0}{46797}{"Audrey 21 Apr 2014"}{{0, 1, 46797, 225884, 194518, 65939, 107735, 219922, 73161, 193969, 105377, 74012, 252139}}

\adfDesignE{13}{1669}{13}{260377}{27}{0}{82}{"Audrey 21 Apr 2014"}{{0, 1, 82, 6724, 30614, 166955, 150706, 120173, 220237, 93421, 109589, 133480, 9526}}

\adfDesignE{13}{1676}{13}{261469}{40}{0}{121}{"Audrey 21 Apr 2014"}{{0, 1, 121, 14641, 202747, 215770, 222739, 20112, 80331, 45698, 38609, 226716, 239860}}

\adfDesignE{13}{1699}{13}{265057}{27}{0}{40024}{"Audrey 21 Apr 2014"}{{0, 1, 40024, 181125, 38050, 160735, 59193, 61166, 41532, 104321, 165840, 22766, 185475}}

\adfDesignE{13}{1750}{13}{273013}{66}{0}{32582}{"Audrey 21 Apr 2014"}{{0, 1, 32582, 112180, 223729, 91178, 107143, 189008, 177428, 181834, 133288, 244838, 144869}}

\adfDesignE{13}{1766}{13}{275509}{23}{0}{73909}{"Audrey 21 Apr 2014"}{{0, 1, 73909, 23338, 201902, 256460, 233958, 105964, 74442, 19048, 243151, 146207, 274674}}

\adfDesignE{13}{1804}{13}{281437}{42}{0}{21343}{"Audrey 21 Apr 2014"}{{0, 1, 21343, 158583, 75607, 201880, 205807, 151542, 86902, 79556, 54287, 252749, 118928}}

\adfDesignE{13}{1805}{13}{281593}{144}{0}{188452}{"Audrey 21 Apr 2014"}{{0, 1, 188452, 210330, 78480, 167007, 279926, 108304, 244768, 114585, 94808, 244552, 239938}}

\adfDesignE{13}{1813}{13}{282841}{27}{0}{19692}{"Audrey 21 Apr 2014"}{{0, 1, 19692, 282694, 216527, 21609, 131564, 217569, 176121, 261231, 131585, 65419, 173034}}

\adfDesignE{13}{1820}{13}{283933}{23}{0}{51874}{"Audrey 21 Apr 2014"}{{0, 1, 51874, 78835, 283724, 231721, 275532, 43681, 122854, 52211, 240460, 161417, 161288}}

\adfDesignE{13}{1841}{13}{287209}{40}{0}{102775}{"Audrey 21 Apr 2014"}{{0, 1, 102775, 15232, 179750, 236161, 275812, 198836, 162341, 51047, 195831, 73141, 232327}}

\adfDesignE{13}{1844}{13}{287677}{43}{0}{121365}{"Audrey 21 Apr 2014"}{{0, 1, 121365, 113148, 233102, 268050, 222382, 110644, 122054, 19626, 231607, 63885, 220198}}

\adfDesignE{13}{1871}{13}{291889}{27}{0}{16858}{"Audrey 21 Apr 2014"}{{0, 1, 16858, 184167, 155882, 273978, 161477, 22452, 207672, 17910, 113554, 85270, 220224}}

\adfDesignE{13}{1885}{13}{294073}{326}{0}{59212}{"Audrey 21 Apr 2014"}{{0, 1, 59212, 122638, 96667, 9532, 82697, 45241, 99135, 284540, 152164, 126194, 98271}}

\adfDesignE{13}{1918}{13}{299221}{23}{0}{64698}{"Audrey 21 Apr 2014"}{{0, 1, 64698, 28635, 150019, 97685, 177389, 92067, 257540, 201535, 57134, 178519, 190883}}

\adfDesignE{13}{1919}{13}{299377}{92}{0}{233913}{"Audrey 21 Apr 2014"}{{0, 1, 233913, 252918, 21033, 229888, 293758, 207260, 295754, 69488, 71083, 138576, 281967}}

\adfDesignE{13}{1939}{13}{302497}{300}{0}{24977}{"Audrey 21 Apr 2014"}{{0, 1, 24977, 101715, 165749, 241328, 94234, 255958, 91368, 61168, 183286, 247321, 45380}}

\adfDesignE{13}{1966}{13}{306709}{10}{0}{10423}{"Audrey 21 Apr 2014"}{{0, 1, 10423, 63943, 305941, 276279, 271925, 283115, 60356, 30429, 24361, 266360, 247121}}

\adfDesignE{13}{1973}{13}{307801}{105}{0}{3832}{"Audrey 21 Apr 2014"}{{0, 1, 3832, 217577, 229956, 264930, 84062, 165738, 114553, 42870, 219907, 232287, 271093}}

\adfDesignE{13}{1993}{13}{310921}{14}{0}{55000}{"Audrey 21 Apr 2014"}{{0, 1, 55000, 49591, 105988, 193092, 242324, 191335, 303755, 117828, 13597, 69995, 212099}}

\adfDesignE{13}{2005}{13}{312793}{105}{0}{6309}{"Audrey 21 Apr 2014"}{{0, 1, 6309, 78770, 244646, 150952, 214276, 288731, 210240, 161840, 92208, 258085, 170700}}

\adfDesignE{13}{2015}{13}{314353}{53}{0}{15357}{"Audrey 21 Apr 2014"}{{0, 1, 15357, 72699, 171040, 241965, 204045, 48361, 178091, 72387, 94951, 193293, 279575}}

\adfDesignE{13}{2059}{13}{321217}{144}{0}{433}{"Audrey 21 Apr 2014"}{{0, 1, 433, 187489, 236053, 63943, 62657, 148253, 271366, 257273, 258127, 306692, 135015}}

\adfDesignE{13}{2065}{13}{322153}{144}{0}{43879}{"Audrey 21 Apr 2014"}{{0, 1, 43879, 180313, 198600, 130750, 278626, 123904, 129588, 191402, 321801, 17936, 316118}}

\adfDesignE{13}{2099}{13}{327457}{27}{0}{112083}{"Audrey 21 Apr 2014"}{{0, 1, 112083, 38541, 305616, 63729, 117966, 251889, 114618, 263727, 97408, 37027, 234680}}

\adfDesignE{13}{2113}{13}{329641}{14}{0}{128814}{"Audrey 21 Apr 2014"}{{0, 1, 128814, 237220, 195662, 313290, 168276, 101427, 226184, 16350, 32551, 320635, 237436}}

\adfDesignE{13}{2161}{13}{337129}{53}{0}{98823}{"Audrey 21 Apr 2014"}{{0, 1, 98823, 32457, 52805, 265853, 265178, 311195, 316105, 71275, 310257, 330606, 305348}}

\adfDesignE{13}{2213}{13}{345241}{105}{0}{172260}{"Audrey 21 Apr 2014"}{{0, 1, 172260, 43650, 145261, 282662, 291685, 318683, 252652, 62578, 226537, 328149, 292569}}

\adfDesignE{13}{2273}{13}{354601}{27}{0}{149796}{"Audrey 21 Apr 2014"}{{0, 1, 149796, 44937, 346670, 235875, 333259, 136384, 150251, 118725, 226147, 173280, 212281}}

\adfDesignE{13}{2333}{13}{363961}{27}{0}{61435}{"Audrey 21 Apr 2014"}{{0, 1, 61435, 347616, 13324, 11651, 231859, 279969, 190538, 352309, 70667, 100337, 160099}}

\adfDesignE{13}{2338}{13}{364741}{358}{0}{263059}{"Audrey 21 Apr 2014"}{{0, 1, 263059, 280738, 88308, 224623, 330275, 140284, 298081, 140117, 136148, 308460, 343093}}

\adfDesignE{13}{2369}{13}{369577}{105}{0}{103848}{"Audrey 21 Apr 2014"}{{0, 1, 103848, 150244, 106703, 235530, 344003, 341147, 152013, 134046, 291303, 247763, 110861}}

\adfDesignE{13}{2389}{13}{372697}{27}{0}{68780}{"Audrey 21 Apr 2014"}{{0, 1, 68780, 45379, 202942, 102716, 334845, 200682, 74565, 269980, 341769, 126636, 95190}}

\adfDesignE{13}{2423}{13}{378001}{27}{0}{55163}{"Audrey 21 Apr 2014"}{{0, 1, 55163, 48519, 206517, 281134, 325816, 174461, 264684, 96866, 375023, 155021, 284801}}

\adfDesignE{13}{2449}{13}{382057}{144}{0}{70789}{"Audrey 21 Apr 2014"}{{0, 1, 70789, 22909, 255293, 258020, 360838, 176333, 252490, 124036, 332487, 182815, 256331}}

\adfDesignE{13}{2501}{13}{390169}{27}{0}{20082}{"Audrey 21 Apr 2014"}{{0, 1, 20082, 242147, 119807, 182120, 279803, 180077, 220022, 208048, 90284, 358114, 50340}}

\adfDesignE{13}{2563}{13}{399841}{27}{0}{55819}{"Audrey 21 Apr 2014"}{{0, 1, 55819, 199689, 72734, 353473, 355242, 338326, 128723, 46367, 388621, 261667, 198384}}

\adfDesignE{13}{2569}{13}{400777}{144}{0}{57165}{"Audrey 21 Apr 2014"}{{0, 1, 57165, 302344, 387412, 271514, 206931, 277460, 251125, 129262, 136681, 221750, 163017}}

\adfDesignE{13}{2603}{13}{406081}{131}{0}{27239}{"Audrey 21 Apr 2014"}{{0, 1, 27239, 53134, 44342, 146844, 391947, 374843, 253894, 259236, 392976, 384185, 107845}}

\adfDesignE{13}{2611}{13}{407329}{27}{0}{147170}{"Audrey 21 Apr 2014"}{{0, 1, 147170, 103983, 234909, 323313, 244404, 156664, 197493, 84015, 15755, 146682, 382256}}

\adfDesignE{13}{2765}{13}{431353}{27}{0}{67682}{"Audrey 21 Apr 2014"}{{0, 1, 67682, 315617, 126528, 16987, 158389, 99542, 330490, 414365, 205282, 16194, 405688}}

\adfDesignE{13}{2801}{13}{436969}{14}{0}{128899}{"Audrey 21 Apr 2014"}{{0, 1, 128899, 79914, 164449, 382430, 371680, 336129, 341683, 54538, 373359, 20926, 367806}}

\adfDesignE{13}{2833}{13}{441961}{144}{0}{4749}{"Audrey 21 Apr 2014"}{{0, 1, 4749, 12990, 256931, 352959, 286179, 33996, 131239, 89001, 151033, 394975, 53791}}

\adfDesignE{13}{2891}{13}{451009}{27}{0}{24814}{"Audrey 21 Apr 2014"}{{0, 1, 24814, 107311, 58018, 37924, 241362, 208157, 252730, 413084, 184833, 135541, 140261}}

\adfDesignE{13}{2935}{13}{457873}{66}{0}{76964}{"Audrey 21 Apr 2014"}{{0, 1, 76964, 412168, 198639, 130399, 368222, 246546, 451351, 327473, 12687, 257032, 265756}}

\adfDesignE{13}{2953}{13}{460681}{79}{0}{251762}{"Audrey 21 Apr 2014"}{{0, 1, 251762, 387897, 262729, 139837, 400774, 389806, 386104, 320843, 268826, 143659, 272529}}

\adfDesignE{13}{2989}{13}{466297}{261}{0}{10053}{"Audrey 21 Apr 2014"}{{0, 1, 10053, 342657, 194882, 235049, 220698, 35868, 133423, 231247, 235546, 87772, 137992}}

\adfDesignE{13}{3005}{13}{468793}{27}{0}{10872}{"Audrey 21 Apr 2014"}{{0, 1, 10872, 64548, 451528, 280913, 368534, 396670, 169433, 187879, 89387, 7575, 316625}}

\adfDesignE{13}{3019}{13}{470977}{27}{0}{11034}{"Audrey 21 Apr 2014"}{{0, 1, 11034, 237090, 244802, 92173, 197539, 434747, 97653, 378803, 262404, 270117, 128522}}

\adfDesignE{13}{3031}{13}{472849}{27}{0}{31126}{"Audrey 21 Apr 2014"}{{0, 1, 31126, 433124, 19785, 178512, 388762, 400102, 150739, 294336, 52961, 112472, 302325}}

\adfDesignE{13}{3059}{13}{477217}{183}{0}{74299}{"Audrey 21 Apr 2014"}{{0, 1, 74299, 372362, 423097, 445779, 165253, 293671, 145955, 31437, 237665, 288401, 385382}}

\adfDesignE{13}{3085}{13}{481273}{79}{0}{81436}{"Audrey 21 Apr 2014"}{{0, 1, 81436, 361429, 119183, 435470, 333915, 296167, 140690, 45802, 65922, 304950, 221400}}

\adfDesignE{13}{3229}{13}{503737}{27}{0}{52920}{"Audrey 21 Apr 2014"}{{0, 1, 52920, 252417, 313611, 174918, 493185, 232493, 257072, 328818, 461369, 18827, 436791}}

\adfDesignE{13}{3235}{13}{504673}{287}{0}{79590}{"Audrey 21 Apr 2014"}{{0, 1, 79590, 417277, 60319, 339634, 174834, 194104, 192157, 165038, 250249, 397965, 252197}}

\adfDesignE{13}{3263}{13}{509041}{27}{0}{63476}{"Audrey 21 Apr 2014"}{{0, 1, 63476, 143061, 157637, 456316, 172475, 78313, 210623, 52724, 273090, 287667, 140781}}

\adfDesignE{13}{3265}{13}{509353}{27}{0}{53895}{"Audrey 21 Apr 2014"}{{0, 1, 53895, 340219, 413711, 26770, 281454, 430990, 181191, 482582, 174004, 247497, 423804}}

\adfDesignE{13}{3323}{13}{518401}{27}{0}{38939}{"Audrey 21 Apr 2014"}{{0, 1, 38939, 441197, 479244, 401319, 280797, 358892, 359831, 117081, 198665, 236713, 197727}}

\adfDesignE{13}{3353}{13}{523081}{27}{0}{151148}{"Audrey 21 Apr 2014"}{{0, 1, 151148, 155229, 277718, 316176, 166807, 40236, 251222, 206904, 205126, 327616, 517222}}

\adfDesignE{13}{3439}{13}{536497}{79}{0}{158214}{"Audrey 21 Apr 2014"}{{0, 1, 158214, 329267, 253941, 370535, 260803, 165075, 502090, 165961, 117480, 42155, 316963}}

\adfDesignE{13}{3449}{13}{538057}{27}{0}{249773}{"Audrey 21 Apr 2014"}{{0, 1, 249773, 456550, 214798, 1270, 296137, 331111, 536618, 536786, 530204, 288453, 324698}}

\adfDesignE{13}{3470}{13}{541333}{973}{0}{187750}{"Audrey 21 Apr 2014"}{{0, 1, 187750, 81539, 50010, 497948, 445234, 41640, 520147, 43384, 449682, 418154, 512509}}

\adfDesignE{13}{3506}{13}{546949}{518}{0}{208602}{"Audrey 21 Apr 2014"}{{0, 1, 208602, 78913, 432522, 247204, 350139, 126218, 296274, 299744, 535157, 341818, 365102}}

\adfDesignE{13}{3538}{13}{551941}{154}{0}{66557}{"Audrey 21 Apr 2014"}{{0, 1, 66557, 507724, 550484, 168067, 399013, 467026, 188185, 383873, 86371, 129132, 365213}}

\adfDesignE{13}{3559}{13}{555217}{27}{0}{160086}{"Audrey 21 Apr 2014"}{{0, 1, 160086, 376327, 308320, 34654, 447197, 298962, 480549, 520562, 503151, 435145, 321565}}

\adfDesignE{13}{3575}{13}{557713}{157}{0}{363126}{"Audrey 21 Apr 2014"}{{0, 1, 363126, 407286, 129557, 192780, 409946, 85801, 494894, 364932, 342354, 64626, 490975}}

\adfDesignE{13}{3649}{13}{569257}{53}{0}{186794}{"Audrey 21 Apr 2014"}{{0, 1, 186794, 529135, 288794, 485345, 233367, 131366, 557619, 83911, 149096, 478013, 292101}}

\adfDesignE{13}{3685}{13}{574873}{105}{0}{48422}{"Audrey 21 Apr 2014"}{{0, 1, 48422, 357990, 446211, 402210, 265126, 442209, 349567, 172662, 261325, 349547, 353968}}

\adfDesignE{13}{3689}{13}{575497}{300}{0}{22247}{"Audrey 21 Apr 2014"}{{0, 1, 22247, 1589, 245166, 222933, 532802, 309882, 66291, 352563, 20448, 264026, 264040}}

\adfDesignE{13}{3725}{13}{581113}{144}{0}{196349}{"Audrey 21 Apr 2014"}{{0, 1, 196349, 150042, 492010, 284144, 474465, 178803, 429465, 296968, 491412, 252268, 240751}}

\adfDesignE{13}{3743}{13}{583921}{14}{0}{56885}{"Audrey 21 Apr 2014"}{{0, 1, 56885, 396964, 488149, 576631, 476181, 44916, 392285, 7289, 50855, 142041, 287408}}

\adfDesignE{13}{3755}{13}{585793}{144}{0}{3949}{"Audrey 21 Apr 2014"}{{0, 1, 3949, 363983, 418638, 93616, 54201, 225304, 491722, 492176, 527643, 582299, 261226}}

\adfDesignE{13}{3799}{13}{592657}{27}{0}{138720}{"Audrey 21 Apr 2014"}{{0, 1, 138720, 258267, 89933, 75910, 498281, 547067, 590704, 516746, 548313, 379980, 504677}}

\adfDesignE{13}{3821}{13}{596089}{27}{0}{147755}{"Audrey 21 Apr 2014"}{{0, 1, 147755, 376489, 510626, 559900, 408724, 45852, 310775, 36188, 39610, 173748, 370777}}

\adfDesignE{13}{3853}{13}{601081}{27}{0}{8259}{"Audrey 21 Apr 2014"}{{0, 1, 8259, 288928, 565863, 57742, 233945, 277421, 500348, 543338, 358877, 34732, 135951}}

\adfDesignE{13}{3865}{13}{602953}{482}{0}{261115}{"Audrey 21 Apr 2014"}{{0, 1, 261115, 323891, 198873, 602176, 308806, 371047, 434600, 776, 33032, 510968, 572433}}

\adfDesignE{13}{3881}{13}{605449}{27}{0}{187191}{"Audrey 21 Apr 2014"}{{0, 1, 187191, 109606, 406483, 156178, 405584, 186291, 557977, 449270, 12674, 309552, 246438}}

\adfDesignE{13}{3905}{13}{609193}{27}{0}{101891}{"Audrey 21 Apr 2014"}{{0, 1, 101891, 517968, 60319, 424245, 239594, 281165, 272997, 184947, 267708, 419253, 275877}}

\adfDesignE{13}{3935}{13}{613873}{547}{0}{72518}{"Audrey 21 Apr 2014"}{{0, 1, 72518, 424206, 166932, 613089, 236077, 141662, 494134, 783, 305278, 48005, 566680}}

\adfDesignE{13}{3941}{13}{614809}{40}{0}{5347}{"Audrey 21 Apr 2014"}{{0, 1, 5347, 309195, 44264, 592952, 559140, 520222, 231118, 21856, 50322, 400201, 339427}}

\adfDesignE{13}{3943}{13}{615121}{27}{0}{13732}{"Audrey 21 Apr 2014"}{{0, 1, 13732, 340798, 612689, 435431, 362372, 378535, 270170, 179689, 239017, 510909, 347383}}

\adfDesignE{13}{3946}{13}{615589}{84}{0}{41539}{"Audrey 21 Apr 2014"}{{0, 1, 41539, 608143, 341873, 39906, 489746, 189411, 100520, 575682, 84304, 433624, 173196}}

\adfDesignE{13}{4033}{13}{629161}{27}{0}{123257}{"Audrey 21 Apr 2014"}{{0, 1, 123257, 566543, 440322, 82572, 268668, 580763, 312316, 546588, 237236, 111016, 505684}}

\adfDesignE{13}{4081}{13}{636649}{391}{0}{148717}{"Audrey 21 Apr 2014"}{{0, 1, 148717, 196478, 612871, 392369, 513127, 48972, 341013, 244279, 611454, 391199, 319414}}

\adfDesignE{13}{4111}{13}{641329}{144}{0}{29286}{"Audrey 21 Apr 2014"}{{0, 1, 29286, 212923, 21111, 15590, 583821, 591995, 118713, 625738, 28222, 477740, 501505}}

\adfDesignE{13}{4133}{13}{644761}{53}{0}{5587}{"Audrey 21 Apr 2014"}{{0, 1, 5587, 266041, 196962, 464428, 240972, 49596, 490383, 180332, 398202, 329124, 602177}}

\adfDesignE{13}{4285}{13}{668473}{27}{0}{110192}{"Audrey 21 Apr 2014"}{{0, 1, 110192, 133292, 23308, 81870, 375905, 462788, 404018, 586602, 182376, 72393, 241147}}

\adfDesignE{13}{4298}{13}{670501}{206}{0}{217663}{"Audrey 21 Apr 2014"}{{0, 1, 217663, 251410, 386216, 199832, 632746, 464192, 298107, 470668, 490593, 625400, 656679}}

\adfDesignE{13}{4303}{13}{671281}{27}{0}{49358}{"Audrey 21 Apr 2014"}{{0, 1, 49358, 133415, 502241, 546510, 556157, 103273, 312101, 124770, 65766, 434593, 528220}}

\adfDesignE{13}{4334}{13}{676117}{266}{0}{34226}{"Audrey 21 Apr 2014"}{{0, 1, 34226, 384432, 332812, 280413, 610640, 312053, 381846, 395703, 31251, 655749, 637576}}

\adfDesignE{13}{4429}{13}{690937}{27}{0}{60935}{"Audrey 21 Apr 2014"}{{0, 1, 60935, 669724, 128972, 191382, 227484, 159446, 576853, 499554, 402518, 552704, 676049}}

\adfDesignE{13}{4481}{13}{699049}{40}{0}{266107}{"Audrey 21 Apr 2014"}{{0, 1, 266107, 669847, 471119, 616073, 366331, 161318, 648034, 82975, 66611, 566933, 278945}}

\adfDesignE{13}{4501}{13}{702169}{131}{0}{15623}{"Audrey 21 Apr 2014"}{{0, 1, 15623, 425486, 636024, 209433, 566388, 648155, 146416, 492735, 120158, 330697, 621898}}

\adfDesignE{13}{4523}{13}{705601}{14}{0}{68208}{"Audrey 21 Apr 2014"}{{0, 1, 68208, 303871, 109394, 520978, 95463, 54276, 474562, 184622, 541930, 347454, 121645}}

\adfDesignE{13}{4535}{13}{707473}{14}{0}{223753}{"Audrey 21 Apr 2014"}{{0, 1, 223753, 370691, 503749, 44164, 552101, 272104, 374878, 663308, 639092, 64678, 536319}}

\adfDesignE{13}{4585}{13}{715273}{53}{0}{253731}{"Audrey 21 Apr 2014"}{{0, 1, 253731, 558723, 382732, 503701, 294064, 165062, 681626, 211571, 167478, 706761, 366188}}

\adfDesignE{13}{4625}{13}{721513}{313}{0}{140391}{"Audrey 21 Apr 2014"}{{0, 1, 140391, 62260, 335178, 339764, 583294, 388506, 692124, 381748, 719341, 270747, 415724}}

\adfDesignE{13}{4658}{13}{726661}{397}{0}{348943}{"Audrey 21 Apr 2014"}{{0, 1, 348943, 446767, 546324, 55487, 644557, 447175, 489512, 671173, 459822, 559380, 417486}}

\adfDesignE{13}{4679}{13}{729937}{27}{0}{423986}{"Audrey 21 Apr 2014"}{{0, 1, 423986, 353395, 364480, 184947, 726580, 56148, 530547, 544989, 309308, 320394, 564847}}

\adfDesignE{13}{4735}{13}{738673}{27}{0}{4183}{"Audrey 21 Apr 2014"}{{0, 1, 4183, 508010, 582282, 280725, 521278, 681851, 166280, 457947, 213212, 287485, 728784}}

\adfDesignE{13}{4775}{13}{744913}{27}{0}{275662}{"Audrey 21 Apr 2014"}{{0, 1, 275662, 218201, 234051, 562006, 416797, 458407, 583853, 182906, 52454, 68305, 671922}}

\adfDesignE{13}{4849}{13}{756457}{27}{0}{23638}{"Audrey 21 Apr 2014"}{{0, 1, 23638, 489778, 554436, 140643, 647176, 116377, 441874, 615813, 85643, 150302, 516604}}

\adfDesignE{13}{4859}{13}{758017}{300}{0}{11710}{"Audrey 21 Apr 2014"}{{0, 1, 11710, 681040, 639560, 39640, 277996, 408162, 279835, 718376, 468311, 426832, 596639}}

\adfDesignE{13}{4891}{13}{763009}{14}{0}{461133}{"Audrey 21 Apr 2014"}{{0, 1, 461133, 665479, 501006, 430706, 744189, 704315, 477955, 332302, 320696, 156224, 547057}}

\adfDesignE{13}{4945}{13}{771433}{183}{0}{24912}{"Audrey 21 Apr 2014"}{{0, 1, 24912, 375612, 535287, 78906, 94988, 356045, 627839, 692526, 651533, 39776, 379740}}

\adfDesignE{13}{4958}{13}{773461}{258}{0}{3059}{"Audrey 21 Apr 2014"}{{0, 1, 3059, 75949, 289691, 551924, 643614, 356981, 651408, 221536, 126788, 340531, 605823}}

~

\adfDesignE{14}{7398}{1}{1346437}{203624}{0}{884546}{"Audrey 19 May 2014"}{{0, 1, 884546, 353231, 283654, 515445, 276719, 555907, 1135074, 559874, 1314234, 186334, 948320, 1065283}}

\adfDesignE{14}{7421}{1}{1350623}{642722}{0}{1339303}{"Audrey 19 May 2014"}{{0, 1, 1339303, 1183838, 1185869, 1155540, 70955, 410085, 1279674, 872618, 420862, 840104, 1109886, 936249}}

\adfDesignE{14}{7494}{1}{1363909}{340445}{0}{274725}{"Audrey 19 May 2014"}{{0, 1, 274725, 557201, 82019, 893095, 569956, 317173, 662051, 604098, 375930, 815861, 791619, 875816}}

\adfDesignE{14}{7673}{1}{1396487}{710243}{0}{1136198}{"Audrey 19 May 2014"}{{0, 1, 1136198, 1193203, 1093133, 937839, 941590, 741964, 705582, 961633, 1044969, 1233436, 1141809, 40539}}

\adfDesignE{14}{8333}{1}{1516607}{689057}{0}{1115641}{"Audrey 19 May 2014"}{{0, 1, 1115641, 1258300, 199318, 897891, 539810, 144759, 145910, 1199179, 1142794, 33548, 696722, 209162}}

\adfDesignE{14}{8396}{1}{1528073}{46983}{0}{1337991}{"Audrey 19 May 2014"}{{0, 1, 1337991, 1408712, 1077771, 823742, 57320, 1188323, 938374, 898476, 963813, 225450, 820385, 956153}}

\adfDesignE{14}{9960}{1}{1812721}{767497}{0}{626762}{"Audrey 19 May 2014"}{{0, 1, 626762, 1274897, 513109, 577727, 872061, 1247841, 845392, 1232404, 12375, 1359312, 739512, 1574933}}

\adfDesignE{14}{10154}{1}{1848029}{925310}{0}{384695}{"Audrey 19 May 2014"}{{0, 1, 384695, 80705, 1770804, 842829, 958192, 111042, 111855, 551989, 1484139, 1533200, 1134418, 276276}}

\adfDesignE{14}{10185}{1}{1853671}{246671}{0}{290726}{"Audrey 19 May 2014"}{{0, 1, 290726, 1624160, 1780001, 1383985, 742179, 1573883, 1144734, 1406557, 1014111, 208365, 1008381, 798614}}

\adfDesignE{14}{10631}{1}{1934843}{760382}{0}{1162296}{"Audrey 19 May 2014"}{{0, 1, 1162296, 1390900, 721023, 1665375, 1132097, 1200859, 438610, 280077, 846571, 1679366, 460018, 631865}}

\adfDesignE{14}{10964}{1}{1995449}{1454356}{0}{1863243}{"Audrey 19 May 2014"}{{0, 1, 1863243, 288645, 365806, 1949377, 884484, 1215145, 148222, 1466897, 1112630, 156704, 1537943, 983597}}

\adfDesignE{14}{13439}{1}{2445899}{1129823}{0}{123084}{"Audrey 19 May 2014"}{{0, 1, 123084, 2218549, 383059, 1284832, 216144, 2270572, 218409, 2223346, 1355348, 1557836, 880018, 1944196}}

\adfDesignE{14}{14505}{1}{2639911}{922070}{0}{2053715}{"Audrey 19 May 2014"}{{0, 1, 2053715, 1735101, 924106, 2496424, 1301081, 289401, 252086, 2493191, 1016651, 518743, 996640, 1762326}}

\adfDesignE{14}{14546}{1}{2647373}{2635246}{0}{2031325}{"Audrey 19 May 2014"}{{0, 1, 2031325, 981889, 2207352, 1593419, 2141304, 208713, 385640, 2502300, 1913770, 1616114, 1708757, 1241027}}

\adfDesignE{14}{16950}{1}{3084901}{2729489}{0}{721202}{"Audrey 19 May 2014"}{{0, 1, 721202, 2591699, 71199, 684053, 3023786, 805258, 2557460, 1468426, 679257, 2313015, 82983, 426166}}

\adfDesignE{14}{17465}{1}{3178631}{2689564}{0}{2555942}{"Audrey 19 May 2014"}{{0, 1, 2555942, 2645448, 2359768, 69373, 2969924, 1224688, 1510333, 165795, 33494, 1831256, 2296187, 1409577}}

\adfDesignE{14}{18158}{1}{3304757}{955634}{0}{3082652}{"Audrey 19 May 2014"}{{0, 1, 3082652, 523286, 561903, 2682290, 1828697, 1802386, 2796665, 2236381, 184609, 2842911, 1953307, 2638211}}

\adfDesignE{14}{18315}{1}{3333331}{2317069}{0}{2148046}{"Audrey 19 May 2014"}{{0, 1, 2148046, 1514655, 891615, 1374951, 63168, 1198042, 1127347, 1575744, 2965563, 199017, 1303363, 2305143}}

\adfDesignE{14}{23570}{1}{4289741}{2941632}{0}{1470660}{"Audrey 19 May 2014"}{{0, 1, 1470660, 610551, 1506504, 2610183, 4127707, 2479851, 4275949, 2842469, 42191, 1802236, 4152277, 4107608}}

\adfDesignE{14}{24975}{1}{4545451}{3605414}{0}{2317309}{"Audrey 19 May 2014"}{{0, 1, 2317309, 2462748, 2556553, 624125, 698641, 158046, 794791, 501278, 3290597, 4173050, 3865245, 1284871}}

\adfDesignE{14}{25326}{1}{4609333}{2386558}{0}{3445679}{"Audrey 19 May 2014"}{{0, 1, 3445679, 2266973, 4388554, 4582378, 4391838, 4279699, 449042, 2906744, 4115482, 3099779, 181348, 2767147}}

\adfDesignE{14}{25403}{1}{4623347}{2930770}{0}{2889923}{"Audrey 19 May 2014"}{{0, 1, 2889923, 3938353, 3192875, 2539700, 1820682, 2994054, 854036, 3080177, 1053514, 3249635, 2840967, 3909512}}

\adfDesignE{14}{25919}{1}{4717259}{4513225}{0}{3182173}{"Audrey 19 May 2014"}{{0, 1, 3182173, 1162982, 2302911, 4066362, 1374316, 1774135, 3075932, 439301, 379977, 1136846, 3022812, 1668547}}

\adfDesignE{14}{26934}{1}{4901989}{4138022}{0}{572763}{"Audrey 19 May 2014"}{{0, 1, 572763, 1644322, 2361083, 4766954, 530737, 4374463, 1423644, 3954134, 3910374, 3671251, 4437062, 2667135}}

\adfDesignE{14}{27111}{1}{4934203}{2523913}{0}{2640072}{"Audrey 19 May 2014"}{{0, 1, 2640072, 3954632, 1178263, 4886631, 1697978, 1538480, 218644, 3230210, 1627897, 462539, 4901759, 3268112}}

\adfDesignE{14}{28503}{1}{5187547}{4792215}{0}{2757766}{"Audrey 19 May 2014"}{{0, 1, 2757766, 2593295, 3233454, 1525849, 3868814, 3919154, 726780, 2576825, 4988966, 4219197, 4238936, 1663792}}

\adfDesignE{14}{30311}{1}{5516603}{2960637}{0}{4138770}{"Audrey 19 May 2014"}{{0, 1, 4138770, 702102, 1163908, 124133, 2215623, 2733372, 5075797, 1129510, 1687294, 2998961, 613959, 4999585}}

\adfDesignE{14}{31640}{1}{5758481}{3816012}{0}{3142251}{"Audrey 19 May 2014"}{{0, 1, 3142251, 2209718, 4341114, 1501903, 4617065, 2497105, 1438831, 444089, 3679052, 2951211, 5189004, 2539542}}

\adfDesignE{14}{31826}{1}{5792333}{4281786}{0}{196886}{"Audrey 19 May 2014"}{{0, 1, 196886, 1804560, 2478606, 4557999, 39424, 307444, 1539534, 5697567, 4798250, 1906532, 3111620, 2523242}}

\adfDesignE{14}{32598}{1}{5932837}{488643}{0}{4180048}{"Audrey 19 May 2014"}{{0, 1, 4180048, 3033604, 1842672, 3452081, 329977, 357603, 1624283, 712925, 2628137, 5395905, 2570038, 3536911}}

\adfDesignE{14}{36338}{1}{6613517}{2020845}{0}{390263}{"Audrey 19 May 2014"}{{0, 1, 390263, 2526176, 2658615, 4064717, 3689985, 4356890, 4354887, 3053104, 2463081, 1138421, 749697, 3621748}}

\adfDesignE{14}{37715}{1}{6864131}{5882021}{0}{4646765}{"Audrey 19 May 2014"}{{0, 1, 4646765, 3583966, 3904087, 5133642, 2737533, 3184366, 143421, 5204275, 2401882, 303564, 4785829, 5155455}}

\adfDesignE{14}{39698}{1}{7225037}{5570528}{0}{6616353}{"Audrey 19 May 2014"}{{0, 1, 6616353, 3539533, 5578606, 224682, 2386885, 2305879, 5709495, 667605, 4499208, 3777319, 986166, 7058490}}

\adfDesignE{14}{40166}{1}{7310213}{4534309}{0}{2728011}{"Audrey 19 May 2014"}{{0, 1, 2728011, 5945092, 5945111, 6605829, 2599169, 5903083, 7140213, 5352933, 548102, 4628315, 3991060, 7094785}}

\adfDesignE{14}{40514}{1}{7373549}{6172672}{0}{6709393}{"Audrey 19 May 2014"}{{0, 1, 6709393, 2744058, 653537, 815662, 460209, 5158093, 1945088, 3919072, 2360866, 1875754, 4697671, 5528341}}

\adfDesignE{14}{41316}{1}{7519513}{1885389}{0}{1463575}{"Audrey 19 May 2014"}{{0, 1, 1463575, 5709880, 4289911, 1124650, 3267076, 5076591, 1033629, 2899309, 4749619, 2315562, 6781641, 6405634}}

\adfDesignE{14}{44549}{1}{8107919}{2018623}{0}{601271}{"Audrey 19 May 2014"}{{0, 1, 601271, 2815150, 2129777, 311788, 5887349, 7216855, 345176, 5916053, 4340088, 2890022, 5325801, 2760264}}

\adfDesignE{14}{44639}{1}{8124299}{7057009}{0}{7656575}{"Audrey 19 May 2014"}{{0, 1, 7656575, 2741003, 5867925, 5708377, 326515, 1670742, 5816705, 3973906, 2163874, 4993747, 19177, 7807247}}

\adfDesignE{14}{44663}{1}{8128667}{711835}{0}{2910401}{"Audrey 19 May 2014"}{{0, 1, 2910401, 5305453, 7486129, 5377614, 4563411, 84514, 4295361, 7783122, 4167695, 1768959, 7511439, 5646570}}

\adfDesignE{14}{45378}{1}{8258797}{4819016}{0}{7624854}{"Audrey 19 May 2014"}{{0, 1, 7624854, 2406432, 7742870, 3431367, 1710546, 7900816, 4325117, 1663684, 586076, 7981768, 5535939, 6902109}}

\adfDesignE{14}{45834}{1}{8341789}{7095269}{0}{6201825}{"Audrey 19 May 2014"}{{0, 1, 6201825, 3963232, 6490542, 2921118, 7404967, 8229405, 3937306, 5958567, 6325188, 4180783, 2084679, 694910}}

\adfDesignE{14}{46266}{1}{8420413}{6626252}{0}{1253567}{"Audrey 19 May 2014"}{{0, 1, 1253567, 4329016, 8033962, 979999, 6672211, 8350846, 3321952, 7835286, 6278421, 4493628, 1143575, 6250427}}

\adfDesignE{14}{47685}{1}{8678671}{3108700}{0}{608157}{"Audrey 19 May 2014"}{{0, 1, 608157, 4693313, 3799648, 8266947, 4183024, 8568564, 2282637, 4850204, 4850561, 4995835, 5525073, 8126733}}

\adfDesignE{14}{48600}{1}{8845201}{6858804}{0}{4062622}{"Audrey 19 May 2014"}{{0, 1, 4062622, 114512, 5623869, 4410262, 7937721, 2325848, 4916789, 8504066, 7762715, 5201697, 7855183, 3201122}}

~

\adfDesignE{16}{15047}{1}{3611281}{3580612}{0}{3354316}{"Audrey 19 May 2014"}{{0, 1, 3354316, 2349421, 545191, 1029799, 1537802, 3132495, 2123382, 813222, 795116, 2073478, 735751, 2749759, 2252868, 1786366}}

\adfDesignE{16}{15454}{1}{3708961}{81852}{0}{1257914}{"Audrey 19 May 2014"}{{0, 1, 1257914, 1017888, 631290, 1434155, 1313309, 2913611, 2911928, 3659397, 384914, 2395651, 3246397, 3488106, 3127235, 1889892}}

\adfDesignE{16}{19843}{1}{4762321}{4741643}{0}{3273462}{"Audrey 19 May 2014"}{{0, 1, 3273462, 2615295, 1496183, 3458800, 3631656, 1247192, 2596466, 4200888, 2468385, 1130664, 241667, 4312881, 3827571, 3597457}}

\adfDesignE{16}{20563}{1}{4935121}{3567244}{0}{2862504}{"Audrey 19 May 2014"}{{0, 1, 2862504, 4635207, 3689383, 892050, 931227, 4495952, 353154, 4422218, 1738146, 4003893, 2511786, 4881881, 1758641, 2304925}}

\adfDesignE{16}{21375}{1}{5130001}{2861556}{0}{761714}{"Audrey 19 May 2014"}{{0, 1, 761714, 5104696, 3370988, 4222901, 3935289, 2541026, 1091267, 3899605, 283948, 1194711, 1827261, 4064039, 2989409, 623152}}

\adfDesignE{16}{24553}{1}{5892721}{2904881}{0}{4275857}{"Audrey 19 May 2014"}{{0, 1, 4275857, 2450056, 3811471, 945740, 4569256, 3780704, 649746, 3383336, 1400347, 1323464, 3728881, 2792919, 4590635, 3546634}}

\adfDesignE{16}{25697}{1}{6167281}{2104735}{0}{1862716}{"Audrey 19 May 2014"}{{0, 1, 1862716, 4773337, 2185592, 3988714, 788185, 397443, 3023948, 5901600, 5214249, 5379095, 3907122, 4537277, 4247370, 3131599}}

\adfDesignE{16}{27753}{1}{6660721}{4654919}{0}{213323}{"Audrey 19 May 2014"}{{0, 1, 213323, 656457, 2378307, 465591, 3258062, 6627281, 111071, 1814336, 5083381, 3402658, 6480838, 5893193, 2468078, 1111749}}

\adfDesignE{16}{27824}{1}{6677761}{1410101}{0}{67175}{"Audrey 19 May 2014"}{{0, 1, 67175, 4991950, 3794874, 3812536, 1615928, 2958345, 3335776, 1804684, 1574506, 5061832, 3652241, 5027796, 1078203, 1290719}}

\adfDesignE{16}{29188}{1}{7005121}{2417497}{0}{1166912}{"Audrey 19 May 2014"}{{0, 1, 1166912, 175280, 812402, 5622815, 4727993, 2529468, 2788498, 4037829, 6627028, 2277127, 3308741, 4041343, 2154890, 1760399}}

\adfDesignE{16}{31353}{1}{7524721}{6811537}{0}{4121839}{"Audrey 19 May 2014"}{{0, 1, 4121839, 877049, 5736128, 344176, 2409534, 5033709, 6640642, 2960994, 4272574, 5115186, 5893894, 7030, 6352320, 2907971}}

\adfDesignE{16}{35469}{1}{8512561}{6971337}{0}{6781767}{"Audrey 19 May 2014"}{{0, 1, 6781767, 1041487, 3441560, 7623427, 7014816, 421005, 2293630, 6242447, 7985112, 1497744, 1309789, 5177444, 7341115, 1416583}}

\adfDesignE{16}{37127}{1}{8910481}{5548083}{0}{3859148}{"Audrey 19 May 2014"}{{0, 1, 3859148, 790099, 4750819, 5951903, 5935983, 6066318, 378929, 8413658, 8172852, 2974497, 7895496, 7741453, 4656485, 3696207}}

\adfDesignE{16}{38023}{1}{9125521}{2389001}{0}{6137734}{"Audrey 19 May 2014"}{{0, 1, 6137734, 3624018, 1819174, 1887998, 4925724, 5600584, 6283840, 7564673, 2488099, 4199796, 6512724, 8343184, 8867195, 4749424}}

\adfDesignE{16}{38338}{1}{9201121}{8427699}{0}{4668136}{"Audrey 19 May 2014"}{{0, 1, 4668136, 391904, 2282514, 3633484, 2584157, 2226455, 1679821, 5114769, 6763876, 6616963, 2306530, 7129396, 1803838, 8004882}}

\adfDesignE{16}{38471}{1}{9233041}{4451261}{0}{6704227}{"Audrey 19 May 2014"}{{0, 1, 6704227, 3719586, 901664, 5093659, 6978264, 1009723, 5063069, 8981026, 6871267, 2254776, 1519091, 450386, 8583392, 6501156}}

\adfDesignE{16}{39163}{1}{9399121}{5507567}{0}{1521262}{"Audrey 19 May 2014"}{{0, 1, 1521262, 5298266, 3702320, 1047494, 4641330, 2868534, 5468512, 9092738, 4195923, 4757790, 5009325, 8031464, 6003184, 4155704}}

\adfDesignE{16}{40168}{1}{9640321}{5965755}{0}{4577590}{"Audrey 19 May 2014"}{{0, 1, 4577590, 3158327, 4170477, 6134130, 5284543, 4065749, 7811619, 9068318, 2733719, 4355777, 996982, 8310696, 6041847, 772472}}

\adfDesignE{16}{40195}{1}{9646801}{7882708}{0}{128327}{"Audrey 19 May 2014"}{{0, 1, 128327, 729622, 7998689, 8843301, 3911389, 4113372, 3031526, 93075, 1295887, 5735411, 5405102, 5885653, 1555037, 9154414}}

\adfDesignE{16}{41843}{1}{10042321}{9958324}{0}{5186217}{"Audrey 19 May 2014"}{{0, 1, 5186217, 6363060, 8355747, 3179010, 4991815, 4264942, 9781765, 4799429, 3524135, 5050505, 591162, 3939817, 6929466, 3339176}}

\adfDesignE{16}{42484}{1}{10196161}{505597}{0}{7053004}{"Audrey 19 May 2014"}{{0, 1, 7053004, 8471953, 3146097, 8760655, 5370561, 3042659, 5648614, 2749453, 7931810, 4825599, 100498, 6271755, 4300611, 3699857}}

\adfDesignE{16}{42987}{1}{10316881}{9702790}{0}{2781838}{"Audrey 19 May 2014"}{{0, 1, 2781838, 2438311, 8345834, 10251327, 580304, 8714920, 1397394, 5504420, 7408712, 9736576, 9137004, 6481176, 6783508, 2973723}}

\adfDesignE{16}{44715}{1}{10731601}{6146866}{0}{1552791}{"Audrey 19 May 2014"}{{0, 1, 1552791, 5240203, 4083151, 7327237, 2825065, 428847, 3188326, 5160137, 3916530, 7906535, 8749963, 2303072, 1488313, 10219435}}

\adfDesignE{16}{45526}{1}{10926241}{1191410}{0}{5335114}{"Audrey 19 May 2014"}{{0, 1, 5335114, 4497223, 2546738, 7234920, 2326421, 3151839, 1033092, 8729966, 10757978, 8599819, 2439288, 5395926, 10575778, 3859584}}

\adfDesignE{16}{46916}{1}{11259841}{8192620}{0}{1053660}{"Audrey 19 May 2014"}{{0, 1, 1053660, 1592682, 1133162, 7712803, 5625481, 6370286, 468409, 2476228, 5817483, 5634359, 3835895, 9198750, 7650451, 8989396}}

\adfDesignE{16}{47191}{1}{11325841}{245201}{0}{8252044}{"Audrey 19 May 2014"}{{0, 1, 8252044, 7569871, 2688012, 1653074, 78421, 8965307, 4223517, 36360, 140068, 11247419, 5434331, 10858294, 8601469, 9532699}}

\adfDesignE{16}{47989}{1}{11517361}{8444741}{0}{1927681}{"Audrey 19 May 2014"}{{0, 1, 1927681, 4202082, 9984293, 8155765, 8702081, 7220159, 8420829, 11335261, 6935819, 2815279, 2369521, 10411811, 1715168, 7943138}}

\adfDesignE{16}{48629}{1}{11670961}{7394792}{0}{10787341}{"Audrey 19 May 2014"}{{0, 1, 10787341, 8684461, 10138290, 434580, 6050183, 7718966, 5252051, 285498, 7750216, 5620777, 4835615, 9405410, 1247173, 3486165}}

~

\adfDesignE{17}{877}{17}{238561}{256}{0}{27838}{"Audrey 18 Apr 2014"}{{0, 1, 27838, 108116, 47632, 57578, 203566, 92314, 58040, 182428, 182657, 116412, 64632, 237115, 63061, 160280, 68257}}

\adfDesignE{17}{975}{17}{265217}{528}{0}{48238}{"Audrey 18 Apr 2014"}{{0, 1, 48238, 155903, 220879, 198661, 188674, 69840, 155586, 46802, 107772, 187319, 215949, 19753, 185750, 117372, 203237}}

\adfDesignE{17}{1041}{17}{283169}{239}{0}{12895}{"Audrey 18 Apr 2014"}{{0, 1, 12895, 60822, 204729, 279037, 236801, 137568, 168744, 83284, 170332, 172376, 195039, 204016, 146310, 195572, 280995}}

\adfDesignE{17}{1503}{17}{408833}{256}{0}{29799}{"Audrey 18 Apr 2014"}{{0, 1, 29799, 403958, 274423, 53311, 298284, 126663, 84481, 264538, 258789, 245365, 62263, 90983, 230794, 41680, 396499}}

\adfDesignE{17}{1831}{17}{498049}{86}{0}{42192}{"Audrey 18 Apr 2014"}{{0, 1, 42192, 137738, 205964, 74136, 198392, 343770, 160862, 175781, 104293, 67341, 379976, 248131, 153172, 447249, 249296}}

\adfDesignE{17}{1881}{17}{511649}{256}{0}{26607}{"Audrey 18 Apr 2014"}{{0, 1, 26607, 321882, 333412, 122722, 431985, 141759, 416934, 300969, 63684, 370349, 27752, 87957, 501022, 189308, 245200}}

\adfDesignE{17}{1987}{17}{540481}{256}{0}{137970}{"Audrey 18 Apr 2014"}{{0, 1, 137970, 520561, 523966, 93346, 366352, 342801, 383103, 381515, 179960, 465022, 207273, 65619, 396680, 293059, 507101}}

\adfDesignE{17}{2035}{17}{553537}{256}{0}{38184}{"Audrey 18 Apr 2014"}{{0, 1, 38184, 1398, 241680, 293793, 211070, 551697, 40639, 195365, 352548, 226529, 214174, 64378, 505272, 327450, 57044}}

\adfDesignE{17}{2437}{17}{662881}{222}{0}{46929}{"Audrey 18 Apr 2014"}{{0, 1, 46929, 240359, 224415, 381088, 212253, 371131, 271305, 116978, 343001, 617487, 204508, 164814, 60698, 96785, 625534}}

\adfDesignE{17}{2766}{17}{752369}{955}{0}{134727}{"Audrey 18 Apr 2014"}{{0, 1, 134727, 462404, 645770, 208368, 403408, 317794, 369555, 265541, 396357, 599764, 726197, 278459, 570246, 124776, 515585}}

\adfDesignE{17}{3127}{17}{850561}{256}{0}{134421}{"Audrey 18 Apr 2014"}{{0, 1, 134421, 537918, 434307, 25890, 514639, 461767, 642371, 50032, 816206, 512775, 816518, 774638, 235856, 188662, 658487}}

\adfDesignE{17}{3513}{17}{955553}{307}{0}{25039}{"Audrey 18 Apr 2014"}{{0, 1, 25039, 108753, 695870, 335528, 63616, 929726, 227128, 562089, 761887, 228501, 540728, 57935, 105011, 644126, 447380}}

\adfDesignE{17}{3607}{17}{981121}{18}{0}{34110}{"Audrey 18 Apr 2014"}{{0, 1, 34110, 863715, 217262, 399907, 302507, 64213, 443358, 923407, 485307, 348258, 648433, 638927, 159197, 686056, 653189}}

\adfDesignE{17}{3781}{17}{1028449}{103}{0}{35113}{"Audrey 18 Apr 2014"}{{0, 1, 35113, 840867, 649079, 681087, 483234, 443840, 466223, 665466, 146378, 611061, 681855, 710344, 363724, 161130, 259741}}

\adfDesignE{17}{3855}{17}{1048577}{375}{0}{177461}{"Audrey 18 Apr 2014"}{{0, 1, 177461, 493480, 497548, 987920, 439605, 711259, 374178, 863533, 192625, 863502, 982796, 245700, 238886, 28913, 242632}}

\adfDesignE{17}{4477}{17}{1217761}{205}{0}{133833}{"Audrey 18 Apr 2014"}{{0, 1, 133833, 443101, 228716, 107932, 1000135, 867140, 542081, 214898, 582497, 1132825, 559247, 894930, 719057, 1210217, 1105478}}

\adfDesignE{17}{4635}{17}{1260737}{256}{0}{109471}{"Audrey 18 Apr 2014"}{{0, 1, 109471, 594656, 692718, 462365, 750576, 292795, 844694, 741609, 780461, 221115, 790502, 56762, 880966, 152171, 193560}}

\adfDesignE{17}{4855}{17}{1320577}{239}{0}{3912}{"Audrey 18 Apr 2014"}{{0, 1, 3912, 777397, 1208810, 1199060, 33216, 524446, 776671, 1009852, 695217, 620861, 267129, 432241, 588232, 718450, 388544}}

\adfDesignE{17}{4911}{17}{1335809}{256}{0}{689454}{"Audrey 18 Apr 2014"}{{0, 1, 689454, 521275, 1065636, 30463, 1248104, 838742, 974959, 942923, 1196394, 578803, 1134520, 362422, 873475, 732798, 832312}}

~

\adfDesignE{18}{44371}{1}{13577527}{10028178}{0}{10633052}{"Audrey 18 May 2014"}{{0, 1, 10633052, 3159775, 13263409, 11458810, 2175304, 2229215, 11083174, 481930, 10820128, 2825065, 6538310, 11056698, 9819523, 3761075, 9966709, 12925565}}

\adfDesignE{18}{47635}{1}{14576311}{9032769}{0}{11463134}{"Audrey 18 May 2014"}{{0, 1, 11463134, 8967874, 12435042, 8820916, 12027251, 10576378, 658774, 6892702, 8801865, 10152819, 12550035, 9456315, 9911505, 2675673, 7922494, 12450332}}

~

\adfDesignE{19}{4379}{19}{1497637}{191}{0}{98502}{"Audrey 18 Apr 2014"}{{0, 1, 98502, 951518, 1307302, 539433, 566243, 1070832, 519754, 87663, 1103521, 532082, 1334349, 426804, 879381, 458456, 584451, 426122, 1094682}}

\adfDesignE{19}{4642}{19}{1587583}{431}{0}{323370}{"Audrey 18 Apr 2014"}{{0, 1, 323370, 415022, 922818, 30482, 1249076, 839260, 542682, 416269, 919126, 10858, 1005447, 748322, 721531, 756292, 1333222, 1546243, 920643}}

\adfDesignE{19}{4720}{19}{1614259}{1572}{0}{268639}{"Audrey 18 Apr 2014"}{{0, 1, 268639, 1463726, 1381881, 830506, 1179203, 957075, 1411477, 1244775, 1359375, 241127, 745260, 657183, 1548402, 520017, 487262, 542626, 1304055}}

\adfDesignE{19}{4831}{19}{1652221}{476}{0}{320661}{"Audrey 18 Apr 2014"}{{0, 1, 320661, 807428, 1030324, 4120, 998741, 682487, 379131, 452190, 782630, 1418519, 670875, 969733, 952429, 392603, 1491488, 229582, 1634826}}

~

\noindent{\bf Part II}\label{app:PART-II}

~

\adfDesignF{11}{1836}{1}{201961}{115}{1}{32}{"Audrey 11 Apr 2014"}{{1, 32, 1024, 32768, 38771, 28906, 117148, 113438, 196679, 32937, 44179}}

\adfDesignF{11}{1862}{1}{204821}{213}{1}{11902}{"Audrey 11 Apr 2014"}{{1, 11902, 126293, 162788, 100937, 77009, 191964, 182094, 71787, 100483, 203668}}

\adfDesignF{11}{1863}{1}{204931}{282}{1}{11138}{"Audrey 11 Apr 2014"}{{1, 11138, 71789, 150051, 55733, 18155, 148424, 173066, 28122, 88268, 74977}}

\adfDesignF{11}{1910}{1}{210101}{23}{1}{36343}{"Audrey 11 Apr 2014"}{{1, 36343, 118763, 98866, 149837, 128373, 167234, 193635, 153911, 68550, 145093}}

\adfDesignF{11}{1949}{1}{214391}{123}{1}{91496}{"Audrey 11 Apr 2014"}{{1, 91496, 192639, 185052, 202958, 154312, 211447, 125263, 149370, 188834, 213756}}

\adfDesignF{11}{2058}{1}{226381}{462}{1}{7324}{"Audrey 11 Apr 2014"}{{1, 7324, 215060, 166823, 33395, 93500, 217856, 44056, 73219, 185748, 94923}}

\adfDesignF{11}{2204}{1}{242441}{33}{1}{74961}{"Audrey 11 Apr 2014"}{{1, 74961, 96464, 235079, 175275, 184162, 134701, 138893, 171869, 157369, 85872}}

\adfDesignF{11}{2414}{1}{265541}{3}{1}{6161}{"Audrey 11 Apr 2014"}{{1, 6161, 251099, 244614, 121679, 42076, 62220, 161757, 9504, 135124, 27929}}

\adfDesignF{11}{2681}{1}{294911}{169}{1}{9601}{"Audrey 11 Apr 2014"}{{1, 9601, 166969, 228084, 120309, 215233, 10656, 269050, 23601, 101553, 34587}}

\adfDesignF{11}{2691}{1}{296011}{76}{1}{10192}{"Audrey 11 Apr 2014"}{{1, 10192, 273014, 55288, 186363, 205120, 153358, 86656, 197139, 214031, 98893}}

\adfDesignF{11}{2844}{1}{312841}{454}{1}{9618}{"Audrey 11 Apr 2014"}{{1, 9618, 217829, 295986, 253089, 307022, 31397, 84781, 160012, 130537, 73933}}

\adfDesignF{11}{3123}{1}{343531}{94}{1}{13552}{"Audrey 21 Apr 2014"}{{1, 13552, 211150, 235101, 182258, 316057, 59956, 72897, 248519, 295095, 83069}}

\adfDesignF{11}{3149}{1}{346391}{38}{1}{16927}{"Audrey 21 Apr 2014"}{{1, 16927, 57972, 312732, 67302, 287346, 229711, 79122, 150488, 297353, 233001}}

\adfDesignF{11}{3152}{1}{346721}{145}{1}{14106}{"Audrey 21 Apr 2014"}{{1, 14106, 308103, 299904, 102903, 175612, 208048, 78544, 168069, 249837, 128478}}

\adfDesignF{11}{3338}{1}{367181}{38}{1}{20870}{"Audrey 21 Apr 2014"}{{1, 20870, 80234, 138220, 77464, 342918, 340970, 76120, 199394, 90507, 102026}}

\adfDesignF{11}{3675}{1}{404251}{227}{1}{8550}{"Audrey 21 Apr 2014"}{{1, 8550, 337320, 159366, 253430, 41140, 48630, 216472, 174522, 72659, 304914}}

\adfDesignF{11}{3711}{1}{408211}{23}{1}{11748}{"Audrey 21 Apr 2014"}{{1, 11748, 40186, 213212, 31880, 196753, 163562, 77199, 297221, 323625, 277457}}

\adfDesignF{11}{3834}{1}{421741}{140}{1}{25530}{"Audrey 21 Apr 2014"}{{1, 25530, 191055, 199485, 329475, 294246, 49688, 359453, 172671, 253698, 233403}}

\adfDesignF{11}{3869}{1}{425591}{48}{1}{42338}{"Audrey 21 Apr 2014"}{{1, 42338, 342543, 146618, 268149, 252437, 236514, 224684, 286751, 54972, 272948}}

\adfDesignF{11}{3906}{1}{429661}{65}{1}{15294}{"Audrey 21 Apr 2014"}{{1, 15294, 170852, 241947, 96886, 303356, 47186, 261865, 93129, 418372, 69756}}

\adfDesignF{11}{3977}{1}{437471}{12}{1}{86861}{"Audrey 21 Apr 2014"}{{1, 86861, 208455, 122536, 367537, 185132, 191634, 186795, 276047, 370428, 191929}}

\adfDesignF{11}{4068}{1}{447481}{698}{1}{19471}{"Audrey 21 Apr 2014"}{{1, 19471, 103434, 298914, 216608, 65943, 153164, 242860, 195333, 187824, 306372}}

\adfDesignF{11}{4199}{1}{461891}{28}{1}{64551}{"Audrey 21 Apr 2014"}{{1, 64551, 112890, 369974, 117519, 333076, 286608, 250894, 174461, 267540, 331941}}

\adfDesignF{11}{4203}{1}{462331}{133}{1}{91232}{"Audrey 21 Apr 2014"}{{1, 91232, 395162, 235197, 248663, 365308, 186990, 382442, 214967, 250655, 403369}}

\adfDesignF{11}{4247}{1}{467171}{245}{1}{51507}{"Audrey 21 Apr 2014"}{{1, 51507, 374111, 400211, 214773, 170802, 201513, 191984, 378502, 456684, 362938}}

\adfDesignF{11}{4613}{1}{507431}{292}{1}{25268}{"Audrey 21 Apr 2014"}{{1, 25268, 123626, 36532, 73587, 169132, 43494, 418277, 250368, 156347, 225661}}

\adfDesignF{11}{4632}{1}{509521}{175}{1}{18179}{"Audrey 21 Apr 2014"}{{1, 18179, 306433, 52414, 29836, 258300, 399685, 104155, 53709, 133675, 172176}}

\adfDesignF{11}{4727}{1}{519971}{69}{1}{26354}{"Audrey 21 Apr 2014"}{{1, 26354, 372031, 451769, 144239, 286596, 372209, 463042, 329440, 105973, 48201}}

\adfDesignF{11}{4817}{1}{529871}{10}{1}{13921}{"Audrey 21 Apr 2014"}{{1, 13921, 391326, 45495, 140050, 240641, 118899, 405846, 297564, 386837, 78904}}

\adfDesignF{11}{4823}{1}{530531}{100}{1}{9861}{"Audrey 21 Apr 2014"}{{1, 9861, 152148, 520291, 354781, 174027, 342993, 118848, 17149, 397431, 34594}}

\adfDesignF{11}{5006}{1}{550661}{225083}{1}{154582}{"Audrey 18 May 2014"}{{1, 154582, 211290, 274887, 275508, 455916, 59027, 58944, 444502, 528584, 290064}}

\adfDesignF{11}{5057}{1}{556271}{377302}{1}{323209}{"Audrey 18 May 2014"}{{1, 323209, 257778, 124306, 144979, 473655, 442069, 247987, 210706, 42108, 514557}}

\adfDesignF{11}{5117}{1}{562871}{246858}{1}{158567}{"Audrey 18 May 2014"}{{1, 158567, 45919, 501688, 39795, 389855, 267339, 202461, 285902, 429223, 493605}}

\adfDesignF{11}{5160}{1}{567601}{482257}{1}{211132}{"Audrey 18 May 2014"}{{1, 211132, 176889, 485351, 145595, 196183, 393782, 57148, 277079, 446363, 448482}}

\adfDesignF{11}{5298}{1}{582781}{489005}{1}{549247}{"Audrey 18 May 2014"}{{1, 549247, 344607, 496092, 117298, 300618, 21726, 499347, 526956, 142978, 497816}}

\adfDesignF{11}{5616}{1}{617761}{461301}{1}{602191}{"Audrey 18 May 2014"}{{1, 602191, 262588, 464899, 445968, 529641, 598980, 219217, 538596, 165855, 496391}}

\adfDesignF{11}{5664}{1}{623041}{384336}{1}{56643}{"Audrey 18 May 2014"}{{1, 56643, 391340, 118922, 402595, 264944, 34425, 439986, 486998, 510480, 408871}}

\adfDesignF{11}{5775}{1}{635251}{59049}{1}{175700}{"Audrey 18 May 2014"}{{1, 175700, 467655, 442905, 161000, 608221, 600727, 144999, 218196, 274601, 82250}}

\adfDesignF{11}{5997}{1}{659671}{165006}{1}{597233}{"Audrey 18 May 2014"}{{1, 597233, 507905, 451264, 505791, 510996, 79338, 423166, 163855, 59049, 659428}}

\adfDesignF{11}{6144}{1}{675841}{84224}{1}{419989}{"Audrey 18 May 2014"}{{1, 419989, 314167, 218210, 408408, 449235, 578127, 287897, 211505, 611610, 555897}}

\adfDesignF{11}{6149}{1}{676391}{420202}{1}{549113}{"Audrey 18 May 2014"}{{1, 549113, 124834, 479129, 155307, 362629, 198805, 284520, 161189, 478270, 588158}}

\adfDesignF{11}{6188}{1}{680681}{546698}{1}{445468}{"Audrey 18 May 2014"}{{1, 445468, 84370, 333745, 415683, 334723, 447547, 486182, 123277, 657599, 74810}}

\adfDesignF{11}{6306}{1}{693661}{132200}{1}{293860}{"Audrey 18 May 2014"}{{1, 293860, 535371, 419938, 688780, 160288, 568797, 50877, 239687, 83880, 426826}}

\adfDesignF{11}{6417}{1}{705871}{75391}{1}{822}{"Audrey 18 May 2014"}{{1, 822, 675684, 597642, 681379, 337835, 293067, 199063, 573585, 670913, 205235}}

\adfDesignF{11}{6653}{1}{731831}{572291}{1}{371355}{"Audrey 18 May 2014"}{{1, 371355, 497878, 432681, 365719, 577758, 232327, 236495, 221570, 636189, 87182}}

\adfDesignF{11}{6686}{1}{735461}{704070}{1}{69379}{"Audrey 18 May 2014"}{{1, 69379, 588857, 186714, 356013, 103703, 530935, 175180, 320195, 209400, 401467}}

\adfDesignF{11}{6726}{1}{739861}{720565}{1}{675756}{"Audrey 18 May 2014"}{{1, 675756, 263031, 569796, 164990, 368906, 197874, 204075, 734188, 395914, 124774}}

\adfDesignF{11}{6737}{1}{741071}{154504}{1}{736699}{"Audrey 18 May 2014"}{{1, 736699, 587609, 266609, 90135, 179552, 533916, 92898, 697923, 411022, 108991}}

\adfDesignF{11}{6761}{1}{743711}{728302}{1}{146060}{"Audrey 18 May 2014"}{{1, 146060, 173565, 27043, 51459, 158174, 255936, 122456, 417521, 302682, 576236}}

\adfDesignF{11}{6774}{1}{745141}{142702}{1}{439732}{"Audrey 18 May 2014"}{{1, 439732, 142324, 24578, 208032, 347418, 513874, 598095, 268885, 600363, 582403}}

\adfDesignF{11}{6803}{1}{748331}{689864}{1}{595796}{"Audrey 18 May 2014"}{{1, 595796, 567104, 113305, 474601, 277405, 404720, 348976, 723194, 572582, 412302}}

\adfDesignF{11}{7037}{1}{774071}{130944}{1}{348496}{"Audrey 18 May 2014"}{{1, 348496, 44329, 344237, 468043, 420350, 453134, 258038, 608707, 119335, 31614}}

\adfDesignF{11}{7154}{1}{786941}{257898}{1}{343342}{"Audrey 18 May 2014"}{{1, 343342, 754105, 525795, 93726, 480920, 139315, 55927, 727634, 301322, 512618}}

\adfDesignF{11}{7173}{1}{789031}{58492}{1}{708864}{"Audrey 18 May 2014"}{{1, 708864, 90394, 633937, 659231, 724803, 538801, 653297, 650088, 681885, 181916}}

\adfDesignF{11}{7359}{1}{809491}{59049}{1}{745766}{"Audrey 18 May 2014"}{{1, 745766, 468769, 341848, 757992, 97261, 315362, 789607, 254485, 286069, 799786}}

\adfDesignF{11}{7544}{1}{829841}{612985}{1}{326112}{"Audrey 18 May 2014"}{{1, 326112, 763189, 8289, 350231, 195678, 660559, 280941, 465628, 82633, 186103}}

\adfDesignF{11}{7608}{1}{836881}{691988}{1}{455026}{"Audrey 18 May 2014"}{{1, 455026, 116871, 577382, 95840, 659811, 101336, 45398, 536625, 720999, 38235}}

\adfDesignF{11}{7671}{1}{843811}{427866}{1}{702203}{"Audrey 18 May 2014"}{{1, 702203, 501060, 274888, 349148, 111750, 137494, 688473, 638356, 312171, 503511}}

\adfDesignF{11}{7743}{1}{851731}{486114}{1}{729493}{"Audrey 18 May 2014"}{{1, 729493, 211711, 717217, 55377, 384262, 694563, 257548, 380529, 481201, 346753}}

\adfDesignF{11}{7791}{1}{857011}{602086}{1}{803046}{"Audrey 18 May 2014"}{{1, 803046, 97847, 588427, 365528, 103667, 178153, 779764, 132851, 449811, 785960}}

\adfDesignF{11}{7859}{1}{864491}{397724}{1}{756066}{"Audrey 18 May 2014"}{{1, 756066, 632007, 249122, 832936, 559988, 826485, 636444, 688884, 639191, 230313}}

\adfDesignF{11}{8258}{1}{908381}{623145}{1}{720725}{"Audrey 18 May 2014"}{{1, 720725, 476490, 276295, 156198, 146220, 404547, 453481, 518306, 756458, 633184}}

\adfDesignF{11}{8276}{1}{910361}{529624}{1}{722859}{"Audrey 18 May 2014"}{{1, 722859, 678906, 456179, 313619, 579057, 769412, 439568, 564360, 823559, 114646}}

\adfDesignF{11}{8499}{1}{934891}{282558}{1}{768342}{"Audrey 18 May 2014"}{{1, 768342, 353431, 897305, 815469, 743544, 87095, 183601, 786770, 435612, 538176}}

\adfDesignF{11}{8567}{1}{942371}{231720}{1}{703673}{"Audrey 18 May 2014"}{{1, 703673, 42173, 738639, 307852, 549142, 939500, 198241, 487176, 644552, 100906}}

\adfDesignF{11}{8729}{1}{960191}{670733}{1}{579361}{"Audrey 18 May 2014"}{{1, 579361, 399496, 281888, 827333, 7586, 238339, 213860, 57011, 339762, 896127}}

\adfDesignF{11}{8832}{1}{971521}{157658}{1}{244055}{"Audrey 18 May 2014"}{{1, 244055, 833557, 170798, 25864, 266583, 95737, 13485, 540048, 18175, 706260}}

\adfDesignF{11}{8949}{1}{984391}{918858}{1}{604829}{"Audrey 18 May 2014"}{{1, 604829, 704603, 791776, 583842, 880325, 810217, 101010, 503048, 479730, 446965}}

\adfDesignF{11}{9015}{1}{991651}{979715}{1}{230218}{"Audrey 18 May 2014"}{{1, 230218, 548178, 953242, 108805, 756881, 666244, 717720, 198387, 779910, 990320}}

\adfDesignF{11}{9060}{1}{996601}{169864}{1}{753769}{"Audrey 18 May 2014"}{{1, 753769, 492256, 801352, 409194, 705297, 190549, 789662, 795626, 606831, 435069}}

\adfDesignF{11}{9149}{1}{1006391}{182762}{1}{74486}{"Audrey 18 May 2014"}{{1, 74486, 937004, 464094, 987616, 408840, 470971, 974819, 263775, 779548, 677192}}

\adfDesignF{11}{9312}{1}{1024321}{543750}{1}{182766}{"Audrey 18 May 2014"}{{1, 182766, 302946, 605623, 190279, 833764, 597659, 201996, 447775, 943676, 815120}}

\adfDesignF{11}{9551}{1}{1050611}{219325}{1}{424969}{"Audrey 18 May 2014"}{{1, 424969, 721283, 852311, 256832, 813251, 621492, 683847, 64589, 59755, 754725}}

\adfDesignF{11}{9564}{1}{1052041}{141908}{1}{79800}{"Audrey 18 May 2014"}{{1, 79800, 35827, 599203, 83909, 749276, 526606, 433096, 461909, 1029724, 208813}}

\adfDesignF{11}{9672}{1}{1063921}{199208}{1}{530002}{"Audrey 18 May 2014"}{{1, 530002, 377979, 749105, 556877, 407381, 617022, 176669, 300049, 170386, 370213}}

\adfDesignF{11}{9677}{1}{1064471}{511164}{1}{733150}{"Audrey 18 May 2014"}{{1, 733150, 33166, 1006318, 385013, 119054, 1011513, 422025, 1036593, 151971, 423551}}

\adfDesignF{11}{9731}{1}{1070411}{938631}{1}{77867}{"Audrey 18 May 2014"}{{1, 77867, 461785, 566283, 247627, 658266, 567887, 978619, 636994, 106880, 1049846}}

\adfDesignF{11}{9744}{1}{1071841}{244110}{1}{18194}{"Audrey 18 May 2014"}{{1, 18194, 894608, 592367, 163943, 917280, 427950, 269276, 894174, 199058, 982354}}

\adfDesignF{11}{10163}{1}{1117931}{689686}{1}{228107}{"Audrey 18 May 2014"}{{1, 228107, 940916, 189184, 940157, 335276, 24891, 957719, 803637, 253572, 916195}}

\adfDesignF{11}{10347}{1}{1138171}{744857}{1}{799353}{"Audrey 18 May 2014"}{{1, 799353, 571893, 279421, 199402, 843724, 979154, 221279, 291890, 368512, 1136226}}

\adfDesignF{11}{10947}{1}{1204171}{734403}{1}{1001031}{"Audrey 18 May 2014"}{{1, 1001031, 123601, 1066552, 1093895, 273527, 991844, 1109902, 1077418, 1020098, 671328}}

\adfDesignF{11}{10982}{1}{1208021}{546501}{1}{986881}{"Audrey 18 May 2014"}{{1, 986881, 1001499, 1041175, 947058, 986629, 1159813, 1139816, 711515, 308150, 173610}}

\adfDesignF{11}{10998}{1}{1209781}{1031128}{1}{44293}{"Audrey 18 May 2014"}{{1, 44293, 814848, 665891, 1059064, 1073258, 681980, 1128132, 766033, 381743, 643443}}

\adfDesignF{11}{11151}{1}{1226611}{69538}{1}{1186909}{"Audrey 18 May 2014"}{{1, 1186909, 53669, 1083280, 278933, 872753, 494533, 413111, 909370, 260434, 580062}}

\adfDesignF{11}{11157}{1}{1227271}{241870}{1}{297193}{"Audrey 18 May 2014"}{{1, 297193, 667192, 752941, 473183, 1055055, 720096, 882632, 57520, 1110872, 119670}}

\adfDesignF{11}{11388}{1}{1252681}{34762}{1}{455274}{"Audrey 18 May 2014"}{{1, 455274, 806092, 1040043, 1140230, 1055896, 651030, 180810, 665387, 60170, 208472}}

\adfDesignF{11}{11447}{1}{1259171}{814922}{1}{795892}{"Audrey 18 May 2014"}{{1, 795892, 475720, 355079, 1232912, 392230, 304011, 141794, 768544, 451210, 119291}}

\adfDesignF{11}{11607}{1}{1276771}{494546}{1}{1246761}{"Audrey 18 May 2014"}{{1, 1246761, 476545, 1273292, 986339, 625474, 612502, 487067, 870509, 36341, 1045795}}

\adfDesignF{11}{11628}{1}{1279081}{621228}{1}{610295}{"Audrey 18 May 2014"}{{1, 610295, 553392, 1265238, 16320, 1089734, 1045580, 479577, 293552, 116656, 925060}}

\adfDesignF{11}{11822}{1}{1300421}{1125595}{1}{1118471}{"Audrey 18 May 2014"}{{1, 1118471, 985103, 136422, 453148, 314063, 637153, 1243379, 131899, 246505, 1235961}}

\adfDesignF{11}{11879}{1}{1306691}{564377}{1}{1106453}{"Audrey 18 May 2014"}{{1, 1106453, 750000, 803321, 885084, 357029, 905090, 710207, 672337, 957524, 693100}}

\adfDesignF{11}{12008}{1}{1320881}{442166}{1}{1089992}{"Audrey 18 May 2014"}{{1, 1089992, 294042, 899181, 1111028, 192475, 680970, 1266624, 109069, 1084805, 1197099}}

\adfDesignF{11}{12057}{1}{1326271}{734947}{1}{587353}{"Audrey 18 May 2014"}{{1, 587353, 565444, 1056080, 114624, 581770, 1241828, 707208, 620850, 1224871, 1257597}}

\adfDesignF{11}{12263}{1}{1348931}{414032}{1}{755196}{"Audrey 18 May 2014"}{{1, 755196, 1065202, 1287742, 641223, 351811, 810196, 561850, 626550, 829068, 1164747}}

\adfDesignF{11}{12381}{1}{1361911}{206536}{1}{1297747}{"Audrey 18 May 2014"}{{1, 1297747, 1323854, 1344836, 623856, 191728, 106671, 526642, 278444, 829593, 286183}}

\adfDesignF{11}{12501}{1}{1375111}{807767}{1}{1221253}{"Audrey 18 May 2014"}{{1, 1221253, 1123410, 336476, 654520, 365703, 459724, 744426, 1225015, 1231345, 888793}}

\adfDesignF{11}{13044}{1}{1434841}{1026139}{1}{213688}{"Audrey 18 May 2014"}{{1, 213688, 181360, 835111, 589357, 1089205, 174907, 788648, 1103533, 1380718, 817677}}

\adfDesignF{11}{13349}{1}{1468391}{1240163}{1}{791513}{"Audrey 18 May 2014"}{{1, 791513, 872237, 870066, 512813, 1110676, 548216, 271571, 91997, 780162, 24312}}

\adfDesignF{11}{13356}{1}{1469161}{235907}{1}{305783}{"Audrey 18 May 2014"}{{1, 305783, 1429566, 1347077, 169238, 376290, 1333872, 1028312, 206949, 314314, 834403}}

\adfDesignF{11}{13845}{1}{1522951}{466092}{1}{1190180}{"Audrey 18 May 2014"}{{1, 1190180, 1248280, 1116125, 211603, 1343474, 694351, 726148, 1070209, 1480407, 56928}}

\adfDesignF{11}{14100}{1}{1551001}{1272531}{1}{1125146}{"Audrey 18 May 2014"}{{1, 1125146, 138099, 706273, 225505, 897142, 1016917, 453178, 1386239, 538272, 1268233}}

\adfDesignF{11}{14226}{1}{1564861}{576124}{1}{1364971}{"Audrey 18 May 2014"}{{1, 1364971, 416187, 1085913, 284601, 62904, 1342436, 1267379, 523841, 531464, 944608}}

\adfDesignF{11}{14358}{1}{1579381}{1224245}{1}{383056}{"Audrey 18 May 2014"}{{1, 383056, 1086712, 419226, 312719, 737319, 80158, 256827, 1060203, 1407552, 573751}}

\adfDesignF{11}{14378}{1}{1581581}{1026300}{1}{1324400}{"Audrey 18 May 2014"}{{1, 1324400, 349341, 1004146, 1181159, 1028310, 691224, 24856, 259466, 340206, 123215}}

\adfDesignF{11}{14430}{1}{1587301}{336099}{1}{930329}{"Audrey 18 May 2014"}{{1, 930329, 1257369, 610548, 1196646, 285271, 744260, 568524, 274380, 273404, 208472}}

\adfDesignF{11}{14625}{1}{1608751}{361377}{1}{91515}{"Audrey 18 May 2014"}{{1, 91515, 1446270, 236778, 471451, 1353947, 457685, 1210490, 1007241, 1054068, 714309}}

\adfDesignF{11}{14681}{1}{1614911}{194052}{1}{431356}{"Audrey 18 May 2014"}{{1, 431356, 1183138, 1426353, 982778, 130180, 238788, 383126, 166760, 1560798, 1571277}}

\adfDesignF{11}{14763}{1}{1623931}{1420273}{1}{1197997}{"Audrey 18 May 2014"}{{1, 1197997, 696760, 1240341, 325150, 1267373, 148052, 108024, 1366538, 848252, 921167}}

\adfDesignF{11}{14967}{1}{1646371}{335817}{1}{391473}{"Audrey 18 May 2014"}{{1, 391473, 311565, 1182452, 1068694, 572339, 635957, 711154, 1192855, 440459, 78535}}

\adfDesignF{11}{15044}{1}{1654841}{1605721}{1}{886332}{"Audrey 18 May 2014"}{{1, 886332, 1604386, 542124, 561567, 469, 324617, 1159220, 1065483, 255204, 219961}}

\adfDesignF{11}{15192}{1}{1671121}{910683}{1}{637065}{"Audrey 18 May 2014"}{{1, 637065, 25923, 618273, 211287, 1440589, 929384, 1518781, 1541096, 1420224, 12982}}

\adfDesignF{11}{15203}{1}{1672331}{290256}{1}{1658166}{"Audrey 18 May 2014"}{{1, 1658166, 1639836, 400650, 684164, 1647416, 59134, 204721, 1621320, 123823, 322424}}

\adfDesignF{11}{15218}{1}{1673981}{1177799}{1}{145216}{"Audrey 18 May 2014"}{{1, 145216, 547999, 514006, 756487, 687048, 1094768, 1528299, 414566, 237553, 769981}}

\adfDesignF{11}{15384}{1}{1692241}{1497167}{1}{29424}{"Audrey 18 May 2014"}{{1, 29424, 1036625, 702216, 1433215, 272440, 128943, 14510, 497508, 790742, 171099}}

\adfDesignF{11}{15431}{1}{1697411}{276722}{1}{1115918}{"Audrey 18 May 2014"}{{1, 1115918, 1350794, 1382219, 594109, 441271, 1623267, 1675003, 768308, 1338411, 1592576}}

\adfDesignF{11}{15480}{1}{1702801}{1368925}{1}{670125}{"Audrey 18 May 2014"}{{1, 670125, 1430303, 656990, 1116797, 630518, 1348615, 1132538, 215948, 1313316, 1701655}}

\adfDesignF{11}{15536}{1}{1708961}{189525}{1}{127195}{"Audrey 18 May 2014"}{{1, 127195, 1543199, 1063228, 365686, 639233, 3938, 168337, 52346, 37414, 1126306}}

\adfDesignF{11}{15555}{1}{1711051}{999497}{1}{912734}{"Audrey 18 May 2014"}{{1, 912734, 1710723, 57173, 107584, 68817, 644519, 1382738, 768092, 1601502, 1301372}}

\adfDesignF{11}{15567}{1}{1712371}{489958}{1}{644871}{"Audrey 18 May 2014"}{{1, 644871, 747436, 1611676, 1247717, 774543, 534334, 1221697, 754552, 1359432, 1377967}}

\adfDesignF{11}{15653}{1}{1721831}{92834}{1}{835674}{"Audrey 18 May 2014"}{{1, 835674, 486310, 1461165, 484588, 159422, 1590465, 1350214, 496133, 596659, 1148524}}

\adfDesignF{11}{15776}{1}{1735361}{1351855}{1}{621056}{"Audrey 18 May 2014"}{{1, 621056, 542471, 149475, 944266, 1110000, 2750, 308776, 1120151, 1511054, 631444}}

\adfDesignF{11}{15786}{1}{1736461}{558386}{1}{1012555}{"Audrey 18 May 2014"}{{1, 1012555, 277490, 605462, 809977, 102786, 1688195, 715215, 1727814, 1409738, 333072}}

\adfDesignF{11}{15809}{1}{1738991}{1237348}{1}{1633435}{"Audrey 18 May 2014"}{{1, 1633435, 353799, 963472, 1160221, 153299, 1421002, 1381593, 1571525, 197581, 1598018}}

\adfDesignF{11}{15933}{1}{1752631}{91610}{1}{35056}{"Audrey 18 May 2014"}{{1, 35056, 328805, 1286624, 1684790, 86171, 1027363, 422909, 1744906, 850205, 1296325}}

\adfDesignF{11}{15938}{1}{1753181}{1151348}{1}{59525}{"Audrey 18 May 2014"}{{1, 59525, 46824, 1393991, 1010726, 1305954, 866310, 789997, 750643, 453609, 335144}}

\adfDesignF{11}{16002}{1}{1760221}{446744}{1}{932713}{"Audrey 18 May 2014"}{{1, 932713, 1275760, 1499996, 1633265, 75484, 1348755, 1297372, 1242460, 1256641, 1759100}}

\adfDesignF{11}{16299}{1}{1792891}{1551334}{1}{1692893}{"Audrey 18 May 2014"}{{1, 1692893, 646897, 893965, 626081, 905882, 1222430, 546131, 1245013, 1333957, 1644096}}

\adfDesignF{11}{16338}{1}{1797181}{1079797}{1}{870009}{"Audrey 18 May 2014"}{{1, 870009, 532673, 225492, 1588649, 705162, 797212, 943121, 1003548, 101598, 421259}}

\adfDesignF{11}{16467}{1}{1811371}{1712744}{1}{1137150}{"Audrey 18 May 2014"}{{1, 1137150, 1347536, 1339869, 939442, 1441114, 950432, 259714, 601776, 785165, 254656}}

\adfDesignF{11}{16578}{1}{1823581}{1508708}{1}{1144864}{"Audrey 18 May 2014"}{{1, 1144864, 146098, 1667771, 1433580, 232243, 1245828, 689728, 1359534, 442446, 755812}}

\adfDesignF{11}{16599}{1}{1825891}{1024}{1}{995761}{"Audrey 18 May 2014"}{{1, 995761, 816917, 1389427, 1179735, 1481210, 137593, 360306, 210821, 989729, 1567955}}

\adfDesignF{11}{16607}{1}{1826771}{1324748}{1}{873028}{"Audrey 18 May 2014"}{{1, 873028, 1531538, 1578721, 1799566, 949802, 1328449, 108176, 269770, 312385, 382419}}

\adfDesignF{11}{16638}{1}{1830181}{1735636}{1}{1436642}{"Audrey 18 May 2014"}{{1, 1436642, 1198120, 754169, 48136, 810227, 40648, 1038849, 65350, 1759943, 168639}}

\adfDesignF{11}{16857}{1}{1854271}{1731245}{1}{709078}{"Audrey 18 May 2014"}{{1, 709078, 465621, 1238804, 1550321, 1368772, 1731125, 1210544, 258467, 826328, 1766565}}

\adfDesignF{11}{16931}{1}{1862411}{597014}{1}{606621}{"Audrey 18 May 2014"}{{1, 606621, 835384, 1306775, 539235, 1131717, 1355437, 715987, 1743028, 1479303, 1460978}}

\adfDesignF{11}{16932}{1}{1862521}{120607}{1}{1031885}{"Audrey 18 May 2014"}{{1, 1031885, 160214, 1334388, 1123895, 28568, 770813, 779455, 579077, 1295362, 346426}}

\adfDesignF{11}{16949}{1}{1864391}{949155}{1}{1290839}{"Audrey 18 May 2014"}{{1, 1290839, 1291100, 745308, 267637, 896361, 1767151, 804106, 1525940, 721423, 12089}}

\adfDesignF{11}{17004}{1}{1870441}{527663}{1}{1254179}{"Audrey 18 May 2014"}{{1, 1254179, 771122, 406701, 981456, 1116493, 735330, 1415933, 1210228, 246163, 1214599}}

\adfDesignF{11}{17207}{1}{1892771}{1711122}{1}{1871670}{"Audrey 18 May 2014"}{{1, 1871670, 451016, 1856743, 1225657, 234587, 1468649, 377034, 1414850, 1819904, 636515}}

\adfDesignF{11}{17309}{1}{1903991}{1891837}{1}{1833197}{"Audrey 18 May 2014"}{{1, 1833197, 486124, 1878860, 796420, 1128003, 1346140, 1826363, 678606, 267748, 1182484}}

\adfDesignF{11}{17774}{1}{1955141}{628343}{1}{135836}{"Audrey 18 May 2014"}{{1, 135836, 753279, 102009, 410257, 285929, 575679, 113208, 537923, 1779176, 1172126}}

\adfDesignF{11}{17906}{1}{1969661}{1154145}{1}{1909498}{"Audrey 18 May 2014"}{{1, 1909498, 1319312, 1600783, 636627, 637605, 888021, 1047202, 732481, 940811, 135964}}

\adfDesignF{11}{17976}{1}{1977361}{1865027}{1}{1550788}{"Audrey 18 May 2014"}{{1, 1550788, 1833026, 324498, 1094090, 1331177, 524032, 781353, 1854252, 222019, 371569}}

\adfDesignF{11}{18284}{1}{2011241}{1228610}{1}{547804}{"Audrey 18 May 2014"}{{1, 547804, 2009011, 1231608, 950418, 869366, 415874, 150144, 1794322, 1056127, 1031530}}

\adfDesignF{11}{18402}{1}{2024221}{1248088}{1}{4819}{"Audrey 18 May 2014"}{{1, 4819, 956330, 1427274, 1754669, 578794, 1855969, 906233, 892130, 1753287, 2015820}}

\adfDesignF{11}{18444}{1}{2028841}{172613}{1}{1265402}{"Audrey 18 May 2014"}{{1, 1265402, 1779605, 1629419, 1558799, 1200245, 22049, 227266, 725705, 343103, 1392611}}

\adfDesignF{11}{18455}{1}{2030051}{1252753}{1}{931812}{"Audrey 18 May 2014"}{{1, 931812, 490134, 2019083, 1192769, 1815387, 524015, 1388303, 1605643, 1707912, 505247}}

\adfDesignF{11}{18497}{1}{2034671}{518621}{1}{372595}{"Audrey 18 May 2014"}{{1, 372595, 1431695, 494429, 626244, 1347571, 419404, 1031238, 746957, 2005351, 1697870}}

\adfDesignF{11}{18525}{1}{2037751}{894069}{1}{1274368}{"Audrey 18 May 2014"}{{1, 1274368, 1649211, 2021766, 624267, 1784603, 801350, 1996403, 1635045, 1760787, 716456}}

\adfDesignF{11}{18531}{1}{2038411}{5144}{1}{1325867}{"Audrey 18 May 2014"}{{1, 1325867, 1732111, 1599841, 1823725, 667689, 1221351, 634930, 1605886, 1657688, 1999788}}

\adfDesignF{11}{18641}{1}{2050511}{547282}{1}{1727992}{"Audrey 18 May 2014"}{{1, 1727992, 183353, 1911033, 194764, 263458, 914627, 1927047, 612707, 146648, 320414}}

\adfDesignF{11}{18935}{1}{2082851}{1825301}{1}{1588367}{"Audrey 18 May 2014"}{{1, 1588367, 215962, 100013, 385852, 1925487, 749667, 1235599, 1659275, 58224, 412957}}

\adfDesignF{11}{18993}{1}{2089231}{1014557}{1}{987913}{"Audrey 18 May 2014"}{{1, 987913, 369305, 890166, 1183345, 167549, 330800, 2018149, 400506, 249505, 1759685}}

\adfDesignF{11}{19095}{1}{2100451}{1340464}{1}{1075283}{"Audrey 18 May 2014"}{{1, 1075283, 368570, 1860179, 1377377, 63571, 1838700, 1933016, 148360, 1832881, 4317}}

\adfDesignF{11}{19245}{1}{2116951}{765651}{1}{535896}{"Audrey 18 May 2014"}{{1, 535896, 1067107, 48389, 938745, 1688782, 147515, 1514198, 2046647, 1861514, 735961}}

\adfDesignF{11}{19592}{1}{2155121}{1781513}{1}{379550}{"Audrey 18 May 2014"}{{1, 379550, 1294376, 1182761, 923008, 1992245, 60085, 1926449, 730433, 1079710, 1206987}}

\adfDesignF{11}{19715}{1}{2168651}{146860}{1}{1107223}{"Audrey 18 May 2014"}{{1, 1107223, 24127, 526303, 913661, 660876, 1730183, 1033100, 1930793, 1297757, 1619231}}

\adfDesignF{11}{19721}{1}{2169311}{1083963}{1}{145117}{"Audrey 18 May 2014"}{{1, 145117, 1441812, 1386054, 1482398, 1425251, 1700005, 1240043, 464648, 1799314, 1931223}}

\adfDesignF{11}{19745}{1}{2171951}{417977}{1}{380297}{"Audrey 18 May 2014"}{{1, 380297, 2106972, 1139715, 2169648, 1644213, 1953969, 1125114, 959907, 1260005, 291865}}

\adfDesignF{11}{19958}{1}{2195381}{575509}{1}{1684936}{"Audrey 18 May 2014"}{{1, 1684936, 1890183, 361969, 194136, 1552439, 1218881, 1458736, 680869, 198643, 1736112}}

\adfDesignF{11}{19992}{1}{2199121}{1776006}{1}{1920026}{"Audrey 18 May 2014"}{{1, 1920026, 1153205, 803601, 732332, 504442, 486430, 354964, 1723470, 1876680, 1440454}}

\adfDesignF{11}{20013}{1}{2201431}{1428494}{1}{446149}{"Audrey 18 May 2014"}{{1, 446149, 2143474, 550933, 1831574, 1332774, 468502, 228210, 1680971, 2031909, 292658}}

\adfDesignF{11}{20030}{1}{2203301}{613378}{1}{1277979}{"Audrey 18 May 2014"}{{1, 1277979, 408676, 63560, 1450574, 731071, 1222867, 1353095, 1457971, 1780444, 1270267}}

\adfDesignF{11}{20250}{1}{2227501}{2142114}{1}{1005196}{"Audrey 18 May 2014"}{{1, 1005196, 42305, 1822690, 1029722, 1725834, 1379654, 607093, 1281268, 2210336, 33406}}

\adfDesignF{11}{20286}{1}{2231461}{1507665}{1}{1436533}{"Audrey 18 May 2014"}{{1, 1436533, 936282, 261322, 1224057, 311998, 2018162, 2014848, 1050799, 2174502, 1960262}}

\adfDesignF{11}{20327}{1}{2235971}{1529063}{1}{908011}{"Audrey 18 May 2014"}{{1, 908011, 973465, 1816279, 528802, 1148340, 523368, 1804563, 1021944, 1720471, 734611}}

\adfDesignF{11}{20379}{1}{2241691}{10109}{1}{428887}{"Audrey 18 May 2014"}{{1, 428887, 2103764, 886550, 867503, 578818, 412535, 952988, 1127708, 1017600, 590410}}

\adfDesignF{11}{20523}{1}{2257531}{544596}{1}{1795798}{"Audrey 18 May 2014"}{{1, 1795798, 650711, 58627, 33630, 1474959, 1162947, 801447, 2211400, 400038, 440566}}

\adfDesignF{11}{21026}{1}{2312861}{1499584}{1}{566051}{"Audrey 18 May 2014"}{{1, 566051, 1535966, 573173, 1634465, 2203356, 1572906, 2233673, 1112453, 971521, 1473601}}

\adfDesignF{11}{21111}{1}{2322211}{284683}{1}{2311112}{"Audrey 18 May 2014"}{{1, 2311112, 110618, 700437, 612165, 370051, 795210, 688221, 1509311, 587365, 1604353}}

\adfDesignF{11}{21161}{1}{2327711}{559011}{1}{74910}{"Audrey 18 May 2014"}{{1, 74910, 1724590, 1076400, 1214960, 1481211, 188062, 417448, 560106, 549685, 2023471}}

\adfDesignF{11}{21264}{1}{2339041}{2147321}{1}{1249957}{"Audrey 18 May 2014"}{{1, 1249957, 1997407, 1549468, 322138, 2096080, 924599, 529348, 1337079, 2019324, 2008845}}

\adfDesignF{11}{21327}{1}{2345971}{876088}{1}{1226530}{"Audrey 18 May 2014"}{{1, 1226530, 823411, 124301, 1288153, 787923, 173595, 1493361, 19515, 2136808, 1310286}}

\adfDesignF{11}{21398}{1}{2353781}{1200788}{1}{1558034}{"Audrey 18 May 2014"}{{1, 1558034, 1477170, 1097600, 900908, 948456, 1799675, 1453795, 2106263, 2019428, 761356}}

\adfDesignF{11}{21561}{1}{2371711}{856197}{1}{559101}{"Audrey 18 May 2014"}{{1, 559101, 46690, 1374424, 353691, 472033, 1980808, 1282158, 1426786, 1971380, 19772}}

\adfDesignF{11}{22058}{1}{2426381}{993764}{1}{1290968}{"Audrey 18 May 2014"}{{1, 1290968, 2191459, 1429456, 310239, 467968, 1666120, 918233, 981994, 1843598, 1031869}}

\adfDesignF{11}{22142}{1}{2435621}{1964137}{1}{1417200}{"Audrey 18 May 2014"}{{1, 1417200, 1357843, 2095541, 1743101, 1444813, 1378836, 2264626, 433016, 950524, 1528225}}

\adfDesignF{11}{22143}{1}{2435731}{59049}{1}{145235}{"Audrey 18 May 2014"}{{1, 145235, 2210496, 2297836, 1835688, 774344, 1714839, 1147415, 1845429, 348768, 2294335}}

\adfDesignF{11}{22293}{1}{2452231}{1989777}{1}{998027}{"Audrey 18 May 2014"}{{1, 998027, 896225, 586363, 795499, 76375, 1615952, 60472, 831603, 2213100, 1735307}}

\adfDesignF{11}{22353}{1}{2458831}{239887}{1}{960258}{"Audrey 18 May 2014"}{{1, 960258, 1836761, 709080, 262120, 2132814, 389365, 1014310, 2238598, 1095365, 1655483}}

\adfDesignF{11}{22419}{1}{2466091}{2455588}{1}{899668}{"Audrey 18 May 2014"}{{1, 899668, 1850932, 1275008, 1397422, 198005, 678955, 608877, 1757279, 1132910, 65307}}

\adfDesignF{11}{22532}{1}{2478521}{1528759}{1}{1720741}{"Audrey 18 May 2014"}{{1, 1720741, 1826078, 323023, 843541, 1981004, 602950, 1981266, 346270, 2059149, 708582}}

\adfDesignF{11}{22560}{1}{2481601}{466789}{1}{236093}{"Audrey 18 May 2014"}{{1, 236093, 664588, 388257, 1863764, 1515939, 1126905, 1938955, 211148, 163876, 1816878}}

\adfDesignF{11}{22797}{1}{2507671}{542094}{1}{1197210}{"Audrey 18 May 2014"}{{1, 1197210, 2270630, 256118, 1559255, 248072, 771906, 1649998, 1352040, 1762281, 1470844}}

\adfDesignF{11}{23016}{1}{2531761}{1217861}{1}{1735247}{"Audrey 18 May 2014"}{{1, 1735247, 563206, 1267706, 1726268, 237587, 167749, 1791550, 2049818, 1148599, 1971074}}

\adfDesignF{11}{23214}{1}{2553541}{891388}{1}{2290780}{"Audrey 18 May 2014"}{{1, 2290780, 701563, 1836429, 722301, 1955305, 2212718, 2309433, 2219350, 1193443, 2433464}}

\adfDesignF{11}{23229}{1}{2555191}{1296054}{1}{1817794}{"Audrey 18 May 2014"}{{1, 1817794, 2025236, 1425977, 1038451, 1762788, 2182684, 2511779, 445716, 2021887, 98833}}

\adfDesignF{11}{24285}{1}{2671351}{2414416}{1}{906536}{"Audrey 18 May 2014"}{{1, 906536, 1111709, 2337711, 1703233, 1152888, 453430, 1826057, 949521, 1564632, 1351037}}

\adfDesignF{11}{24489}{1}{2693791}{1566596}{1}{2361976}{"Audrey 18 May 2014"}{{1, 2361976, 568473, 2152489, 914414, 1561666, 1514343, 1887849, 236996, 1012923, 1257825}}

\adfDesignF{11}{24807}{1}{2728771}{2702182}{1}{2074927}{"Audrey 18 May 2014"}{{1, 2074927, 881308, 2081660, 57279, 851499, 906203, 1754295, 901099, 1687138, 2448446}}

\adfDesignF{11}{24999}{1}{2749891}{2130037}{1}{869195}{"Audrey 18 May 2014"}{{1, 869195, 394467, 1334621, 1631854, 59948, 1567192, 1195007, 2530954, 1761158, 2405058}}

\adfDesignF{11}{25016}{1}{2751761}{1402597}{1}{2732251}{"Audrey 18 May 2014"}{{1, 2732251, 897082, 1881901, 858513, 380577, 1945669, 540805, 1897885, 2711427, 2664455}}

\adfDesignF{11}{25194}{1}{2771341}{1879612}{1}{680817}{"Audrey 18 May 2014"}{{1, 680817, 2233898, 2592640, 1962524, 779188, 2656399, 2659544, 1582416, 1842191, 2409769}}

\adfDesignF{11}{25401}{1}{2794111}{260698}{1}{1256076}{"Audrey 18 May 2014"}{{1, 1256076, 1406405, 34029, 1494237, 1022537, 1210887, 158895, 847290, 518806, 433170}}

\adfDesignF{11}{25733}{1}{2830631}{2500215}{1}{2427481}{"Audrey 18 May 2014"}{{1, 2427481, 751742, 2381977, 369831, 459013, 1410675, 2601115, 1709272, 1465102, 576946}}

\adfDesignF{11}{25814}{1}{2839541}{1395332}{1}{977335}{"Audrey 18 May 2014"}{{1, 977335, 1863399, 2724987, 2628499, 2685629, 1099955, 2692735, 525179, 386805, 1452722}}

\adfDesignF{11}{25817}{1}{2839871}{2633896}{1}{570169}{"Audrey 18 May 2014"}{{1, 570169, 1295707, 2242801, 1571166, 1359717, 738399, 1343681, 1892935, 2722336, 462443}}

\adfDesignF{11}{25880}{1}{2846801}{1098032}{1}{837399}{"Audrey 18 May 2014"}{{1, 837399, 1675677, 2103616, 2484397, 779807, 1848210, 220931, 2341882, 2452844, 2336042}}

\adfDesignF{11}{26183}{1}{2880131}{460355}{1}{711328}{"Audrey 18 May 2014"}{{1, 711328, 349242, 2794102, 2186976, 586574, 1932302, 1199271, 400605, 1392300, 2847954}}

\adfDesignF{11}{26207}{1}{2882771}{65947}{1}{34172}{"Audrey 18 May 2014"}{{1, 34172, 203329, 680478, 863330, 2317117, 2333838, 52421, 1129621, 1105122, 2811655}}

\adfDesignF{11}{26220}{1}{2884201}{1098927}{1}{748974}{"Audrey 18 May 2014"}{{1, 748974, 2263382, 2058710, 451131, 1242444, 2526017, 568598, 1303998, 918628, 2339122}}

\adfDesignF{11}{26577}{1}{2923471}{2614431}{1}{357925}{"Audrey 18 May 2014"}{{1, 357925, 882934, 2783792, 2594967, 2209420, 1901058, 1232871, 792993, 990548, 870846}}

\adfDesignF{11}{26589}{1}{2924791}{1204488}{1}{1607820}{"Audrey 18 May 2014"}{{1, 1607820, 2777468, 1182857, 2117109, 2855551, 863033, 1898303, 1903693, 1897760, 445151}}

\adfDesignF{11}{26636}{1}{2929961}{2273282}{1}{1633374}{"Audrey 18 May 2014"}{{1, 1633374, 2405755, 981947, 2608090, 1829281, 1975280, 490155, 449642, 737965, 1538315}}

\adfDesignF{11}{26762}{1}{2943821}{2533905}{1}{1713546}{"Audrey 18 May 2014"}{{1, 1713546, 2177012, 233352, 179762, 801896, 1298867, 1316795, 53527, 345845, 710860}}

\adfDesignF{11}{26894}{1}{2958341}{1052402}{1}{2201719}{"Audrey 18 May 2014"}{{1, 2201719, 2367292, 100872, 265075, 2109786, 1425685, 1159442, 1198534, 824628, 180330}}

\adfDesignF{11}{26960}{1}{2965601}{467076}{1}{1771750}{"Audrey 18 May 2014"}{{1, 1771750, 507197, 1732134, 703665, 562557, 1524660, 819317, 1758063, 284324, 2198736}}

\adfDesignF{11}{27272}{1}{2999921}{2119342}{1}{1869638}{"Audrey 18 May 2014"}{{1, 1869638, 2302792, 2807410, 1640641, 535063, 2461008, 2243013, 2774184, 412600, 953176}}

\adfDesignF{11}{27285}{1}{3001351}{1828321}{1}{831841}{"Audrey 18 May 2014"}{{1, 831841, 977582, 745820, 392112, 16916, 1088868, 2334453, 1612867, 2982233, 1022711}}

\adfDesignF{11}{27473}{1}{3022031}{700267}{1}{940980}{"Audrey 18 May 2014"}{{1, 940980, 365524, 1337286, 782035, 2657676, 1671081, 409150, 1261662, 2896503, 2788257}}

\adfDesignF{11}{27569}{1}{3032591}{2401338}{1}{2485323}{"Audrey 18 May 2014"}{{1, 2485323, 544073, 1204771, 1189228, 957997, 2327657, 2821429, 2292770, 2167209, 2205088}}

\adfDesignF{11}{27656}{1}{3042161}{1877108}{1}{119335}{"Audrey 18 May 2014"}{{1, 119335, 486584, 774633, 1724909, 275772, 2196083, 2605460, 1546256, 184305, 2255306}}

\adfDesignF{11}{27860}{1}{3064601}{2487891}{1}{97065}{"Audrey 18 May 2014"}{{1, 97065, 1030751, 2881569, 2555518, 2549730, 1559493, 2350852, 1388122, 2878965, 1095540}}

\adfDesignF{11}{28124}{1}{3093641}{1048744}{1}{1941205}{"Audrey 18 May 2014"}{{1, 1941205, 2465514, 1297064, 2210476, 1812786, 452399, 1142843, 2367382, 1006220, 772315}}

\adfDesignF{11}{28223}{1}{3104531}{615239}{1}{2307698}{"Audrey 18 May 2014"}{{1, 2307698, 1045238, 340895, 763372, 71547, 596633, 1700458, 819029, 2067132, 2706121}}

\adfDesignF{11}{28716}{1}{3158761}{352910}{1}{2617238}{"Audrey 18 May 2014"}{{1, 2617238, 423333, 2523416, 1682315, 2835504, 2042074, 1647222, 237206, 1670088, 115408}}

\adfDesignF{11}{28755}{1}{3163051}{1552019}{1}{929462}{"Audrey 18 May 2014"}{{1, 929462, 794222, 82, 302660, 1865184, 6724, 2676763, 1118640, 551368, 1244047}}

\adfDesignF{11}{28934}{1}{3182741}{447321}{1}{3084819}{"Audrey 18 May 2014"}{{1, 3084819, 2302192, 1482747, 3092686, 2173140, 3030921, 3117570, 278957, 1438649, 2277505}}

\adfDesignF{11}{29139}{1}{3205291}{439111}{1}{2486250}{"Audrey 18 May 2014"}{{1, 2486250, 110799, 1689337, 153871, 677027, 3015391, 489300, 2005315, 2864017, 2535147}}

\adfDesignF{11}{29156}{1}{3207161}{1689116}{1}{1910012}{"Audrey 18 May 2014"}{{1, 1910012, 202644, 2660765, 101292, 155340, 385648, 433745, 361225, 367414, 3043397}}

\adfDesignF{11}{29307}{1}{3223771}{924198}{1}{1602658}{"Audrey 18 May 2014"}{{1, 1602658, 2134653, 1342138, 1149787, 1382704, 784458, 2003471, 2909918, 2744085, 64982}}

\adfDesignF{11}{29393}{1}{3233231}{333682}{1}{424892}{"Audrey 18 May 2014"}{{1, 424892, 2525548, 1637764, 679513, 1808989, 1881482, 585501, 198059, 2381391, 809784}}

\adfDesignF{11}{29567}{1}{3252371}{2276567}{1}{1194147}{"Audrey 18 May 2014"}{{1, 1194147, 1254514, 2775619, 2563522, 2706517, 522169, 1976723, 1773014, 1118365, 377264}}

\adfDesignF{11}{29765}{1}{3274151}{3244913}{1}{2352691}{"Audrey 18 May 2014"}{{1, 2352691, 2952770, 2600763, 2653866, 2550181, 1129950, 2484208, 1448913, 1551045, 3194669}}

\adfDesignF{11}{29841}{1}{3282511}{587470}{1}{832779}{"Audrey 18 May 2014"}{{1, 832779, 1786294, 100980, 2656622, 2706159, 1481234, 2478085, 2177581, 830583, 1362237}}

\adfDesignF{11}{29900}{1}{3289001}{1490694}{1}{1395243}{"Audrey 18 May 2014"}{{1, 1395243, 2539167, 977429, 2884608, 1619051, 1462568, 116581, 1279728, 1960025, 2210604}}

\adfDesignF{11}{30038}{1}{3304181}{206890}{1}{2651754}{"Audrey 18 May 2014"}{{1, 2651754, 3177185, 3280717, 276555, 2767063, 2065250, 324964, 790418, 92546, 1094452}}

\adfDesignF{11}{30224}{1}{3324641}{2578429}{1}{652699}{"Audrey 18 May 2014"}{{1, 652699, 3136143, 2533385, 1057637, 128946, 2962980, 387243, 611473, 1286782, 541275}}

\adfDesignF{11}{30425}{1}{3346751}{749665}{1}{456216}{"Audrey 18 May 2014"}{{1, 456216, 1940717, 3169822, 2446705, 841005, 1708938, 2479403, 1722566, 1527693, 440689}}

\adfDesignF{11}{30960}{1}{3405601}{379122}{1}{1345697}{"Audrey 18 May 2014"}{{1, 1345697, 2734468, 3286893, 1526631, 1631971, 2617927, 2338866, 3005018, 2429737, 2921998}}

\adfDesignF{11}{30971}{1}{3406811}{2653888}{1}{3350137}{"Audrey 18 May 2014"}{{1, 3350137, 2726314, 1386458, 1975023, 1728714, 3340313, 774686, 2425804, 1842009, 891407}}

\adfDesignF{11}{31365}{1}{3450151}{441075}{1}{2840256}{"Audrey 18 May 2014"}{{1, 2840256, 781262, 2216667, 649083, 1299606, 688766, 1645586, 451826, 592251, 2635300}}

\adfDesignF{11}{31443}{1}{3458731}{1431016}{1}{2575}{"Audrey 18 May 2014"}{{1, 2575, 3171894, 1563159, 2630272, 755102, 580828, 1460308, 652503, 2710690, 307592}}

\adfDesignF{11}{31472}{1}{3461921}{3065333}{1}{2928810}{"Audrey 18 May 2014"}{{1, 2928810, 933826, 1712877, 793744, 408167, 726918, 2916763, 1458588, 3154305, 2275606}}

\adfDesignF{11}{31551}{1}{3470611}{1924533}{1}{944998}{"Audrey 18 May 2014"}{{1, 944998, 1774205, 2708600, 852357, 2376962, 3249544, 2085668, 2513597, 2439630, 1878104}}

\adfDesignF{11}{31848}{1}{3503281}{3373981}{1}{2997680}{"Audrey 18 May 2014"}{{1, 2997680, 1459912, 1333026, 2939840, 31964, 3108170, 924248, 1939542, 539497, 2242525}}

\adfDesignF{11}{32067}{1}{3527371}{2009701}{1}{1031276}{"Audrey 18 May 2014"}{{1, 1031276, 3140079, 2801109, 696231, 2898964, 1652272, 2514328, 2755170, 1027968, 2646828}}

\adfDesignF{11}{32366}{1}{3560261}{1841010}{1}{322417}{"Audrey 18 May 2014"}{{1, 322417, 221211, 3038635, 2079337, 2210185, 736951, 831949, 876732, 2818888, 1104738}}

\adfDesignF{11}{32460}{1}{3570601}{1199914}{1}{3241588}{"Audrey 18 May 2014"}{{1, 3241588, 3214253, 2440689, 2613741, 3059611, 604785, 325323, 409978, 2072864, 3440773}}

\adfDesignF{11}{32466}{1}{3571261}{2508908}{1}{2962086}{"Audrey 18 May 2014"}{{1, 2962086, 878854, 2564843, 2738019, 3298159, 3288426, 524180, 3379554, 2877025, 2487680}}

\adfDesignF{11}{32811}{1}{3609211}{2895967}{1}{29323}{"Audrey 18 May 2014"}{{1, 29323, 846111, 796439, 2385627, 12919, 3465893, 2217101, 2944091, 862484, 876855}}

\adfDesignF{11}{32955}{1}{3625051}{2280760}{1}{3400031}{"Audrey 18 May 2014"}{{1, 3400031, 2913083, 1535466, 301592, 347931, 2417878, 2121877, 1579823, 2479906, 1027667}}

\adfDesignF{11}{33300}{1}{3663001}{2418626}{1}{2009383}{"Audrey 18 May 2014"}{{1, 2009383, 265418, 797496, 3542493, 3280543, 332388, 1509269, 2042100, 1370082, 3165832}}

\adfDesignF{11}{33351}{1}{3668611}{457058}{1}{1925348}{"Audrey 18 May 2014"}{{1, 1925348, 2261710, 3603078, 707639, 2580192, 2299219, 1205064, 2227265, 2472265, 2729885}}

\adfDesignF{11}{33552}{1}{3690721}{1771454}{1}{2935081}{"Audrey 18 May 2014"}{{1, 2935081, 363690, 3446223, 2356902, 2754154, 663967, 279581, 1664642, 1450340, 2539024}}

\adfDesignF{11}{33656}{1}{3702161}{3387857}{1}{2268914}{"Audrey 18 May 2014"}{{1, 2268914, 1101905, 2825455, 588016, 2570732, 360904, 1359632, 3191822, 2189802, 2053622}}

\adfDesignF{11}{33818}{1}{3719981}{1633132}{1}{2838251}{"Audrey 18 May 2014"}{{1, 2838251, 1523748, 53787, 426259, 2932265, 2616132, 949930, 1703098, 1890578, 3665856}}

\adfDesignF{11}{33963}{1}{3735931}{2571738}{1}{2012991}{"Audrey 18 May 2014"}{{1, 2012991, 38103, 2332643, 2297381, 2897739, 508982, 844343, 523325, 1899488, 1588728}}

\adfDesignF{11}{34014}{1}{3741541}{3334844}{1}{716451}{"Audrey 18 May 2014"}{{1, 716451, 25611, 509497, 1153646, 1974200, 2819970, 1792667, 3027288, 3286467, 3401907}}

\adfDesignF{11}{34130}{1}{3754301}{1959207}{1}{1908076}{"Audrey 18 May 2014"}{{1, 1908076, 1855521, 301051, 1363371, 252480, 2878461, 1491795, 3532735, 3314293, 1873721}}

\adfDesignF{11}{34212}{1}{3763321}{484835}{1}{476355}{"Audrey 18 May 2014"}{{1, 476355, 883009, 3127346, 1232696, 2400808, 1276150, 3665478, 728320, 2073931, 2952511}}

\adfDesignF{11}{34436}{1}{3787961}{3025574}{1}{1591449}{"Audrey 18 May 2014"}{{1, 1591449, 3435781, 1564623, 1545377, 1970569, 3364620, 1132351, 1876381, 1370939, 1087714}}

\adfDesignF{11}{34440}{1}{3788401}{217305}{1}{3356222}{"Audrey 18 May 2014"}{{1, 3356222, 2941939, 3734935, 1424715, 698946, 2158802, 3355118, 2729029, 2493936, 3625164}}

\adfDesignF{11}{34536}{1}{3798961}{3313567}{1}{1774604}{"Audrey 18 May 2014"}{{1, 1774604, 2254568, 2237819, 2466287, 2379195, 1697990, 3758941, 1923815, 1311351, 2989195}}

\adfDesignF{11}{34601}{1}{3806111}{2654845}{1}{3675164}{"Audrey 18 May 2014"}{{1, 3675164, 586754, 286719, 2292122, 3413926, 3399583, 1353570, 852369, 2847743, 322604}}

\adfDesignF{11}{34637}{1}{3810071}{1140918}{1}{1545496}{"Audrey 18 May 2014"}{{1, 1545496, 1515690, 220304, 3386082, 1757391, 1168018, 1227980, 3712270, 2082416, 2434707}}

\adfDesignF{11}{34919}{1}{3841091}{3321057}{1}{1812309}{"Audrey 18 May 2014"}{{1, 1812309, 772655, 680890, 1862532, 1875226, 3031673, 967829, 1730739, 2796842, 3674759}}

\adfDesignF{11}{35016}{1}{3851761}{625387}{1}{797574}{"Audrey 18 May 2014"}{{1, 797574, 2104565, 2806164, 3792432, 3467800, 799452, 1613508, 3272448, 707615, 3749007}}

\adfDesignF{11}{35631}{1}{3919411}{711695}{1}{912966}{"Audrey 18 May 2014"}{{1, 912966, 1054485, 3625635, 1714525, 21258, 2827367, 1130621, 2450926, 769561, 1170299}}

\adfDesignF{11}{35862}{1}{3944821}{3203079}{1}{2680680}{"Audrey 18 May 2014"}{{1, 2680680, 1535960, 2444408, 536297, 1911183, 832647, 1541740, 1517920, 2537847, 240601}}

\adfDesignF{11}{35870}{1}{3945701}{1231103}{1}{2833302}{"Audrey 18 May 2014"}{{1, 2833302, 516086, 1999485, 2050494, 2815984, 2129686, 1434902, 2441240, 3669892, 3783134}}

\adfDesignF{11}{35874}{1}{3946141}{2403679}{1}{1938157}{"Audrey 18 May 2014"}{{1, 1938157, 2554519, 1506705, 3633724, 2410276, 285699, 3065482, 530395, 3266951, 538796}}

\adfDesignF{11}{35898}{1}{3948781}{753388}{1}{739883}{"Audrey 18 May 2014"}{{1, 739883, 3394878, 822736, 96052, 1030259, 2384238, 2432900, 1634288, 1985608, 1274281}}

\adfDesignF{11}{36021}{1}{3962311}{3300431}{1}{2739371}{"Audrey 18 May 2014"}{{1, 2739371, 32028, 3084226, 3516546, 1177098, 3207424, 2667890, 500886, 3906516, 2941880}}

\adfDesignF{11}{36084}{1}{3969241}{993975}{1}{2941683}{"Audrey 18 May 2014"}{{1, 2941683, 1767990, 1398798, 3414118, 1455524, 2463095, 1643917, 1581612, 3694713, 3453995}}

\adfDesignF{11}{36287}{1}{3991571}{3341933}{1}{3564017}{"Audrey 18 May 2014"}{{1, 3564017, 445829, 1500839, 3219296, 2520659, 2340343, 3965113, 111518, 3340194, 2941617}}

\adfDesignF{11}{36602}{1}{4026221}{877386}{1}{3122585}{"Audrey 18 May 2014"}{{1, 3122585, 139486, 303530, 1644324, 2471765, 2471978, 3256118, 756868, 989222, 949007}}

\adfDesignF{11}{36861}{1}{4054711}{59049}{1}{2933101}{"Audrey 18 May 2014"}{{1, 2933101, 2466662, 1971966, 3189864, 367007, 3648872, 3114508, 4014083, 1928862, 693340}}

\adfDesignF{11}{36882}{1}{4057021}{911464}{1}{2657807}{"Audrey 18 May 2014"}{{1, 2657807, 3193826, 3148946, 1245207, 943278, 2573333, 3147448, 3844943, 3876710, 3767648}}

\adfDesignF{11}{36932}{1}{4062521}{2768526}{1}{194320}{"Audrey 18 May 2014"}{{1, 194320, 3192226, 2962309, 1034306, 1240487, 1750305, 946959, 1184185, 1514718, 2230268}}

\adfDesignF{11}{36944}{1}{4063841}{2540223}{1}{2451153}{"Audrey 18 May 2014"}{{1, 2451153, 1877528, 3521652, 1600631, 886344, 2845304, 915814, 3550439, 2470759, 199580}}

\adfDesignF{11}{37016}{1}{4071761}{3914380}{1}{2137613}{"Audrey 18 May 2014"}{{1, 2137613, 2138915, 226517, 9323, 1767665, 1690928, 2645793, 1411348, 1524028, 2734913}}

\adfDesignF{11}{37835}{1}{4161851}{126109}{1}{2436362}{"Audrey 18 May 2014"}{{1, 2436362, 3158890, 642705, 3057119, 242630, 1843424, 3324242, 3280882, 2888495, 4096356}}

\adfDesignF{11}{37857}{1}{4164271}{1964586}{1}{1995198}{"Audrey 18 May 2014"}{{1, 1995198, 1018109, 1811053, 584187, 1773939, 4036808, 2271567, 4113977, 4114346, 3266441}}

\adfDesignF{11}{38069}{1}{4187591}{3581872}{1}{3902414}{"Audrey 18 May 2014"}{{1, 3902414, 2904109, 3154959, 3121562, 408306, 675184, 2674203, 1946034, 546848, 1604335}}

\adfDesignF{11}{38207}{1}{4202771}{1617065}{1}{3908275}{"Audrey 18 May 2014"}{{1, 3908275, 3714431, 3758562, 2323518, 3200666, 1537231, 1503231, 3366709, 1778488, 125514}}

\adfDesignF{11}{38241}{1}{4206511}{1664613}{1}{340252}{"Audrey 18 May 2014"}{{1, 340252, 4034273, 787276, 1613072, 2244908, 3549903, 4027016, 800469, 2210471, 1424914}}

\adfDesignF{11}{38379}{1}{4221691}{3393935}{1}{3340256}{"Audrey 18 May 2014"}{{1, 3340256, 1421113, 3421146, 2283262, 605786, 2716461, 3131098, 3578564, 2866529, 1965930}}

\adfDesignF{11}{38724}{1}{4259641}{3015005}{1}{3844611}{"Audrey 18 May 2014"}{{1, 3844611, 2797783, 1479187, 999592, 2807794, 4023473, 2465630, 177094, 782635, 1920405}}

\adfDesignF{11}{38835}{1}{4271851}{1879047}{1}{2835535}{"Audrey 18 May 2014"}{{1, 2835535, 2920277, 208497, 1993601, 2070639, 643233, 181695, 213821, 1660507, 87747}}

\adfDesignF{11}{38909}{1}{4279991}{1832646}{1}{3750727}{"Audrey 18 May 2014"}{{1, 3750727, 3530728, 3926309, 1463672, 913610, 3940158, 3311119, 3348498, 2107844, 3667271}}

\adfDesignF{11}{38946}{1}{4284061}{193100}{1}{3890316}{"Audrey 18 May 2014"}{{1, 3890316, 3525557, 2412987, 346021, 1953338, 400520, 2106132, 3879674, 3948189, 3241631}}

\adfDesignF{11}{39186}{1}{4310461}{1308077}{1}{569428}{"Audrey 18 May 2014"}{{1, 569428, 2439381, 1476357, 1184044, 2739056, 2282189, 2983307, 4305991, 2139291, 1433260}}

\adfDesignF{11}{39251}{1}{4317611}{197634}{1}{286272}{"Audrey 18 May 2014"}{{1, 286272, 3401204, 697267, 544483, 162665, 1000245, 2492731, 1613196, 1172752, 1582017}}

\adfDesignF{11}{39845}{1}{4382951}{614176}{1}{2889264}{"Audrey 18 May 2014"}{{1, 2889264, 3476929, 2651697, 1737596, 1239512, 3200676, 861662, 233356, 2120605, 3503457}}

\adfDesignF{11}{39888}{1}{4387681}{4157931}{1}{424287}{"Audrey 18 May 2014"}{{1, 424287, 1682301, 3662350, 3431343, 2182512, 169256, 4333226, 979561, 690644, 4382924}}

\adfDesignF{11}{40388}{1}{4442681}{1746373}{1}{2530334}{"Audrey 18 May 2014"}{{1, 2530334, 2657682, 3536303, 1699697, 1569214, 2721769, 2752180, 1061853, 2576403, 1107969}}

\adfDesignF{11}{40412}{1}{4445321}{3145655}{1}{1361421}{"Audrey 18 May 2014"}{{1, 1361421, 3884254, 1865, 771874, 2705801, 3478225, 3706327, 879530, 1166286, 4271021}}

\adfDesignF{11}{40698}{1}{4476781}{4224036}{1}{4258561}{"Audrey 18 May 2014"}{{1, 4258561, 448903, 1309182, 560256, 1824790, 3996150, 1271552, 1762502, 499613, 1975614}}

\adfDesignF{11}{40857}{1}{4494271}{1173151}{1}{3511146}{"Audrey 18 May 2014"}{{1, 3511146, 1338636, 1870888, 4090189, 1019747, 2756866, 2642907, 903023, 323652, 4014300}}

\adfDesignF{11}{40880}{1}{4496801}{2090228}{1}{3596277}{"Audrey 18 May 2014"}{{1, 3596277, 3872639, 4213495, 2344410, 1133849, 3385388, 1081842, 2754641, 3965358, 632906}}

\adfDesignF{11}{40902}{1}{4499221}{926114}{1}{3399219}{"Audrey 18 May 2014"}{{1, 3399219, 1901148, 2813051, 12753, 245572, 3537696, 2550370, 667053, 321900, 2548121}}

\adfDesignF{11}{41111}{1}{4522211}{3379184}{1}{3566871}{"Audrey 18 May 2014"}{{1, 3566871, 1891580, 611666, 2064558, 2965719, 3735104, 1556300, 3323947, 1559431, 1335878}}

\adfDesignF{11}{41135}{1}{4524851}{780756}{1}{1931314}{"Audrey 18 May 2014"}{{1, 1931314, 3341766, 860929, 4382842, 823837, 1000135, 1807659, 3820025, 2278923, 2376824}}

\adfDesignF{11}{41267}{1}{4539371}{4470238}{1}{1883328}{"Audrey 18 May 2014"}{{1, 1883328, 576685, 1641591, 2191023, 4162527, 720276, 2104485, 1761076, 1398520, 1717972}}

\adfDesignF{11}{41298}{1}{4542781}{3950549}{1}{3303522}{"Audrey 18 May 2014"}{{1, 3303522, 3067535, 1468731, 4324598, 3464058, 3935825, 178748, 19390, 2079480, 872017}}

\adfDesignF{11}{41351}{1}{4548611}{1271297}{1}{1999622}{"Audrey 18 May 2014"}{{1, 1999622, 4351668, 2636223, 539252, 1890473, 3611603, 2032144, 18274, 2100265, 3563530}}

\adfDesignF{11}{41540}{1}{4569401}{2780330}{1}{550227}{"Audrey 18 May 2014"}{{1, 550227, 4088274, 3750507, 1904870, 4551115, 370080, 1791397, 368607, 91003, 811523}}

\adfDesignF{11}{41594}{1}{4575341}{3721893}{1}{3128985}{"Audrey 18 May 2014"}{{1, 3128985, 3691375, 814197, 851812, 1550703, 172660, 3435502, 3230859, 1901795, 4098816}}

\adfDesignF{11}{41699}{1}{4586891}{3766683}{1}{3657535}{"Audrey 18 May 2014"}{{1, 3657535, 173218, 338928, 1621493, 812395, 2677871, 208121, 1519612, 1927009, 824490}}

\adfDesignF{11}{41886}{1}{4607461}{1283249}{1}{3214567}{"Audrey 18 May 2014"}{{1, 3214567, 2550207, 103502, 139902, 3906607, 317179, 3494342, 75276, 501233, 4127028}}

\adfDesignF{11}{42243}{1}{4646731}{4335734}{1}{3262514}{"Audrey 18 May 2014"}{{1, 3262514, 408894, 1690188, 275125, 4077173, 4250971, 662137, 3164366, 2064763, 3377523}}

\adfDesignF{11}{42435}{1}{4667851}{2385762}{1}{741540}{"Audrey 18 May 2014"}{{1, 741540, 4055949, 2650928, 1725341, 2752401, 2255641, 2974757, 955157, 1414593, 3812947}}

\adfDesignF{11}{42441}{1}{4668511}{249014}{1}{2966773}{"Audrey 18 May 2014"}{{1, 2966773, 2167767, 1879956, 2637953, 4596978, 3668540, 2384565, 2235318, 522671, 282033}}

\adfDesignF{11}{42474}{1}{4672141}{870068}{1}{450179}{"Audrey 18 May 2014"}{{1, 450179, 2344025, 4424920, 1368202, 2787987, 1946420, 725335, 3994757, 1919193, 3399686}}

\adfDesignF{11}{42554}{1}{4680941}{2953373}{1}{3731674}{"Audrey 18 May 2014"}{{1, 3731674, 3290084, 2476182, 69261, 1315599, 1619103, 2525144, 3802537, 504833, 4070287}}

\adfDesignF{11}{42866}{1}{4715261}{1247232}{1}{2795264}{"Audrey 18 May 2014"}{{1, 2795264, 2145470, 3198620, 3732700, 4259871, 1301861, 72205, 5276, 3191717, 2873320}}

\adfDesignF{11}{42906}{1}{4719661}{4556576}{1}{4113417}{"Audrey 18 May 2014"}{{1, 4113417, 2346144, 194468, 2073588, 302522, 3879092, 3395005, 273931, 1686443, 614033}}

\adfDesignF{11}{42933}{1}{4722631}{59049}{1}{4059619}{"Audrey 18 May 2014"}{{1, 4059619, 2418664, 3881130, 3079934, 4107406, 3195599, 2906604, 3474612, 4118349, 1816499}}

\adfDesignF{11}{43067}{1}{4737371}{3102729}{1}{1913305}{"Audrey 18 May 2014"}{{1, 1913305, 3643340, 1868637, 2308440, 4063480, 1770573, 4282746, 1878427, 1999456, 4695821}}

\adfDesignF{11}{43205}{1}{4752551}{2290887}{1}{3019914}{"Audrey 18 May 2014"}{{1, 3019914, 1341252, 4299007, 2312780, 225810, 2247554, 2296543, 2703410, 527512, 36421}}

\adfDesignF{11}{43263}{1}{4758931}{1024}{1}{594650}{"Audrey 18 May 2014"}{{1, 594650, 1013476, 1999422, 4006984, 3873210, 695775, 1142610, 1421906, 2855337, 1432353}}

\adfDesignF{11}{43710}{1}{4808101}{2127140}{1}{662241}{"Audrey 18 May 2014"}{{1, 662241, 1825568, 2638145, 1779282, 1795094, 1105808, 3956721, 3221185, 221117, 2027242}}

\adfDesignF{11}{43977}{1}{4837471}{146593}{1}{3556624}{"Audrey 18 May 2014"}{{1, 3556624, 3634882, 4346947, 297819, 3040483, 2824507, 4257044, 1125976, 1922500, 4018043}}

\adfDesignF{11}{44156}{1}{4857161}{3410918}{1}{93390}{"Audrey 18 May 2014"}{{1, 93390, 3088105, 4191575, 2869938, 508679, 2497230, 4581446, 3643772, 4024581, 3644249}}

\adfDesignF{11}{44204}{1}{4862441}{1316606}{1}{1552347}{"Audrey 18 May 2014"}{{1, 1552347, 4073219, 3955208, 2398066, 4080594, 1465778, 2227693, 3956035, 2504611, 2961094}}

\adfDesignF{11}{44387}{1}{4882571}{3361228}{1}{510973}{"Audrey 18 May 2014"}{{1, 510973, 2805075, 4692928, 1933998, 3037367, 4332034, 4668235, 1041173, 1472598, 4801044}}

\adfDesignF{11}{44441}{1}{4888511}{831792}{1}{1373748}{"Audrey 18 May 2014"}{{1, 1373748, 3227020, 290187, 404359, 974091, 3892994, 795600, 173464, 465866, 3068203}}

\adfDesignF{11}{44477}{1}{4892471}{1855707}{1}{3826908}{"Audrey 18 May 2014"}{{1, 3826908, 4299644, 3123436, 2982586, 1198269, 1078662, 1002853, 679639, 63076, 1314810}}

\adfDesignF{11}{44603}{1}{4906331}{1621397}{1}{1352992}{"Audrey 18 May 2014"}{{1, 1352992, 911647, 4390755, 2125526, 3196128, 3410858, 2992853, 3252594, 70798, 2827503}}

\adfDesignF{11}{44681}{1}{4914911}{4205287}{1}{3817956}{"Audrey 18 May 2014"}{{1, 3817956, 2441717, 801380, 3183560, 4057407, 3058585, 4396009, 1555767, 3071956, 3105128}}

\adfDesignF{11}{44733}{1}{4920631}{2513466}{1}{508908}{"Audrey 18 May 2014"}{{1, 508908, 4701672, 2622854, 1335848, 4116917, 1126301, 3887273, 3364430, 2579680, 359271}}

\adfDesignF{11}{44961}{1}{4945711}{428345}{1}{831885}{"Audrey 18 May 2014"}{{1, 831885, 4041550, 592528, 868465, 3435067, 3298316, 3437103, 637303, 2369799, 270827}}

\adfDesignF{11}{45221}{1}{4974311}{4106306}{1}{1315042}{"Audrey 18 May 2014"}{{1, 1315042, 1319681, 3294233, 1917551, 2781046, 4232067, 2141727, 1095023, 2867509, 3907675}}

\adfDesignF{11}{45794}{1}{5037341}{4988828}{1}{3294634}{"Audrey 18 May 2014"}{{1, 3294634, 4724267, 1847608, 3785639, 4877379, 446194, 3736307, 3800256, 525256, 3186505}}

\adfDesignF{11}{46026}{1}{5062861}{3072287}{1}{1331682}{"Audrey 18 May 2014"}{{1, 1331682, 3563793, 4219924, 3385164, 657170, 621785, 4164403, 28025, 2039619, 239878}}

\adfDesignF{11}{46418}{1}{5105981}{3856783}{1}{1510406}{"Audrey 18 May 2014"}{{1, 1510406, 4609922, 1161986, 1989148, 581916, 1166499, 1966791, 3393308, 1366830, 2677117}}

\adfDesignF{11}{46452}{1}{5109721}{5038486}{1}{1014178}{"Audrey 18 May 2014"}{{1, 1014178, 836710, 1507910, 749890, 1286222, 2491147, 701763, 5106329, 3860178, 2884556}}

\adfDesignF{11}{46488}{1}{5113681}{2713980}{1}{4321822}{"Audrey 18 May 2014"}{{1, 4321822, 1111661, 506903, 2687418, 2218088, 3522202, 1196459, 192632, 3721342, 976196}}

\adfDesignF{11}{46626}{1}{5128861}{2657856}{1}{485269}{"Audrey 18 May 2014"}{{1, 485269, 4607268, 1505694, 3954765, 2518944, 3993806, 1893439, 2059663, 3816972, 808484}}

\adfDesignF{11}{46758}{1}{5143381}{4621875}{1}{3153909}{"Audrey 18 May 2014"}{{1, 3153909, 2571092, 2577981, 4215738, 2664362, 2280878, 4970834, 3033863, 2973450, 2418178}}

\adfDesignF{11}{47150}{1}{5186501}{2891864}{1}{145163}{"Audrey 18 May 2014"}{{1, 145163, 4729507, 1914269, 3866770, 2862785, 2066330, 3749457, 638549, 542615, 230558}}

\adfDesignF{11}{47337}{1}{5207071}{693914}{1}{4902032}{"Audrey 18 May 2014"}{{1, 4902032, 3639822, 465059, 310023, 1717605, 4179496, 397438, 2144011, 1146171, 1926626}}

\adfDesignF{11}{47708}{1}{5247881}{4616606}{1}{2625749}{"Audrey 18 May 2014"}{{1, 2625749, 1958702, 5240773, 2888825, 174477, 3292735, 1190253, 3565281, 788880, 4513729}}

\adfDesignF{11}{47792}{1}{5257121}{2762931}{1}{2899452}{"Audrey 18 May 2014"}{{1, 2899452, 1995574, 4796033, 4375008, 2892755, 4170504, 4038537, 1094196, 1399512, 3881154}}

\adfDesignF{11}{47955}{1}{5275051}{4328146}{1}{2431988}{"Audrey 18 May 2014"}{{1, 2431988, 4374261, 2384627, 3729178, 87329, 4251741, 1745653, 3902956, 4837026, 3905546}}

\adfDesignF{11}{48162}{1}{5297821}{3095200}{1}{912965}{"Audrey 18 May 2014"}{{1, 912965, 4211116, 4606166, 1073557, 1890221, 3998367, 928883, 4864983, 3821541, 181305}}

\adfDesignF{11}{48285}{1}{5311351}{1597440}{1}{4535905}{"Audrey 18 May 2014"}{{1, 4535905, 2518153, 5278508, 44933, 4658793, 458396, 821059, 661109, 2580257, 4999641}}

\adfDesignF{11}{48555}{1}{5341051}{1611599}{1}{4299441}{"Audrey 18 May 2014"}{{1, 4299441, 2338266, 736148, 3527484, 1264139, 161342, 729095, 716638, 2639529, 4952121}}

\adfDesignF{11}{48761}{1}{5363711}{4644268}{1}{620996}{"Audrey 18 May 2014"}{{1, 620996, 1302249, 4712534, 2587420, 2748316, 3312224, 1960824, 4916906, 652250, 4004835}}

\adfDesignF{11}{49040}{1}{5394401}{295726}{1}{1027046}{"Audrey 18 May 2014"}{{1, 1027046, 2314576, 3756222, 1510661, 2299390, 2248957, 2276441, 708272, 3738464, 1697574}}

\adfDesignF{11}{49085}{1}{5399351}{4213965}{1}{262991}{"Audrey 18 May 2014"}{{1, 262991, 3979122, 3459188, 4060669, 3364093, 4725356, 674934, 3302820, 2141197, 1026384}}

\adfDesignF{11}{49106}{1}{5401661}{2233287}{1}{949301}{"Audrey 18 May 2014"}{{1, 949301, 2480649, 1455094, 4036713, 3939332, 4085005, 3686317, 1548855, 2877816, 1949222}}

\adfDesignF{11}{49203}{1}{5412331}{1024}{1}{3930415}{"Audrey 18 May 2014"}{{1, 3930415, 78482, 1849347, 191646, 3383158, 5305854, 4284289, 111150, 3918254, 4009059}}

\adfDesignF{11}{49776}{1}{5475361}{1607651}{1}{4128283}{"Audrey 18 May 2014"}{{1, 4128283, 2372269, 2838697, 2042146, 924232, 2366889, 740373, 2301417, 4991562, 4670936}}

\adfDesignF{11}{49811}{1}{5479211}{1912936}{1}{3803555}{"Audrey 18 May 2014"}{{1, 3803555, 1353386, 3186418, 2740595, 2139321, 2247963, 2414286, 2348913, 4622289, 2539328}}

\adfDesignF{11}{49985}{1}{5498351}{3165595}{1}{1673230}{"Audrey 18 May 2014"}{{1, 1673230, 4283912, 807802, 907534, 529044, 5272875, 1644736, 3469813, 3909825, 4992983}}

~

\adfDesignF{13}{1517}{1}{236653}{91}{1}{67452}{"Audrey 12 Apr 2014"}{{1, 67452, 118379, 228088, 180246, 141970, 233448, 117182, 186717, 235730, 218196, 69869, 95946}}

\adfDesignF{13}{1726}{1}{269257}{13}{1}{16602}{"Audrey 12 Apr 2014"}{{1, 16602, 176493, 82112, 244490, 242962, 185264, 30217, 36843, 184839, 244306, 150021, 21392}}

\adfDesignF{13}{1925}{1}{300301}{76}{1}{13805}{"Audrey 12 Apr 2014"}{{1, 13805, 187191, 81650, 148597, 25454, 40300, 184048, 236180, 96943, 156859, 268285, 62192}}

\adfDesignF{13}{2416}{1}{376897}{141}{1}{70110}{"Audrey 12 Apr 2014"}{{1, 70110, 298323, 280309, 300616, 107520, 287200, 246672, 255075, 299394, 365616, 195893, 308447}}

\adfDesignF{13}{2745}{1}{428221}{181}{1}{41933}{"Audrey 12 Apr 2014"}{{1, 41933, 101063, 199763, 230898, 169024, 197621, 336822, 371904, 98054, 348561, 169241, 304441}}

\adfDesignF{13}{2842}{1}{443353}{99}{1}{5406}{"Audrey 12 Apr 2014"}{{1, 5406, 406891, 178513, 305150, 367740, 7588, 232252, 421969, 113229, 288834, 390691, 385207}}

\adfDesignF{13}{2897}{1}{451933}{213}{1}{159965}{"Audrey 12 Apr 2014"}{{1, 159965, 354765, 304482, 287921, 338802, 204637, 346649, 432051, 280324, 332534, 383344, 189989}}

\adfDesignF{13}{3777}{1}{589213}{354}{1}{24669}{"Audrey 21 Apr 2014"}{{1, 24669, 491745, 140161, 129825, 280270, 155288, 325959, 92760, 382361, 341805, 349515, 231706}}

\adfDesignF{13}{3961}{1}{617917}{118}{1}{94949}{"Audrey 21 Apr 2014"}{{1, 94949, 521488, 456985, 137025, 144290, 353403, 514596, 442580, 464918, 126619, 154279, 296369}}

\adfDesignF{13}{4457}{1}{695293}{11}{1}{66638}{"Audrey 21 Apr 2014"}{{1, 66638, 481946, 333878, 281457, 202891, 278073, 670124, 530187, 678097, 631109, 349144, 363506}}

\adfDesignF{13}{4621}{1}{720877}{233}{1}{15240}{"Audrey 21 Apr 2014"}{{1, 15240, 135206, 272974, 663470, 261998, 632694, 526685, 434882, 579419, 323187, 338216, 141290}}

\adfDesignF{13}{4953}{1}{772669}{4096}{1}{366166}{"Audrey 21 Apr 2014"}{{1, 366166, 151331, 309611, 707739, 639419, 109843, 259812, 222836, 347707, 401549, 89117, 258214}}

\adfDesignF{13}{4963}{1}{774229}{541}{1}{29760}{"Audrey 21 Apr 2014"}{{1, 29760, 713853, 195749, 191244, 64061, 303562, 301148, 463805, 656417, 398021, 175489, 378035}}

\adfDesignF{13}{5180}{1}{808081}{495980}{1}{703093}{"Audrey 18 May 2014"}{{1, 703093, 255304, 190418, 318956, 324112, 420254, 403729, 380522, 472823, 482787, 39169, 49237}}

\adfDesignF{13}{6510}{1}{1015561}{30488}{1}{254773}{"Audrey 18 May 2014"}{{1, 254773, 715775, 933110, 590662, 931868, 2067, 555193, 850109, 188031, 194032, 667500, 210245}}

\adfDesignF{13}{6713}{1}{1047229}{742790}{1}{574865}{"Audrey 18 May 2014"}{{1, 574865, 948840, 445805, 857674, 878520, 25634, 530151, 671035, 402522, 89690, 369264, 489373}}

\adfDesignF{13}{7018}{1}{1094809}{125283}{1}{439857}{"Audrey 18 May 2014"}{{1, 439857, 628778, 660357, 1062777, 692406, 184277, 249265, 212991, 586539, 449264, 1080366, 321976}}

\adfDesignF{13}{7510}{1}{1171561}{719103}{1}{1170557}{"Audrey 18 May 2014"}{{1, 1170557, 1008016, 180640, 229395, 484937, 492628, 971791, 232149, 63043, 1140683, 540926, 514600}}

\adfDesignF{13}{7725}{1}{1205101}{240252}{1}{237497}{"Audrey 18 May 2014"}{{1, 237497, 72704, 294760, 298630, 1125058, 495904, 1191558, 1197799, 1142346, 564933, 1177967, 636550}}

\adfDesignF{13}{7740}{1}{1207441}{1013096}{1}{1015747}{"Audrey 18 May 2014"}{{1, 1015747, 537683, 181081, 580095, 962847, 1060765, 438018, 24648, 1050362, 1145609, 582552, 872679}}

\adfDesignF{13}{7846}{1}{1223977}{1103072}{1}{336459}{"Audrey 18 May 2014"}{{1, 336459, 249928, 857098, 786943, 902243, 673928, 57840, 778237, 667150, 1031866, 750421, 251748}}

\adfDesignF{13}{7858}{1}{1225849}{488001}{1}{651344}{"Audrey 18 May 2014"}{{1, 651344, 1055171, 928729, 1129897, 925928, 778665, 913896, 660314, 1214517, 1032470, 959223, 782486}}

\adfDesignF{13}{7971}{1}{1243477}{176954}{1}{424197}{"Audrey 18 May 2014"}{{1, 424197, 781616, 942026, 834402, 71052, 649718, 453735, 194873, 677975, 1113561, 923188, 394518}}

\adfDesignF{13}{8097}{1}{1263133}{797816}{1}{1239112}{"Audrey 18 May 2014"}{{1, 1239112, 1019793, 753749, 1206826, 998137, 541729, 1186990, 14419, 1002776, 264014, 310499, 303886}}

\adfDesignF{13}{8965}{1}{1398541}{550017}{1}{1017406}{"Audrey 18 May 2014"}{{1, 1017406, 231637, 729712, 674304, 633284, 501145, 438959, 625742, 1019560, 210414, 396673, 513868}}

\adfDesignF{13}{9070}{1}{1414921}{165182}{1}{1141833}{"Audrey 18 May 2014"}{{1, 1141833, 814597, 43926, 36750, 50653, 956153, 1270560, 727866, 876635, 474636, 687000, 1408916}}

\adfDesignF{13}{9211}{1}{1436917}{25009}{1}{1159043}{"Audrey 18 May 2014"}{{1, 1159043, 1224881, 1383713, 1006200, 1364894, 1376043, 1351869, 1090970, 1366295, 42159, 294335, 1271850}}

\adfDesignF{13}{9545}{1}{1489021}{119708}{1}{378534}{"Audrey 18 May 2014"}{{1, 378534, 987347, 138298, 983835, 422643, 1351080, 140913, 651280, 372634, 1168247, 1330171, 1009143}}

\adfDesignF{13}{9748}{1}{1520689}{1043855}{1}{576910}{"Audrey 18 May 2014"}{{1, 576910, 1070804, 439725, 410770, 750385, 948586, 1440208, 903527, 1110284, 1487372, 598490, 907761}}

\adfDesignF{13}{9811}{1}{1530517}{34694}{1}{158460}{"Audrey 18 May 2014"}{{1, 158460, 1440215, 1079030, 1387145, 267228, 135041, 438683, 687074, 419245, 1472315, 206522, 1492143}}

\adfDesignF{13}{9870}{1}{1539721}{1032466}{1}{1166325}{"Audrey 18 May 2014"}{{1, 1166325, 1296545, 518884, 42250, 366, 372233, 301402, 525861, 66290, 133956, 741830, 992941}}

\adfDesignF{13}{10526}{1}{1642057}{608345}{1}{219756}{"Audrey 18 May 2014"}{{1, 219756, 1445223, 1255047, 930698, 60053, 1437016, 696141, 563248, 512885, 405637, 458270, 226310}}

\adfDesignF{13}{10710}{1}{1670761}{867354}{1}{774273}{"Audrey 18 May 2014"}{{1, 774273, 1228792, 604483, 1245407, 290678, 925067, 831352, 987387, 276271, 1645153, 1008564, 207138}}

\adfDesignF{13}{10856}{1}{1693537}{911226}{1}{284858}{"Audrey 18 May 2014"}{{1, 284858, 1641883, 1087861, 814941, 879103, 1286795, 1314756, 1524783, 120813, 184177, 209143, 812108}}

\adfDesignF{13}{11552}{1}{1802113}{996983}{1}{1687707}{"Audrey 18 May 2014"}{{1, 1687707, 1788230, 636945, 1713711, 261056, 46013, 1610908, 951636, 1788679, 1529928, 886583, 1515507}}

\adfDesignF{13}{11745}{1}{1832221}{1512906}{1}{1337381}{"Audrey 18 May 2014"}{{1, 1337381, 1282276, 495333, 279218, 1382711, 234558, 687409, 1687974, 1351983, 75999, 823086, 1355397}}

\adfDesignF{13}{12170}{1}{1898521}{13683}{1}{1027637}{"Audrey 18 May 2014"}{{1, 1027637, 685687, 1262469, 1733361, 1351959, 611251, 1384348, 1774393, 1201933, 1537536, 1866871, 750722}}

\adfDesignF{13}{12310}{1}{1920361}{250582}{1}{703474}{"Audrey 18 May 2014"}{{1, 703474, 559337, 908560, 346893, 332207, 655223, 1536399, 1092467, 1339602, 264540, 588533, 1274569}}

\adfDesignF{13}{12348}{1}{1926289}{826001}{1}{1285881}{"Audrey 18 May 2014"}{{1, 1285881, 68052, 1243409, 275948, 372365, 1348124, 1777474, 1094534, 1269182, 1411005, 1353571, 58188}}

\adfDesignF{13}{12415}{1}{1936741}{582729}{1}{690199}{"Audrey 18 May 2014"}{{1, 690199, 286054, 870465, 1520407, 951704, 71536, 737351, 1490279, 1424349, 334074, 777712, 529574}}

\adfDesignF{13}{12823}{1}{2000389}{339102}{1}{131192}{"Audrey 18 May 2014"}{{1, 131192, 1994297, 933936, 1105462, 1568593, 835259, 1990086, 591788, 753817, 1528871, 639980, 1929441}}

\adfDesignF{13}{13098}{1}{2043289}{223464}{1}{1374839}{"Audrey 18 May 2014"}{{1, 1374839, 1007269, 1758097, 1815278, 939862, 1358219, 1283976, 1614094, 1275838, 181587, 1796184, 1941068}}

\adfDesignF{13}{13111}{1}{2045317}{1293688}{1}{571020}{"Audrey 18 May 2014"}{{1, 571020, 1449577, 1759274, 696443, 1670965, 1781898, 1230751, 242918, 1728054, 435015, 560967, 145019}}

\adfDesignF{13}{13790}{1}{2151241}{2050215}{1}{1117727}{"Audrey 18 May 2014"}{{1, 1117727, 1948189, 1373937, 1580939, 882120, 1809915, 397902, 395655, 1512574, 1707326, 1501722, 830680}}

\adfDesignF{13}{14035}{1}{2189461}{1502366}{1}{392155}{"Audrey 18 May 2014"}{{1, 392155, 2182307, 1412132, 822113, 1970187, 1685305, 1032120, 689757, 1265473, 523516, 227393, 934307}}

\adfDesignF{13}{14167}{1}{2210053}{1420652}{1}{1777783}{"Audrey 18 May 2014"}{{1, 1777783, 1791856, 522002, 606760, 524734, 2023475, 607931, 438701, 736501, 1897645, 1527648, 805281}}

\adfDesignF{13}{14575}{1}{2273701}{124863}{1}{1485502}{"Audrey 18 May 2014"}{{1, 1485502, 697165, 1028144, 1342960, 1771510, 110620, 1312568, 1001782, 491559, 1432963, 968113, 1999319}}

\adfDesignF{13}{14868}{1}{2319409}{1638710}{1}{1070241}{"Audrey 18 May 2014"}{{1, 1070241, 1176930, 1933318, 1573055, 294787, 163260, 1826872, 1111422, 1361733, 397775, 1508279, 1498781}}

\adfDesignF{13}{15818}{1}{2467609}{2261893}{1}{216255}{"Audrey 18 May 2014"}{{1, 216255, 99257, 1559453, 1256921, 916678, 831875, 1029198, 752126, 1028504, 1195305, 1237198, 2215274}}

\adfDesignF{13}{16018}{1}{2498809}{1858260}{1}{1755090}{"Audrey 18 May 2014"}{{1, 1755090, 1581193, 1226723, 452344, 327143, 1569895, 269318, 119371, 1704212, 1051788, 2445024, 2490752}}

\adfDesignF{13}{16151}{1}{2519557}{893163}{1}{2232876}{"Audrey 18 May 2014"}{{1, 2232876, 565978, 1812225, 2178175, 580591, 46506, 1109058, 2080446, 154200, 1936922, 1392234, 1028130}}

\adfDesignF{13}{16503}{1}{2574469}{1275896}{1}{1735974}{"Audrey 18 May 2014"}{{1, 1735974, 1827939, 130752, 1438594, 701575, 1611344, 614141, 1429461, 2205604, 2475922, 1009741, 265766}}

\adfDesignF{13}{16605}{1}{2590381}{369129}{1}{2263084}{"Audrey 18 May 2014"}{{1, 2263084, 710335, 1361017, 677997, 1194637, 1363475, 1790462, 1281273, 2561390, 101724, 206965, 2029926}}

\adfDesignF{13}{16737}{1}{2610973}{1746603}{1}{12323}{"Audrey 18 May 2014"}{{1, 12323, 419895, 2028572, 637254, 1685231, 2033344, 2001204, 196907, 891044, 1193747, 322399, 1632944}}

\adfDesignF{13}{16795}{1}{2620021}{1896984}{1}{1044858}{"Audrey 18 May 2014"}{{1, 1044858, 2169758, 553190, 2164210, 1787479, 723300, 733950, 1242463, 579943, 2246235, 137977, 2336762}}

\adfDesignF{13}{16840}{1}{2627041}{1839078}{1}{2035108}{"Audrey 18 May 2014"}{{1, 2035108, 456073, 582215, 1356072, 1745579, 1951913, 2438463, 2569184, 1361105, 1611284, 1333488, 948802}}

\adfDesignF{13}{16876}{1}{2632657}{2148336}{1}{2562561}{"Audrey 18 May 2014"}{{1, 2562561, 911254, 896407, 1712204, 1548389, 674795, 517999, 2580097, 1158617, 418561, 1510409, 1304648}}

\adfDesignF{13}{16961}{1}{2645917}{1240632}{1}{1528657}{"Audrey 18 May 2014"}{{1, 1528657, 352676, 1817797, 1094640, 915257, 736755, 533317, 1465146, 450430, 699766, 285834, 703392}}

\adfDesignF{13}{17255}{1}{2691781}{1917664}{1}{1257392}{"Audrey 18 May 2014"}{{1, 1257392, 920628, 630250, 214476, 1635326, 2399235, 332523, 109067, 1806657, 1319214, 2157134, 677002}}

\adfDesignF{13}{17326}{1}{2702857}{343922}{1}{445985}{"Audrey 18 May 2014"}{{1, 445985, 2076452, 65595, 1364764, 2498996, 2458538, 367883, 1474141, 1837205, 178089, 1569720, 1879773}}

\adfDesignF{13}{17723}{1}{2764789}{1900971}{1}{1886569}{"Audrey 18 May 2014"}{{1, 1886569, 2064171, 1242377, 2178075, 1302306, 1855099, 51938, 498562, 1455134, 1930944, 752626, 1370932}}

\adfDesignF{13}{17898}{1}{2792089}{794578}{1}{2143827}{"Audrey 18 May 2014"}{{1, 2143827, 721076, 2578379, 2200018, 2616217, 1764327, 2396197, 1295091, 543746, 198442, 300782, 2786520}}

\adfDesignF{13}{18715}{1}{2919541}{2753363}{1}{506902}{"Audrey 18 May 2014"}{{1, 506902, 834194, 2886253, 1193204, 2024120, 1584905, 980553, 1180179, 707771, 2739657, 2464685, 414822}}

\adfDesignF{13}{18823}{1}{2936389}{2662473}{1}{1561062}{"Audrey 18 May 2014"}{{1, 1561062, 2400355, 465222, 912728, 962666, 2021650, 1677882, 1224350, 2405156, 1459156, 2626647, 2837848}}

\adfDesignF{13}{18922}{1}{2951833}{1416907}{1}{2552560}{"Audrey 18 May 2014"}{{1, 2552560, 2235531, 2852742, 943144, 2256997, 1229723, 902509, 489351, 136280, 1165082, 2432583, 514495}}

\adfDesignF{13}{19287}{1}{3008773}{2903595}{1}{1934294}{"Audrey 18 May 2014"}{{1, 1934294, 2816065, 2958818, 2096898, 1663313, 1229208, 2901178, 2483026, 962517, 1239647, 111868, 664578}}

\adfDesignF{13}{19933}{1}{3109549}{805063}{1}{2098761}{"Audrey 18 May 2014"}{{1, 2098761, 304210, 2154483, 436311, 82955, 2179794, 1750415, 698941, 2138194, 100088, 1427371, 2176221}}

\adfDesignF{13}{19991}{1}{3118597}{1752975}{1}{3032127}{"Audrey 18 May 2014"}{{1, 3032127, 1783891, 2108641, 823529, 2609865, 2245355, 1765176, 1842848, 3002946, 2119988, 2006494, 1607915}}

\adfDesignF{13}{20888}{1}{3258529}{990328}{1}{878200}{"Audrey 18 May 2014"}{{1, 878200, 79222, 3166250, 198430, 1612138, 877564, 2011010, 1658993, 234352, 2493289, 1994431, 1088765}}

\adfDesignF{13}{21053}{1}{3284269}{2175508}{1}{1808134}{"Audrey 18 May 2014"}{{1, 1808134, 3280292, 1600192, 2679453, 949738, 1268324, 3085593, 512636, 1909892, 777677, 873713, 959969}}

\adfDesignF{13}{21333}{1}{3327949}{2892221}{1}{532112}{"Audrey 18 May 2014"}{{1, 532112, 1279624, 1592539, 146702, 1522880, 452856, 3308629, 2958570, 997441, 2362974, 1129908, 355509}}

\adfDesignF{13}{21683}{1}{3382549}{2082753}{1}{1455744}{"Audrey 18 May 2014"}{{1, 1455744, 3349742, 2998472, 648667, 419114, 2180039, 188687, 76583, 3192810, 775426, 877213, 750247}}

\adfDesignF{13}{21685}{1}{3382861}{1328745}{1}{528470}{"Audrey 18 May 2014"}{{1, 528470, 1685323, 3001730, 2627231, 2040645, 787821, 912017, 503015, 3119670, 1165106, 2271488, 1654649}}

\adfDesignF{13}{22693}{1}{3540109}{2895411}{1}{2563062}{"Audrey 18 May 2014"}{{1, 2563062, 2127487, 1630968, 1653437, 600303, 1013879, 25503, 1197610, 1590427, 1836263, 1365512, 2096093}}

\adfDesignF{13}{22857}{1}{3565693}{1414064}{1}{1510244}{"Audrey 18 May 2014"}{{1, 1510244, 2189463, 2915273, 2033625, 3077266, 1661801, 1274315, 837198, 3478363, 1809857, 917642, 3254803}}

\adfDesignF{13}{23027}{1}{3592213}{863299}{1}{930234}{"Audrey 18 May 2014"}{{1, 930234, 3512973, 468600, 3381489, 841781, 1163736, 1076757, 1259283, 18496, 2500007, 3592077, 2807844}}

\adfDesignF{13}{23048}{1}{3595489}{2503663}{1}{407423}{"Audrey 18 May 2014"}{{1, 407423, 560266, 2039864, 1014589, 714795, 100752, 2579672, 2338221, 2226488, 2118658, 2976659, 900057}}

\adfDesignF{13}{23172}{1}{3614833}{1342441}{1}{1001961}{"Audrey 18 May 2014"}{{1, 1001961, 3580262, 2150908, 2259151, 1689175, 972577, 1158190, 2217099, 1718651, 1589403, 1526300, 1825320}}

\adfDesignF{13}{23392}{1}{3649153}{967838}{1}{3447944}{"Audrey 18 May 2014"}{{1, 3447944, 1358299, 1330444, 908131, 3556943, 1188038, 1127629, 633467, 2003534, 157610, 2238233, 295492}}

\adfDesignF{13}{23967}{1}{3738853}{332163}{1}{1618179}{"Audrey 18 May 2014"}{{1, 1618179, 1316344, 3319534, 1320045, 2557507, 1383730, 2281883, 189404, 839494, 2884377, 19550, 964217}}

\adfDesignF{13}{24016}{1}{3746497}{57609}{1}{2009978}{"Audrey 18 May 2014"}{{1, 2009978, 2746013, 167883, 1644778, 2444629, 3551255, 2196583, 1636542, 289567, 1244079, 1969588, 2578086}}

\adfDesignF{13}{24548}{1}{3829489}{1800909}{1}{3693804}{"Audrey 18 May 2014"}{{1, 3693804, 2065602, 1433562, 2204296, 847118, 906505, 226266, 111103, 1687638, 962214, 761887, 417960}}

\adfDesignF{13}{24731}{1}{3858037}{3028551}{1}{699627}{"Audrey 18 May 2014"}{{1, 699627, 1068865, 3501645, 3465526, 3780263, 963550, 3074766, 1893600, 2361770, 3251097, 3231025, 3572561}}

\adfDesignF{13}{24836}{1}{3874417}{15903}{1}{2387215}{"Audrey 18 May 2014"}{{1, 2387215, 728099, 3397413, 3298942, 478901, 3803274, 1432450, 1808133, 3351486, 53486, 1169255, 1337847}}

\adfDesignF{13}{25045}{1}{3907021}{131674}{1}{982993}{"Audrey 18 May 2014"}{{1, 982993, 2525392, 3562297, 2466440, 3094433, 2592461, 940439, 808096, 643734, 987681, 511796, 419342}}

\adfDesignF{13}{25138}{1}{3921529}{658119}{1}{1381076}{"Audrey 18 May 2014"}{{1, 1381076, 1956640, 1733675, 345002, 365594, 554278, 1896212, 131796, 2523961, 1499929, 3543115, 3676366}}

\adfDesignF{13}{25187}{1}{3929173}{2050547}{1}{904026}{"Audrey 18 May 2014"}{{1, 904026, 2883022, 308001, 67381, 206887, 2592262, 3452768, 2004346, 3476316, 1749280, 698105, 1303470}}

\adfDesignF{13}{25282}{1}{3943993}{2182053}{1}{197790}{"Audrey 18 May 2014"}{{1, 197790, 417533, 582643, 1427503, 3247486, 1555960, 3554610, 2231734, 2971300, 2974063, 1252806, 3250529}}

\adfDesignF{13}{25401}{1}{3962557}{1257544}{1}{2694752}{"Audrey 18 May 2014"}{{1, 2694752, 1484672, 3390952, 3327751, 3006459, 3601147, 2976583, 2460964, 3778969, 1111274, 482666, 3384266}}

\adfDesignF{13}{25578}{1}{3990169}{1452207}{1}{1992704}{"Audrey 18 May 2014"}{{1, 1992704, 679069, 1489775, 2845938, 551553, 1788769, 1340803, 2348743, 1039142, 227249, 3691824, 1955275}}

\adfDesignF{13}{26385}{1}{4116061}{2696324}{1}{1724676}{"Audrey 18 May 2014"}{{1, 1724676, 2894838, 3227379, 885416, 1210277, 2956993, 1771475, 50752, 2719187, 136842, 1407574, 1594895}}

\adfDesignF{13}{26472}{1}{4129633}{4051739}{1}{2457306}{"Audrey 18 May 2014"}{{1, 2457306, 3405036, 786396, 3402422, 3097460, 1997433, 4015550, 3295807, 1335790, 1862057, 2215176, 1036997}}

\adfDesignF{13}{26605}{1}{4150381}{686812}{1}{3418987}{"Audrey 18 May 2014"}{{1, 3418987, 2876908, 2759247, 1158646, 2210437, 1702133, 1819834, 3217342, 1031203, 3999481, 423048, 285019}}

\adfDesignF{13}{26725}{1}{4169101}{393839}{1}{3904505}{"Audrey 18 May 2014"}{{1, 3904505, 3499224, 1614778, 2347196, 3964952, 2136248, 4037772, 3840350, 2076332, 2480605, 86254, 3364591}}

\adfDesignF{13}{27528}{1}{4294369}{1654042}{1}{2987420}{"Audrey 18 May 2014"}{{1, 2987420, 358268, 2813952, 1564783, 2888696, 3274739, 1587004, 3992514, 2587841, 107387, 3529764, 73845}}

\adfDesignF{13}{27580}{1}{4302481}{861596}{1}{3289919}{"Audrey 18 May 2014"}{{1, 3289919, 581544, 1101375, 1207412, 583973, 2607409, 2564020, 3983147, 1119115, 1215707, 4127576, 3433688}}

\adfDesignF{13}{28273}{1}{4410589}{1063515}{1}{806526}{"Audrey 18 May 2014"}{{1, 806526, 1701778, 1422907, 2696816, 1129989, 1092555, 79976, 2269840, 3852555, 803443, 2754316, 3442243}}

\adfDesignF{13}{28355}{1}{4423381}{337722}{1}{795965}{"Audrey 18 May 2014"}{{1, 795965, 3843976, 25616, 2076411, 2828156, 1519068, 2611032, 4333459, 4386612, 2674992, 1139549, 305449}}

\adfDesignF{13}{28538}{1}{4451929}{1852600}{1}{1671875}{"Audrey 18 May 2014"}{{1, 1671875, 133330, 3008720, 336403, 2671197, 3879244, 827939, 3395158, 3527315, 4275420, 4033248, 3403653}}

\adfDesignF{13}{29052}{1}{4532113}{2835217}{1}{3052811}{"Audrey 18 May 2014"}{{1, 3052811, 3645154, 1784327, 3028915, 2564346, 4396203, 2870027, 1456116, 951947, 1050366, 3756953, 3167625}}

\adfDesignF{13}{29581}{1}{4614637}{3838352}{1}{4595336}{"Audrey 18 May 2014"}{{1, 4595336, 3357641, 2133087, 1079127, 2241191, 379747, 3161346, 2205905, 3083194, 1631358, 3442130, 377759}}

\adfDesignF{13}{29593}{1}{4616509}{2290477}{1}{1823794}{"Audrey 18 May 2014"}{{1, 1823794, 2120882, 3620951, 3365193, 912192, 3564627, 1933696, 1141308, 638596, 813177, 4584270, 3180367}}

\adfDesignF{13}{29603}{1}{4618069}{837917}{1}{2453728}{"Audrey 18 May 2014"}{{1, 2453728, 1347648, 2560432, 2518205, 4560102, 1274224, 143588, 4175916, 4371855, 2836926, 3071047, 3012811}}

\adfDesignF{13}{29971}{1}{4675477}{720028}{1}{4586469}{"Audrey 18 May 2014"}{{1, 4586469, 2166026, 4327364, 455825, 1717806, 3647883, 2405678, 2908222, 2760329, 277241, 500678, 2299340}}

\adfDesignF{13}{30128}{1}{4699969}{3364159}{1}{2889542}{"Audrey 18 May 2014"}{{1, 2889542, 340985, 3077617, 2937103, 2754487, 1072183, 1644704, 1831652, 292484, 1076717, 4014529, 1567841}}

\adfDesignF{13}{30232}{1}{4716193}{4659018}{1}{424238}{"Audrey 18 May 2014"}{{1, 424238, 3239571, 3319768, 3318352, 1553855, 3177108, 4430041, 2971837, 459095, 922289, 520923, 3960080}}

\adfDesignF{13}{30245}{1}{4718221}{1678962}{1}{3462212}{"Audrey 18 May 2014"}{{1, 3462212, 2543847, 4324820, 4220605, 1593337, 1979780, 3642726, 8934, 3463353, 3776762, 430111, 3581059}}

\adfDesignF{13}{30428}{1}{4746769}{1029066}{1}{2509210}{"Audrey 18 May 2014"}{{1, 2509210, 1435424, 3139375, 2546408, 3310157, 1636846, 1006720, 2819546, 1669072, 1598265, 3529465, 3280125}}

\adfDesignF{13}{30456}{1}{4751137}{3492588}{1}{2098995}{"Audrey 18 May 2014"}{{1, 2098995, 3158555, 2319918, 3708603, 274993, 2800179, 1404460, 2289899, 3556592, 2053557, 3347883, 1493187}}

\adfDesignF{13}{30510}{1}{4759561}{2943665}{1}{271238}{"Audrey 18 May 2014"}{{1, 271238, 1518267, 208143, 3138013, 836025, 1984227, 883829, 3201415, 1573808, 1027336, 4063223, 332719}}

\adfDesignF{13}{30738}{1}{4795129}{3284147}{1}{2397255}{"Audrey 18 May 2014"}{{1, 2397255, 3692137, 921865, 4788087, 2179499, 3955813, 831336, 1638474, 1175971, 2864344, 584397, 3741595}}

\adfDesignF{13}{30853}{1}{4813069}{3914313}{1}{4774288}{"Audrey 18 May 2014"}{{1, 4774288, 2288433, 459118, 3300142, 1510877, 981069, 473556, 1696068, 187846, 2130740, 3353721, 2909486}}

\adfDesignF{13}{31396}{1}{4897777}{1253400}{1}{2988138}{"Audrey 18 May 2014"}{{1, 2988138, 2880539, 1688819, 4594849, 1783745, 1827682, 1836726, 623312, 942165, 3352515, 1529026, 441368}}

\adfDesignF{13}{31585}{1}{4927261}{2990339}{1}{2965757}{"Audrey 18 May 2014"}{{1, 2965757, 1844817, 3971720, 1900091, 3250668, 2648554, 1988093, 3710273, 1911499, 1105237, 3989159, 277697}}

\adfDesignF{13}{31686}{1}{4943017}{3804720}{1}{3680276}{"Audrey 18 May 2014"}{{1, 3680276, 1352238, 4257173, 2044919, 4302736, 4494616, 1615825, 194901, 3332789, 758032, 2327287, 1297309}}

\adfDesignF{13}{31752}{1}{4953313}{3145144}{1}{4869435}{"Audrey 18 May 2014"}{{1, 4869435, 1814424, 285653, 4125960, 799804, 1711360, 1556660, 4756513, 2751484, 737657, 3592150, 2718777}}

\adfDesignF{13}{31990}{1}{4990441}{2964731}{1}{3850812}{"Audrey 18 May 2014"}{{1, 3850812, 1968273, 2429963, 4282906, 1890381, 447523, 2900351, 4658633, 2323780, 305204, 4427502, 457317}}

\adfDesignF{13}{32718}{1}{5104009}{3615298}{1}{1897782}{"Audrey 18 May 2014"}{{1, 1897782, 4024800, 2693028, 3147962, 61955, 1132486, 124305, 1599539, 2697811, 207257, 3462416, 4470703}}

\adfDesignF{13}{32813}{1}{5118829}{3606216}{1}{4050922}{"Audrey 18 May 2014"}{{1, 4050922, 1447739, 4667084, 2727439, 3451459, 1189282, 328074, 610958, 4738263, 4785736, 4919141, 2915705}}

\adfDesignF{13}{32888}{1}{5130529}{4383827}{1}{4028086}{"Audrey 18 May 2014"}{{1, 4028086, 4422910, 2417709, 410348, 4244940, 2508401, 1398944, 1519324, 824627, 4952923, 4113231, 5072259}}

\adfDesignF{13}{33011}{1}{5149717}{966786}{1}{1051775}{"Audrey 18 May 2014"}{{1, 1051775, 4492704, 376721, 1354078, 253798, 2810755, 4400803, 1390536, 1073966, 764568, 4599782, 3179098}}

\adfDesignF{13}{33441}{1}{5216797}{212988}{1}{3340372}{"Audrey 18 May 2014"}{{1, 3340372, 3198212, 78414, 1769435, 2730384, 3368530, 2543078, 2982096, 3045716, 852967, 367016, 1807764}}

\adfDesignF{13}{33497}{1}{5225533}{148605}{1}{4690183}{"Audrey 18 May 2014"}{{1, 4690183, 39582, 4538148, 4300357, 1277261, 218832, 4738660, 3100043, 358618, 5161653, 2270048, 659812}}

\adfDesignF{13}{33782}{1}{5269993}{1551822}{1}{2151515}{"Audrey 18 May 2014"}{{1, 2151515, 2503829, 1916384, 2338406, 827794, 1263581, 3266277, 3154029, 2597527, 1526625, 1829667, 2974330}}

\adfDesignF{13}{34098}{1}{5319289}{2235936}{1}{4756539}{"Audrey 18 May 2014"}{{1, 4756539, 3691885, 4933159, 1701850, 619794, 2025519, 4303271, 4393468, 2687356, 1508723, 4445485, 2167973}}

\adfDesignF{13}{34686}{1}{5411017}{4967432}{1}{2499471}{"Audrey 18 May 2014"}{{1, 2499471, 670287, 3834637, 2509842, 3912598, 2224269, 555219, 1081193, 2561644, 1283547, 1935354, 3987023}}

\adfDesignF{13}{34731}{1}{5418037}{4086301}{1}{685127}{"Audrey 18 May 2014"}{{1, 685127, 1952597, 3991112, 333768, 5216951, 396914, 5221048, 819067, 2570308, 769265, 4672480, 461547}}

\adfDesignF{13}{34915}{1}{5446741}{1746770}{1}{990929}{"Audrey 18 May 2014"}{{1, 990929, 1815561, 813423, 3026341, 614045, 3220472, 2653106, 1282853, 1358447, 616300, 4601557, 793929}}

\adfDesignF{13}{35001}{1}{5460157}{3386767}{1}{3921644}{"Audrey 18 May 2014"}{{1, 3921644, 510413, 2743871, 959628, 4459008, 3222679, 248622, 3119549, 392049, 140339, 3072501, 4510481}}

\adfDesignF{13}{35068}{1}{5470609}{156801}{1}{847499}{"Audrey 18 May 2014"}{{1, 847499, 1887564, 4059874, 5095185, 4625473, 3980288, 3178132, 3880109, 5427491, 1206238, 3736150, 370259}}

\adfDesignF{13}{35125}{1}{5479501}{4890882}{1}{5455624}{"Audrey 18 May 2014"}{{1, 5455624, 243025, 83634, 3088847, 1683641, 2802680, 1555353, 2894197, 2724843, 2518563, 1794724, 2552373}}

\adfDesignF{13}{35176}{1}{5487457}{4691468}{1}{2097693}{"Audrey 18 May 2014"}{{1, 2097693, 978347, 1147670, 2190270, 3487435, 3100104, 1893426, 5089618, 4837504, 732877, 1462012, 420328}}

\adfDesignF{13}{35540}{1}{5544241}{1010103}{1}{1170896}{"Audrey 18 May 2014"}{{1, 1170896, 895613, 4215303, 4034853, 1331922, 2591222, 4397349, 1534342, 4601033, 3244350, 4477702, 770860}}

\adfDesignF{13}{35938}{1}{5606329}{7300}{1}{3938247}{"Audrey 18 May 2014"}{{1, 3938247, 3593747, 3513918, 3305488, 5143813, 5255306, 4540797, 3447767, 4830492, 354603, 72357, 1247767}}

\adfDesignF{13}{36461}{1}{5687917}{12746}{1}{291611}{"Audrey 18 May 2014"}{{1, 291611, 2616171, 998022, 544303, 3418248, 2639112, 1855581, 4910947, 5175025, 4516420, 571270, 902874}}

\adfDesignF{13}{36922}{1}{5759833}{371861}{1}{2230174}{"Audrey 18 May 2014"}{{1, 2230174, 2676446, 2304372, 2884641, 1722506, 2035692, 921144, 3601443, 5095401, 4465577, 3270913, 3350688}}

\adfDesignF{13}{37932}{1}{5917393}{5721191}{1}{619916}{"Audrey 18 May 2014"}{{1, 619916, 2593457, 5315870, 2622613, 1951151, 5007151, 3099415, 5359433, 1041669, 5851286, 3076906, 4882883}}

\adfDesignF{13}{37935}{1}{5917861}{1340346}{1}{2284866}{"Audrey 18 May 2014"}{{1, 2284866, 5856698, 1028357, 824417, 3545378, 2276610, 2226009, 2971900, 2839560, 2698776, 1376348, 1660385}}

\adfDesignF{13}{37965}{1}{5922541}{1961752}{1}{410785}{"Audrey 18 May 2014"}{{1, 410785, 5200594, 319639, 172645, 3470391, 5258071, 3760658, 4069713, 3639012, 2196020, 243285, 871891}}

\adfDesignF{13}{39136}{1}{6105217}{2213186}{1}{2128429}{"Audrey 18 May 2014"}{{1, 2128429, 4679267, 4352924, 3091084, 41760, 3445954, 2614618, 1375282, 1064809, 3910755, 3074350, 746852}}

\adfDesignF{13}{39900}{1}{6224401}{7479}{1}{1506039}{"Audrey 18 May 2014"}{{1, 1506039, 418324, 3287020, 2159262, 705169, 4716971, 2304564, 6142391, 530653, 2147072, 4072709, 3131830}}

\adfDesignF{13}{40021}{1}{6243277}{2933239}{1}{5496355}{"Audrey 18 May 2014"}{{1, 5496355, 5727918, 3731563, 5895301, 3308362, 674836, 1113803, 5182438, 4730380, 1645865, 5923432, 272685}}

\adfDesignF{13}{40630}{1}{6338281}{884257}{1}{5191662}{"Audrey 18 May 2014"}{{1, 5191662, 4518174, 1513049, 4529146, 2297166, 6238573, 3512855, 1567164, 674670, 4093601, 632969, 3260656}}

\adfDesignF{13}{40706}{1}{6350137}{2106325}{1}{950418}{"Audrey 18 May 2014"}{{1, 950418, 86748, 3031993, 303159, 3004361, 2519615, 6335411, 6196617, 5276626, 5034466, 6177014, 5534530}}

\adfDesignF{13}{41732}{1}{6510193}{1229535}{1}{608091}{"Audrey 18 May 2014"}{{1, 608091, 2212074, 6213074, 1996500, 339895, 1483081, 4192467, 3361504, 3889952, 5236240, 1062312, 1955774}}

\adfDesignF{13}{41963}{1}{6546229}{4419252}{1}{2563873}{"Audrey 18 May 2014"}{{1, 2563873, 3086176, 5402768, 595823, 6141726, 5101664, 4966001, 3048659, 322124, 6229383, 1383197, 435979}}

\adfDesignF{13}{41987}{1}{6549973}{3672901}{1}{4212275}{"Audrey 18 May 2014"}{{1, 4212275, 6066060, 3165417, 3706846, 3349005, 4287301, 12960, 3609018, 3361654, 1323421, 3222178, 2983702}}

\adfDesignF{13}{42208}{1}{6584449}{5342916}{1}{351060}{"Audrey 18 May 2014"}{{1, 351060, 1991667, 5146608, 1983329, 1503684, 1444261, 6524762, 4582947, 5598466, 5291950, 2850548, 2237411}}

\adfDesignF{13}{42253}{1}{6591469}{3316064}{1}{2000970}{"Audrey 18 May 2014"}{{1, 2000970, 5151823, 1014857, 5233239, 3832042, 4516261, 5365639, 1346525, 894934, 5336874, 151783, 4703866}}

\adfDesignF{13}{42797}{1}{6676333}{6501062}{1}{3911760}{"Audrey 18 May 2014"}{{1, 3911760, 4820252, 2106605, 2116230, 4964443, 6051926, 6101730, 5323497, 5572091, 4354416, 1804265, 6283448}}

\adfDesignF{13}{43321}{1}{6758077}{6200822}{1}{2454404}{"Audrey 18 May 2014"}{{1, 2454404, 3222032, 1685068, 2644704, 5327531, 806535, 331454, 4988387, 3518526, 1255207, 793869, 4590}}

\adfDesignF{13}{44218}{1}{6898009}{3976243}{1}{649254}{"Audrey 18 May 2014"}{{1, 649254, 324535, 5961985, 4164813, 1971502, 4111459, 3474784, 2873659, 2514120, 5902783, 2735653, 6703506}}

\adfDesignF{13}{44520}{1}{6945121}{2122690}{1}{6069021}{"Audrey 18 May 2014"}{{1, 6069021, 4217564, 4905230, 2188896, 1970041, 6880094, 6272258, 348941, 3126078, 2914703, 3845859, 5877161}}

\adfDesignF{13}{44760}{1}{6982561}{3962562}{1}{4761515}{"Audrey 18 May 2014"}{{1, 4761515, 5638836, 2778852, 4356879, 2364343, 4228004, 3544520, 2292018, 581588, 5171147, 874381, 5303282}}

\adfDesignF{13}{45397}{1}{7081933}{1027168}{1}{3006905}{"Audrey 18 May 2014"}{{1, 3006905, 2145657, 4159925, 6793143, 1283711, 3606671, 6045638, 2821092, 2128994, 4777885, 3308336, 2413640}}

\adfDesignF{13}{45437}{1}{7088173}{1629864}{1}{4752237}{"Audrey 18 May 2014"}{{1, 4752237, 921755, 6536357, 1335207, 1049227, 73949, 5732919, 2988927, 370577, 6060726, 2145765, 3473218}}

\adfDesignF{13}{45583}{1}{7110949}{906292}{1}{2565761}{"Audrey 18 May 2014"}{{1, 2565761, 4920544, 4598506, 2998541, 4725029, 1733796, 3128591, 661254, 4180486, 912042, 3586093, 1544101}}

\adfDesignF{13}{46005}{1}{7176781}{5598337}{1}{1787412}{"Audrey 18 May 2014"}{{1, 1787412, 4297441, 3261954, 2954743, 194683, 5326830, 5461128, 2347016, 4880757, 890028, 5567171, 6091522}}

\adfDesignF{13}{46311}{1}{7224517}{6716104}{1}{684696}{"Audrey 18 May 2014"}{{1, 684696, 2479769, 3603435, 5505573, 6432480, 3802853, 2841201, 6042789, 7020278, 3349225, 1998977, 6810342}}

\adfDesignF{13}{46576}{1}{7265857}{2812755}{1}{1388330}{"Audrey 18 May 2014"}{{1, 1388330, 2707368, 5199056, 1350682, 5434786, 3650588, 2212117, 1426136, 1360380, 5826105, 898254, 4875482}}

\adfDesignF{13}{47133}{1}{7352749}{165624}{1}{3212426}{"Audrey 18 May 2014"}{{1, 3212426, 1998239, 2400097, 4640177, 4544451, 5874602, 4213323, 3875904, 5212241, 1137651, 1935617, 5071765}}

\adfDesignF{13}{47297}{1}{7378333}{3744924}{1}{3888511}{"Audrey 18 May 2014"}{{1, 3888511, 4061892, 2809708, 3200043, 6124799, 4273581, 1471975, 4884477, 843482, 6045145, 3859394, 2806990}}

\adfDesignF{13}{47321}{1}{7382077}{7040555}{1}{5147249}{"Audrey 18 May 2014"}{{1, 5147249, 1264079, 3611256, 859129, 6244195, 51413, 2816541, 5669796, 6268178, 3432586, 1032959, 513003}}

\adfDesignF{13}{48316}{1}{7537297}{2474478}{1}{4799702}{"Audrey 18 May 2014"}{{1, 4799702, 1529361, 1486983, 1212469, 3168914, 7143557, 6885924, 6669081, 5007243, 7310623, 3484317, 4062904}}

\adfDesignF{13}{48358}{1}{7543849}{3532933}{1}{5147971}{"Audrey 18 May 2014"}{{1, 5147971, 3529049, 5274574, 1267207, 1464247, 4189698, 6179234, 7251162, 3860391, 4584266, 4194418, 5864725}}

\adfDesignF{13}{48747}{1}{7604533}{201626}{1}{5350844}{"Audrey 18 May 2014"}{{1, 5350844, 890823, 4343151, 2180647, 5339797, 152164, 3687172, 390247, 6892932, 7381719, 3333257, 5684444}}

\adfDesignF{13}{48792}{1}{7611553}{2999755}{1}{2083720}{"Audrey 18 May 2014"}{{1, 2083720, 414398, 4382028, 1455171, 216828, 2277186, 6317932, 3817747, 111079, 5430256, 3854451, 85415}}

\adfDesignF{13}{49231}{1}{7680037}{4096}{1}{1402863}{"Audrey 18 May 2014"}{{1, 1402863, 7435482, 5451899, 2699906, 3663440, 309171, 2147035, 650360, 625191, 4776470, 591591, 966739}}

\adfDesignF{13}{49247}{1}{7682533}{1198764}{1}{3077488}{"Audrey 18 May 2014"}{{1, 3077488, 5580673, 5699797, 5588082, 5050445, 5401733, 3962517, 6330934, 1091279, 1773801, 4421139, 5799842}}

\adfDesignF{13}{49330}{1}{7695481}{1497966}{1}{7453922}{"Audrey 18 May 2014"}{{1, 7453922, 3613539, 151567, 2725645, 5882843, 1544704, 1329192, 263435, 6532705, 1646765, 3897217, 3435870}}

\adfDesignF{13}{49358}{1}{7699849}{2831833}{1}{3720546}{"Audrey 18 May 2014"}{{1, 3720546, 5099423, 2013035, 93753, 1489589, 2580109, 4771365, 4097300, 1276204, 2202742, 7055190, 4099988}}

\adfDesignF{13}{49896}{1}{7783777}{3296528}{1}{1062721}{"Audrey 18 May 2014"}{{1, 1062721, 4367580, 6057418, 6059838, 1414471, 588905, 2688374, 2860466, 1018406, 1737315, 364823, 2914790}}

~

\adfDesignF{15}{3471}{1}{728911}{400308}{1}{78969}{"Mononoko 13 Apr 2014"}{{1, 78969, 269356, 422073, 498351, 375229, 497840, 112175, 621103, 190528, 353681, 152102, 347380, 414646, 40032}}

\adfDesignF{15}{3627}{1}{761671}{212943}{1}{583695}{"Mononoko 13 Apr 2014"}{{1, 583695, 606370, 333528, 88386, 204427, 400576, 253095, 387220, 125360, 557243, 539071, 663877, 40923, 547925}}

\adfDesignF{15}{3733}{1}{783931}{293}{1}{509}{"Mononoko 13 Apr 2014"}{{1, 509, 259081, 171821, 440548, 34666, 398512, 588410, 39048, 277157, 749264, 384910, 720371, 573062, 66226}}

\adfDesignF{15}{3813}{1}{800731}{779314}{1}{394840}{"Mononoko 13 Apr 2014"}{{1, 394840, 303555, 638658, 717469, 444587, 477605, 603314, 632377, 591336, 356143, 729017, 694593, 330427, 292657}}

\adfDesignF{15}{3877}{1}{814171}{39}{1}{64725}{"Mononoko 13 Apr 2014"}{{1, 64725, 415830, 546003, 137749, 631575, 794307, 692780, 531846, 582470, 182595, 769310, 519732, 550493, 93952}}

\adfDesignF{15}{4048}{1}{850081}{310}{1}{270807}{"Mononoko 13 Apr 2014"}{{1, 270807, 793460, 397931, 282190, 145754, 242486, 699195, 708506, 2156, 704326, 336788, 207507, 593725, 565735}}

\adfDesignF{15}{4823}{1}{1012831}{114}{1}{261892}{"Mononoko 13 Apr 2014"}{{1, 261892, 530006, 906957, 720679, 20649, 303199, 454939, 509903, 787619, 992181, 447740, 27886, 608802, 517364}}

\adfDesignF{15}{5183}{1}{1088431}{543}{1}{133989}{"Mononoko 13 Apr 2014"}{{1, 133989, 471207, 1026137, 466573, 526781, 285921, 762962, 899036, 1010541, 561649, 668521, 942693, 251689, 699748}}

\adfDesignF{15}{5919}{1}{1242991}{1192582}{1}{987292}{"Mononoko 13 Apr 2014"}{{1, 987292, 652001, 164176, 1139010, 260229, 742132, 53738, 512643, 1042421, 982761, 756558, 537252, 566172, 304551}}

\adfDesignF{15}{6427}{1}{1349671}{78}{1}{342782}{"Mononoko 13 Apr 2014"}{{1, 342782, 1191277, 1302551, 974688, 1194321, 85605, 655899, 826167, 459019, 155349, 921284, 852166, 570624, 1265635}}

\adfDesignF{15}{6427}{1}{1349671}{923357}{1}{1265635}{"Mononoko 13 Apr 2014"}{{1, 1265635, 570624, 852166, 921284, 155349, 459019, 826167, 655899, 85605, 1194321, 974688, 1302551, 1191277, 342782}}

\adfDesignF{15}{7163}{1}{1504231}{8}{1}{107579}{"Mononoko 13 Apr 2014"}{{1, 107579, 1192158, 430422, 1129696, 230801, 503893, 332500, 908551, 590342, 1273429, 892759, 1483804, 165258, 1288424}}

\adfDesignF{15}{7696}{1}{1616161}{19}{1}{2187}{"Mononoko 13 Apr 2014"}{{1, 2187, 1550647, 559211, 1176741, 604255, 1102148, 701625, 717086, 590912, 1011905, 511826, 980050, 339864, 1464669}}

\adfDesignF{15}{8155}{1}{1712551}{1357867}{1}{124431}{"Mononoko 13 Apr 2014"}{{1, 124431, 1612721, 898224, 694631, 1180991, 1314913, 529514, 981911, 1641648, 531559, 273207, 1282867, 1544967, 1088823}}

\adfDesignF{15}{8516}{1}{1788361}{119}{1}{337416}{"Mononoko 13 Apr 2014"}{{1, 337416, 707435, 191846, 395180, 1647081, 418336, 1502768, 396436, 1600020, 141279, 1032609, 1366519, 1200079, 1581522}}

\adfDesignF{15}{8734}{1}{1834141}{128}{1}{315590}{"Mononoko 13 Apr 2014"}{{1, 315590, 1357659, 929646, 1455062, 139256, 1782680, 741565, 1443314, 1221038, 1694884, 1570012, 1569058, 1295322, 992182}}

\adfDesignF{15}{9104}{1}{1911841}{122}{1}{227143}{"Mononoko 13 Apr 2014"}{{1, 227143, 1001223, 1573416, 433153, 510337, 733679, 604650, 1093033, 1210618, 1401503, 951019, 305968, 1157233, 268070}}

\adfDesignF{15}{9133}{1}{1917931}{95}{1}{181780}{"Mononoko 13 Apr 2014"}{{1, 181780, 1853132, 769982, 559442, 911347, 1449604, 639168, 1616991, 272713, 1006583, 286547, 1343562, 1448889, 1085776}}

\adfDesignF{15}{10255}{1}{2153551}{475867}{1}{183912}{"Audrey 18 May 2014"}{{1, 183912, 2105289, 976278, 1232013, 413493, 131704, 957951, 984104, 1955257, 1740057, 1837935, 1243862, 193169, 1119832}}

\adfDesignF{15}{10436}{1}{2191561}{874080}{1}{1645157}{"Audrey 18 May 2014"}{{1, 1645157, 976186, 1138841, 169454, 917273, 1789725, 1067398, 825894, 1306578, 1274287, 948240, 147977, 226826, 715529}}

\adfDesignF{15}{10594}{1}{2224741}{1701199}{1}{279268}{"Audrey 18 May 2014"}{{1, 279268, 95328, 809098, 1585340, 148415, 635390, 977101, 1884195, 1852681, 2076325, 1310083, 1152312, 1756189, 1011461}}

\adfDesignF{15}{11624}{1}{2441041}{699660}{1}{1675998}{"Audrey 18 May 2014"}{{1, 1675998, 2391279, 2035371, 1041070, 1982511, 487803, 959633, 2099859, 827737, 458529, 277240, 1531170, 746852, 572234}}

\adfDesignF{15}{11877}{1}{2494171}{16384}{1}{123648}{"Audrey 18 May 2014"}{{1, 123648, 2053845, 2323682, 109420, 1180656, 1924458, 1492700, 715600, 1792575, 1313514, 446065, 1441797, 1949060, 592176}}

\adfDesignF{15}{12699}{1}{2666791}{394954}{1}{1903941}{"Audrey 18 May 2014"}{{1, 1903941, 990853, 363199, 2560586, 1373670, 696786, 1582020, 1642886, 1512887, 1293120, 66064, 93918, 660706, 1260109}}

\adfDesignF{15}{13719}{1}{2880991}{1208896}{1}{2343110}{"Audrey 18 May 2014"}{{1, 2343110, 1091959, 755300, 1426565, 934175, 1538126, 2553473, 1752681, 1352005, 1946815, 1880746, 2116550, 373010, 102421}}

\adfDesignF{15}{14911}{1}{3131311}{417363}{1}{899888}{"Audrey 18 May 2014"}{{1, 899888, 680901, 2833919, 1133830, 1109556, 1253780, 122364, 1343917, 2978187, 2021754, 977643, 2328046, 2084786, 2150605}}

\adfDesignF{15}{14975}{1}{3144751}{1638218}{1}{731126}{"Audrey 18 May 2014"}{{1, 731126, 452896, 629102, 1547592, 194041, 2413054, 326041, 1781615, 881531, 2950709, 571, 2365814, 734034, 715628}}

\adfDesignF{15}{15626}{1}{3281461}{399752}{1}{1581984}{"Audrey 18 May 2014"}{{1, 1581984, 1515386, 2413281, 2069508, 1778328, 1328005, 2730273, 2190155, 504372, 1503132, 371472, 2317263, 1959486, 707581}}

\adfDesignF{15}{16316}{1}{3426361}{184507}{1}{431490}{"Audrey 18 May 2014"}{{1, 431490, 2016082, 427890, 793615, 2991649, 2252065, 136362, 1368288, 2898849, 434711, 742806, 1273917, 1630183, 3160258}}

\adfDesignF{15}{16629}{1}{3492091}{16384}{1}{2507626}{"Audrey 18 May 2014"}{{1, 2507626, 844722, 322828, 2334890, 2527626, 3445871, 26570, 2018631, 594683, 964464, 1030685, 2620799, 1150632, 562518}}

\adfDesignF{15}{17977}{1}{3775171}{2981099}{1}{1282911}{"Audrey 18 May 2014"}{{1, 1282911, 3108222, 2523440, 120013, 3198950, 2400376, 1387100, 842804, 3346676, 576220, 91884, 3055020, 408927, 308482}}

\adfDesignF{15}{17993}{1}{3778531}{881513}{1}{3402745}{"Audrey 18 May 2014"}{{1, 3402745, 78733, 2917123, 2094449, 3395386, 3381746, 1636219, 820903, 2973144, 383144, 772571, 2063579, 40505, 2489469}}

\adfDesignF{15}{18247}{1}{3831871}{2172186}{1}{448206}{"Audrey 18 May 2014"}{{1, 448206, 2781261, 1257788, 3266808, 2886767, 560013, 2140565, 1711023, 3104024, 945103, 2823652, 2741916, 863060, 1292910}}

\adfDesignF{15}{18351}{1}{3853711}{2804622}{1}{2487409}{"Audrey 18 May 2014"}{{1, 2487409, 1155983, 2299929, 3138484, 3606307, 319643, 636311, 355967, 3026432, 247403, 1046659, 2061417, 1197815, 1542506}}

\adfDesignF{15}{18971}{1}{3983911}{2965399}{1}{967533}{"Audrey 18 May 2014"}{{1, 967533, 618864, 1470945, 3350422, 3194713, 1193281, 1337973, 990269, 3276521, 789197, 1823097, 2027074, 1522697, 1340879}}

\adfDesignF{15}{20187}{1}{4239271}{1239661}{1}{1507791}{"Audrey 18 May 2014"}{{1, 1507791, 1687072, 2850028, 3537223, 73461, 161963, 3147778, 1030031, 1822858, 4165809, 2569517, 3643692, 359212, 3118461}}

\adfDesignF{15}{21575}{1}{4530751}{4530402}{1}{1679089}{"Audrey 18 May 2014"}{{1, 1679089, 3037404, 66549, 4312699, 1273682, 2225674, 852905, 973710, 1598085, 3257068, 625988, 640442, 3490492, 3150718}}

\adfDesignF{15}{22583}{1}{4742431}{2888421}{1}{1918638}{"Audrey 18 May 2014"}{{1, 1918638, 1984224, 2111938, 2378131, 3761582, 2048189, 3162190, 1585438, 986286, 980848, 775604, 4338448, 1045055, 1378014}}

\adfDesignF{15}{22843}{1}{4797031}{295184}{1}{2453452}{"Audrey 18 May 2014"}{{1, 2453452, 1885791, 4474311, 719296, 2460357, 4655390, 1922001, 2957111, 4272152, 2336673, 2485220, 989239, 2162640, 4602614}}

\adfDesignF{15}{23147}{1}{4860871}{2527800}{1}{2096364}{"Audrey 18 May 2014"}{{1, 2096364, 4245041, 4608512, 1433480, 3336229, 2793910, 4162242, 4608215, 1004860, 1524641, 4831468, 1314459, 505015, 2422531}}

\adfDesignF{15}{24068}{1}{5054281}{2725430}{1}{1195915}{"Audrey 18 May 2014"}{{1, 1195915, 2792655, 96583, 4629033, 1292300, 3127444, 4414103, 3795886, 3982730, 3761980, 730922, 2901804, 1161812, 1496799}}

\adfDesignF{15}{24087}{1}{5058271}{4782969}{1}{2966827}{"Audrey 18 May 2014"}{{1, 2966827, 3216157, 641027, 2422478, 3495601, 1911773, 4840990, 1035395, 446075, 1562669, 179671, 2059395, 3381849, 2189718}}

\adfDesignF{15}{24503}{1}{5145631}{787481}{1}{777743}{"Audrey 18 May 2014"}{{1, 777743, 4958737, 3429877, 824008, 3841049, 1439047, 3114635, 2591090, 3204447, 1304581, 2928841, 2217890, 4270295, 1117176}}

\adfDesignF{15}{24814}{1}{5210941}{2221765}{1}{5042978}{"Audrey 18 May 2014"}{{1, 5042978, 4745736, 4378061, 101354, 422345, 3417539, 1524680, 1868605, 3685756, 4788595, 1961365, 4151466, 4175216, 1423831}}

\adfDesignF{15}{24981}{1}{5246011}{16384}{1}{2353839}{"Audrey 18 May 2014"}{{1, 2353839, 4504315, 2767779, 504923, 2673303, 2863860, 4175672, 1593351, 3048347, 2572707, 28312, 1812035, 884881, 1692741}}

\adfDesignF{15}{25115}{1}{5274151}{4173306}{1}{2421442}{"Audrey 18 May 2014"}{{1, 2421442, 2209644, 4080168, 2415090, 634376, 2937291, 2365968, 2487955, 2832303, 4639774, 5189569, 698539, 3980179, 26758}}

\adfDesignF{15}{26079}{1}{5476591}{4782969}{1}{1389576}{"Audrey 18 May 2014"}{{1, 1389576, 1434769, 473940, 4628508, 4451527, 2220326, 2986243, 4078859, 1726745, 1025063, 1866689, 1055579, 923792, 4597929}}

\adfDesignF{15}{26873}{1}{5643331}{2725323}{1}{989894}{"Audrey 18 May 2014"}{{1, 989894, 4709720, 4173981, 96378, 3493377, 3351837, 733814, 5439389, 3561046, 2149953, 1301600, 199797, 1673292, 1985907}}

\adfDesignF{15}{30349}{1}{6373291}{4708691}{1}{3602018}{"Audrey 18 May 2014"}{{1, 3602018, 4546418, 6274368, 1899905, 5091474, 2758244, 4049275, 4558646, 5936371, 1281816, 13029, 4150889, 1913568, 4910306}}

\adfDesignF{15}{32264}{1}{6775441}{6159154}{1}{6683498}{"Audrey 18 May 2014"}{{1, 6683498, 4540322, 4420687, 705308, 6387808, 1321259, 3140893, 6496644, 1939268, 387632, 5546125, 5869667, 2633551, 4130865}}

\adfDesignF{15}{35744}{1}{7506241}{7436542}{1}{1098999}{"Audrey 18 May 2014"}{{1, 1098999, 7093896, 38997, 4534134, 5656498, 4505327, 610602, 548239, 3160173, 1849742, 1901915, 7307984, 6919005, 7318175}}

\adfDesignF{15}{36325}{1}{7628251}{16384}{1}{2008050}{"Audrey 18 May 2014"}{{1, 2008050, 1836904, 75656, 4412135, 1332306, 2642086, 7478302, 4862273, 6597463, 6295944, 2978115, 5941296, 2690322, 4246904}}

\adfDesignF{15}{37181}{1}{7808011}{1831}{1}{894977}{"Audrey 18 May 2014"}{{1, 894977, 6830105, 4382839, 6189600, 75030, 1229710, 7400198, 2218894, 809742, 7732980, 5683324, 1385719, 1206278, 808669}}

\adfDesignF{15}{37305}{1}{7834051}{7637863}{1}{3857130}{"Audrey 18 May 2014"}{{1, 3857130, 1434075, 1813027, 7373309, 4844839, 945792, 7172096, 3910617, 5847351, 2989211, 3031129, 7061931, 2110407, 2447442}}

\adfDesignF{15}{37607}{1}{7897471}{3795634}{1}{6038380}{"Audrey 18 May 2014"}{{1, 6038380, 1727725, 6124848, 2869342, 2591241, 1820946, 1906361, 2097993, 3585762, 5306229, 38145, 4263385, 7572101, 1442367}}

\adfDesignF{15}{38822}{1}{8152621}{569017}{1}{4745128}{"Audrey 18 May 2014"}{{1, 4745128, 4953744, 7029530, 2897253, 6320116, 636266, 3477118, 5922473, 765065, 1832504, 2771227, 7874380, 3353239, 4490303}}

\adfDesignF{15}{40011}{1}{8402311}{4211710}{1}{4260094}{"Audrey 18 May 2014"}{{1, 4260094, 5692917, 6704597, 2405110, 6589165, 4888088, 6537331, 8311083, 637562, 1813145, 7656440, 4574374, 1788942, 5359639}}

\adfDesignF{15}{40213}{1}{8444731}{8136940}{1}{7319538}{"Audrey 18 May 2014"}{{1, 7319538, 8326267, 3229448, 7021105, 6765352, 7152394, 4980958, 882069, 2925382, 1679378, 2417530, 3582237, 4333214, 6942975}}

\adfDesignF{15}{40515}{1}{8508151}{4782969}{1}{1425116}{"Audrey 18 May 2014"}{{1, 1425116, 412699, 993907, 4297883, 6465283, 3744643, 118160, 6890119, 4300459, 2042867, 3338392, 7977292, 624125, 8417960}}

\adfDesignF{15}{41087}{1}{8628271}{3033550}{1}{838314}{"Audrey 18 May 2014"}{{1, 838314, 6317917, 3144756, 2431773, 5018094, 7670924, 1234107, 4969614, 5364614, 3610176, 119033, 1076247, 513901, 831884}}

\adfDesignF{15}{43541}{1}{9143611}{5884218}{1}{594538}{"Audrey 18 May 2014"}{{1, 594538, 1719406, 5638239, 8109872, 8983005, 359645, 8419386, 2502551, 4140907, 160605, 8189428, 8148430, 1002821, 6036443}}

\adfDesignF{15}{44081}{1}{9257011}{1262883}{1}{5477250}{"Audrey 18 May 2014"}{{1, 5477250, 7458535, 6729485, 5195044, 2570782, 4048433, 7728817, 9100843, 4409433, 6686228, 8988339, 3326670, 2683694, 8909545}}

\adfDesignF{15}{44675}{1}{9381751}{4782969}{1}{5097201}{"Audrey 18 May 2014"}{{1, 5097201, 2703290, 253815, 2608915, 1191965, 6951859, 8393394, 5271978, 3891009, 8189785, 6714442, 7666818, 3855958, 2881827}}

\adfDesignF{15}{45459}{1}{9546391}{9270645}{1}{6579381}{"Audrey 18 May 2014"}{{1, 6579381, 1157796, 120262, 6446178, 4749817, 166279, 6031090, 4640178, 6902344, 4796573, 2800731, 2357505, 4785951, 5744260}}

\adfDesignF{15}{45631}{1}{9582511}{7268550}{1}{8263007}{"Audrey 18 May 2014"}{{1, 8263007, 6052382, 8651471, 4346427, 5021292, 1536540, 6583731, 5520401, 2448101, 4561218, 9365475, 6528909, 4993150, 2787983}}

\adfDesignF{15}{45752}{1}{9607921}{2321712}{1}{3218218}{"Audrey 18 May 2014"}{{1, 3218218, 1398127, 1602097, 8364837, 2017905, 6741864, 7933574, 6288705, 987739, 7590015, 9255760, 276220, 1717119, 255345}}

\adfDesignF{15}{46563}{1}{9778231}{4782969}{1}{5029921}{"Audrey 18 May 2014"}{{1, 5029921, 8159151, 4086825, 3432303, 2962007, 878911, 7100390, 8513781, 515573, 6816223, 3869399, 4296921, 6955856, 5830355}}

\adfDesignF{15}{46847}{1}{9837871}{5258672}{1}{7410376}{"Audrey 18 May 2014"}{{1, 7410376, 4651961, 4592817, 2467723, 2624596, 8350000, 2157173, 4817858, 3704187, 7213274, 3915366, 3028737, 427196, 3665961}}

\adfDesignF{15}{47241}{1}{9920611}{16384}{1}{1666079}{"Audrey 18 May 2014"}{{1, 1666079, 2514608, 2033066, 4751429, 5323331, 4231483, 1964117, 1025227, 8134786, 4597279, 4023049, 5441886, 6862318, 6955007}}

\adfDesignF{15}{47672}{1}{10011121}{3955732}{1}{5859876}{"Audrey 18 May 2014"}{{1, 5859876, 1705376, 699877, 6573029, 7570890, 3686841, 1495513, 5657650, 8825291, 2440230, 464404, 6810232, 3653594, 4623922}}

\adfDesignF{15}{48063}{1}{10093231}{8415804}{1}{3557875}{"Audrey 18 May 2014"}{{1, 3557875, 8484051, 2873478, 1486426, 8030604, 7264624, 9381203, 3521421, 7727151, 2062626, 9363963, 2321208, 3698332, 879654}}

~

\adfDesignF{17}{17303}{1}{4706417}{1335756}{1}{2584947}{"Audrey 18 May 2014"}{{1, 2584947, 1337808, 518584, 2226606, 2306736, 4697676, 486290, 1666517, 4043244, 2977666, 3817469, 1102409, 3940078, 4041935, 3195866, 3413921}}

\adfDesignF{17}{18696}{1}{5085313}{4875883}{1}{2430680}{"Audrey 18 May 2014"}{{1, 2430680, 2168679, 410302, 1621052, 518257, 3529652, 2619808, 4605832, 2412199, 4426641, 3490769, 1031473, 4520441, 1266414, 431047, 113944}}

\adfDesignF{17}{19425}{1}{5283601}{3650885}{1}{2723639}{"Audrey 18 May 2014"}{{1, 2723639, 1896715, 38549, 3124340, 1995897, 1333520, 497865, 35691, 1708351, 2159053, 2116099, 419835, 2051145, 5268113, 549552, 499640}}

\adfDesignF{17}{20856}{1}{5672833}{3915286}{1}{2905756}{"Audrey 18 May 2014"}{{1, 2905756, 1658501, 2795930, 4993293, 5259901, 3962570, 4628993, 332567, 2799768, 4495143, 5668113, 1741874, 384820, 3892791, 2278155, 3257321}}

\adfDesignF{17}{23838}{1}{6483937}{225977}{1}{1898010}{"Audrey 18 May 2014"}{{1, 1898010, 5466522, 1725001, 5159860, 4238997, 630150, 3982480, 1144373, 3277848, 5343421, 3465786, 6236683, 3827646, 2624621, 6462543, 2871491}}

\adfDesignF{17}{27549}{1}{7493329}{3048848}{1}{2831669}{"Audrey 18 May 2014"}{{1, 2831669, 229176, 5803357, 896015, 6672951, 5439053, 4169411, 117836, 2101243, 6718749, 3170912, 1817930, 1375421, 4326609, 6598711, 184259}}

\adfDesignF{17}{31059}{1}{8448049}{2314190}{1}{3652653}{"Audrey 18 May 2014"}{{1, 3652653, 5321493, 550014, 5098599, 1490460, 7661604, 4757032, 535627, 1244668, 5838756, 1984050, 5622686, 6884018, 3891272, 4064321, 539089}}

\adfDesignF{17}{32060}{1}{8720321}{6755027}{1}{436219}{"Audrey 18 May 2014"}{{1, 436219, 891420, 6507269, 7805917, 4524706, 5271474, 3350390, 4136773, 8075473, 5165106, 4435839, 1624446, 925214, 2029344, 3744342, 2118314}}

\adfDesignF{17}{35378}{1}{9622817}{2287808}{1}{6286281}{"Audrey 18 May 2014"}{{1, 6286281, 8702702, 645679, 5796382, 4165678, 2447333, 9540385, 7602875, 9498903, 8812316, 3994294, 5055749, 3804549, 1325822, 4256393, 4670011}}

\adfDesignF{17}{39734}{1}{10807649}{7129953}{1}{2719644}{"Audrey 18 May 2014"}{{1, 2719644, 317659, 9970581, 7029217, 9212184, 1336656, 75771, 702041, 705766, 4612553, 9858587, 4383199, 1282997, 2382822, 9907882, 348334}}

\adfDesignF{17}{39995}{1}{10878641}{5888417}{1}{3101860}{"Audrey 18 May 2014"}{{1, 3101860, 8456278, 2928956, 1456138, 623967, 6622387, 3806519, 3715016, 8841767, 10012981, 7951789, 8551343, 4806910, 9468872, 3009712, 3674632}}

\adfDesignF{17}{42020}{1}{11429441}{7533392}{1}{4023933}{"Audrey 18 May 2014"}{{1, 4023933, 4870994, 10828005, 246198, 3169736, 6313328, 7497268, 3129581, 3291130, 4757031, 3383328, 3676905, 2262927, 8128309, 2356305, 641667}}

\adfDesignF{17}{44211}{1}{12025393}{8831646}{1}{7894392}{"Audrey 18 May 2014"}{{1, 7894392, 6208059, 11504850, 6928069, 10179460, 8859773, 9437661, 10499347, 1728870, 1381581, 11177577, 5184924, 6309882, 4786381, 3311188, 2836522}}

\adfDesignF{17}{45989}{1}{12509009}{10935723}{1}{5741865}{"Audrey 18 May 2014"}{{1, 5741865, 6868636, 1800670, 7750690, 1736469, 2262046, 8530910, 5708554, 2351231, 2950493, 6330475, 4447639, 400803, 11786820, 12232006, 6663755}}

\adfDesignF{17}{46764}{1}{12719809}{9545478}{1}{9652259}{"Audrey 18 May 2014"}{{1, 9652259, 2700480, 7225913, 7735475, 12710767, 7603480, 4293120, 4953487, 2064741, 5438910, 12621294, 2466028, 9737844, 3291490, 8877183, 386000}}

\adfDesignF{17}{49535}{1}{13473521}{10116348}{1}{659333}{"Audrey 18 May 2014"}{{1, 659333, 10323345, 9110668, 10302930, 8874952, 12030237, 3071616, 12850618, 12836944, 10903051, 11563038, 8463972, 5334728, 5244727, 2001878, 11163172}}

% 1437 designs.
% Designs go before here.
%%%%%%%%%%%%%%%%%%%%%%%%%%%%%%%%%%%%%%%%%%%%%%%%%%%%%%%%%%
%%%%%%%%%%%%%%%%%%%%%%%%%%%%%%%%%%%%%%%%%%%%%%%%%%%%%%%%%%
%%%%%%%%%%%%%%%%%%%%%%%%%%%%%%%%%%%%%%%%%%%%%%%%%%%%%%%%%%
%%%%%%%%%%%%%%%%%%%%%%%%%%%%%%%%%%%%%%%%%%%%%%%%%%%%%%%%%%

\end{document}